\documentclass[11pt]{article}

\usepackage{graphicx}
    \graphicspath{ {figures/} }
    \usepackage[margin=0.1cm]{caption}
    \usepackage{float}
\usepackage{amsmath}
\usepackage{amsthm}
\usepackage{amssymb}
\usepackage[mathscr]{eucal}
\usepackage[all]{xy}
\usepackage{epsfig}
\usepackage{epstopdf}
\usepackage{amscd}
\usepackage[normalem]{ulem}
\usepackage{enumitem}
    \setlist[enumerate,1]{label=\textnormal{(\alph*)}}
\usepackage{thmtools, thm-restate}
\usepackage{array,ragged2e}

\usepackage[left=1.25in,right=1.25in,top=1.25in,bottom=1.25in,
            ]{geometry}
\usepackage{hyperref}
\usepackage{color}

\theoremstyle{plain}
\newtheorem{theorem}{Theorem}

\newtheorem{corollary}[theorem]{Corollary}

\newtheorem*{conjecture}{Conjecture}

\theoremstyle{remark}

\theoremstyle{definition}

\newtheorem*{question*}{Question}

\newcommand{\ka}[2]{$#1\textnormal{a}_{#2}$}
\newcommand{\kn}[2]{$#1\textnormal{n}_{#2}$}

\newcommand{\knotmosaictable}{below in Section \ref{sec:table}}

\makeatletter
\newcommand\xleftrightarrow[2][]{%
  \ext@arrow 9999{\longleftrightarrowfill@}{#1}{#2}}
\newcommand\longleftrightarrowfill@{%
  \arrowfill@\leftarrow\relbar\rightarrow}
\makeatother

\title{Tabulating Knot Mosaics: Crossing Number 10 or Less 
}
\author{Aaron Heap, Douglas Baldwin, James Canning and Greg Vinal}

\begin{document}

\maketitle
\begin{abstract} The study of knot mosaics is based upon representing knot diagrams using a set of tiles on a square grid. This branch of knot theory has many unanswered questions, especially regarding the efficiency with which we draw knots as mosaics. While any knot or link can be displayed as a mosaic, for most of them it is still unknown what size of mosaic (mosaic number) is necessary and how many non-blank tiles (tile number) are necessary to depict a given knot or link. We implement an algorithmic programming approach to find the mosaic number and tile number of all prime knots with crossing number 10 or less. We also introduce an online repository which includes a table of knot mosaics and a tool that allows users can create and identify their own knot mosaics.
\end{abstract}


\section{Introduction}

\begin{figure}[ht]
  \centering
  \includegraphics{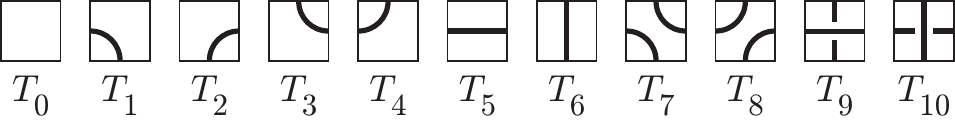}\\
  \caption{Possible tiles used for constructing knot mosaics.}
  \label{fig:tiles}
\end{figure}

\begin{figure}[ht]
  \centering
  \includegraphics{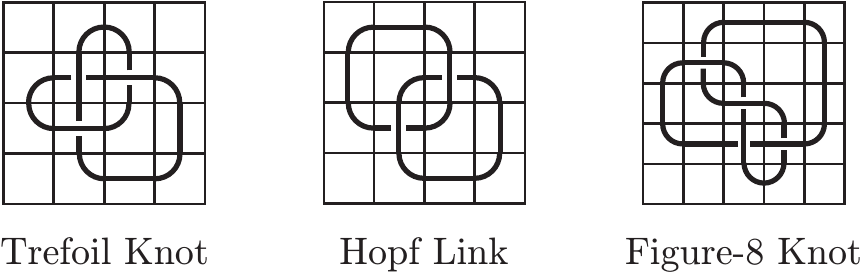}\\
  \caption{Examples of knot mosaics.}
  \label{fig:examples}
\end{figure}

The study of mosaic knot theory was first introduced by Lomonaco and Kauffman \cite{Lom-Kauff} in 2008. A \emph{knot mosaic} is a two-dimensional representation of a knot (or link), made up of a finite set of tiles arranged onto a square array. The possible tiles for constructing these mosaics are shown in Figure \ref{fig:tiles}. Knot mosaics provide a rigid construction of a knot diagram. Some examples are shown in Figure \ref{fig:examples}.

An \emph{$n\times n$ knot mosaic}, or \emph{$n$-mosaic}, is a knot mosaic with $n$ rows and $n$ columns. The first two examples in Figure \ref{fig:examples} are 4-mosaics, whereas the third example is a 5-mosaic. While a knot mosaic could depict a knot or a link, the primary focus of this paper is prime knots.

Let $K$ be a knot. There are many questions regarding the efficiency with which we can create a knot mosaic for $K$. What is the smallest mosaic on which $K$ can be depicted? How many non-blank tiles (tiles other than $T_0$) are necessary to create a knot diagram of $K$? How many crossing tiles ($T_9$ or $T_{10}$) are needed to create a knot diagram for $K$ on the smallest possible mosaic? To discuss these questions, we consider the following knot invariants related to knot mosaics. The \emph{mosaic number} of a knot $K$ is the smallest integer $n$ such that $K$ can be represented on an $n$-mosaic. The \emph{tile number} of a knot $K$ is the fewest number of non-blank tiles that are needed to represent $K$ on any $n$-mosaic. Examining the interplay between these numbers and the crossing number of $K$ have led to interesting results. For example, Ludwig, Evans, and Paat \cite{Ludwig} found an infinite family of knots whose mosaic number can be realized only when the crossing number is not. Similarly, Heap and Knowles \cite{Heap2} found knots whose tile number is realized only when the crossing number is not, and, perhaps more surprisingly, knots whose tile number is realized only when the mosaic number is not. For example, the mosaic number of the knot $9_{10}$ is 6, and to achieve a representation of this knot on a 6-mosaic, at least 32 non-blank tiles are needed. However, on a 7-mosaic, one can represent $9_{10}$ with only 27 non-blank tiles. That is, the tile number of $9_{10}$ is 27, but this cannot be realized on a 6-mosaic.

While it is fairly simple to create a knot mosaic for a given knot with few crossings and determine the mosaic number and tile number of that knot, it becomes much more difficult to do so as the number of crossings increases. The original paper on knot mosaics \cite{Lom-Kauff} provided example mosaics of a few different knots (the unknot, trefoil, and figure eight) and links (Hopf link, Solomon's knot, Borromean rings, and various unlinks). In a 2008 preprint (published in 2014) Kuriya and Shehab \cite{Kuriya} provided a few more knot mosaic examples, including the knots $4_1$, $5_1$, $5_2$, $6_2$, $6_3$, and $7_4$. They determined the mosaic number of all but one of these, and they correctly conjectured that the mosaic number of $6_3$ is 6. However, the first significant attempt at knot mosaic tabulation came in 2018, when Lee, Ludwig, Paat, and Peiffer \cite{Lee2} provided mosaics for all 36 prime knots with crossing number 8 or less, and they successfully provided the mosaic number of all of these knots. In 2018 and 2019, the tile numbers for all of these knots were provided in \cite{Heap1} and \cite{Heap2}, along with the tile numbers for every prime knot with mosaic number 6 and crossing number 9 or larger. Heap and LaCourt were able to determine the tile numbers for many knots with mosaic number 7 and crossing number 9 or more. In this paper, we complete this task for all 250 prime knots with crossing number 10 or less.

\begin{theorem}\label{thm:main-thm}
  Every prime knot with crossing number 10 or less has:
  \begin{enumerate}
    \item Mosaic number 7 or less, and
    \item Tile number 31 or less.
  \end{enumerate}
\end{theorem}

The specific details for individual knots are given below in Theorems \ref{thm:6-29}, \ref{thm:7-29}, and \ref{thm:7-31}. These results answer Question 3.8 posed in \cite{Lee2}, asking for the determination of the tile number for every knot with 10 or fewer crossings.

For a deeper overview of background material, Adams et al. \cite{Encyclopedia} provide a quality summary of the primary topics in this paper, including knots, knot mosaics, and the Dowker-Thistlethwaite codes we will use below. Additionally, we have created an online repository for information related to knot mosaics. This website, Knot Mosaic Space \cite{KnotSpace}, provides a table of knot mosaics and an interactive tool that allows users to build and identify their own knot mosaics.

\section{Determining Tile Numbers}
There are 49 prime knots with crossing number 9, and the mosaic number and tile number have been determined previously for 29 of them. There are 165 prime knots with crossing number 10, and the mosaic number and tile number have been determined previously for 58 of them. These results are summarized in Table \ref{table:summary}.

\begin{table}[h]
\centering
{\renewcommand{\arraystretch}{1.3} 
\begin{tabular}{|>{\raggedright\arraybackslash}m{.2\linewidth}|m{.73\linewidth}|}
\hline
\multicolumn{2}{|>{\centering\arraybackslash}m{.96\linewidth}|}{\text{Knots with Mosaic Number $m$, Tile Number $t$, and Crossing Number 9 or 10}} \\
\hline
$m=6$, $t=22$ &  $9_5$, $9_{20}$ \hfill \cite{Heap1} \\
\hline
$m=6$, $t=24$  &  $9_8$, $9_{11}$, $9_{12}$, $9_{14}$, $9_{17}$, $9_{19}$, $9_{21}$, $9_{23}$, $9_{26}$, $9_{27}$, $9_{31}$, $10_{41}$, $10_{44}$, $10_{85}$, $10_{100}$, $10_{116}$, $10_{124}$, $10_{125}$, $10_{126}$, $10_{127}$, $10_{141}$, $10_{143}$, $10_{148}$, $10_{155}$, $10_{159}$ \hfill \cite{Heap2} \\
\hline
$m=6$, $t=27$  &  $9_{1}$, $9_{2}$, $9_{3}$, $9_{4}$, $9_{7}$, $9_{9}$, $9_{13}$, $9_{24}$, $9_{28}$, $9_{37}$, $9_{46}$, $9_{48}$, $10_{1}$, $10_{2}$, $10_{3}$, $10_{4}$, $10_{12}$, $10_{22}$, $10_{28}$, $10_{34}$, $10_{63}$, $10_{65}$, $10_{66}$, $10_{75}$, $10_{78}$, $10_{140}$, $10_{142}$, $10_{144}$ \hfill \cite{Heap2} \\
\hline
$m=6$, $t=27^*$  &  $9_{10}$, $10_{11}$, $10_{20}$, $10_{21}$ \hfill \cite{Heap2}, \cite{Heap3} \\
\hline
$m=6$, $t= \; ?^{\dag}$  &  $9_{16}$, $9_{35}$, $10_{61}$, $10_{62}$, $10_{64}$, $10_{74}$, $10_{76}$, $10_{77}$, $10_{139}$ \hfill \cite{Heap2} \\
\hline
$m=7$, $t=27$  &  $9_6$, $9_{15}$, $9_{18}$, $10_{5}$, $10_{6}$, $10_{7}$, $10_{8}$, $10_{9}$, $10_{10}$, $10_{13}$, $10_{14}$, $10_{15}$, $10_{16}$, $10_{17}$, $10_{18}$, $10_{19}$, $10_{24}$, $10_{25}$, $10_{26}$, $10_{29}$, $10_{30}$, $10_{31}$, $10_{32}$, $10_{33}$, $10_{35}$, $10_{36}$, $10_{38}$, $10_{39}$ \hfill \cite{Heap3} \\
\hline
\multicolumn{2}{|>{\raggedright\arraybackslash}m{.96\linewidth}|}{$^*$Tile number realized on a 7-mosaic, not on a 6-mosaic.} \\
\multicolumn{2}{|>{\raggedright\arraybackslash}m{.96\linewidth}|}{$^{\dag}$32 non-blank tiles needed on a 6-mosaic; tile number previously unknown.} \\
\hline
\end{tabular}
}
\caption{Prime knots with crossing number 9 or 10 whose mosaic number and tile number have been determined previously.}
\label{table:summary}
\end{table}

\begin{figure}[ht]
  \centering
  \includegraphics{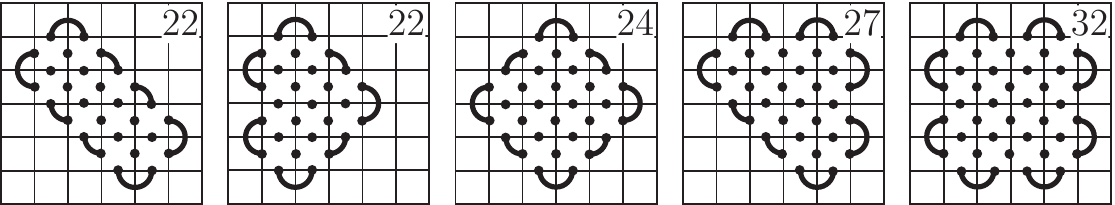}\\
    \caption{All Space-efficient 6-mosaics, with number of non-blank tiles indicated}
    \label{fig:layouts-6}
  \end{figure}

Knot mosaics for every prime knot with mosaic number 6 or less were given in \cite{Heap1} and \cite{Heap2}. The authors accomplished this by first providing a list of possible space-efficient layouts (up to rotation) for a 6-mosaic, shown in Figure \ref{fig:layouts-6}. That is, if the number of non-blank tiles within the mosaic is to be minimized, it will have an outer shell as depicted, and the interior of this shell must be populated with tiles $T_7$, $T_8$, $T_9$, and $T_{10}$ (from Figure \ref{fig:tiles}). Going through these layouts and considering every possible option on the interior, they were able to determine every prime knot with mosaic number 6, and they were able to determine the tile number for most of these. These included all prime knots with crossing number 8 or less and some prime knots with crossing number 9 through 13. The tile number exceptions were the knots with mosaic number 6 needing 32 non-blank tiles, as increasing to a 7-mosaic could result in a decrease in non-blank tiles, such as the knot $9_{10}$ mentioned above.

\begin{figure}[ht]
    \centering
    \includegraphics[width=\linewidth]{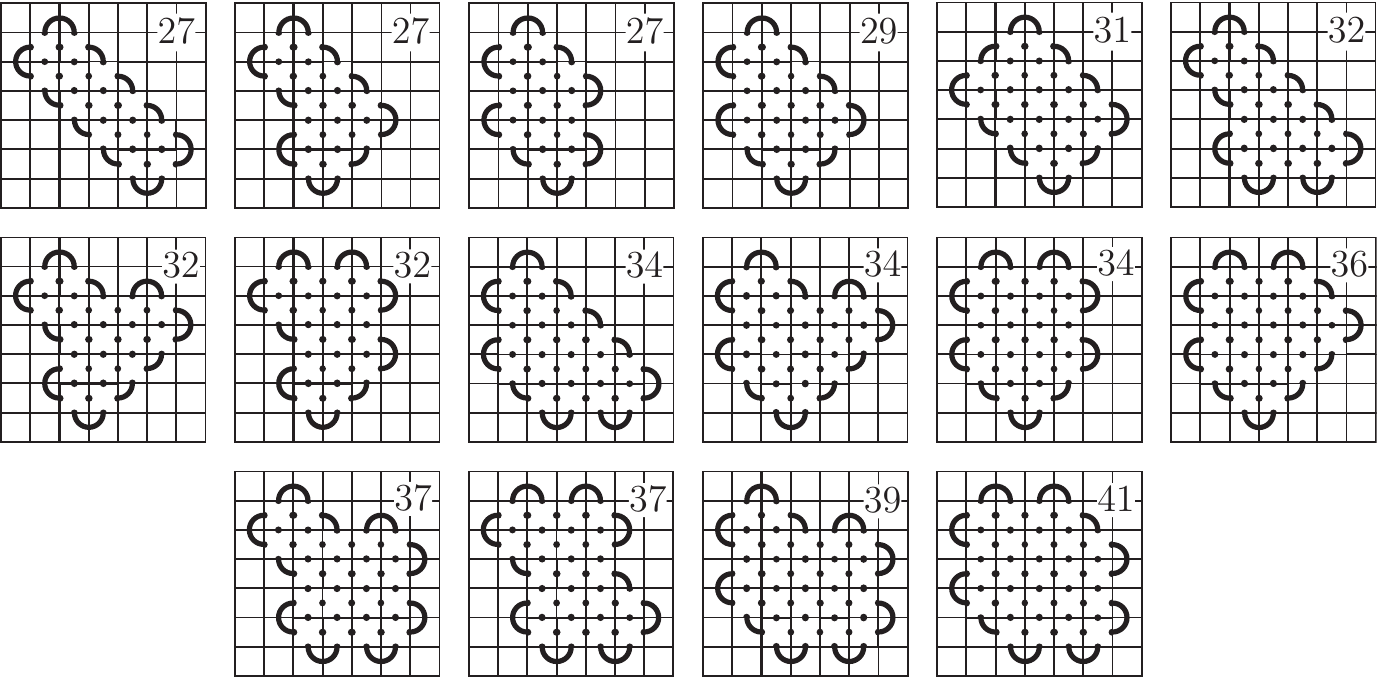}
    \caption{All Space-efficient 7-mosaics, with number of non-blank tiles indicated}
    \label{fig:layouts-7}
  \end{figure}

Heap and LaCourt \cite{Heap3} extended this work to 7-mosaics. Again, they first determined the possible space efficient layouts, shown in Figure \ref{fig:layouts-7}. Because of the vast number of options for completing the interior of these layouts, they limited themselves to determining which prime knots could be obtained on the first three layouts, each of which use 27 non-blank tiles.

\begin{figure}[ht]
    \centering
    \begin{minipage}{0.18\linewidth}
        \captionsetup{skip=3pt}
        \includegraphics[width=\linewidth]{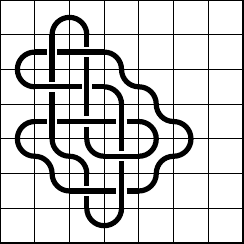}
        \caption*{$10_{74}$}
    \end{minipage} \hspace{.11in}
    \begin{minipage}{0.18\linewidth}
        \captionsetup{skip=3pt}
        \includegraphics[width=\linewidth]{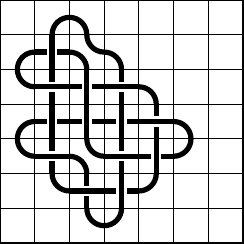}
        \caption*{$10_{76}$}
    \end{minipage} \hspace{.11in}
    \begin{minipage}{0.18\linewidth}
        \captionsetup{skip=3pt}
        \includegraphics[width=\linewidth]{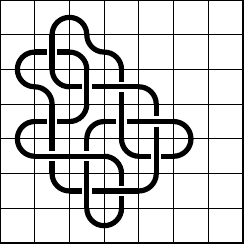}
        \caption*{$10_{77}$}
    \end{minipage} \hspace{.11in}
    \begin{minipage}{0.18\linewidth}
        \captionsetup{skip=3pt}
        \includegraphics[width=\linewidth]{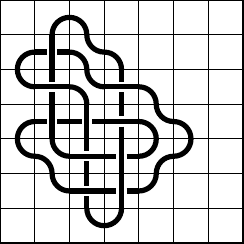}
        \caption*{$10_{139}$}
    \end{minipage}
    \caption{Prime knots with mosaic number 6, tile number 29, and crossing number 9 or 10.}
    \label{fig:6-32and7-29}
\end{figure}

Using these layouts as a starting point, we are now ready to complete the story for all prime knots with crossing number 10 or less by determining the mosaic number and tile number for each of the remaining knots. We first see that every knot with crossing number 9 or 10 whose mosaic number is 6 and whose tile number was previously unknown has tile number 29. Knot mosaics for these are provided in Figure \ref{fig:6-32and7-29}.

\begin{theorem}\label{thm:6-29}
  The following prime knots have mosaic number 6, needing 32 non-blank tiles on a 6-mosaic, and tile number 29 realized on a 7-mosaic: $9_{16}$, $9_{35}$, $10_{61}$, $10_{62}$, $10_{64}$, $10_{74}$, $10_{76}$, $10_{77}$, and $10_{139}$.
\end{theorem}

\begin{figure}[ht]
    \centering
    \begin{minipage}{0.18\linewidth}
        \captionsetup{skip=3pt}
        \centering
        \includegraphics[width=\linewidth]{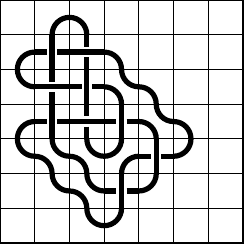}
        \caption*{$9_{16}$}
    \end{minipage} \hfill
    \begin{minipage}{0.18\linewidth}
        \captionsetup{skip=3pt}
        \includegraphics[width=\linewidth]{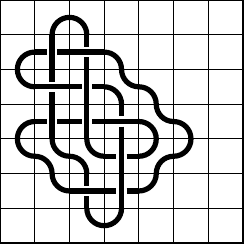}
        \caption*{$9_{35}$}
    \end{minipage} \hfill
    \begin{minipage}{0.18\linewidth}
        \captionsetup{skip=3pt}
        \includegraphics[width=\linewidth]{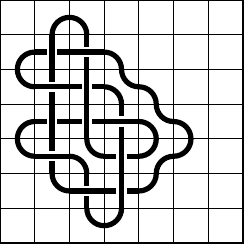}
        \caption*{$10_{61}$}
    \end{minipage} \hfill
    \begin{minipage}{0.18\linewidth}
        \captionsetup{skip=3pt}
        \includegraphics[width=\linewidth]{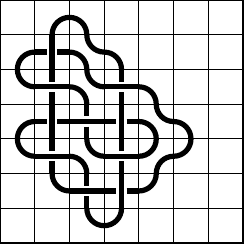}
        \caption*{$10_{62}$}
    \end{minipage} \hfill
    \begin{minipage}{0.18\linewidth}
        \captionsetup{skip=3pt}
        \includegraphics[width=\linewidth]{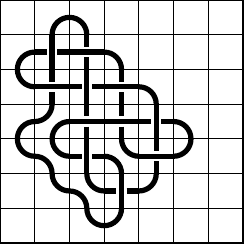}
        \caption*{$10_{64}$}
    \end{minipage}
\end{figure}

Viewing this more broadly, this means that every knot with crossing number 9 or 10 and mosaic number 6 that needs 32 non-blank tiles to be realized on a 6-mosaic (the last layout of Figure \ref{fig:layouts-6}) does not have its tile number realized on a 6-mosaic. There were thirteen such knots identified in \cite{Heap2}, and we now know that some have tile number 27 (see \cite{Heap3}), and the rest have tile number 29. None of them have tile number 32. We will say more about this for higher crossing numbers in Section \ref{sec:additional}.

The prime knots with mosaic number 7 and tile number 27 (those that fit on the first three layouts of Figure \ref{fig:layouts-7}) were given in \cite{Heap3}. We next provide the complete list of knots with crossing number 9 or 10 that fit on the fourth layout of Figure \ref{fig:layouts-7}. Knot mosaics for these are provided \knotmosaictable.

\begin{theorem}\label{thm:7-29}
  The following prime knots have mosaic number 7 and tile number 29:
  \begin{enumerate}
    \item $9_{22}$, $9_{25}$, $9_{29}$, $9_{30}$, $9_{32}$, $9_{33}$, $9_{34}$, $9_{36}$, $9_{38}$, $9_{39}$, $9_{42}$, $9_{43}$, $9_{44}$, $9_{45}$, $9_{47}$, $9_{49}$;
    \item $10_{23}$, $10_{27}$, $10_{37}$, $10_{40}$, $10_{42}$, $10_{43}$, $10_{45}$, $10_{46}$, $10_{47}$, $10_{48}$, $10_{49}$, $10_{50}$, $10_{51}$, $10_{52}$, $10_{53}$, $10_{54}$, $10_{55}$, $10_{56}$, $10_{57}$, $10_{67}$, $10_{68}$, $10_{69}$, $10_{70}$, $10_{71}$, $10_{72}$, $10_{73}$, $10_{79}$, $10_{82}$, $10_{83}$, $10_{84}$, $10_{86}$, $10_{87}$, $10_{90}$, $10_{91}$, $10_{92}$, $10_{93}$, $10_{94}$, $10_{95}$, $10_{101}$, $10_{102}$, $10_{103}$, $10_{106}$, $10_{107}$, $10_{112}$, $10_{113}$, $10_{114}$, $10_{117}$, $10_{128}$, $10_{129}$, $10_{130}$, $10_{131}$, $10_{132}$, $10_{133}$, $10_{134}$, $10_{135}$, $10_{136}$, $10_{145}$, $10_{146}$, $10_{147}$, $10_{149}$, $10_{150}$, $10_{151}$, $10_{152}$, $10_{153}$, $10_{156}$, $10_{158}$, $10_{160}$, $10_{161}$, $10_{162}$, $10_{163}$, $10_{164}$.
  \end{enumerate}
\end{theorem}

The prime knots with crossing number 9 or 10 that fit on the fifth layout of Figure \ref{fig:layouts-7} have tile number 31. Again, knot mosaics for these are provided \knotmosaictable.

\begin{theorem}\label{thm:7-31}
The following prime knots have mosaic number 7 and tile number 31:
\begin{enumerate}
    \item $9_{40}$, $9_{41}$;
    \item $10_{58}$, $10_{59}$, $10_{60}$, $10_{80}$, $10_{81}$, $10_{88}$, $10_{89}$, $10_{96}$, $10_{97}$, $10_{98}$, $10_{99}$, $10_{104}$, $10_{105}$, $10_{108}$, $10_{109}$, $10_{110}$, $10_{111}$, $10_{115}$, $10_{118}$, $10_{119}$, $10_{120}$, $10_{121}$, $10_{122}$, $10_{123}$, $10_{137}$, $10_{138}$, $10_{154}$, $10_{157}$, $10_{165}$.
\end{enumerate}
\end{theorem}

With this theorem, the mosaic number and tile number of every prime knot with crossing number 10 or less have been determined.

\section{Algorithm for Finding Knot Mosaics}

In order to prove the results of this paper, we simply need to find knot mosaics with the stated characteristics. As with the results from these previous papers, this requires a certain amount of trust from the reader that all possibilities were considered and are correct as stated. The authors of \cite{Heap2} and \cite{Heap3} exerted a significant amount of time and mental effort to fill the layouts of Figures \ref{fig:layouts-6} and \ref{fig:layouts-7}. Although they were able to simplify their workload by making use of symmetries, knot mosaic ambient isotopies, and the fact that they needed at least 8 or 9 crossing tiles in each mosaic, the number of possibilities that arise in the case of 7-mosaics increases significantly as the number of interior tile locations needing to be filled increases, limiting the scope of \cite{Heap3} to the first three layouts of Figure \ref{fig:layouts-7}, which use 27 non-blank tiles. In order to extend their work to the other layouts, the authors of this paper set out to automate the process by piecing together several computer programs into one pipeline. This pipeline creates every possible knot mosaic that can fit within a given layout and identifies the knot that is depicted. In this section, we describe each program within the pipeline.

\begin{figure}[ht]
    \centering
    \begin{minipage}{0.2\linewidth}
        \centering
        \includegraphics[width=\linewidth]{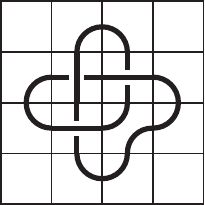}
    \end{minipage}
    \begin{minipage}{0.5in}
        \centering
        $\xleftrightarrow{\hspace{.2in}}$
    \end{minipage}
    \begin{minipage}{1.2in}
        \centering
        $\left[
          \begin{array}{cccc}
            0 & 2 & 1 & 0 \\
            2 & 10 & 9 & 1 \\
            3 & 9 & 8 & 4 \\
            0 & 3 & 4 & 0 \\
          \end{array}
        \right]$
    \end{minipage}
    \caption{Writing a knot mosaic as a matrix.}
    \label{fig:mosaic-matrix}
\end{figure}

To make coding easier, we rewrite a knot mosaic as a matrix with numerical entries. To do this, we simply replace each tile $T_{i}$ in the mosaic with the number $i$, as shown in Figure \ref{fig:mosaic-matrix}. By representing a knot as a matrix, computers can efficiently store and read the knot diagram.

\subsection{Step 1: Generate All Possible Knot Mosaics}
We begin the pipeline by selecting one of the identified layouts from Figure \ref{fig:layouts-7}. Each layout has a given number of interior positions that need to be filled. In Python, we create a matrix with numerical entries that represents the layout as described above, leaving the empty inner tiles as a null value. For a layout with $n$ interior tile positions to fill, the pipeline generates a vector of length $n$ for every possible combination of the entries 7, 8, 9, and 10, which correspond to the tiles with four connection points. (For 7-mosaics, $n\geq 13$.) Replacing the null values of the layout matrix with each of these vectors will create all possible mosaics for this layout. However, this requires $4^n$ vectors, and as $n$ increases it becomes infeasible to create so many vectors at once. To help alleviate this, we split the process in half.

Let $k = floor(\frac{n}{2})$. Using the \textit{product} function from the \textit{itertools} package, we create $4^k$ vectors of length $k$ and $4^{n-k}$ vectors of length $n-k$, including all possible configurations of integers 7, 8, 9, and 10 for each length. By concatenating every pair of vectors of the two lengths we can create all $4^n$ possibilities more quickly and using less processing memory.

We can obtain further optimization by restricting the number of crossing tiles, $T_9$ and $T_{10}$, in the mosaic. By only looking for knots with a minimum crossing number $m$ (in our case $m=9$), we only need to check mosaics with at least $m$ crossings. With this minimum, we only need to concatenate those vectors that together contain at least $m$ entries that are 9 or 10. Therefore, we count the number of entries that equal 9 or 10 in each of the smaller vectors, and only concatenate pairs of vectors that will result in a total of at least $m$ such entries. Doing so reduces the number of length $n$ vectors created from $4^n=2^{2n}$ to
\[\sum_{i=m}^n \binom{n}{i}2^i \cdot 2^{n-i} = 2^n\sum_{i=m}^n \binom{n}{i}.\]

\begin{figure}[ht]
  \centering
  \includegraphics[width=.95\linewidth]{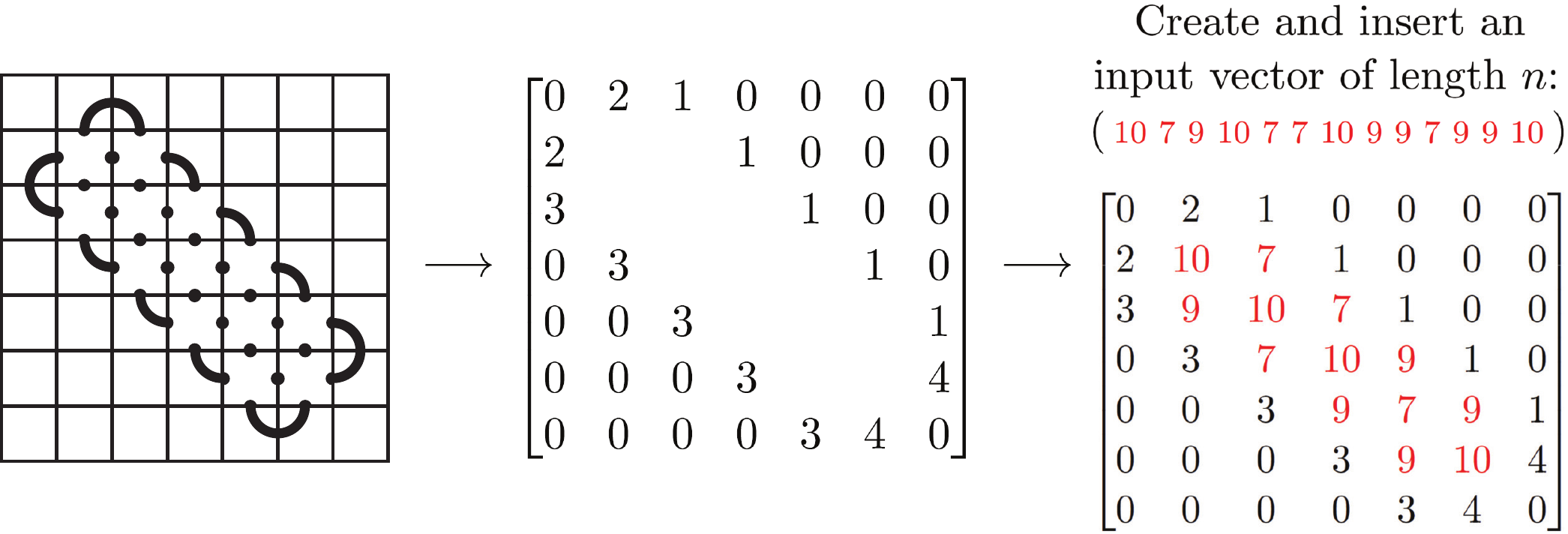}
  \caption{Generating the knot mosaic matrices.}
  \label{fig:pipline-step1}
\end{figure}

We then fill in the layout matrix by putting the $n$ values from the vectors into the $n$ empty spots of the matrix, one vector at a time. The process is illustrated in Figure \ref{fig:pipline-step1}. We send this matrix representation of a knot mosaic on to the next step in the pipeline.

\subsection{Step 2: Get the Reduced DT Code}
Now that we have a method to create all the possible knot mosaics for a particular layout, we need a way to identify which knot is represented. To do so, we will use the Dowker-Thistlethwaite code (DT code) \cite{DT}. The DT code is a representation of a knot diagram as a list of even integers and is determined using the following method:
\begin{enumerate}
    \item Start at an arbitrary base point and orientation on a knot with $n$ crossings.
    \item Travel along the knot in the direction of the orientation, labeling each crossing sequentially with integers 1 through 2$n$. When assigning an even number, if it is an undercrossing, assign the negative even number instead.
    \item  Each crossing will have two numbers, one even and one odd, which can be listed as an ordered pair. Arrange the $n$ odd numbers in ascending order with their associated even numbers, for example\\ $[(1,4),(3,-6,),(5,2)]$.
    \item Take the even numbers in the order determined by the odds ($4,-6,2$). This sequence of even numbers is the DT code.
\end{enumerate}

Although it is not a knot invariant, the DT code does uniquely identify knot diagrams, and every DT code of a knot can be reduced to the representation of the simplest knot diagram of that knot. That is, given any knot diagram, the DT code of that diagram can be reduced to a unique DT code for that knot type.

A Python function is used to produce the DT code of a given knot mosaic matrix. The function identifies the first nonzero entry of the matrix, starting at the top left entry and continuing along the row and on to the next row, if necessary. (Note that the first nonzero entry must be 2.) If this first nonzero entry is in the $(i,j)$ entry of the matrix, then the function moves to the $(i,j+1)$ entry. Depending on the tile corresponding to the number in this position, the function proceeds to the $(i,j+2)$ or $(i+1,j+1)$ entry. The function continues this process by incrementing the row or column index by $\pm 1$, depending on the tile corresponding to the numerical entry in the current position. This allows the function to trace along the knot (or link component) until it reaches the initial starting location again. As this process traverses along the knot, at each crossing (entry 9 or 10), the function records the associated label required to obtain the DT code. If the program ever returns to the initial position before it has passed the expected number of nonzero positions, then the knot diagram is a link and may be disregarded since we are only interested in finding knots. Finally, the function obtains the DT code from the crossing labels.

Now that the DT code has been determined, we need to reduce the DT code in order to identify the knot. We do so using \textit{KnotScape}, a program developed by Thistlethwaite and Hoste \cite{Thistle}. We adapted this C program to take multiple DT codes in sequence and output the reduced DT code of each one. If the program finds a composite knot, we eliminate it since our goal is to identify prime knots.

\begin{figure}[ht]
    \centering
    \begin{minipage}{.22\linewidth}
        \centering
        \includegraphics[width=\linewidth]{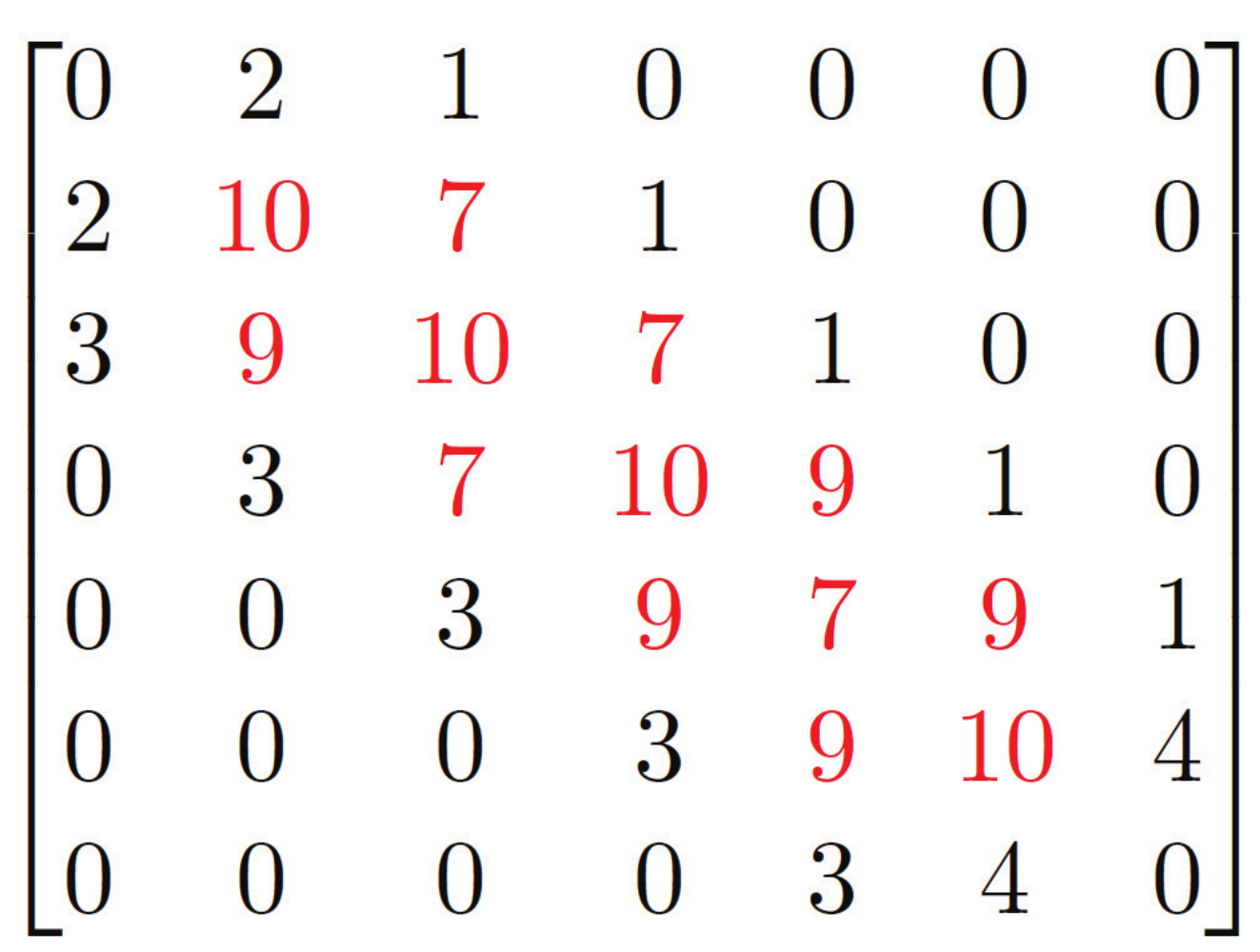}
    \end{minipage}
    \begin{minipage}{.15in}
        \centering
        $\rightarrow$
    \end{minipage}
    \begin{minipage}{.25\linewidth}
        \centering
        Get the DT code: \newline
        $\left[
           \begin{smallmatrix}
             10 & 14 & 16 & 12 & 18 & 8 & 2 & 4 & 6 \\
           \end{smallmatrix}
         \right]$
    \end{minipage}
    \begin{minipage}{.15in}
        \centering
        $\rightarrow$
    \end{minipage}
    \begin{minipage}{.27\linewidth}
        \centering
        Reduce the DT code: \newline
        $\left[
           \begin{smallmatrix}
             4 & 12 & 14 & 16 & 18 & 2 & 10 & 6 & 8 \\
           \end{smallmatrix}
         \right]$
    \end{minipage}
    \begin{minipage}{.15in}
        \centering
        $\rightarrow$
    \end{minipage}
    \begin{minipage}{.1\linewidth}
        \centering
        Identify the knot: $9_6$
    \end{minipage}
    \caption{Identifying the knot mosaic.}
    \label{fig:pipline-step3}
\end{figure}

\subsection{Step 3: Identify the Knot}
The final step in the pipeline is to identify the name of the knot. Using the reduced DT code from the previous step, we reference a table that contains every prime knot with crossing number 16 or less and their reduced DT codes. If the DT code is found in the table, the program outputs the name of the knot and the original vector that created the knot mosaic, so that we can recreate the knot mosaic that has been identified. If the DT code is not found in the table, it is from a knot with crossing number greater than 16. Steps 2 and 3 are depicted in Figure \ref{fig:pipline-step3}.

\subsection{Pipeline Output and Testing}
Once the pipeline has completed its task, we run a bash script to count the number of crossing in each knot mosaic and identify the knot mosaic with the fewest number of crossings for each unique knot. We compare all the knots found on the current layout with all the knots found on previous layouts and smaller mosaics. We identify any knots that were not previously found, any knots whose tile number was not realized on a smaller mosaic, and any knots whose crossing number had not been realized previously.

To test the pipeline, we ran it using every possible layout for a 4-mosaic, 5-mosaic, and 6-mosaic (see Figure \ref{fig:layouts-6}). We also ran the pipeline on the first three layouts for 7-mosaics (see Figure \ref{fig:layouts-7}). The results matched and confirmed the findings from \cite{Heap1}, \cite{Heap2}, \cite{Heap3}, \cite{Kuriya}, and \cite{Lee2}.

\section{Additional Results}\label{sec:additional}

After successful testing of the pipeline, we began running the larger 7-mosaic layouts through it. We have successfully found all of the knots that can be realized on the fourth layout (29 non-blank tiles) listed in Figure \ref{fig:layouts-7}, which includes knots with crossing number up to 14. We have also completed the fifth layout (31 non-blank tiles) for all knots with crossing number 16 or less. This layout allows for 17 crossing tiles, but our database of DT codes does not include information for knots with crossing number 17 or larger. The result is that, when combined with the findings of \cite{Heap1}, \cite{Heap2}, \cite{Heap3}, \cite{Kuriya}, and \cite{Lee2}, we have completed the search for knot mosaics for every knot with crossing number 10 or less, and the mosaic number and tile number has been determined for all of these knots. This was summarized above in Theorem \ref{thm:main-thm}, with specific details for crossing numbers 9 and 10 provided in Table \ref{table:summary} and Theorems \ref{thm:6-29}, \ref{thm:7-29}, and \ref{thm:7-31}. For the rest of this section, we provide a few additional results that were achieved using our pipeline algorithm.

\subsection{6-Mosaics with Tile Number 32}
We discussed above how the tile number need not be realized on a minimal mosaic. In \cite{Heap3} and in Theorem \ref{thm:6-29} we see examples of this, where 32 non-blank tiles are necessary for the knot to be depicted on a 6-mosaic but the tile number is first realized on a 7-mosaic. Because of this, there are several knots with mosaic number 6 whose tile numbers were not previously known. It also seems natural to ask if there is any knot with mosaic number 6 and tile number 32. We know from above that there are no such knots with crossing number 10 or less. In fact, there are also no such knots with crossing number 11 and 12. One knot, \ka{11}{341}, has tile number 27 and was given in \cite{Heap3}. The rest have tile number 29 or 31. However, there are six knots with crossing number 13, mosaic number 6, and tile number 32. We summarize the results here, with mosaics given \knotmosaictable, and we now know the tile number of every prime knot with mosaic number 6.

\begin{theorem}\label{thm:6mosaic7tile} The following prime knots have crossing number 11 or larger, mosaic number 6, needing 32 non-blank tiles on a 6-mosaic, and tile number strictly less than 32 realized on a 7-mosaic:
  \begin{enumerate}
    \item Tile number 27: \ka{11}{341};
    \item Tile number 29: \ka{11}{46}, \ka{11}{58}, \ka{11}{59}, \ka{11}{106}, \ka{11}{139}, \ka{11}{165}, \ka{11}{166}, \ka{11}{179}, \ka{11}{181}, \ka{11}{246}, \ka{11}{339}, \ka{11}{340}, \ka{11}{364}, \ka{12}{165}, \ka{12}{373}, \ka{12}{376}, \ka{12}{380}, \ka{12}{444}, \ka{12}{503}, \ka{13}{1236};
    \item Tile number 31: \ka{11}{43}, \ka{11}{44}, \ka{11}{47}, \ka{11}{247}, \ka{11}{342}, \ka{11}{367}, \kn{11}{71}, \kn{11}{72}, \kn{11}{73}, \kn{11}{74}, \kn{11}{75}, \kn{11}{76}, \kn{11}{77}, \kn{11}{78}, \ka{12}{119}, \ka{12}{169}, \ka{12}{379}, \ka{12}{722}, \ka{12}{803}, \ka{12}{1148}, \ka{12}{1149}, \ka{12}{1166}, \ka{13}{1461}, \ka{13}{4573}.
  \end{enumerate}
\end{theorem}

Space-efficient minimal mosaics were given for each of the following knots in \cite{Heap2}. Our pipeline did not find mosaics depicting these knots using the first five layouts of Figure \ref{fig:layouts-7}. Therefore, their tile number, 32, was realized on a 6-mosaic.

\begin{theorem}
The following knots are the only prime knots with mosaic number 6 and tile number 32: \ka{13}{1230}, \kn{13}{2399}, \kn{13}{2400}, \kn{13}{2401}, \kn{13}{2402}, and \kn{13}{2403}.
\end{theorem}

\subsection{Crossing Numbers Realized}

Many knots that are represented on as a knot mosaic with mosaic number or tile number realized necessarily have more crossings than the crossing number of the knot. Mosaics in which the crossing number is realized have been found for every knot with crossing number 10 or less.

\begin{theorem}\label{thm:CrossingRealized}
The knots $7_3$, $8_1$, $8_7$, and $8_8$ have mosaic number 6 and tile number 22, and their crossing number is realized only when at least 24 non-blank tiles are used. The knots $10_1$, $10_{34}$, and $10_{78}$ have mosaic number 6 and tile number 27, and their crossing number is realized with 32 non-blank tiles on a 6-mosaic or 29 non-blank tiles on a 7-mosaic. The crossing number for the following knots is first realized on a 7-mosaic in which the mosaic number or tile number is not realized, using the given number of non-blank tiles:
\begin{enumerate}
    \item 27 tiles: $8_6$, $8_9$, $9_9$, $9_{10}$, $9_{13}$, $9_{21}$, $9_{26}$;

    \item 29 tiles: $8_3$, $9_3$, $9_7$, $9_{12}$, $9_{15}$, $9_{16}$, $9_{19}$, $9_{24}$, $9_{37}$, $9_{46}$, $9_{48}$, $10_5$, $10_{11}$, $10_{13}$, $10_{14}$, $10_{15}$, $10_{16}$, $10_{18}$, $10_{21}$, $10_{22}$, $10_{24}$, $10_{31}$, $10_{33}$, $10_{35}$, $10_{36}$, $10_{38}$, $10_{39}$, $10_{62}$, $10_{63}$, $10_{65}$, $10_{74}$, $10_{139}$, $10_{140}$, $10_{142}$, $10_{144}$;

    \item 31 tiles: $9_4$, $9_{29}$, $9_{35}$, $10_{6}$, $10_{7}$, $10_{9}$, $10_{12}$, $10_{17}$,  $10_{20}$, $10_{37}$,  $10_{48}$,  $10_{50}$, $10_{51}$, $10_{56}$,  $10_{61}$, $10_{64}$, $10_{67}$, $10_{68}$, $10_{70}$, $10_{72}$, $10_{77}$, $10_{79}$, $10_{84}$, $10_{90}$, $10_{91}$, $10_{92}$, $10_{93}$, $10_{103}$, $10_{114}$, $10_{152}$, $10_{153}$, $10_{158}$, $10_{163}$;

    \item 34 tiles: $10_{3}$ and $10_{76}$.
\end{enumerate}
\end{theorem}

\begin{corollary}\label{cor:xNum}
The prime knots with crossing number 9 or 10 not listed in Theorem \ref{thm:CrossingRealized} have mosaics in which the crossing number, mosaic number, and tile number are realized.
\end{corollary}

\begin{figure}[h]
  \centering
  \includegraphics[width=.8\linewidth]{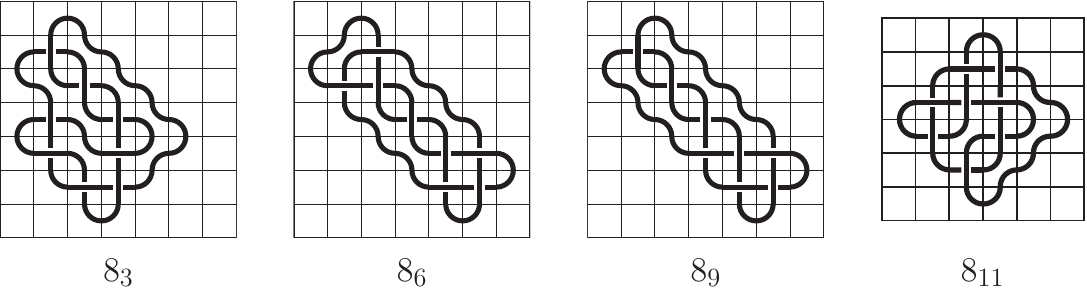}
  \caption{Simplest knot mosaics with crossing number realized for $8_3$, $8_6$, $8_9$, and $8_{11}$.}
  \label{fig:8-crossing}
\end{figure}

Additionally, we were able to answer Question 3.10 posed in \cite{Lee2}: Does there exist a representation of $8_3$, $8_6$, $8_9$, or $8_{11}$ on a 6-mosaic with only 8 crossing tiles? A positive answer for $8_{11}$ was provide in \cite{Heap2}, where the authors provided a mosaic of $8_{11}$ with mosaic number, tile number, and crossing number realized. The answer is no for the remaining knots, as they are all listed in Theorem \ref{thm:CrossingRealized}. Knot mosaics of each of these, with crossing number realized, are given in Figure \ref{fig:8-crossing}. A minimal mosaic for $7_3$, with crossing number realized, was given in \cite{Heap2}. Knot Mosaics for the remaining knots with crossing number 8, 9, or 10 realized can be found on the table of knots given \knotmosaictable


\subsection{Online Resources}

We have created an online repository for information related to knot mosaics. This website, Knot Mosaic Space \cite{KnotSpace}, provides access to the source code and output of the pipeline algorithm. It also includes a table of knot mosaics and a tool for building knot mosaics.

At the time of this writing, the pipeline is currently running for various layouts, searching for knots with crossing number greater than 10. While we limit the scope of this paper to tabulating knot mosaics for crossing number 10 or less, our interest and search continues for knots with larger crossing numbers. Results, as they become available, will be posted on the table of knot mosaics found at Knot Mosaic Space.

Knot Mosaic Space is also home to an interactive web-based tool that we developed for creating and identifying knot mosaics. Once the user builds and submits their knot mosaic, it is run through the pipeline algorithm described above. The knot mosaic is converted to its numerical matrix representation, the DT code is determined and reduced, and the name of the knot depicted is returned. Additionally, the results of this paper and the output of the pipeline are stored in a database for this tool so users can easily search for and display the knot mosaic representation of the knots that we have found. Users are able to search for mosaics based on whether they want the mosaic number, tile number, or crossing number realized.

\subsection{Different Layouts Producing the Same Knots}

It was previously known that the first two 6-mosaic layouts (both with 22 non-blank tiles) produce the same set of knots \cite{Heap2}. It was also found in \cite{Heap3} that the same is true for the first three 7-mosaic layouts (all with 27 non-blank tiles). That is, these three layouts all produce the same set of knots. In analyzing our results, we confirm these prior results, and we have reached the same conclusion for the sixth, seventh, and eighth 7-mosaic layouts (all with 32 non-blank tiles) and the ninth, tenth, and eleventh 7-mosaic layouts (all with 34 non-blank tiles). This provides more evidence to support the following conjecture, first proposed in \cite{Heap3}.

\begin{conjecture}
Space-efficient layouts of the same mosaic size and the same number of non-blank tiles produce the same set of prime knots.
\end{conjecture}


\section{Mosaics from Theorems}\label{sec:table}
In this section we present the mosaics for the knots referenced in the theorems. Mosaics that are marked with an $\ast$ are space-efficient mosaics that have more crossings than the crossing number of the knot they represent. These images were created using a program that takes the matrix representation of a mosaic and draws the mosaic using the Python PyCairo package.

\subsection{Mosaics from Theorems \ref{thm:7-29} and \ref{thm:7-31}}

These are the prime knots with crossing number 9 or 10, mosaic number 7, and tile number 29 or 31.

\begin{figure}[H]
    \centering
    \begin{minipage}{0.155\linewidth}
        \captionsetup{skip=3pt}
        \centering
        \includegraphics[width=\linewidth]{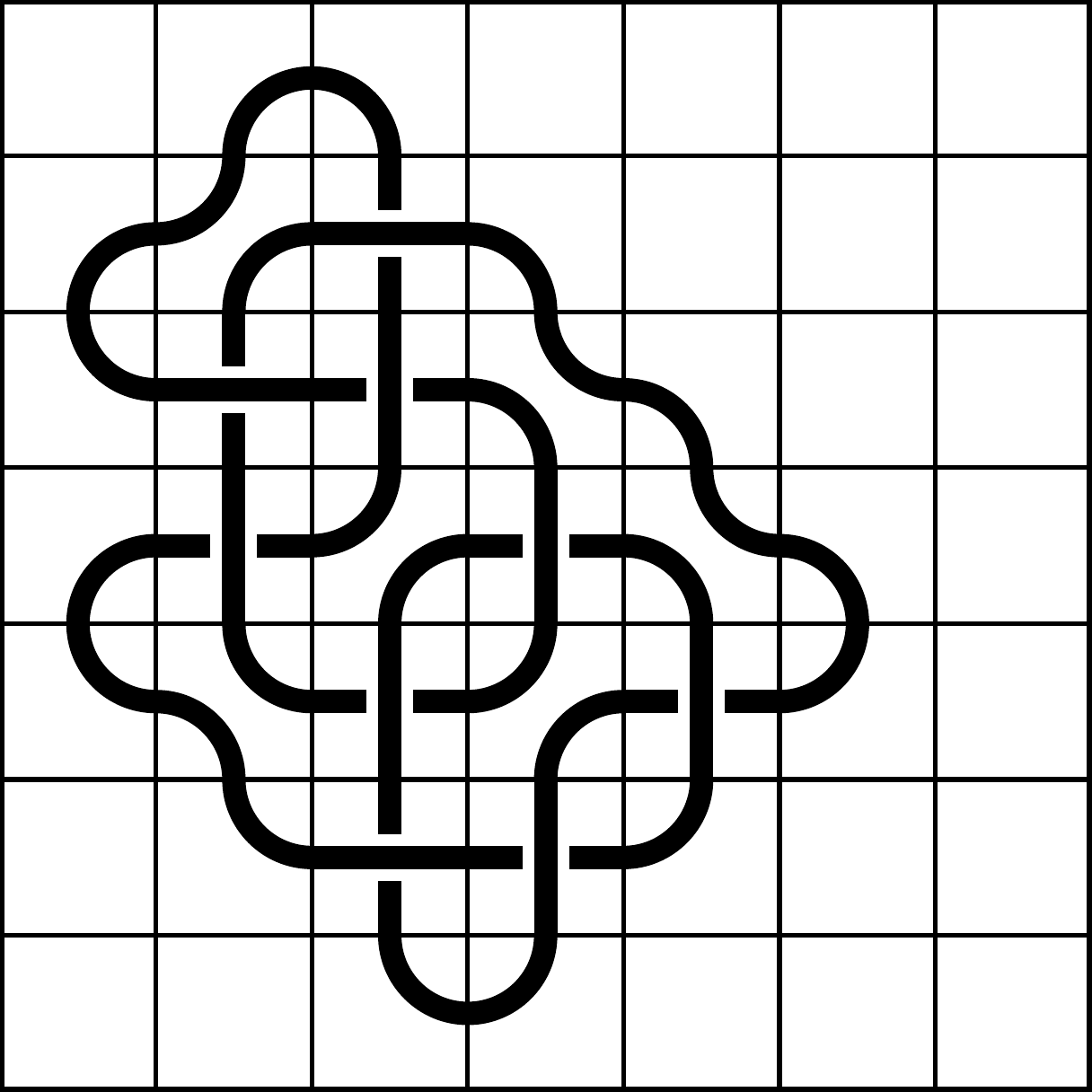}
        \caption*{$9_{22}$}
    \end{minipage} \hfill
     \begin{minipage}{0.155\linewidth}
        \captionsetup{skip=3pt}
        \centering
        \includegraphics[width=\linewidth]{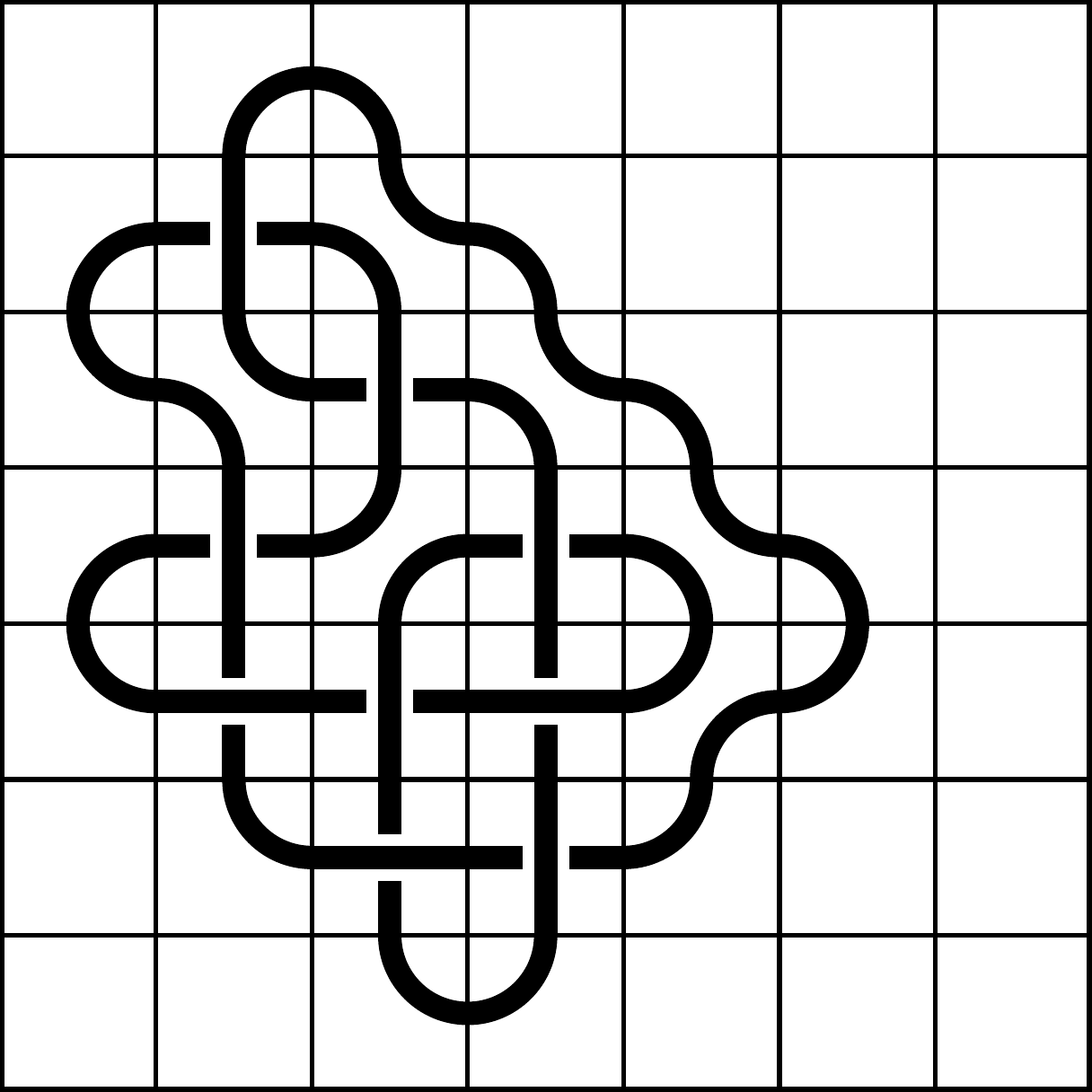}
        \caption*{ $9_{25}$ }
    \end{minipage} \hfill
    \begin{minipage}{0.155\linewidth}
        \captionsetup{skip=3pt}
        \centering
        \includegraphics[width=\linewidth]{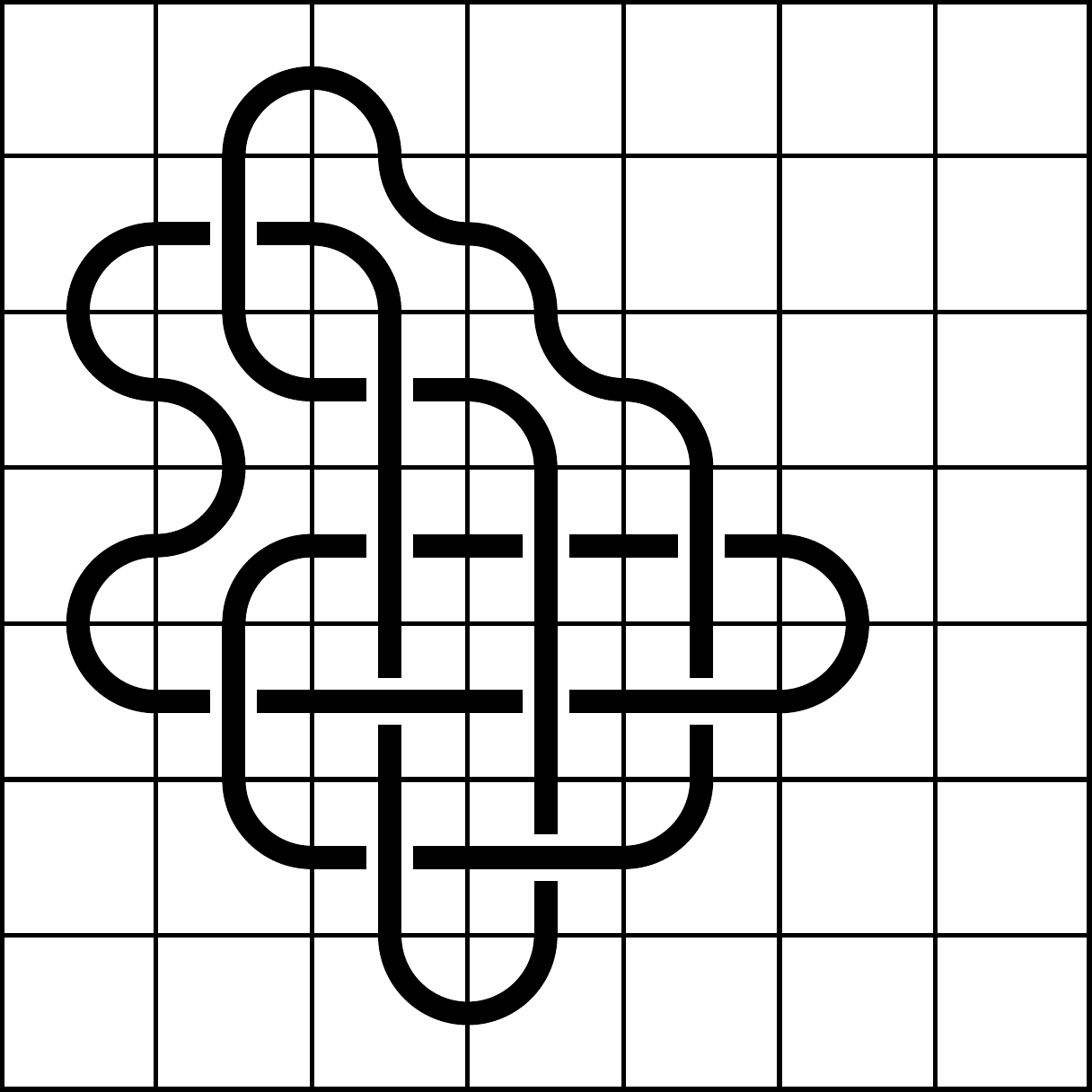}
        \caption*{\phantom{{\Large $\ast$}} $9_{29}$ {\Large $\ast$}}
    \end{minipage} \hfill
     \begin{minipage}{0.155\linewidth}
        \captionsetup{skip=3pt}
        \centering
        \includegraphics[width=\linewidth]{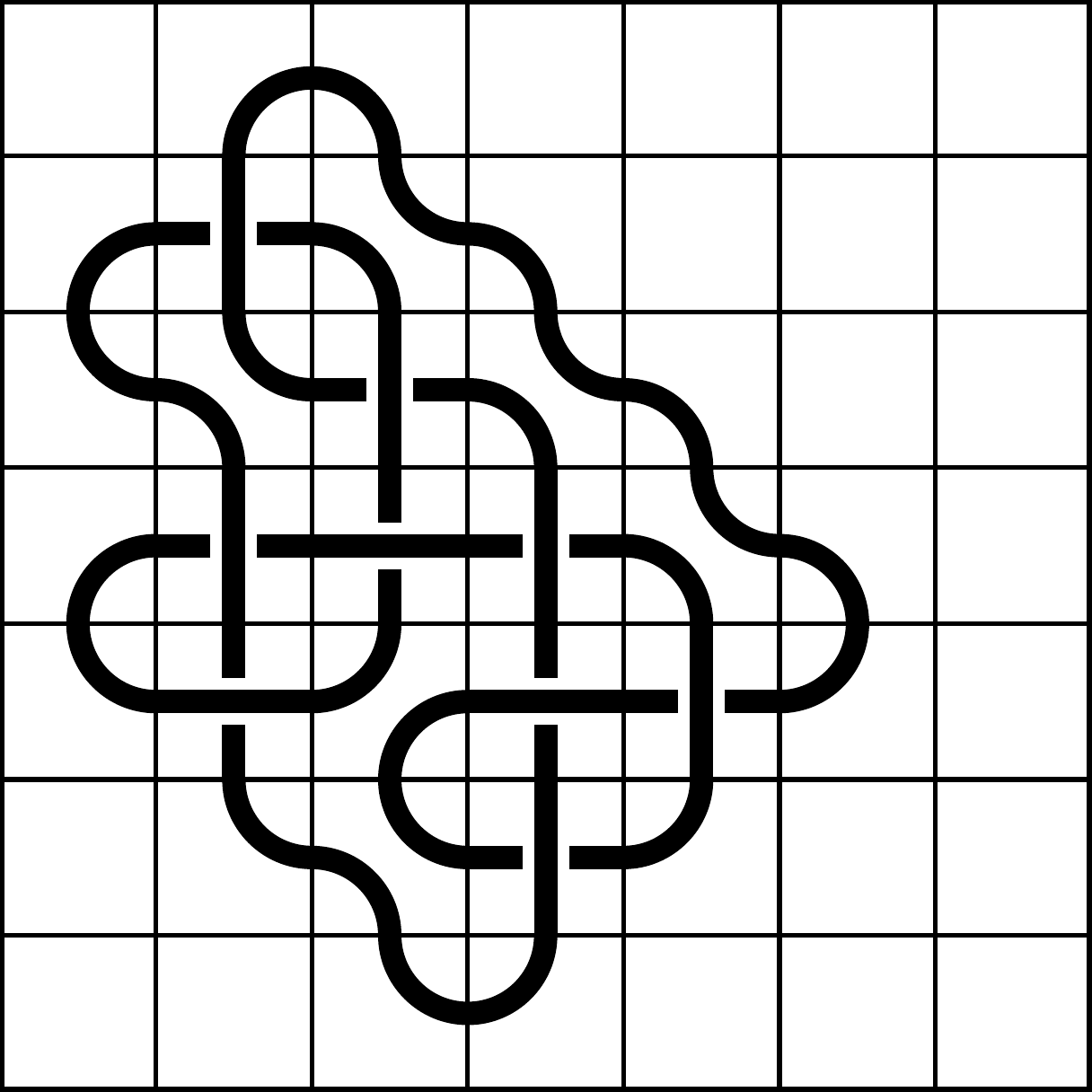}
        \caption*{ $9_{30}$ }
    \end{minipage} \hfill
    \begin{minipage}{0.155\linewidth}
        \captionsetup{skip=3pt}
        \centering
        \includegraphics[width=\linewidth]{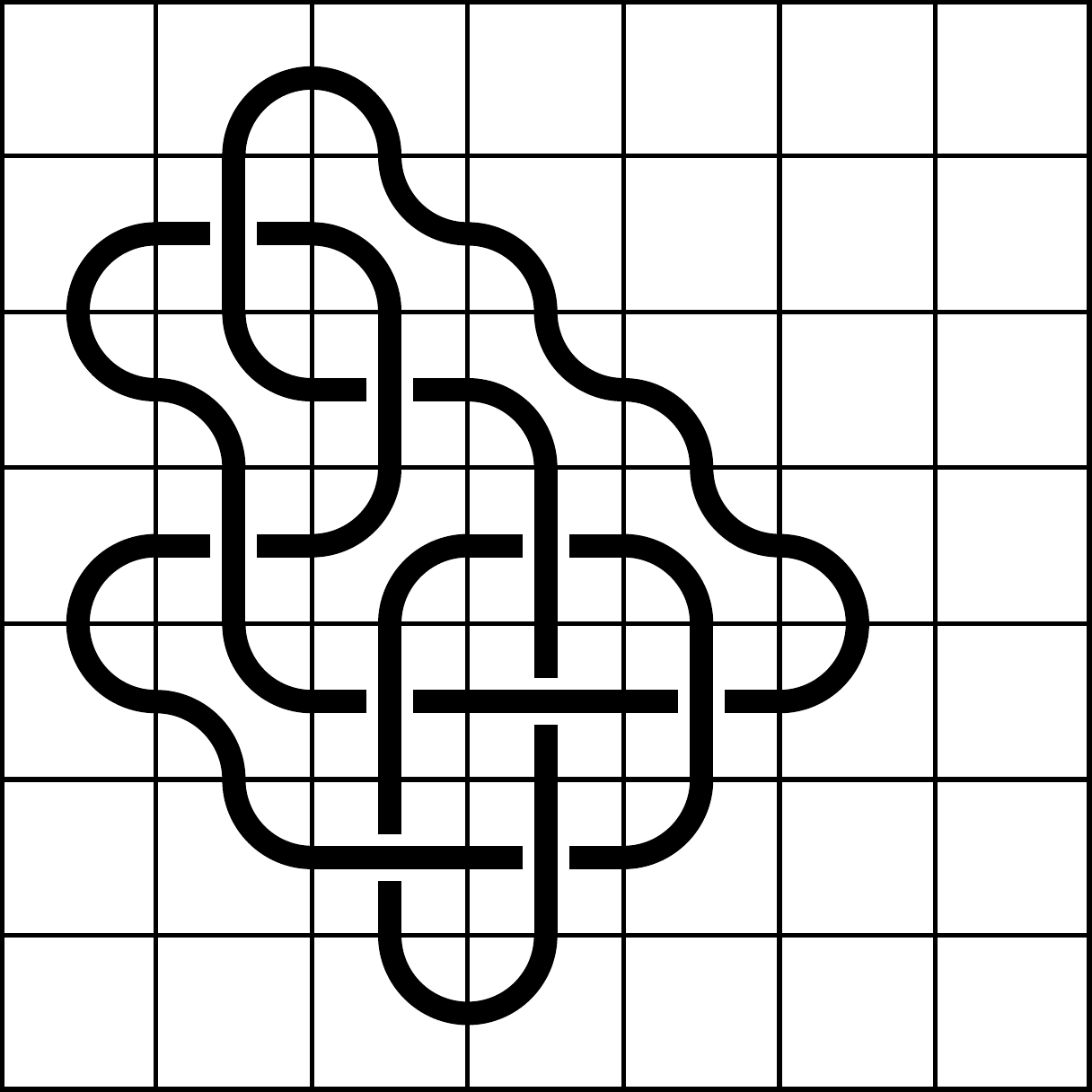}
        \caption*{ $9_{32}$ }
    \end{minipage}  \hfill
    \begin{minipage}{0.155\linewidth}
        \captionsetup{skip=3pt}
        \centering
        \includegraphics[width=\linewidth]{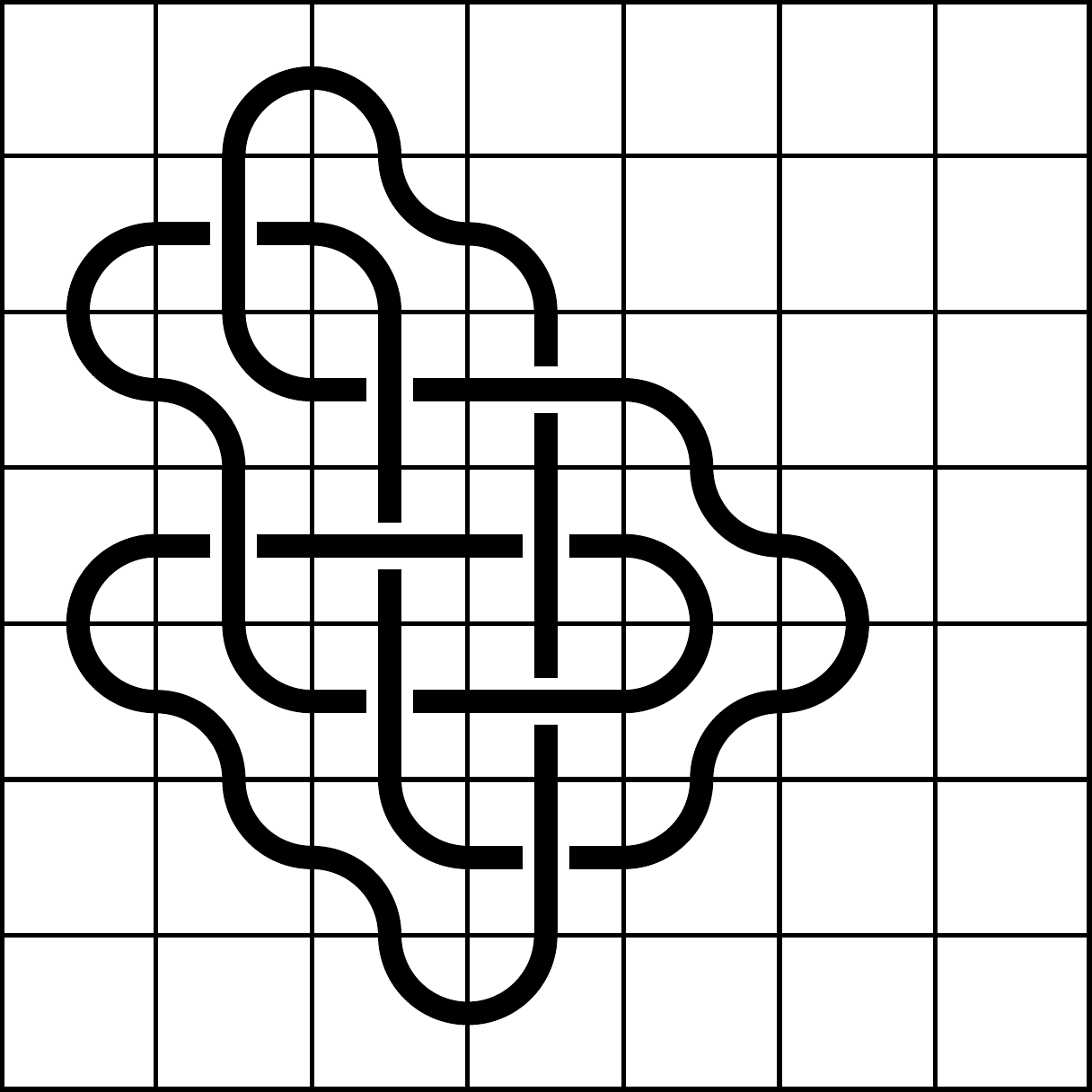}
        \caption*{ $9_{33}$ }
    \end{minipage} \newline
\end{figure}
\unskip

\begin{figure}[H]
     \centering
     \begin{minipage}{0.155\linewidth}
        \captionsetup{skip=3pt}
        \centering
        \includegraphics[width=\linewidth]{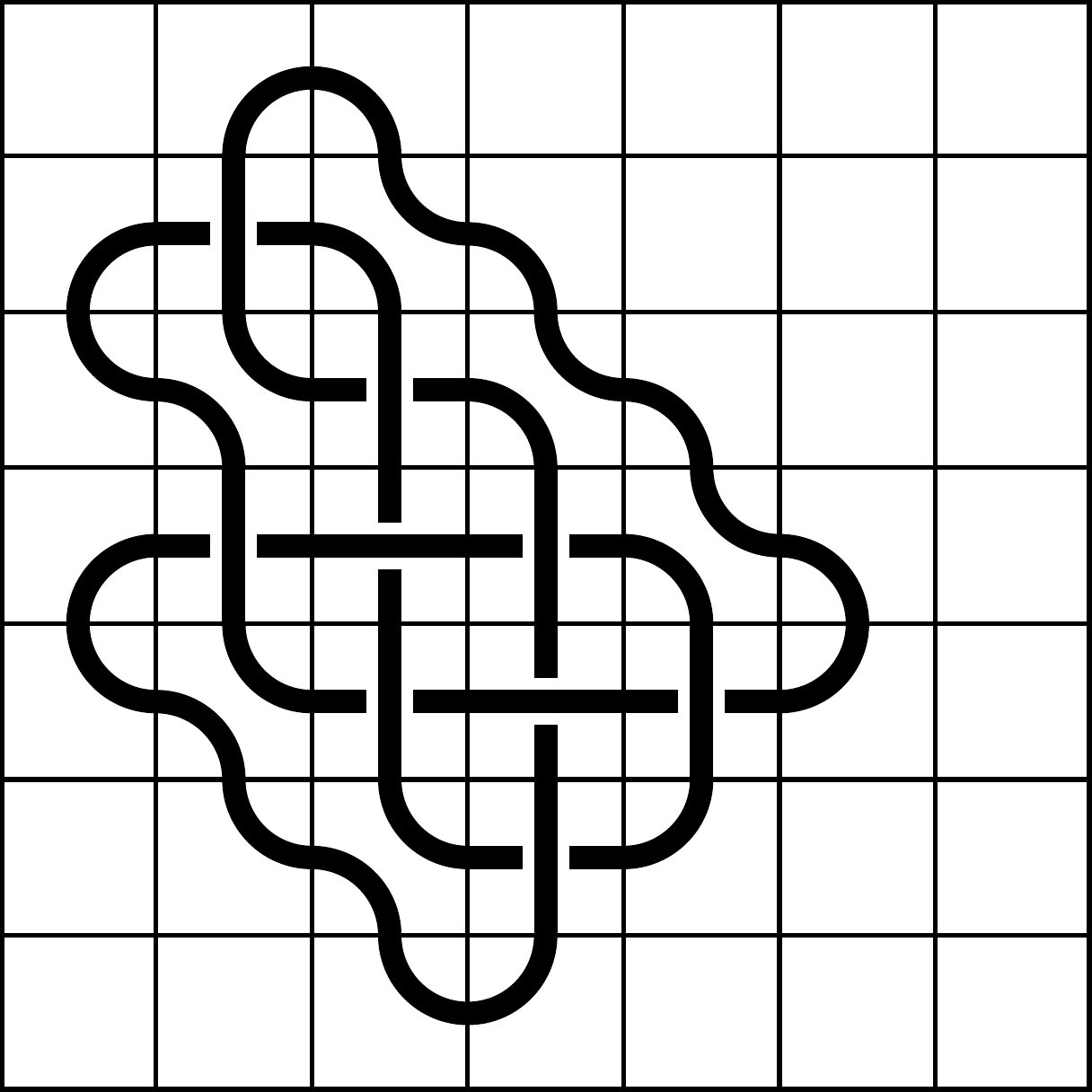}
        \caption*{$9_{34}$}
    \end{minipage} \hfill
    \begin{minipage}{0.155\linewidth}
        \captionsetup{skip=3pt}
        \centering
        \includegraphics[width=\linewidth]{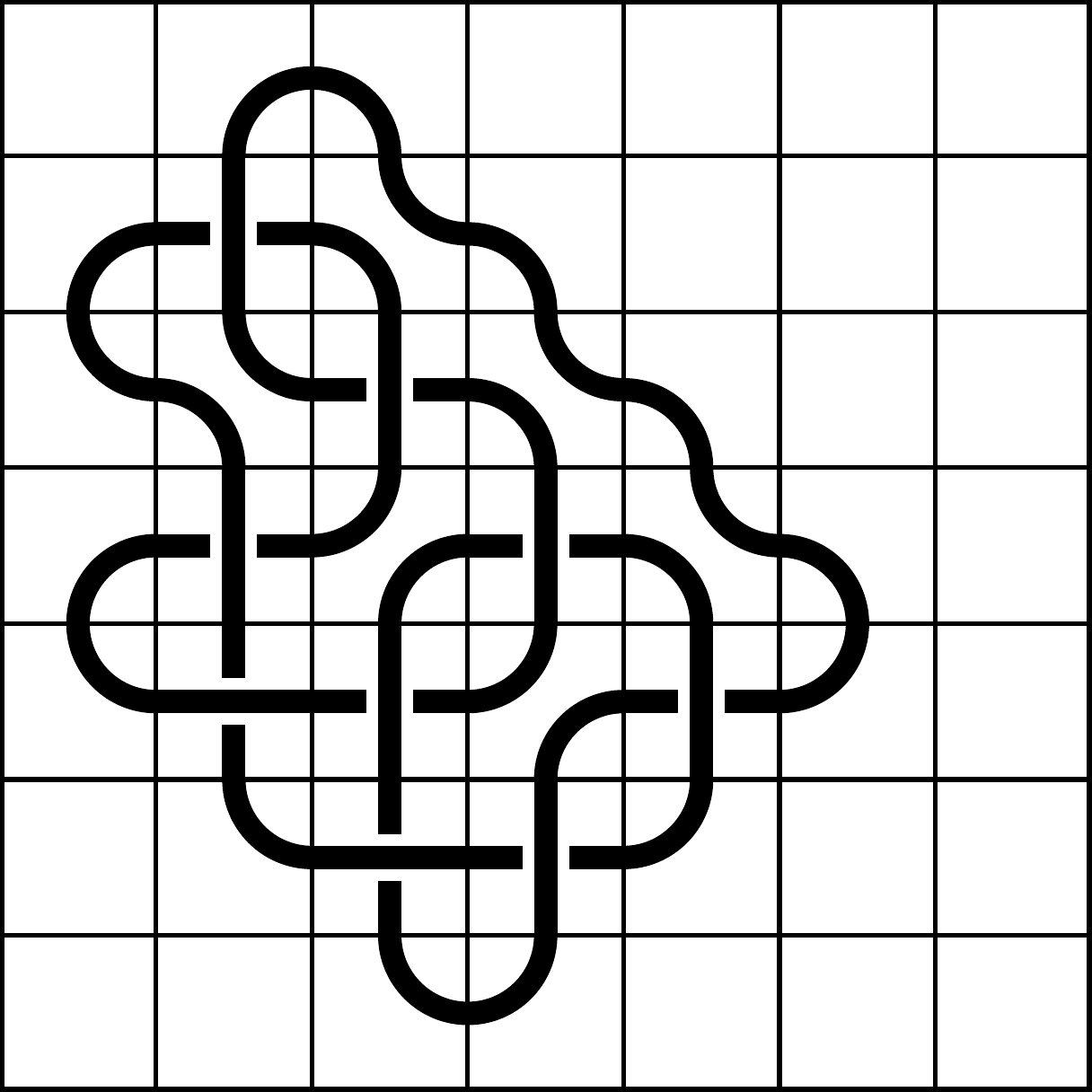}
        \caption*{ $9_{36}$ }
    \end{minipage} \hfill
     \begin{minipage}{0.155\linewidth}
        \captionsetup{skip=3pt}
        \centering
        \includegraphics[width=\linewidth]{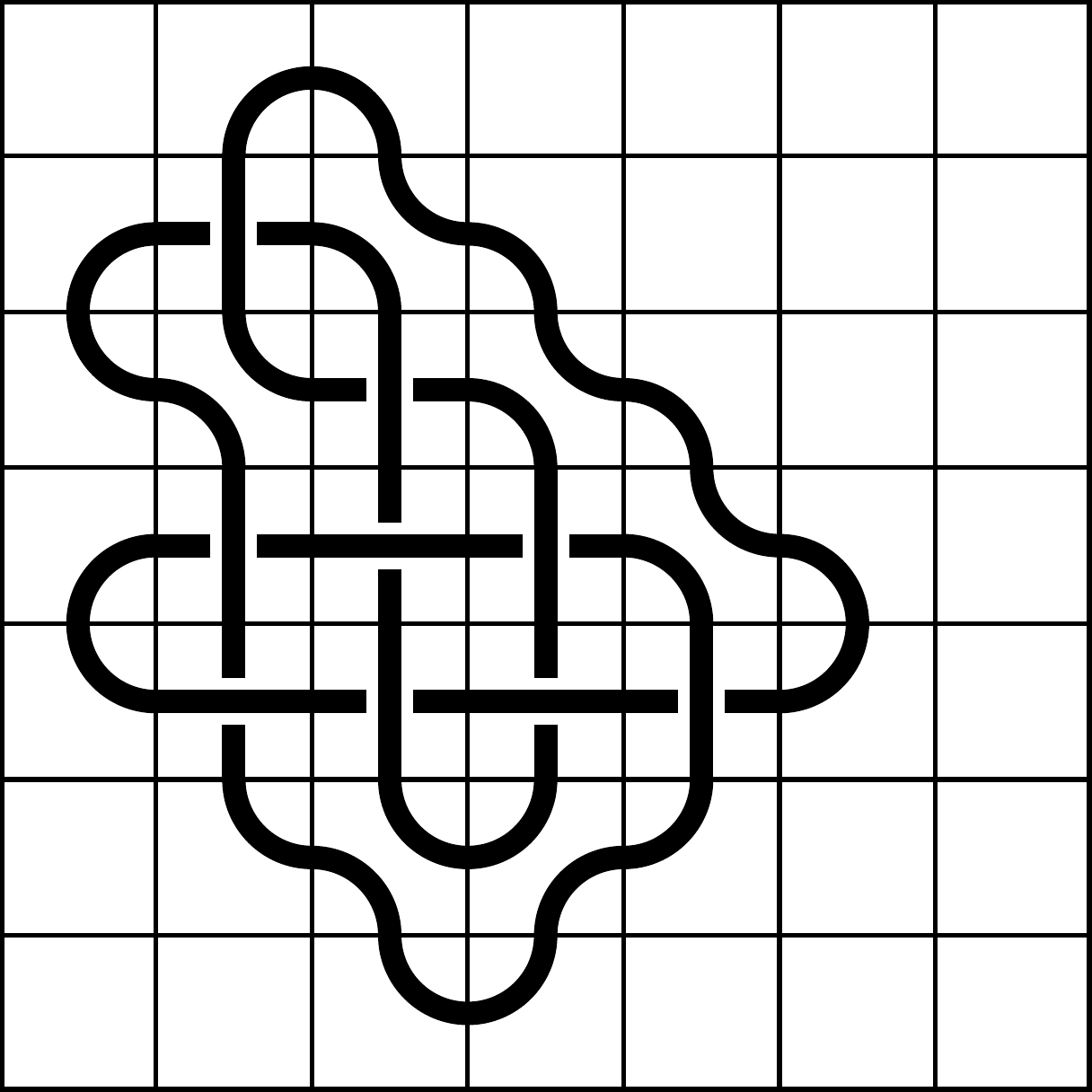}
        \caption*{$9_{38}$}
    \end{minipage} \hfill
    \begin{minipage}{0.155\linewidth}
        \captionsetup{skip=3pt}
        \centering
        \includegraphics[width=\linewidth]{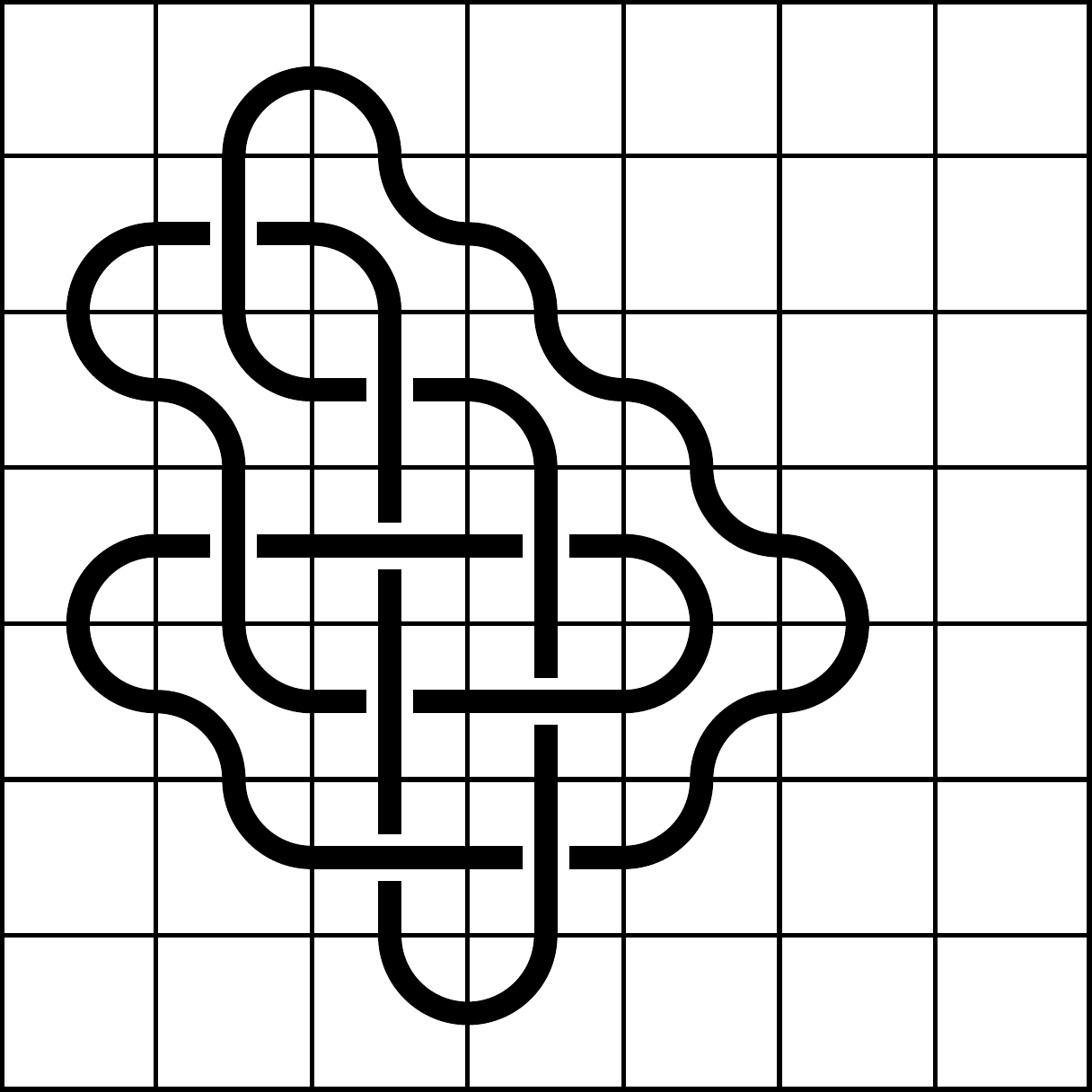}
        \caption*{$9_{39}$ }
    \end{minipage}   \hfill
    \begin{minipage}{0.155\linewidth}
        \captionsetup{skip=3pt}
        \centering
        \includegraphics[width=\linewidth]{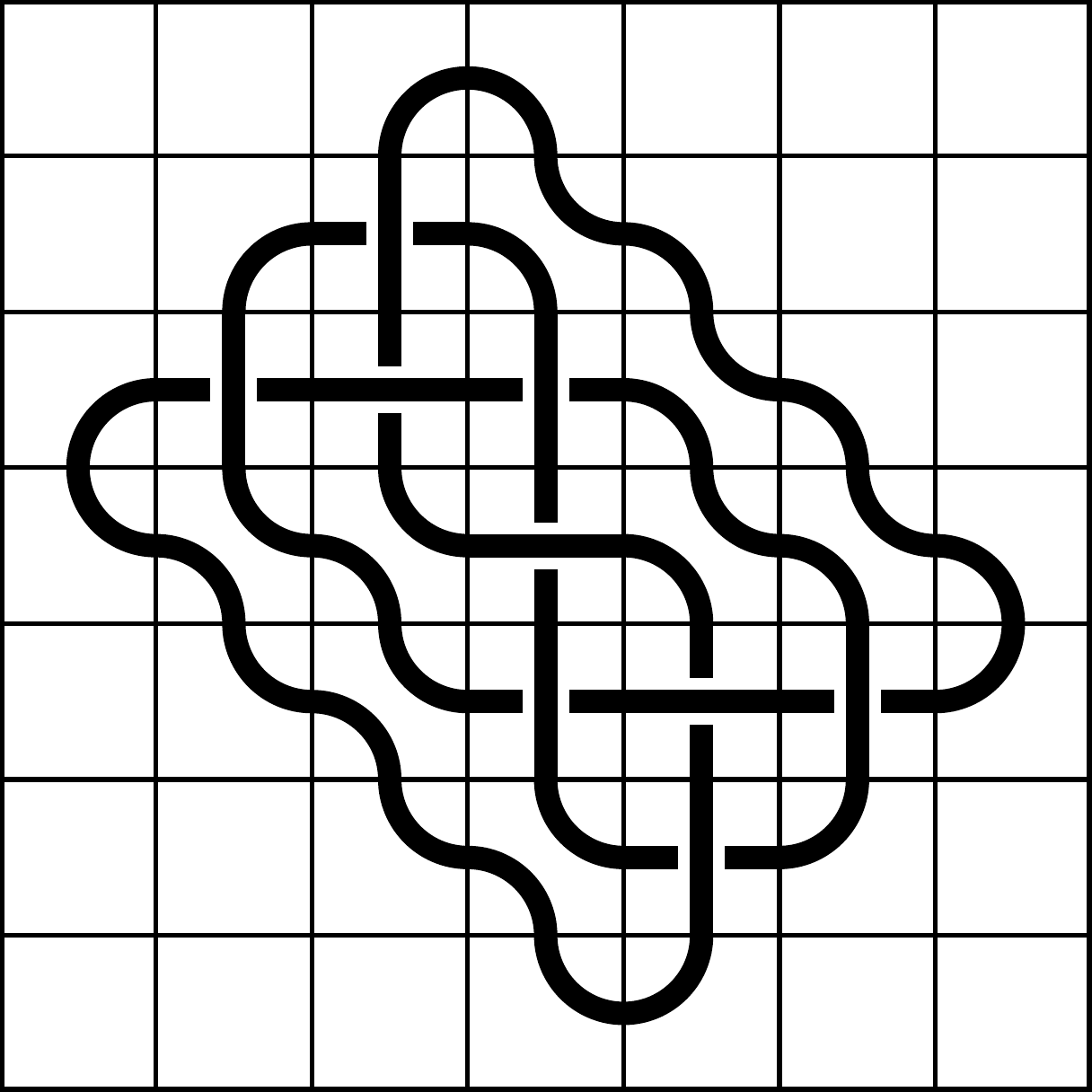}
        \caption*{$9_{40}$}
    \end{minipage} \hfill
     \begin{minipage}{0.155\linewidth}
        \captionsetup{skip=3pt}
        \centering
        \includegraphics[width=\linewidth]{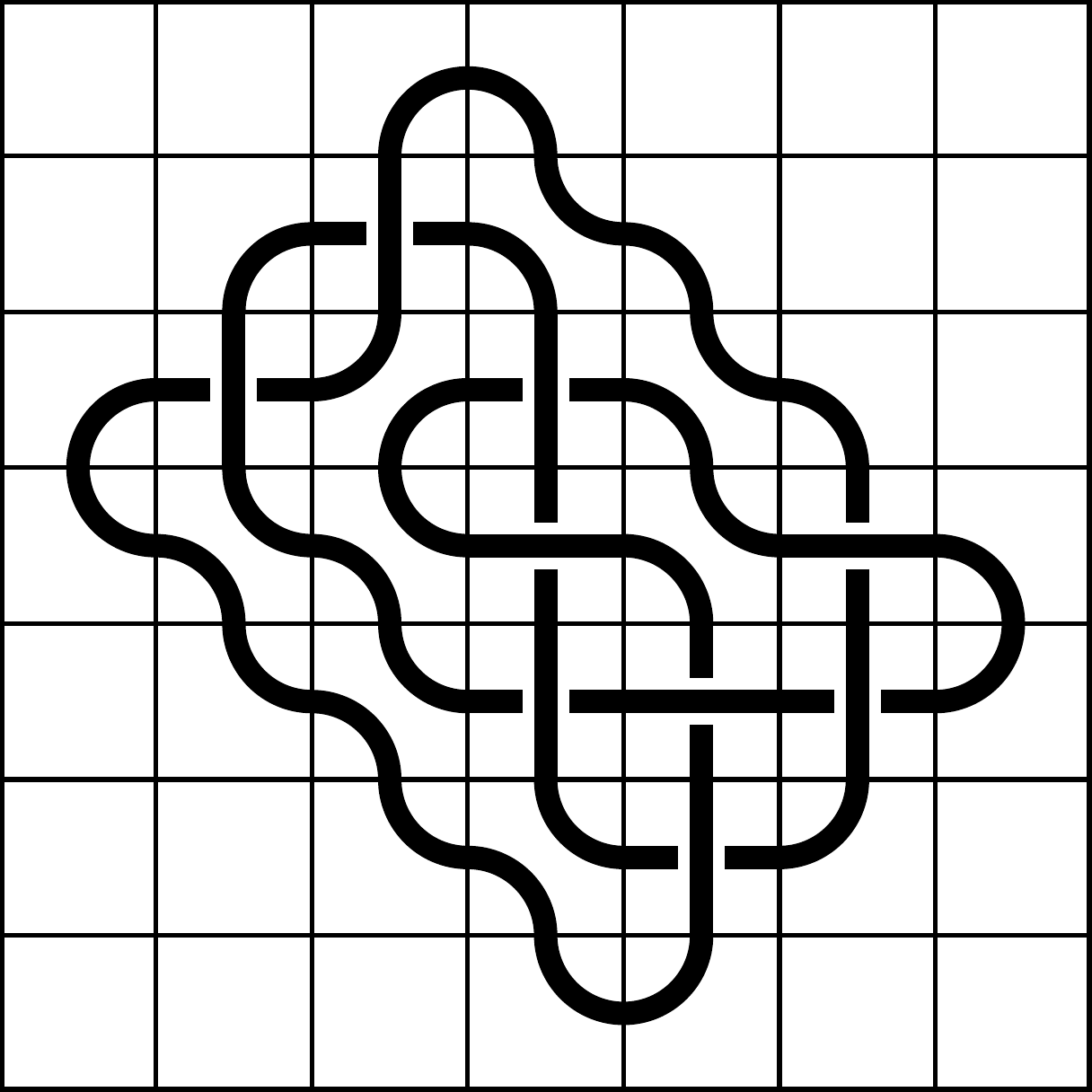}
        \caption*{ $9_{41}$ }
    \end{minipage} \newline
\end{figure}
\unskip

\begin{figure}[H]
    \centering
    \begin{minipage}{0.155\linewidth}
        \captionsetup{skip=3pt}
        \centering
        \includegraphics[width=\linewidth]{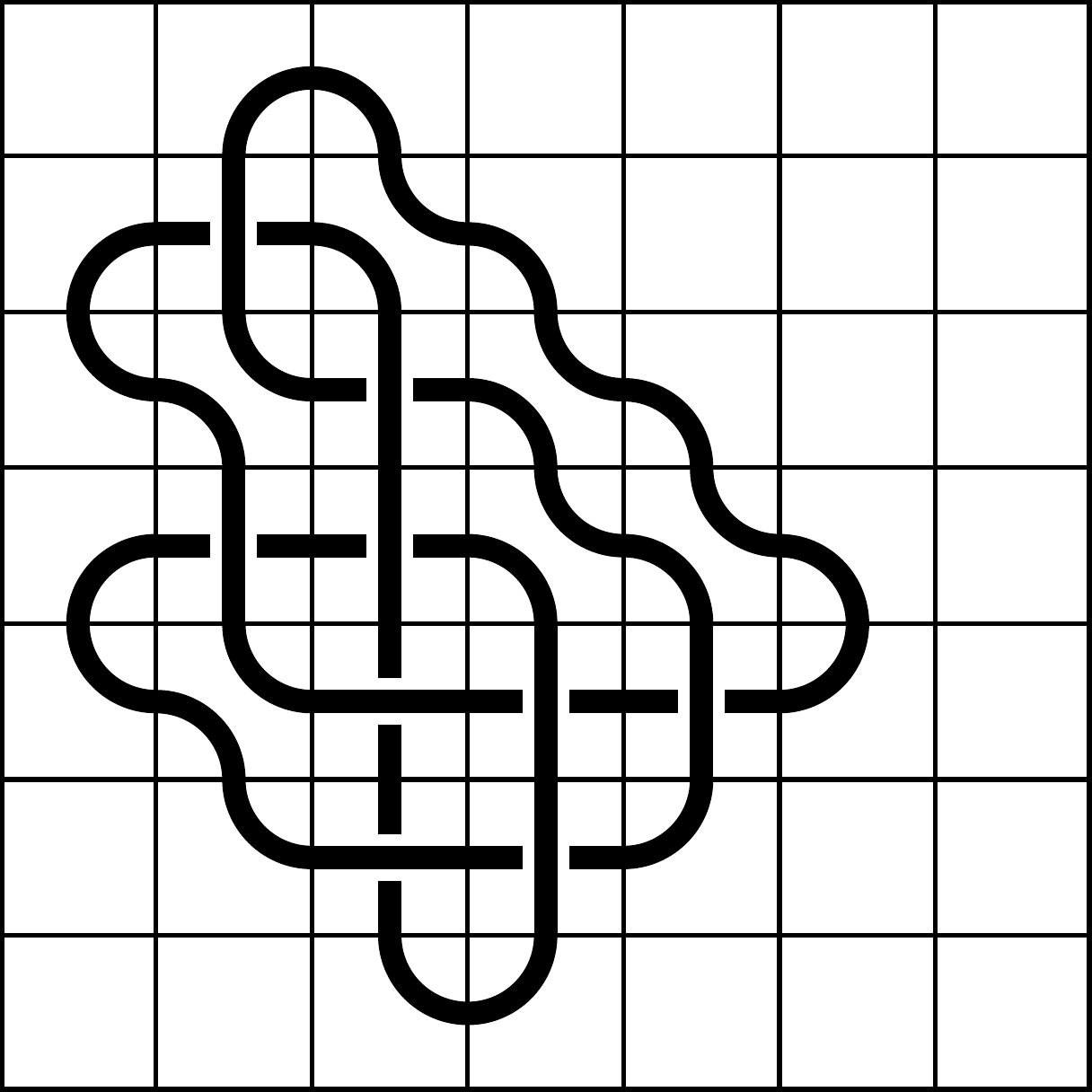}
        \caption*{$9_{42}$ }
    \end{minipage} \hfill
     \begin{minipage}{0.155\linewidth}
        \captionsetup{skip=3pt}
        \centering
        \includegraphics[width=\linewidth]{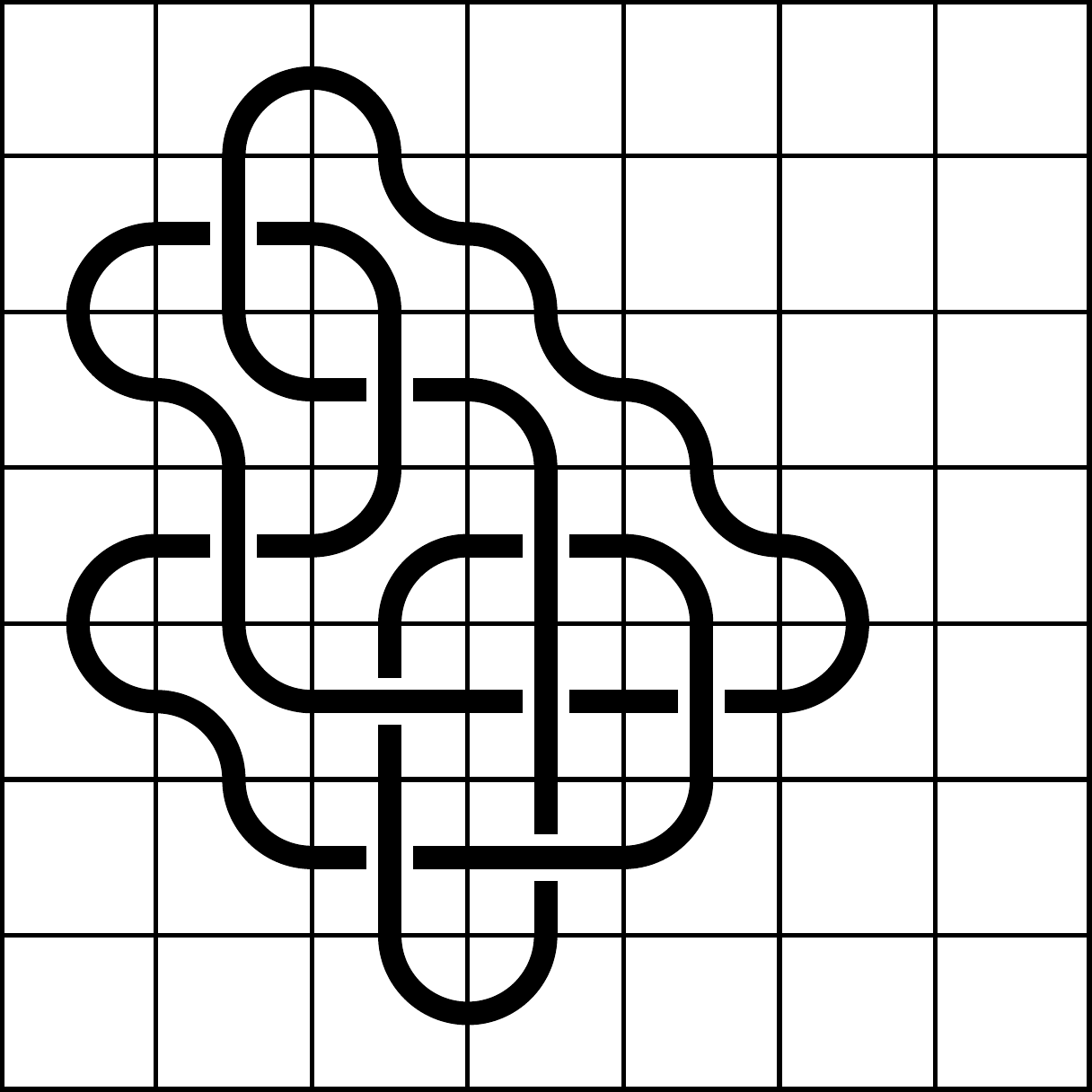}
        \caption*{$9_{43}$ }
    \end{minipage} \hfill
    \begin{minipage}{0.155\linewidth}
        \captionsetup{skip=3pt}
        \centering
        \includegraphics[width=\linewidth]{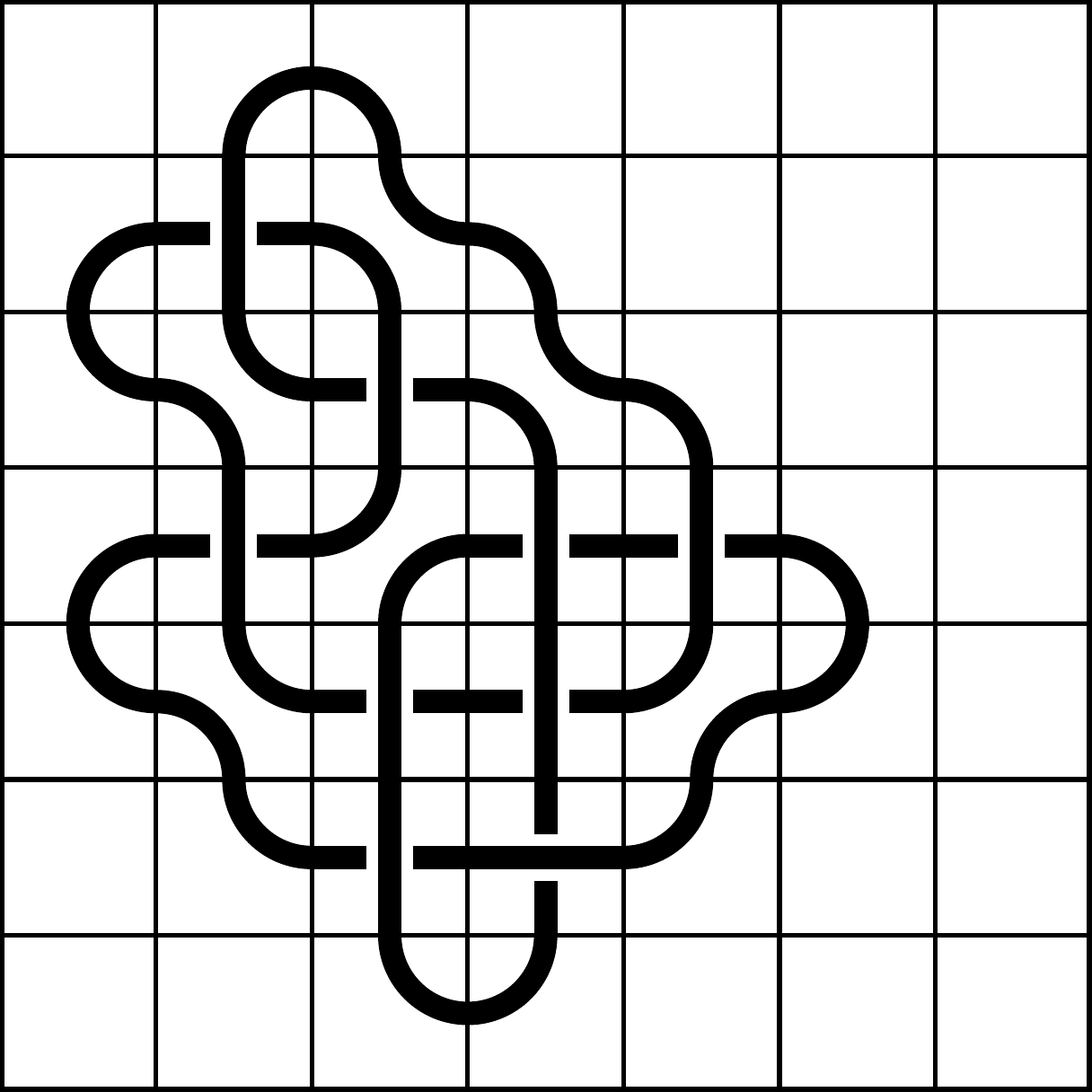}
        \caption*{$9_{44}$ }
    \end{minipage} \hfill
     \begin{minipage}{0.155\linewidth}
        \captionsetup{skip=3pt}
        \centering
        \includegraphics[width=\linewidth]{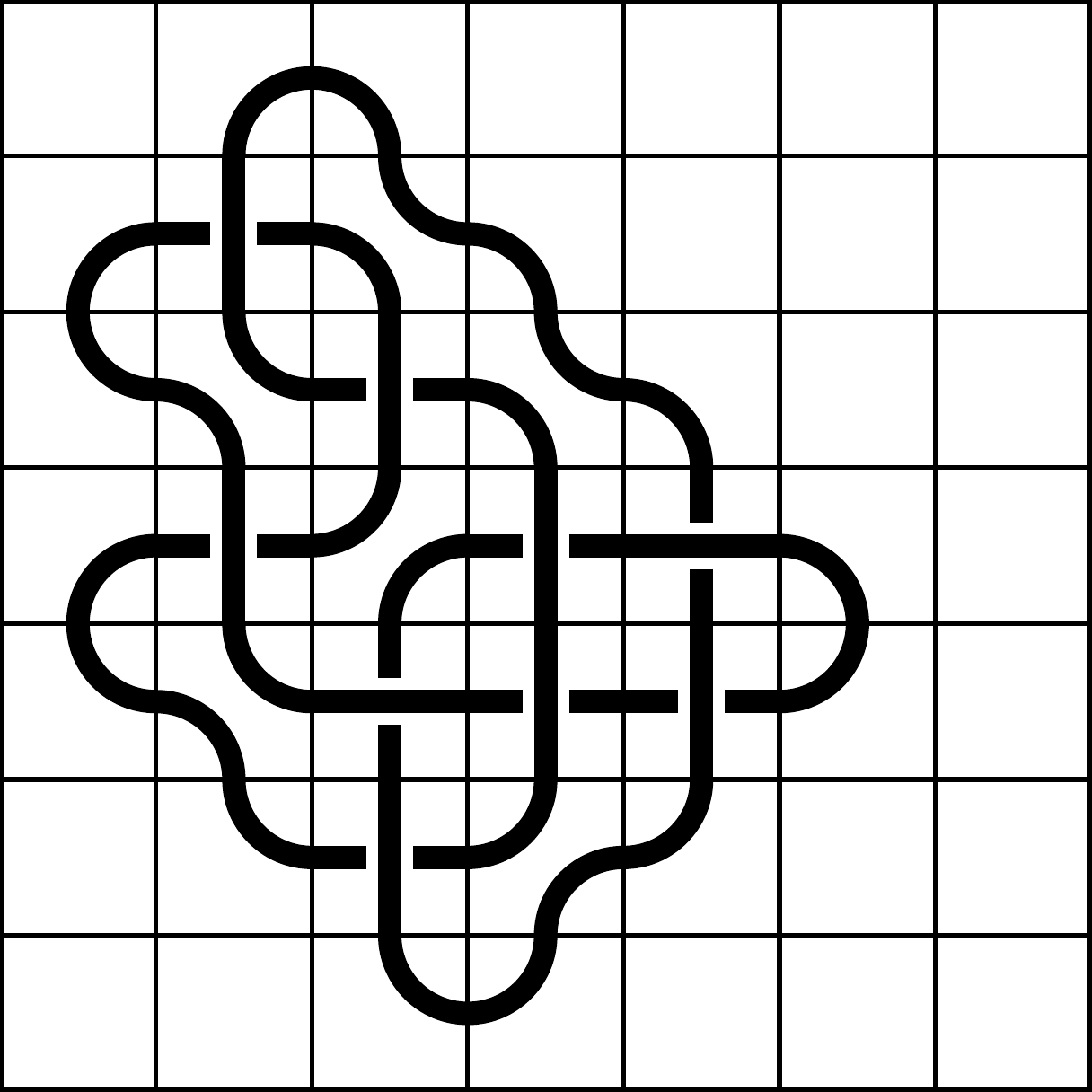}
        \caption*{$9_{45}$ }
    \end{minipage} \hfill
    \begin{minipage}{0.155\linewidth}
        \captionsetup{skip=3pt}
        \centering
        \includegraphics[width=\linewidth]{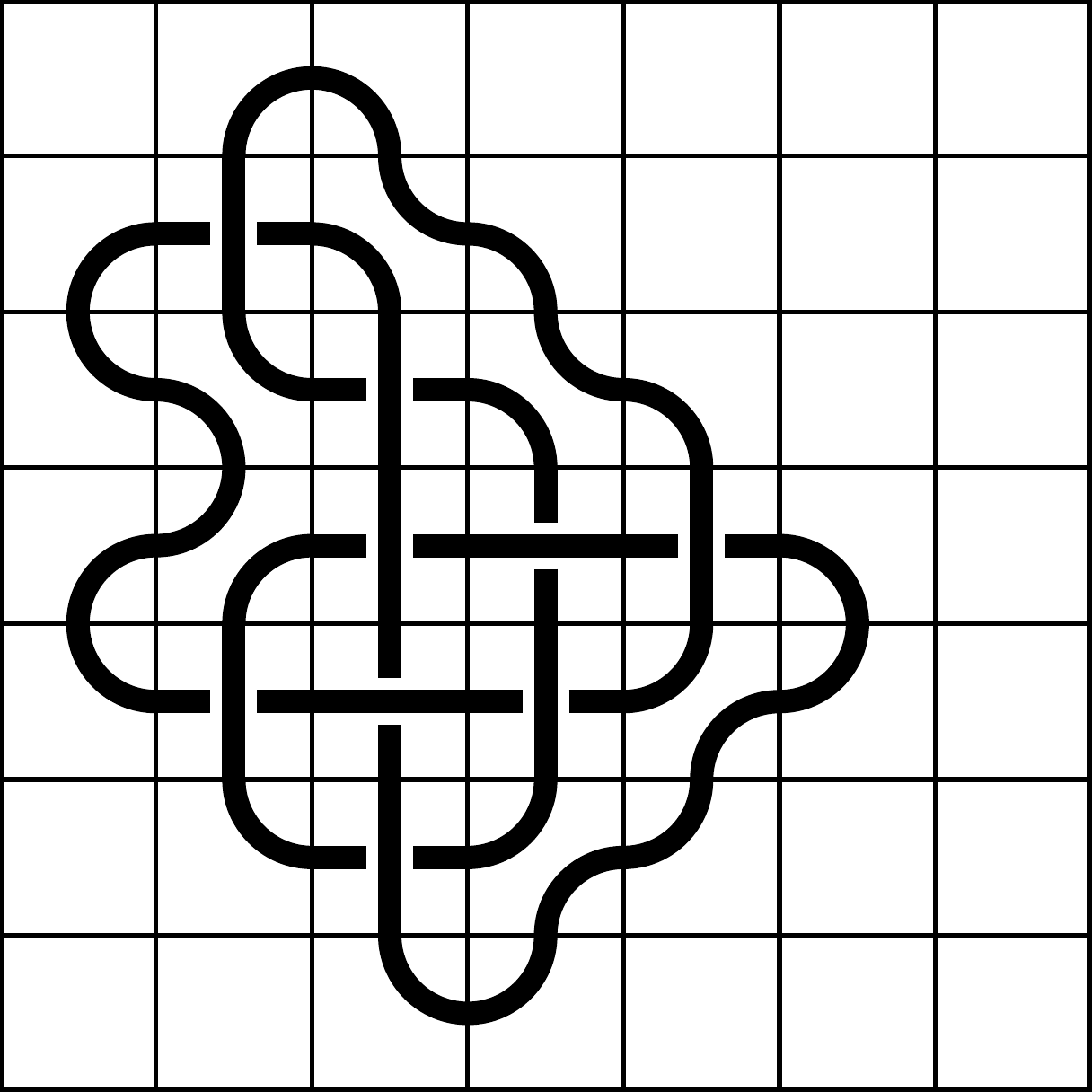}
        \caption*{$9_{47}$ }
    \end{minipage}   \hfill
    \begin{minipage}{0.155\linewidth}
        \captionsetup{skip=3pt}
        \centering
        \includegraphics[width=\linewidth]{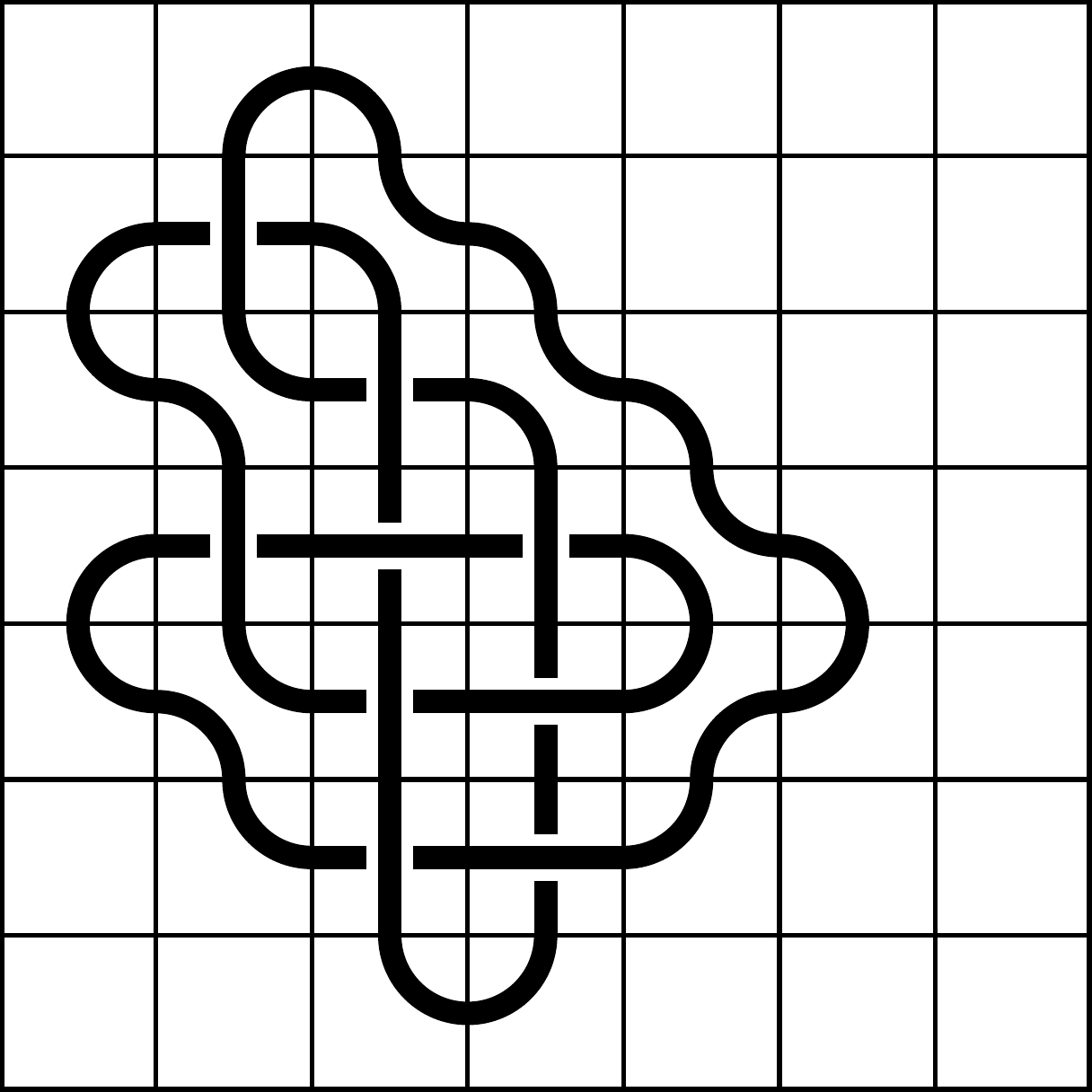}
        \caption*{$9_{49}$}
    \end{minipage} \newline
\end{figure}
\unskip

\begin{figure}[H]
    \centering
    \begin{minipage}{0.155\linewidth}
        \captionsetup{skip=3pt}
        \centering
        \includegraphics[width=\linewidth]{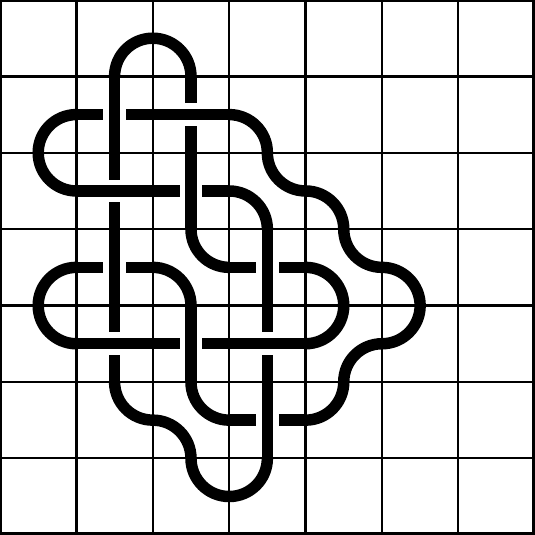}
        \caption*{$10_{23}$ }
    \end{minipage} \hfill
    \begin{minipage}{0.155\linewidth}
        \captionsetup{skip=3pt}
        \centering
        \includegraphics[width=\linewidth]{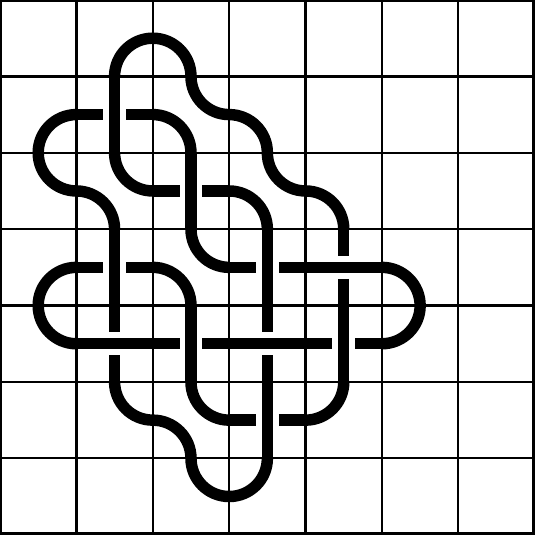}
        \caption*{$10_{27}$}
    \end{minipage} \hfill
    \begin{minipage}{0.155\linewidth}
        \captionsetup{skip=3pt}
        \centering
        \includegraphics[width=\linewidth]{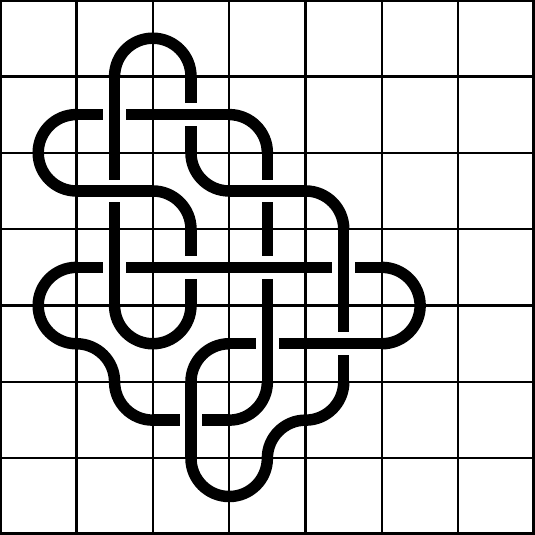}
        \caption*{\phantom{{\Large $\ast$}} $10_{37}$ {\Large $\ast$}}
    \end{minipage} \hfill
    \begin{minipage}{0.155\linewidth}
        \captionsetup{skip=3pt}
        \centering
        \includegraphics[width=\linewidth]{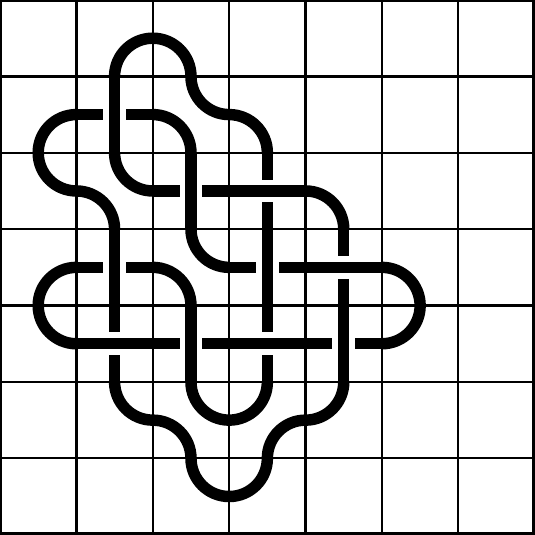}
        \caption*{$10_{40}$}
    \end{minipage} \hfill
        \begin{minipage}{0.155\linewidth}
        \captionsetup{skip=3pt}
        \centering
        \includegraphics[width=\linewidth]{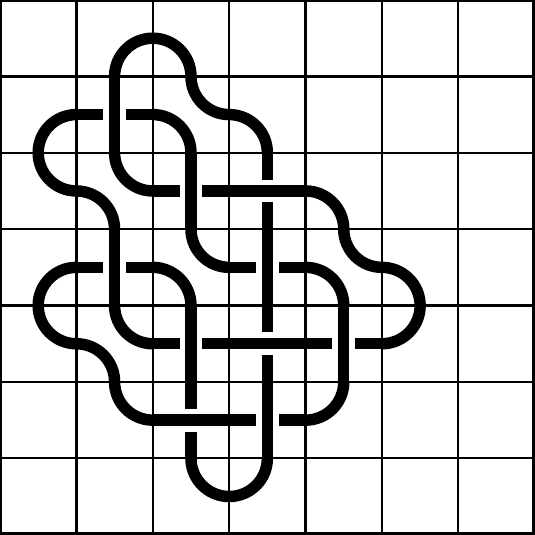}
        \caption*{$10_{42}$}
    \end{minipage} \hfill
    \begin{minipage}{0.155\linewidth}
        \captionsetup{skip=3pt}
        \centering
        \includegraphics[width=\linewidth]{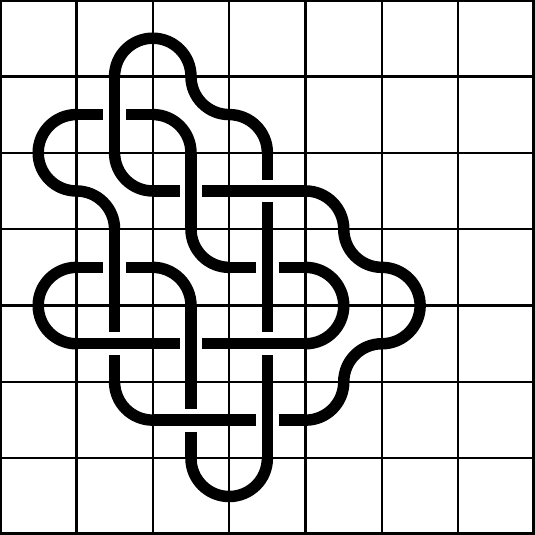}
        \caption*{$10_{43}$ }
    \end{minipage}  \newline
\end{figure}
\unskip

\begin{figure}[H]
    \centering
    \begin{minipage}{0.155\linewidth}
        \captionsetup{skip=3pt}
        \centering
        \includegraphics[width=\linewidth]{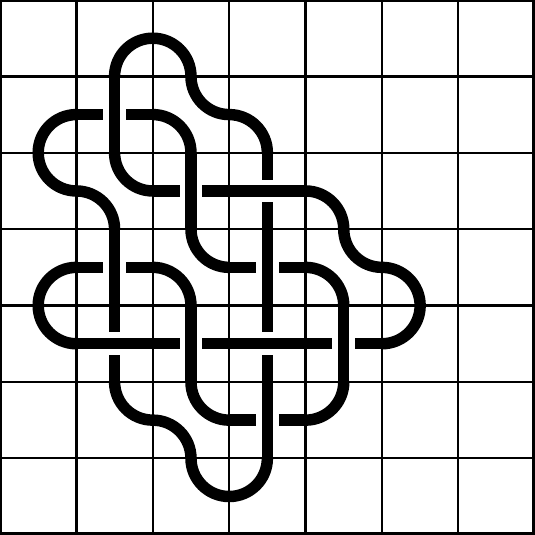}
        \caption*{$10_{45}$}
    \end{minipage} \hfill
     \begin{minipage}{0.155\linewidth}
        \captionsetup{skip=3pt}
        \centering
        \includegraphics[width=\linewidth]{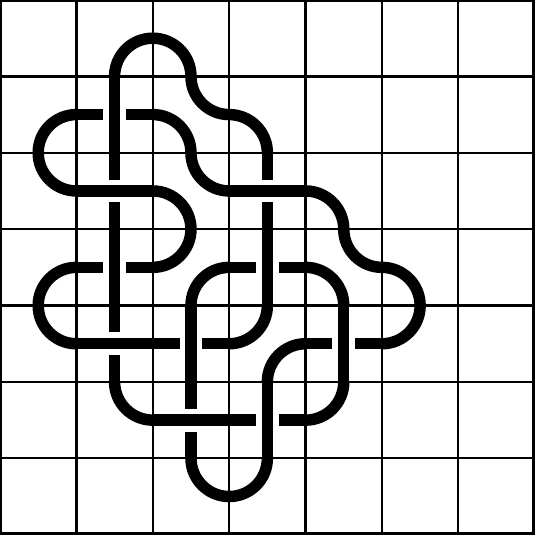}
        \caption*{$10_{46}$ }
    \end{minipage} \hfill
    \begin{minipage}{0.155\linewidth}
        \captionsetup{skip=3pt}
        \centering
        \includegraphics[width=\linewidth]{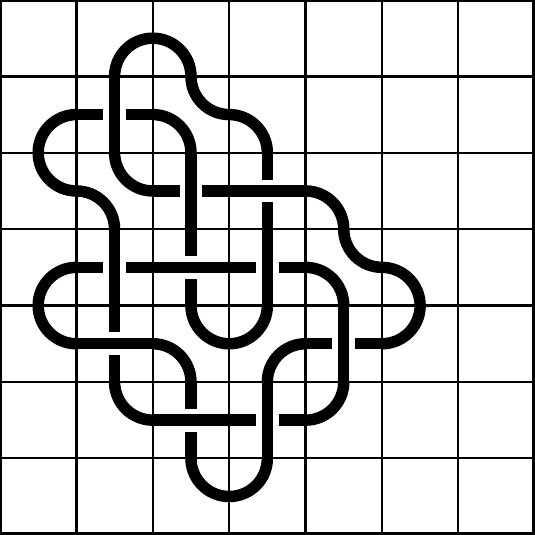}
        \caption*{$10_{47}$ }
    \end{minipage}   \hfill
        \begin{minipage}{0.155\linewidth}
        \captionsetup{skip=3pt}
        \centering
        \includegraphics[width=\linewidth]{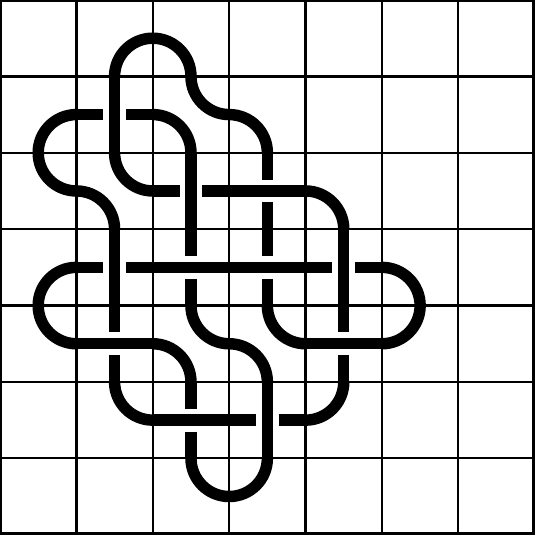}
        \caption*{\phantom{{\Large $\ast$}} $10_{48}$ {\Large $\ast$}}
    \end{minipage} \hfill
     \begin{minipage}{0.155\linewidth}
        \captionsetup{skip=3pt}
        \centering
        \includegraphics[width=\linewidth]{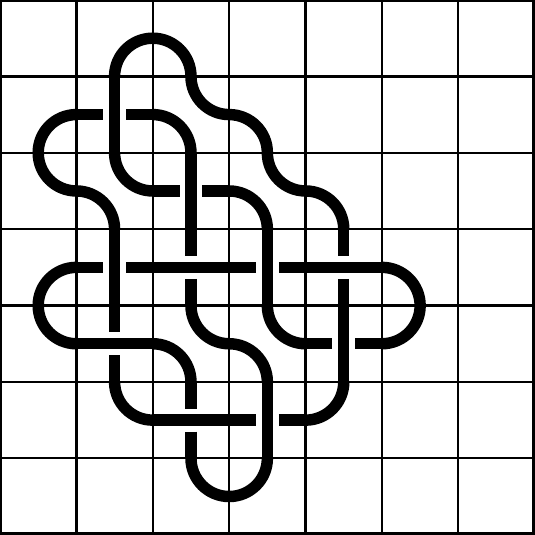}
        \caption*{$10_{49}$ }
    \end{minipage} \hfill
    \begin{minipage}{0.155\linewidth}
        \captionsetup{skip=3pt}
        \centering
        \includegraphics[width=\linewidth]{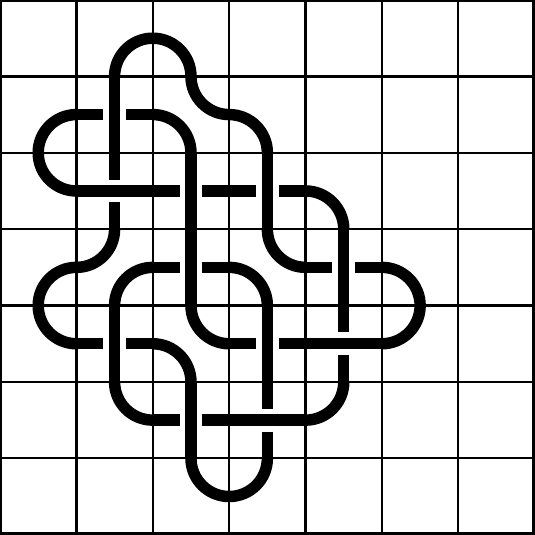}
        \caption*{\phantom{{\Large $\ast$}} $10_{50}$ {\Large $\ast$}}
    \end{minipage} \newline
\end{figure}
\unskip

\begin{figure}[H]
    \centering
    \begin{minipage}{0.155\linewidth}
        \captionsetup{skip=3pt}
        \centering
        \includegraphics[width=\linewidth]{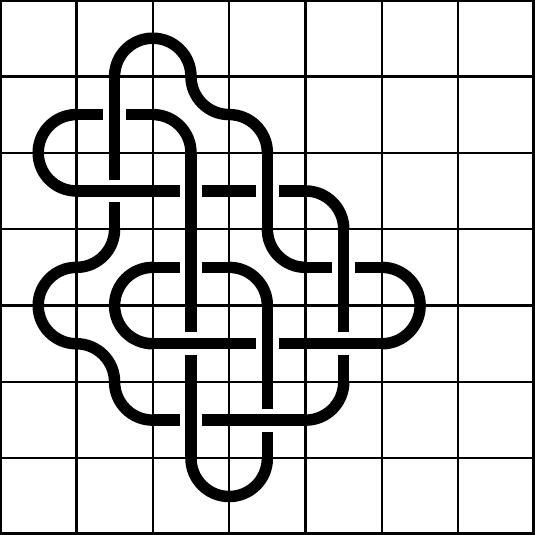}
        \caption*{\phantom{{\Large $\ast$}} $10_{51}$ {\Large $\ast$}}
    \end{minipage} \hfill
    \begin{minipage}{0.155\linewidth}
        \captionsetup{skip=3pt}
        \centering
        \includegraphics[width=\linewidth]{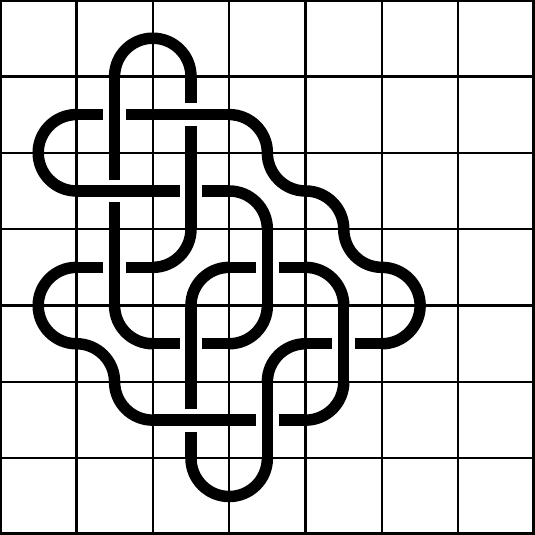}
        \caption*{$10_{52}$ }
    \end{minipage}    \hfill
    \begin{minipage}{0.155\linewidth}
        \captionsetup{skip=3pt}
        \centering
        \includegraphics[width=\linewidth]{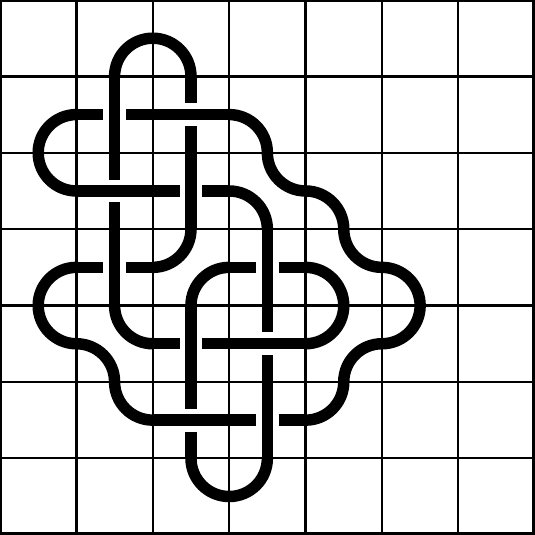}
        \caption*{$10_{53}$ }
    \end{minipage} \hfill
    \begin{minipage}{0.155\linewidth}
        \captionsetup{skip=3pt}
        \centering
        \includegraphics[width=\linewidth]{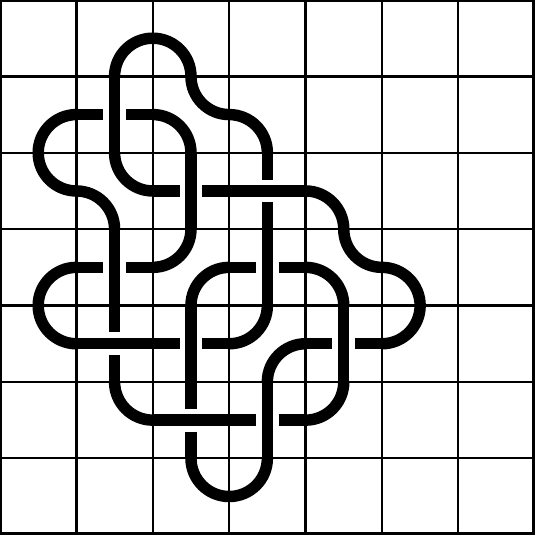}
        \caption*{$10_{54}$ }
    \end{minipage} \hfill
    \begin{minipage}{0.155\linewidth}
        \captionsetup{skip=3pt}
        \centering
        \includegraphics[width=\linewidth]{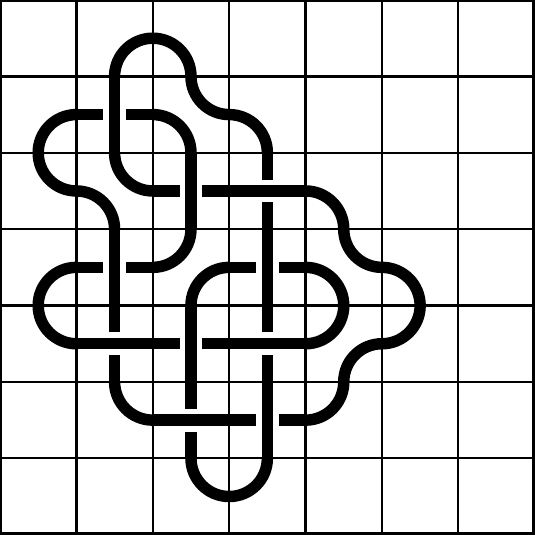}
        \caption*{$10_{55}$ }
    \end{minipage} \hfill
    \begin{minipage}{0.155\linewidth}
        \captionsetup{skip=3pt}
        \centering
        \includegraphics[width=\linewidth]{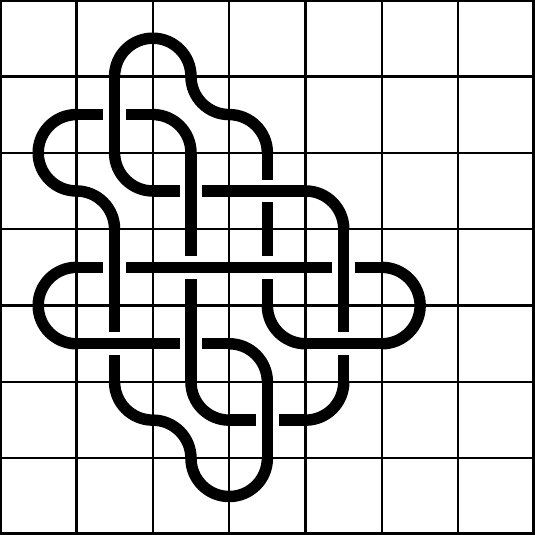}
        \caption*{\phantom{{\Large $\ast$}} $10_{56}$ {\Large $\ast$}}
    \end{minipage} \newline
\end{figure}
\unskip

\begin{figure}[H]
    \centering
    \begin{minipage}{0.155\linewidth}
        \captionsetup{skip=3pt}
        \centering
        \includegraphics[width=\linewidth]{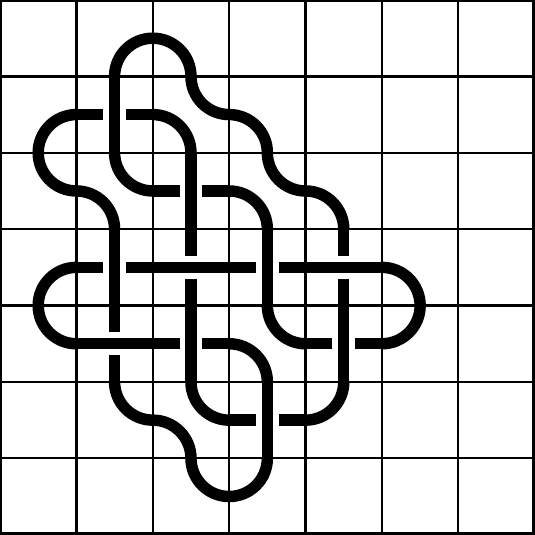}
        \caption*{$10_{57}$ }
    \end{minipage}    \hfill
    \begin{minipage}{0.155\linewidth}
        \captionsetup{skip=3pt}
        \centering
        \includegraphics[width=\linewidth]{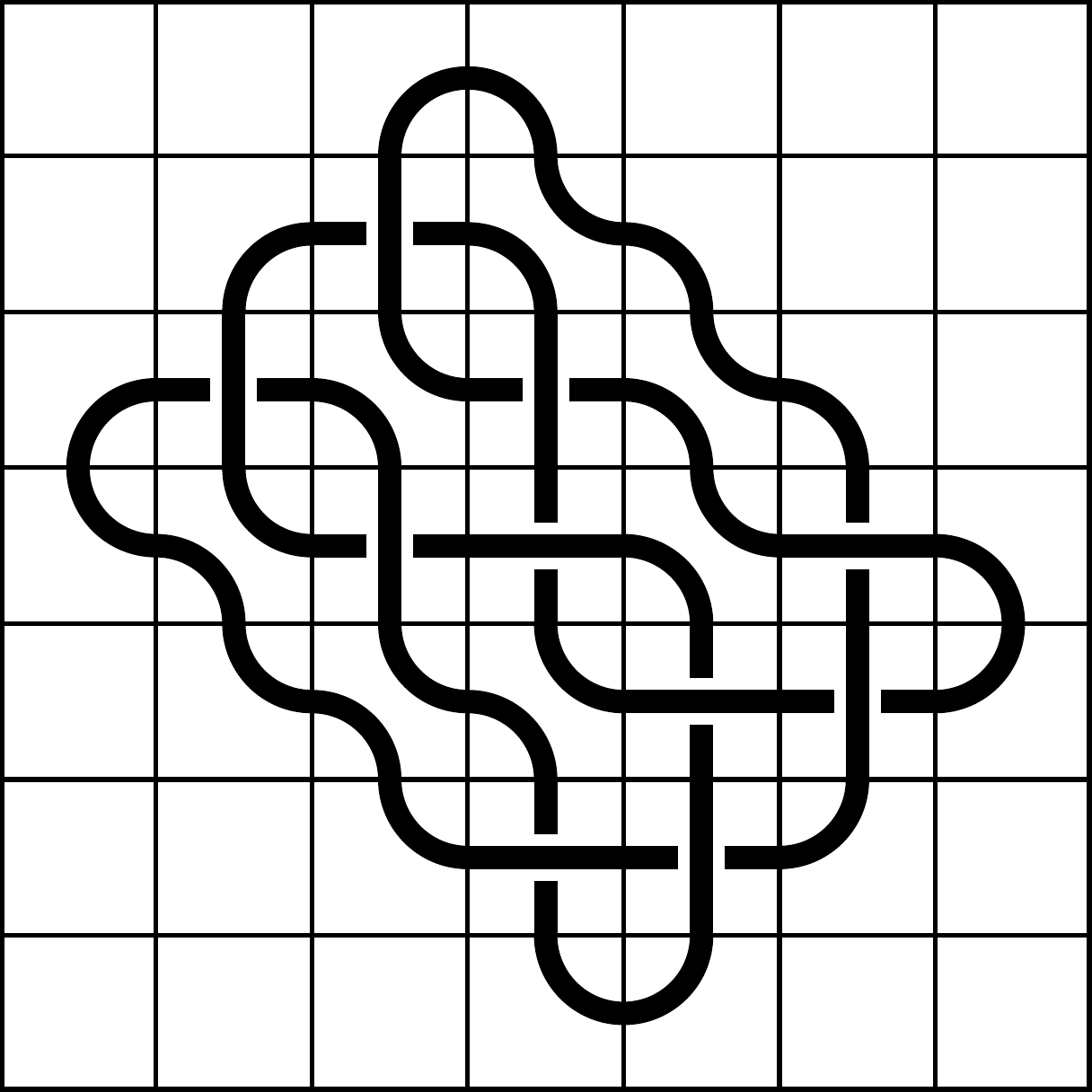}
        \caption*{$10_{58}$}
    \end{minipage} \hfill
    \begin{minipage}{0.155\linewidth}
        \captionsetup{skip=3pt}
        \centering
        \includegraphics[width=\linewidth]{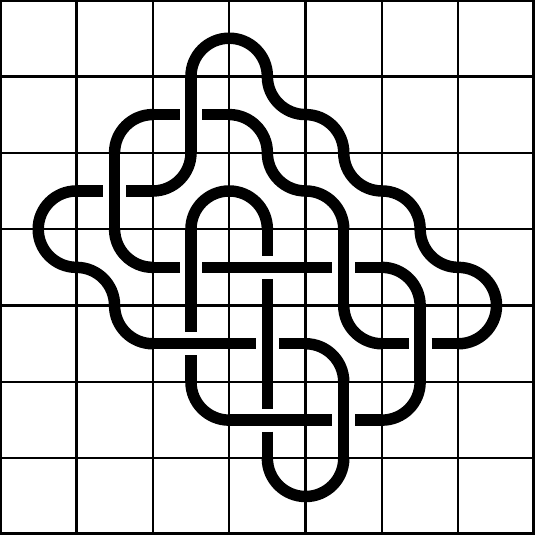}
        \caption*{ $10_{59}$ }
    \end{minipage} \hfill
    \begin{minipage}{0.155\linewidth}
        \captionsetup{skip=3pt}
        \centering
        \includegraphics[width=\linewidth]{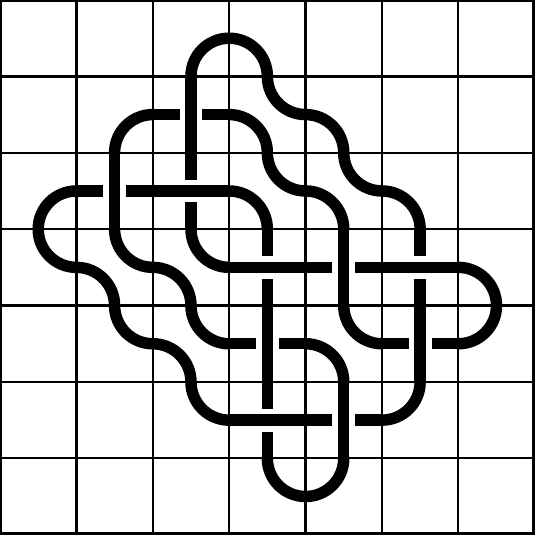}
        \caption*{ $10_{60}$ }
    \end{minipage} \hfill
    \begin{minipage}{0.155\linewidth}
        \captionsetup{skip=3pt}
        \centering
        \includegraphics[width=\linewidth]{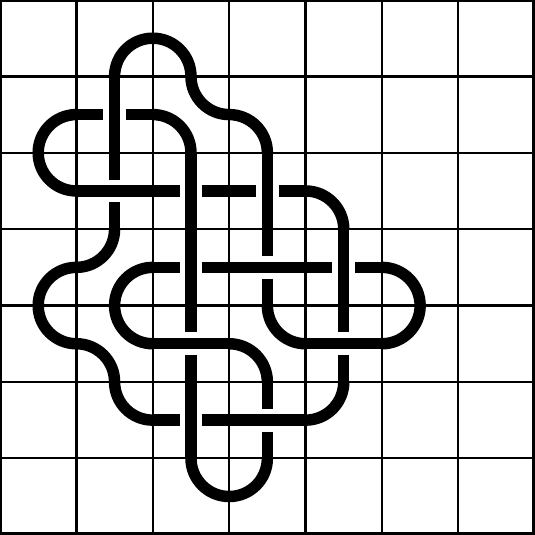}
        \caption*{\phantom{{\Large $\ast$}} $10_{67}$ {\Large $\ast$} }
    \end{minipage} \hfill
    \begin{minipage}{0.155\linewidth}
        \captionsetup{skip=3pt}
        \centering
        \includegraphics[width=\linewidth]{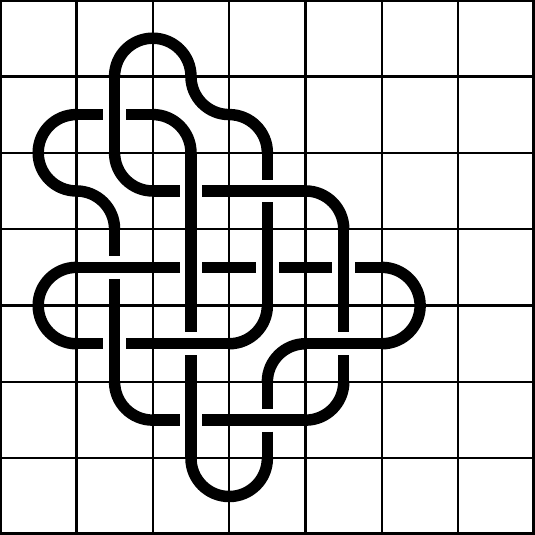}
        \caption*{\phantom{{\Large $\ast$}} $10_{68}$ {\Large $\ast$} }
    \end{minipage} \newline
\end{figure}
\unskip

\begin{figure}[H]
    \centering
    \begin{minipage}{0.155\linewidth}
        \captionsetup{skip=3pt}
        \centering
        \includegraphics[width=\linewidth]{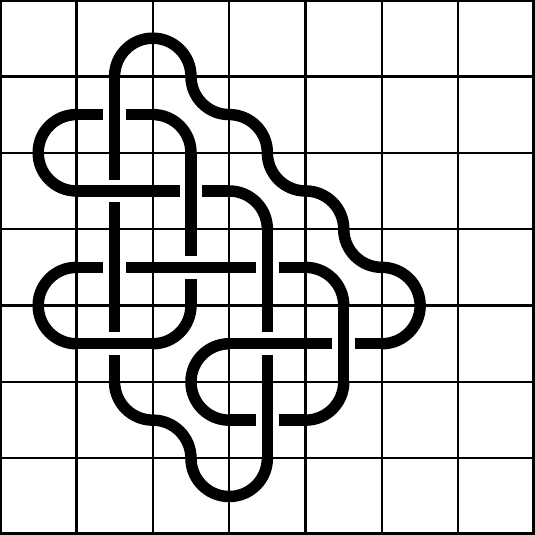}
        \caption*{$10_{69}$ }
    \end{minipage} \hfill
    \begin{minipage}{0.155\linewidth}
        \captionsetup{skip=3pt}
        \centering
        \includegraphics[width=\linewidth]{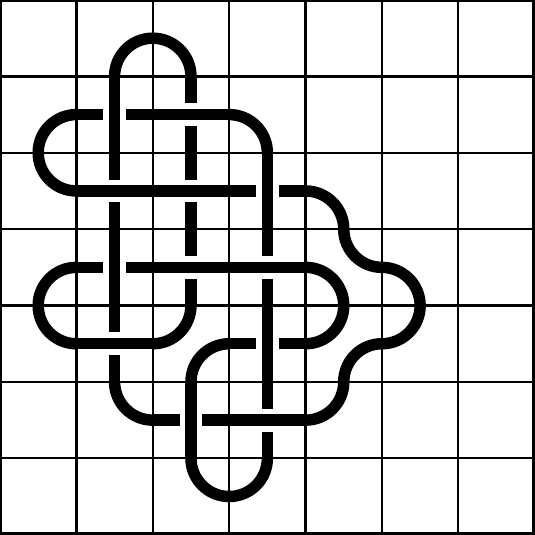}
        \caption*{\phantom{{\Large $\ast$}} $10_{70}$ {\Large $\ast$} }
    \end{minipage} \hfill
    \begin{minipage}{0.155\linewidth}
        \captionsetup{skip=3pt}
        \centering
        \includegraphics[width=\linewidth]{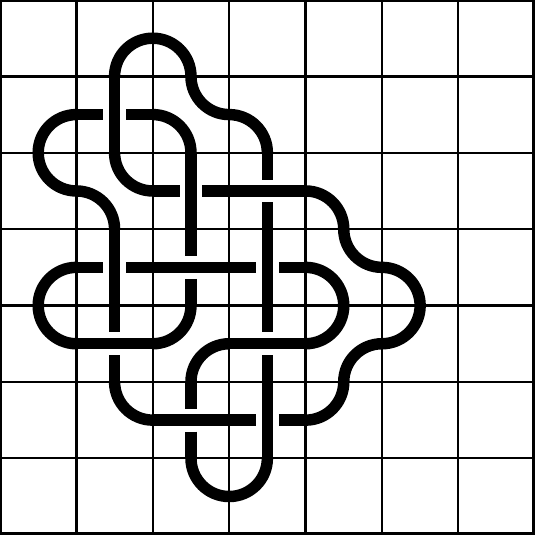}
        \caption*{$10_{71}$ }
    \end{minipage}  \hfill
    \begin{minipage}{0.155\linewidth}
        \captionsetup{skip=3pt}
        \centering
        \includegraphics[width=\linewidth]{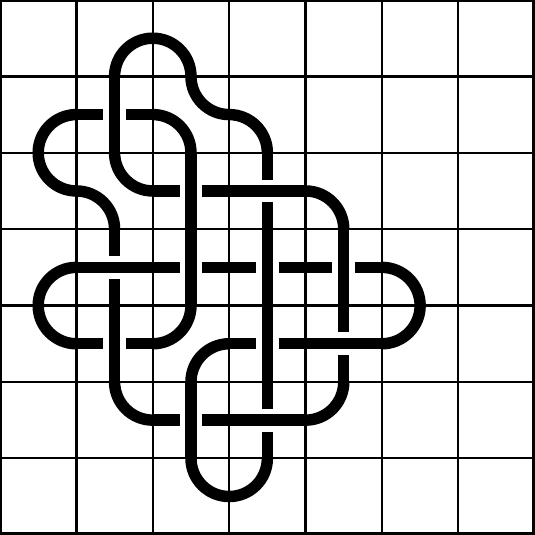}
        \caption*{\phantom{{\Large $\ast$}} $10_{72}$ {\Large $\ast$} }
    \end{minipage} \hfill
    \begin{minipage}{0.155\linewidth}
        \captionsetup{skip=3pt}
        \centering
        \includegraphics[width=\linewidth]{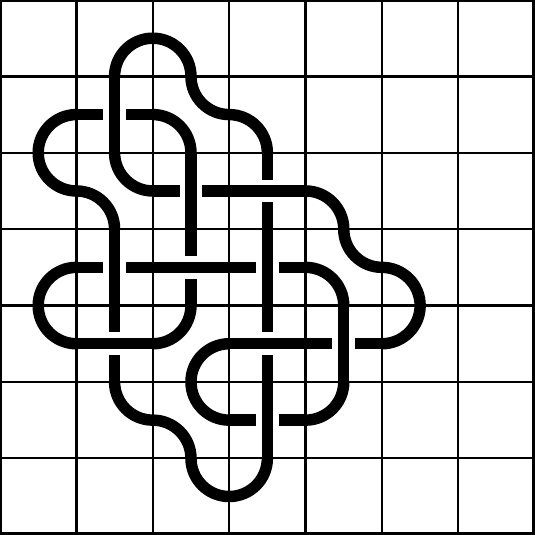}
        \caption*{$10_{73}$ }
    \end{minipage} \hfill
    \begin{minipage}{0.155\linewidth}
        \captionsetup{skip=3pt}
        \centering
        \includegraphics[width=\linewidth]{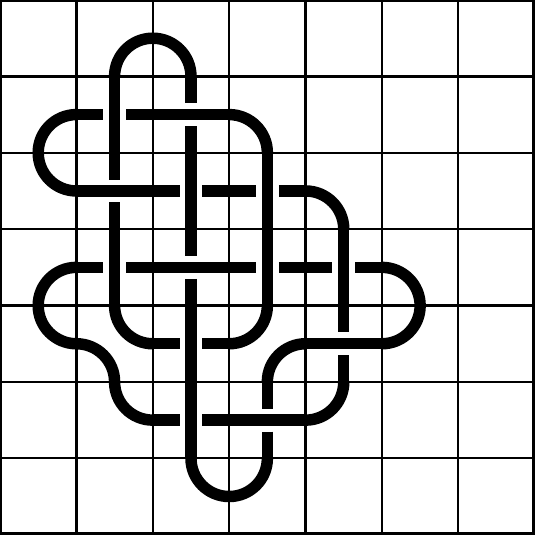}
        \caption*{\phantom{{\Large $\ast$}} $10_{79}$ {\Large $\ast$} }
    \end{minipage} \newline
\end{figure}
\unskip

\begin{figure}[H]
    \centering
    \begin{minipage}{0.155\linewidth}
        \captionsetup{skip=3pt}
        \centering
        \includegraphics[width=\linewidth]{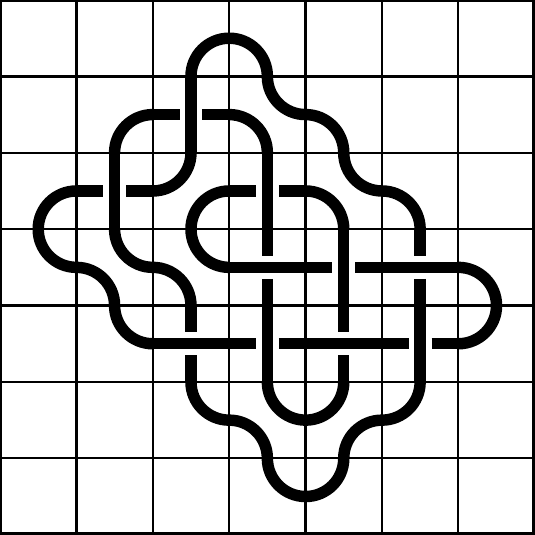}
        \caption*{$10_{80}$}
    \end{minipage} \hfill
    \begin{minipage}{0.155\linewidth}
        \captionsetup{skip=3pt}
        \centering
        \includegraphics[width=\linewidth]{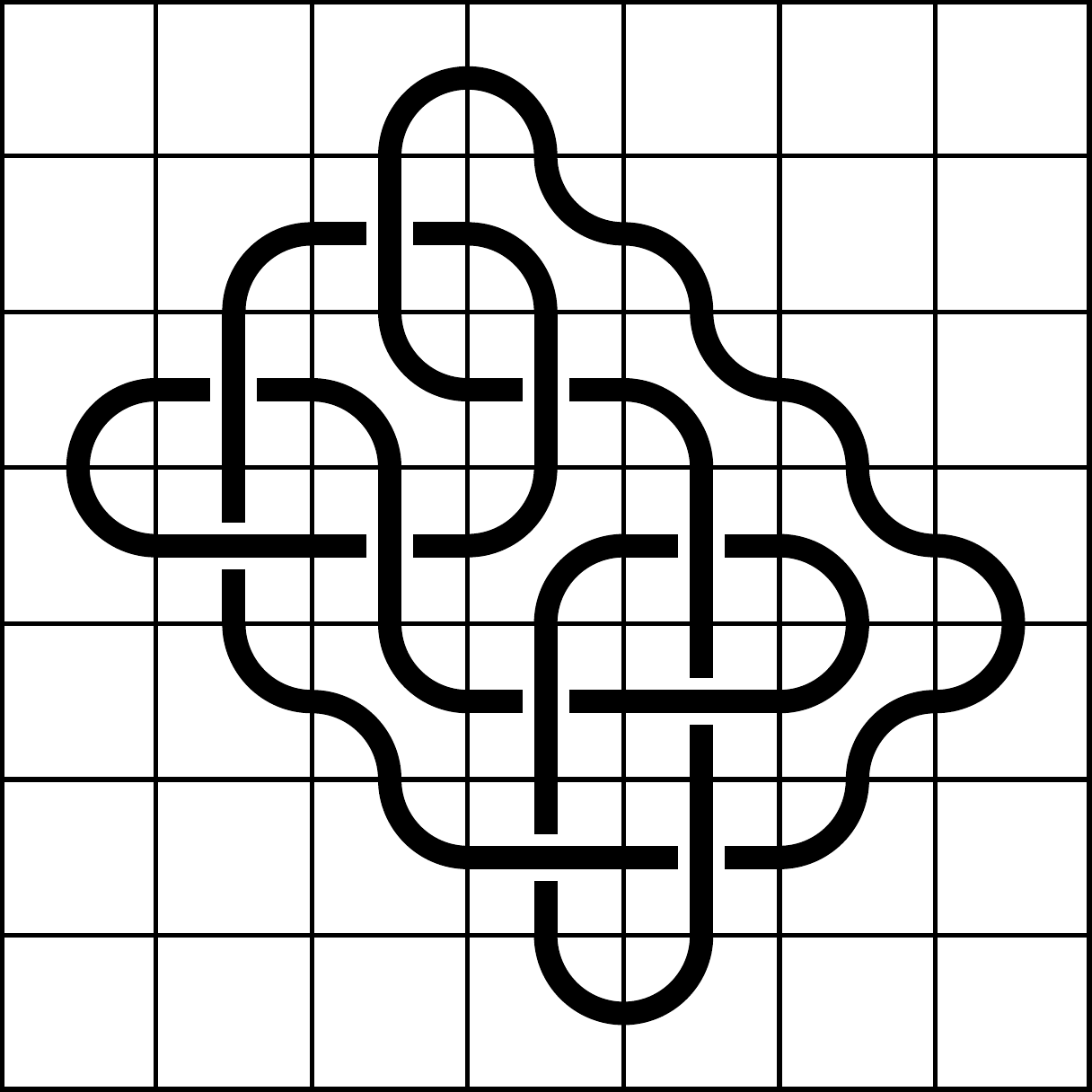}
        \caption*{ $10_{81}$ }
    \end{minipage} \hfill
    \begin{minipage}{0.155\linewidth}
        \captionsetup{skip=3pt}
        \centering
        \includegraphics[width=\linewidth]{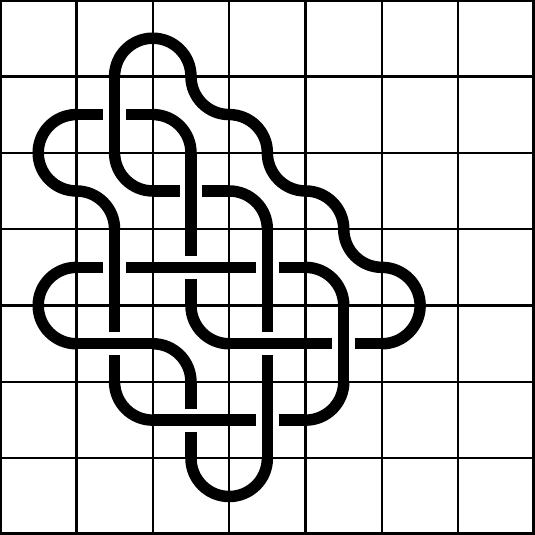}
        \caption*{$10_{82}$ }
    \end{minipage} \hfill
    \begin{minipage}{0.155\linewidth}
        \captionsetup{skip=3pt}
        \centering
        \includegraphics[width=\linewidth]{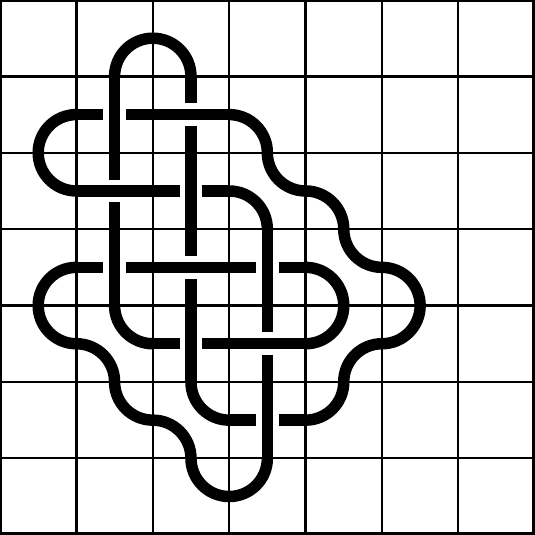}
        \caption*{$10_{83}$ }
    \end{minipage}  \hfill
    \begin{minipage}{0.155\linewidth}
        \captionsetup{skip=3pt}
        \centering
        \includegraphics[width=\linewidth]{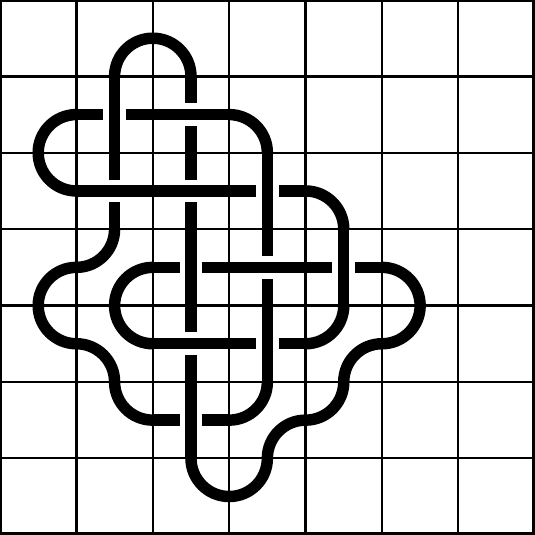}
        \caption*{\phantom{{\Large $\ast$}} $10_{84}$ {\Large $\ast$} }
    \end{minipage} \hfill
    \begin{minipage}{0.155\linewidth}
        \captionsetup{skip=3pt}
        \centering
        \includegraphics[width=\linewidth]{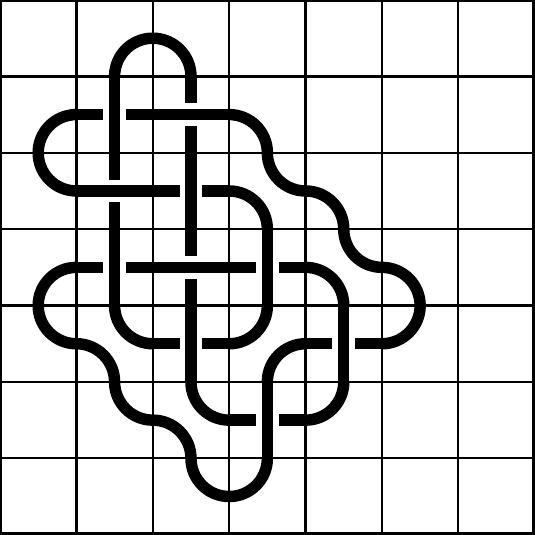}
        \caption*{$10_{86}$ }
    \end{minipage} \newline
\end{figure}
\unskip

\begin{figure}[H]
    \centering
    \begin{minipage}{0.155\linewidth}
        \captionsetup{skip=3pt}
        \centering
        \includegraphics[width=\linewidth]{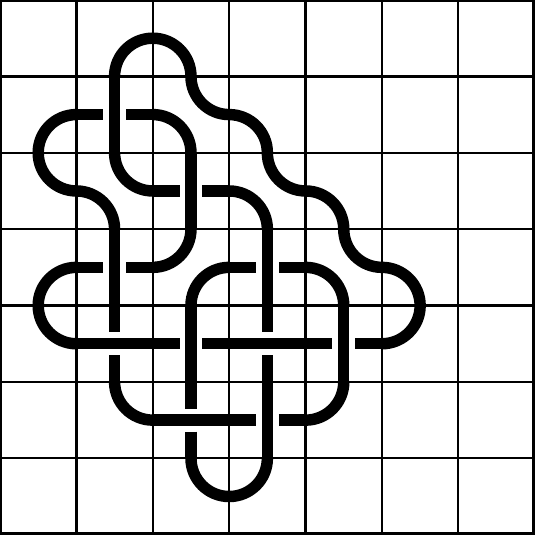}
        \caption*{$10_{87}$ }
    \end{minipage} \hfill
    \begin{minipage}{0.155\linewidth}
        \captionsetup{skip=3pt}
        \centering
        \includegraphics[width=\linewidth]{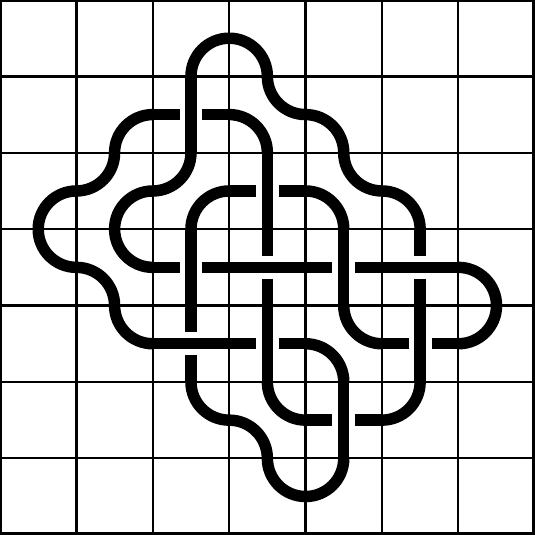}
        \caption*{$10_{88}$}
    \end{minipage} \hfill
    \begin{minipage}{0.155\linewidth}
        \captionsetup{skip=3pt}
        \centering
        \includegraphics[width=\linewidth]{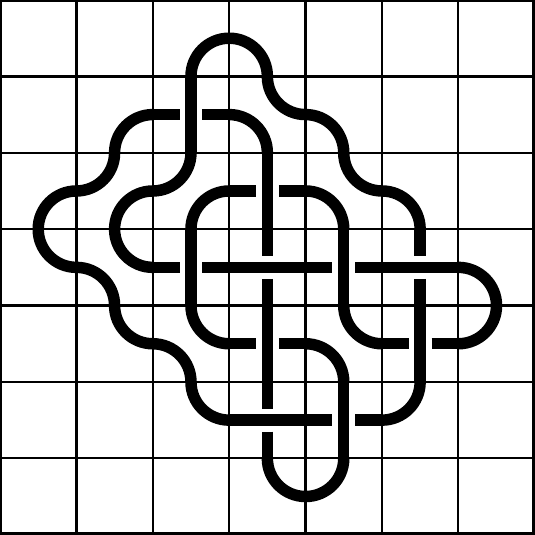}
        \caption*{ $10_{89}$ }
    \end{minipage} \hfill
    \begin{minipage}{0.155\linewidth}
        \captionsetup{skip=3pt}
        \centering
        \includegraphics[width=\linewidth]{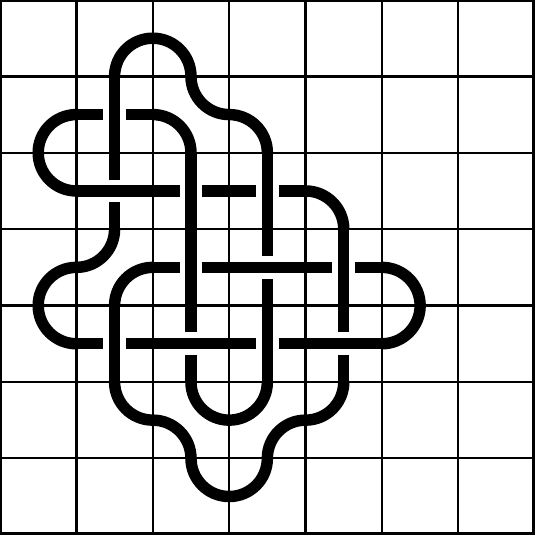}
        \caption*{\phantom{{\Large $\ast$}} $10_{90}$ {\Large $\ast$} }
    \end{minipage} \hfill
    \begin{minipage}{0.155\linewidth}
        \captionsetup{skip=3pt}
        \centering
        \includegraphics[width=\linewidth]{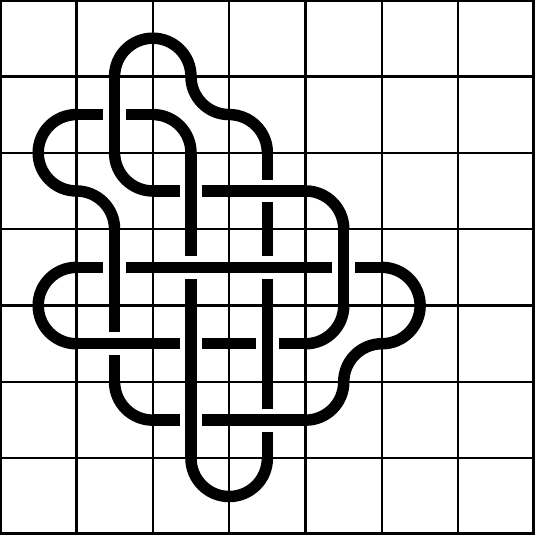}
        \caption*{\phantom{{\Large $\ast$}} $10_{91}$ {\Large $\ast$} }
    \end{minipage}   \hfill
    \begin{minipage}{0.155\linewidth}
        \captionsetup{skip=3pt}
        \centering
        \includegraphics[width=\linewidth]{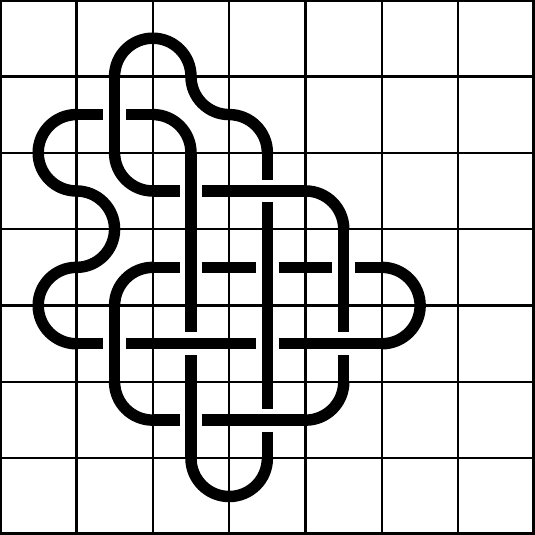}
        \caption*{\phantom{{\Large $\ast$}} $10_{92}$ {\Large $\ast$} }
    \end{minipage}  \newline
\end{figure}
\unskip

\begin{figure}[H]
    \centering
    \begin{minipage}{0.155\linewidth}
        \captionsetup{skip=3pt}
        \centering
        \includegraphics[width=\linewidth]{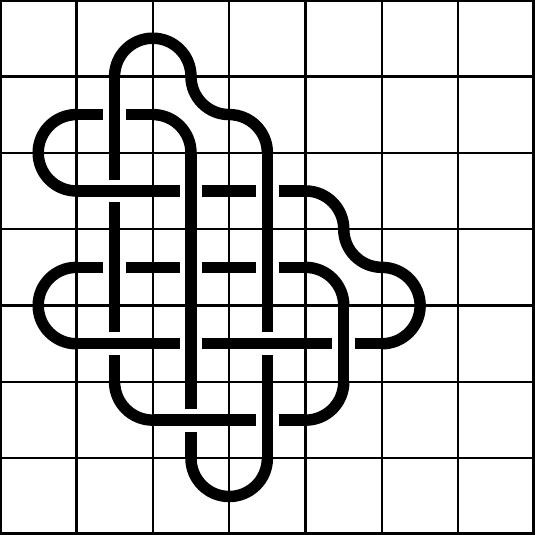}
        \caption*{\phantom{{\Large $\ast$}} $10_{93}$ {\Large $\ast$} }
    \end{minipage} \hfill
    \begin{minipage}{0.155\linewidth}
        \captionsetup{skip=3pt}
        \centering
        \includegraphics[width=\linewidth]{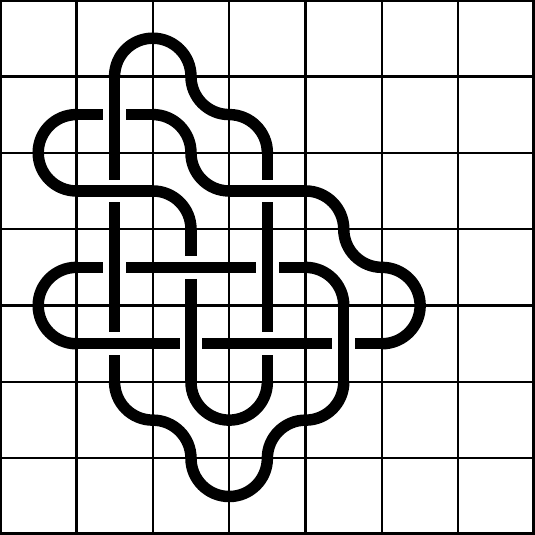}
        \caption*{$10_{94}$ }
    \end{minipage} \hfill
    \begin{minipage}{0.155\linewidth}
        \captionsetup{skip=3pt}
        \centering
        \includegraphics[width=\linewidth]{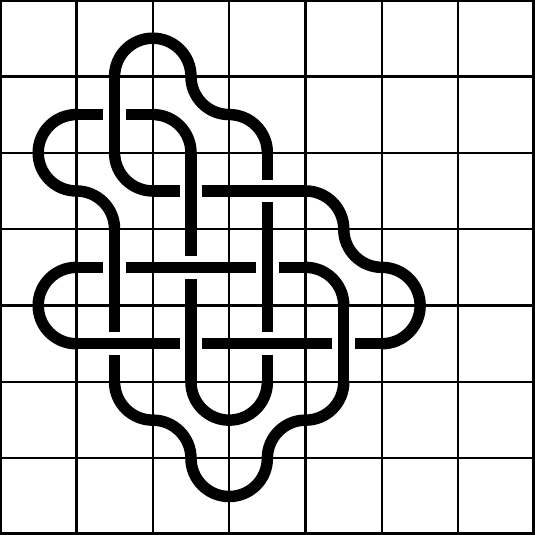}
        \caption*{$10_{95}$ }
    \end{minipage} \hfill
    \begin{minipage}{0.155\linewidth}
        \captionsetup{skip=3pt}
        \centering
        \includegraphics[width=\linewidth]{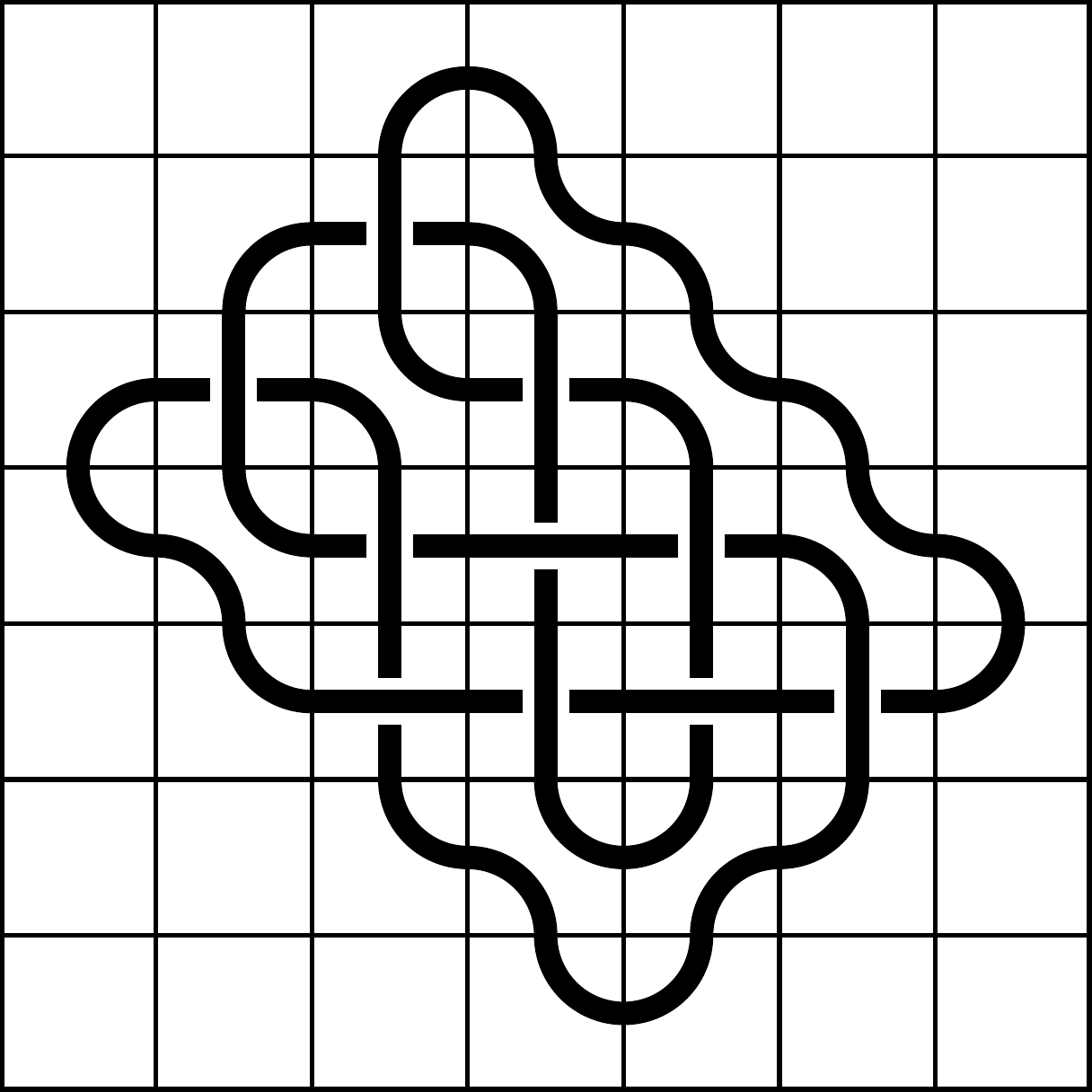}
        \caption*{ $10_{96}$ }
    \end{minipage} \hfill
    \begin{minipage}{0.155\linewidth}
        \captionsetup{skip=3pt}
        \centering
        \includegraphics[width=\linewidth]{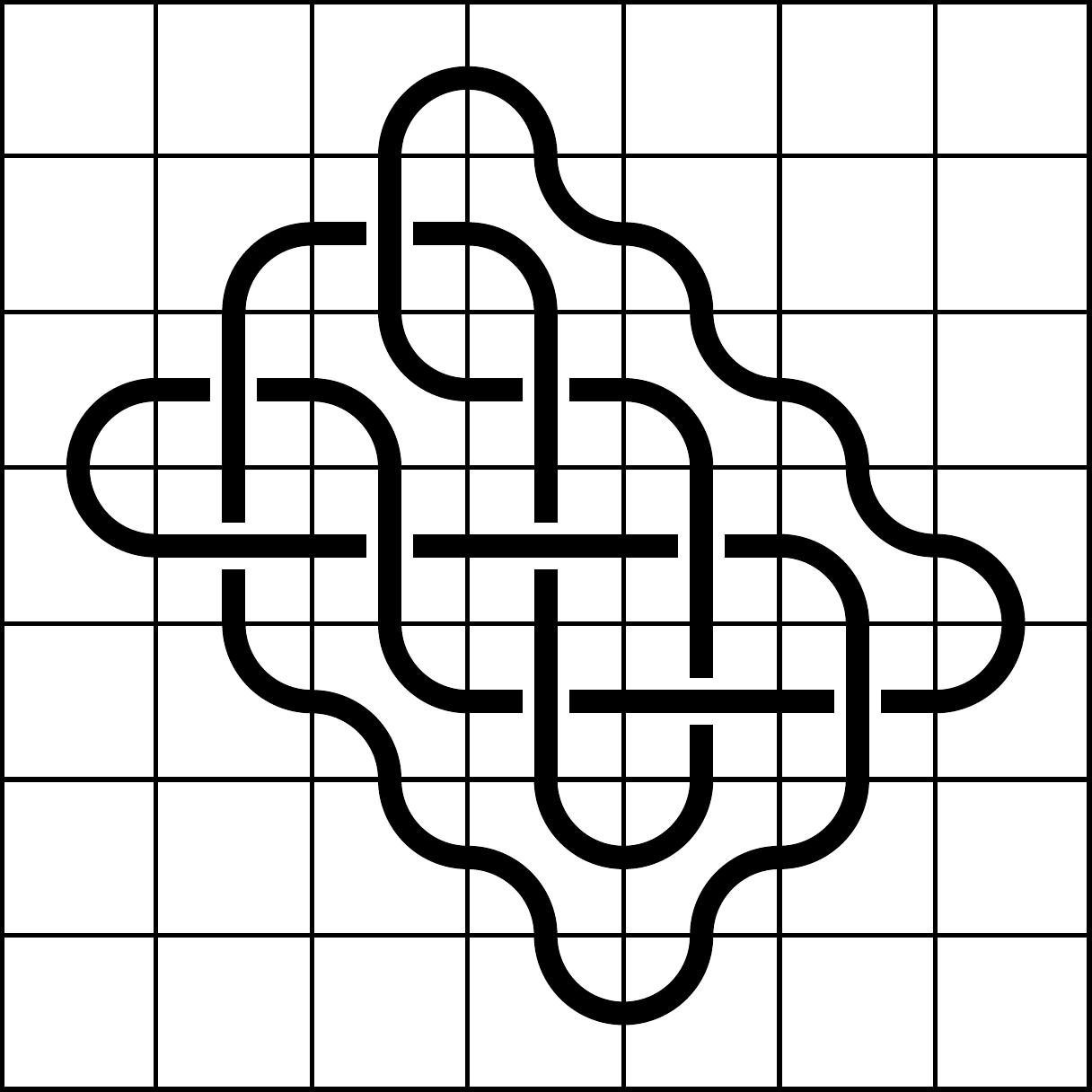}
        \caption*{$10_{97}$}
    \end{minipage} \hfill
    \begin{minipage}{0.155\linewidth}
        \captionsetup{skip=3pt}
        \centering
        \includegraphics[width=\linewidth]{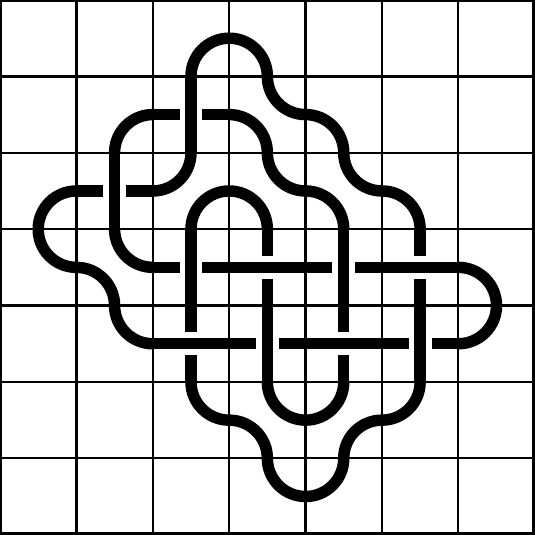}
        \caption*{ $10_{98}$ }
    \end{minipage}  \newline
\end{figure}
\unskip

\begin{figure}[H]
    \centering
    \begin{minipage}{0.155\linewidth}
        \captionsetup{skip=3pt}
        \centering
        \includegraphics[width=\linewidth]{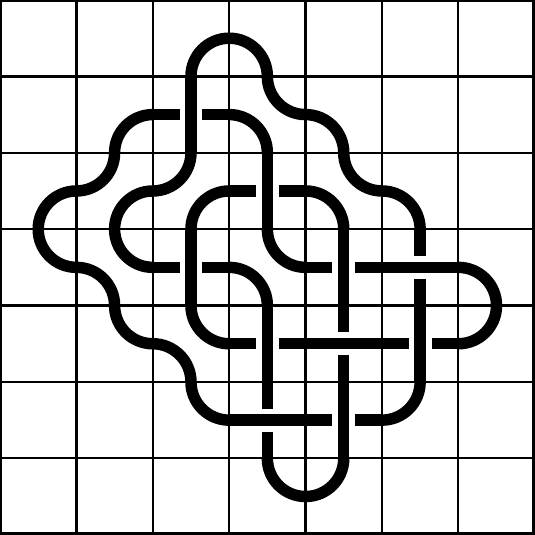}
        \caption*{$10_{99}$}
    \end{minipage} \hfill
    \begin{minipage}{0.155\linewidth}
        \captionsetup{skip=3pt}
        \centering
        \includegraphics[width=\linewidth]{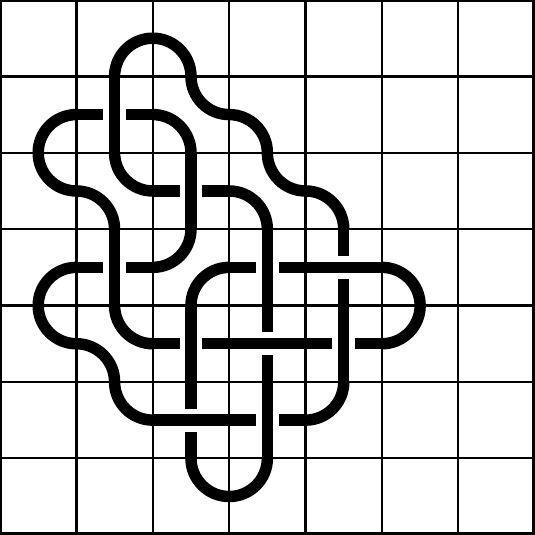}
        \caption*{$10_{101}$ }
    \end{minipage} \hfill
    \begin{minipage}{0.155\linewidth}
        \captionsetup{skip=3pt}
        \centering
        \includegraphics[width=\linewidth]{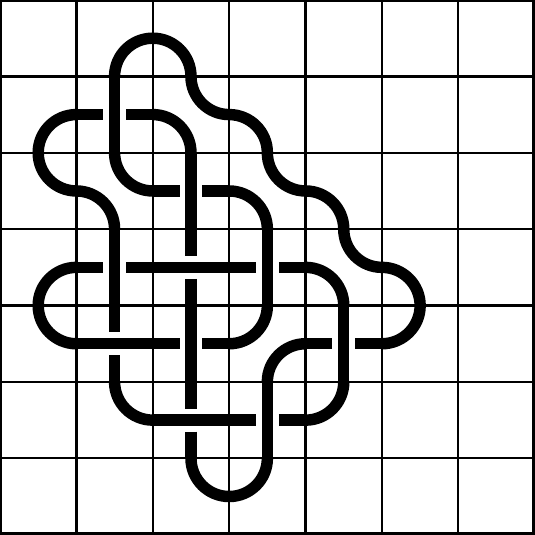}
        \caption*{$10_{102}$ }
    \end{minipage} \hfill
    \begin{minipage}{0.155\linewidth}
        \captionsetup{skip=3pt}
        \centering
        \includegraphics[width=\linewidth]{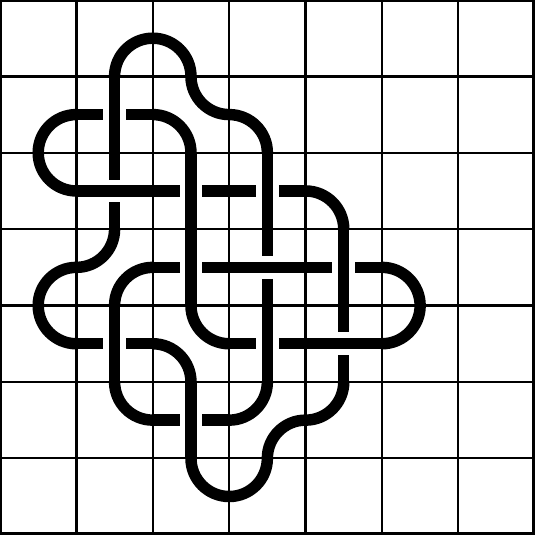}
        \caption*{\phantom{{\Large $\ast$}} $10_{103}$ {\Large $\ast$}}
    \end{minipage} \hfill
    \begin{minipage}{0.155\linewidth}
        \captionsetup{skip=3pt}
        \centering
        \includegraphics[width=\linewidth]{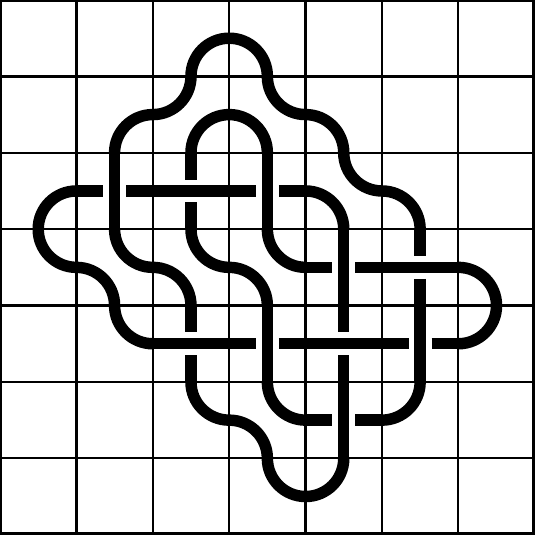}
        \caption*{ $10_{104}$ }
    \end{minipage} \hfill
    \begin{minipage}{0.155\linewidth}
        \captionsetup{skip=3pt}
        \centering
        \includegraphics[width=\linewidth]{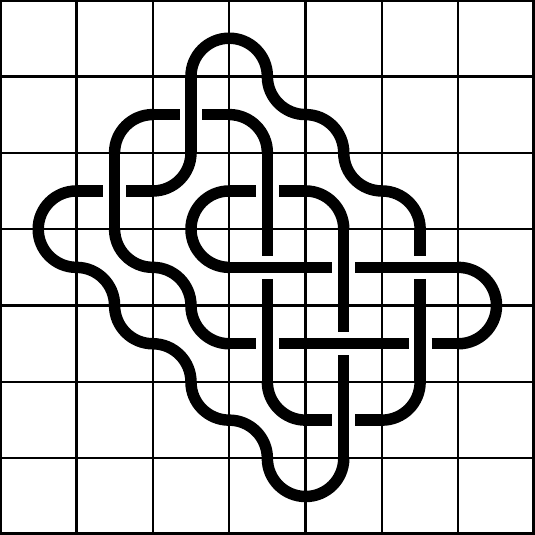}
        \caption*{ $10_{105}$ }
    \end{minipage} \newline
\end{figure}
\unskip

\begin{figure}[H]
    \centering
    \begin{minipage}{0.155\linewidth}
        \captionsetup{skip=3pt}
        \centering
        \includegraphics[width=\linewidth]{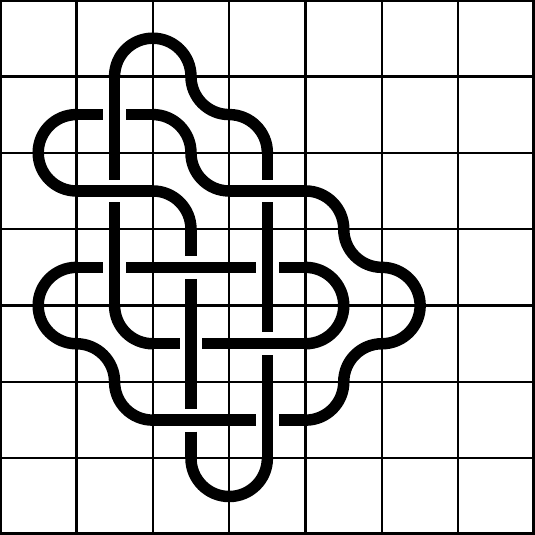}
        \caption*{$10_{106}$ }
    \end{minipage} \hfill
    \begin{minipage}{0.155\linewidth}
        \captionsetup{skip=3pt}
        \centering
        \includegraphics[width=\linewidth]{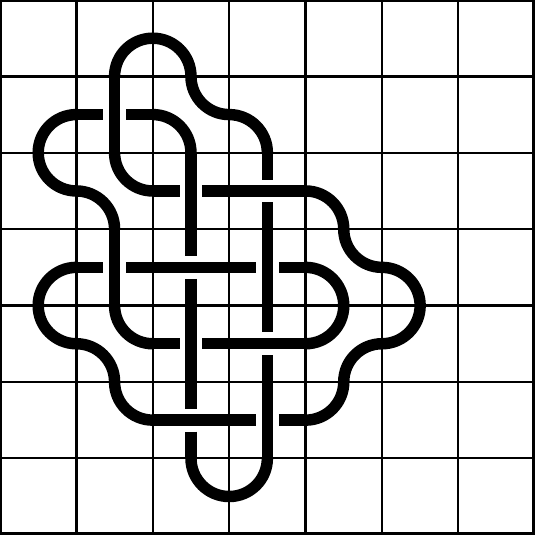}
        \caption*{$10_{107}$ }
    \end{minipage} \hfill
    \begin{minipage}{0.155\linewidth}
        \captionsetup{skip=3pt}
        \centering
        \includegraphics[width=\linewidth]{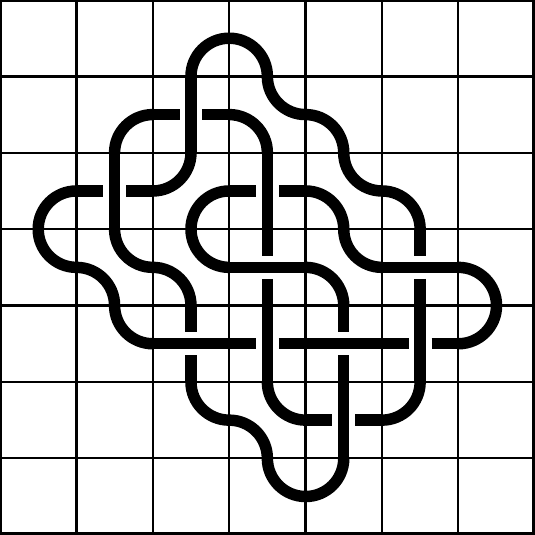}
        \caption*{$10_{108}$}
    \end{minipage} \hfill
    \begin{minipage}{0.155\linewidth}
        \captionsetup{skip=3pt}
        \centering
        \includegraphics[width=\linewidth]{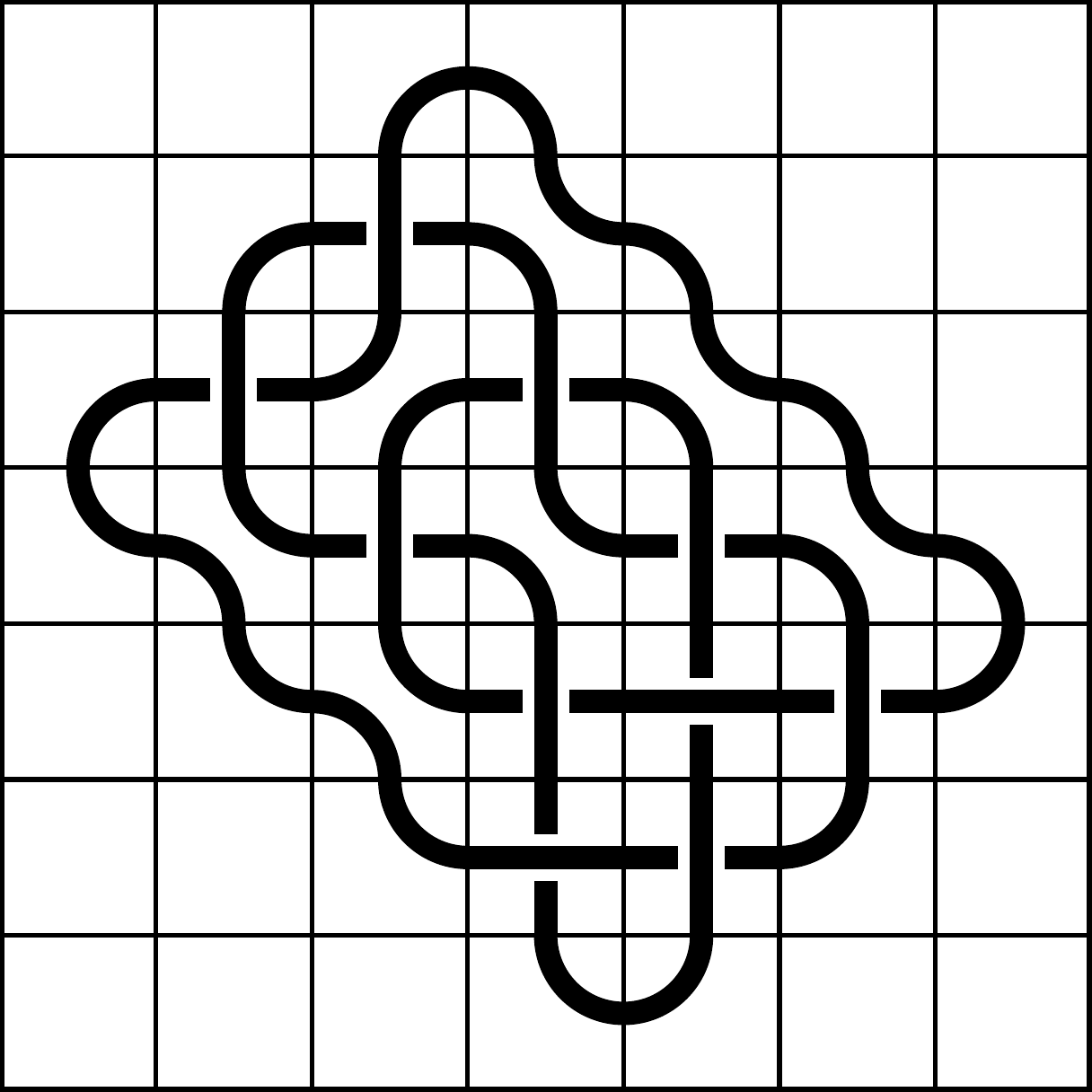}
        \caption*{ $10_{109}$ }
    \end{minipage} \hfill
    \begin{minipage}{0.155\linewidth}
        \captionsetup{skip=3pt}
        \centering
        \includegraphics[width=\linewidth]{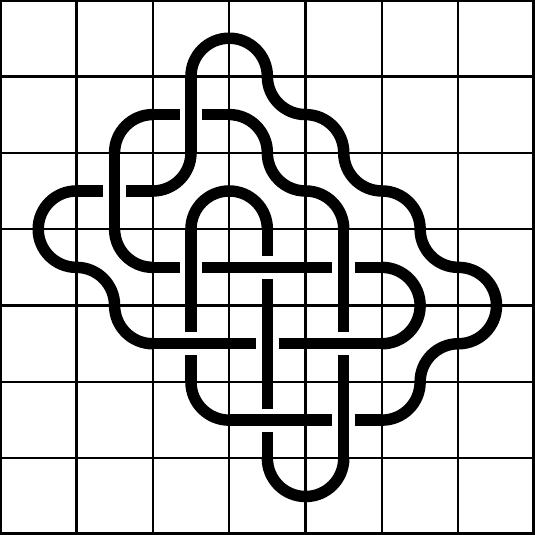}
        \caption*{$10_{110}$}
    \end{minipage} \hfill
    \begin{minipage}{0.155\linewidth}
        \captionsetup{skip=3pt}
        \centering
        \includegraphics[width=\linewidth]{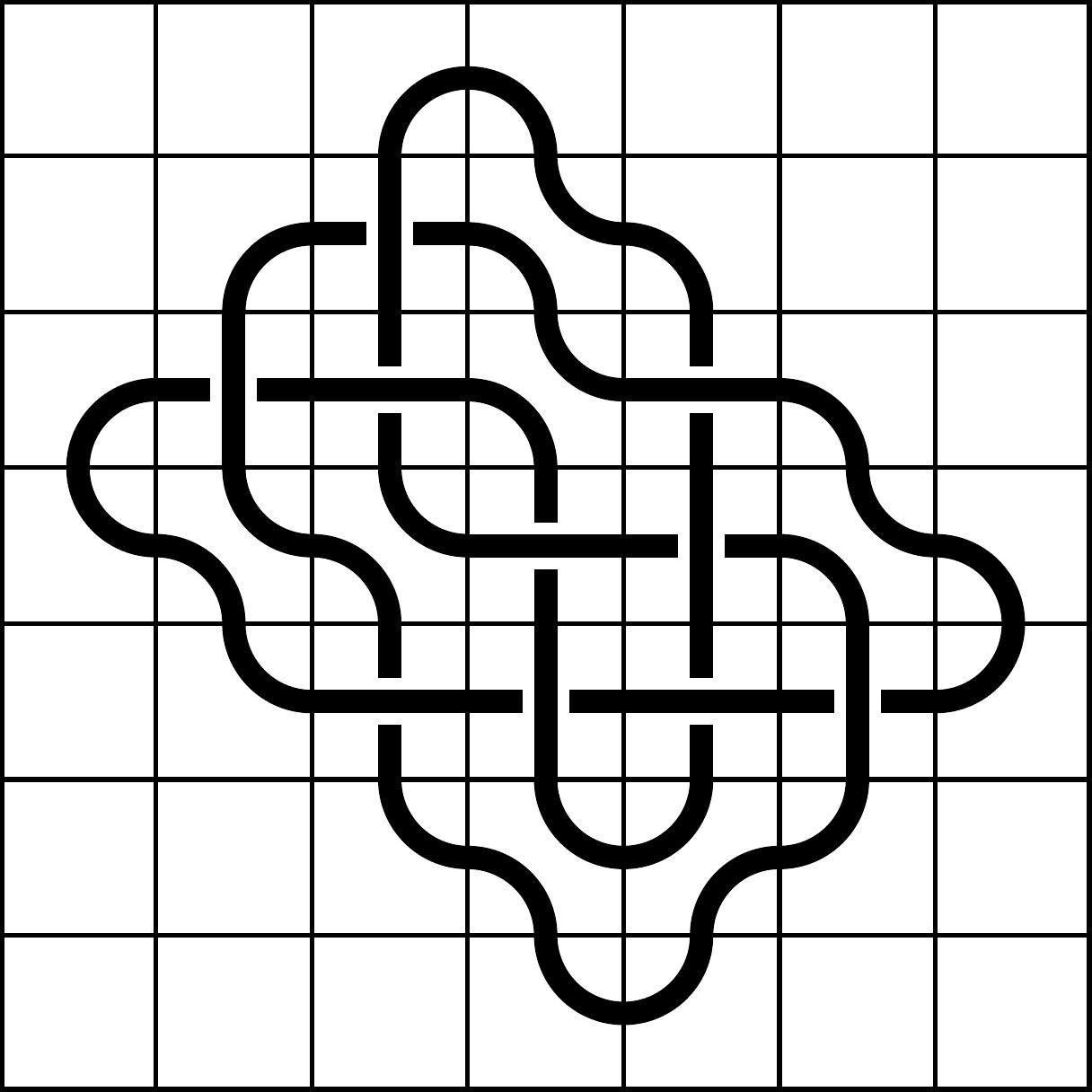}
        \caption*{ $10_{111}$ }
    \end{minipage}  \newline
\end{figure}
\unskip

\begin{figure}[H]
    \centering
    \begin{minipage}{0.155\linewidth}
        \captionsetup{skip=3pt}
        \centering
        \includegraphics[width=\linewidth]{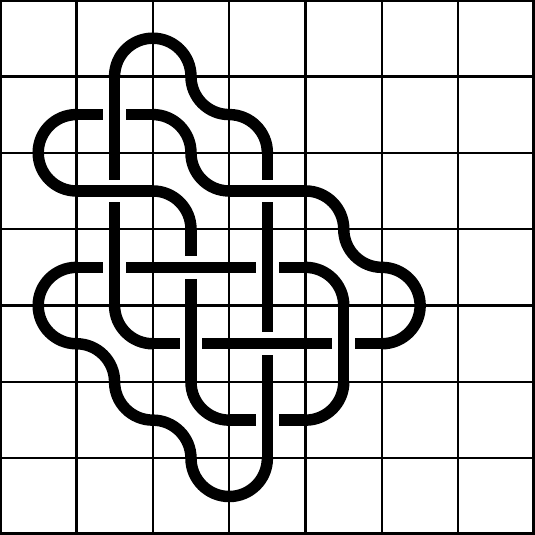}
        \caption*{$10_{112}$ }
    \end{minipage}  \hfill
    \begin{minipage}{0.155\linewidth}
        \captionsetup{skip=3pt}
        \centering
        \includegraphics[width=\linewidth]{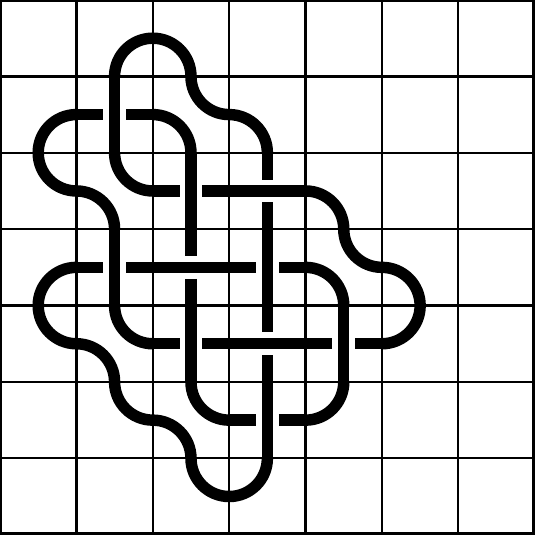}
        \caption*{$10_{113}$ }
    \end{minipage} \hfill
    \begin{minipage}{0.155\linewidth}
        \captionsetup{skip=3pt}
        \centering
        \includegraphics[width=\linewidth]{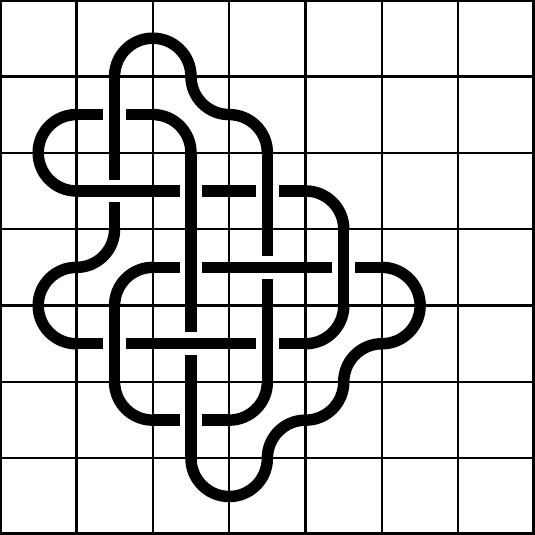}
        \caption*{\phantom{{\Large $\ast$}} $10_{114}$ {\Large $\ast$} }
    \end{minipage} \hfill
    \begin{minipage}{0.155\linewidth}
        \captionsetup{skip=3pt}
        \centering
        \includegraphics[width=\linewidth]{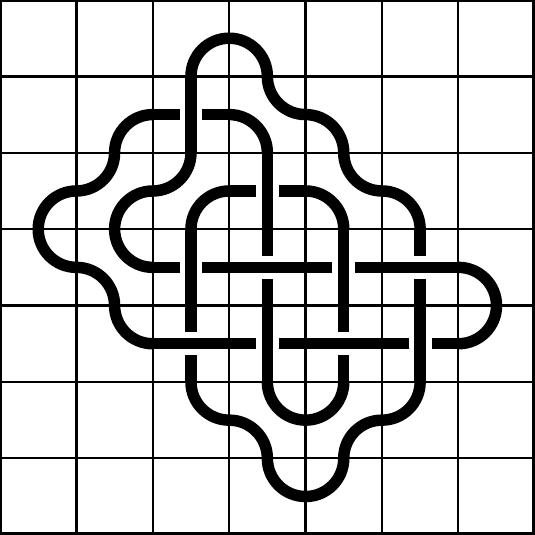}
        \caption*{ $10_{115}$ }
    \end{minipage}\hfill
    \begin{minipage}{0.155\linewidth}
        \captionsetup{skip=3pt}
        \centering
        \includegraphics[width=\linewidth]{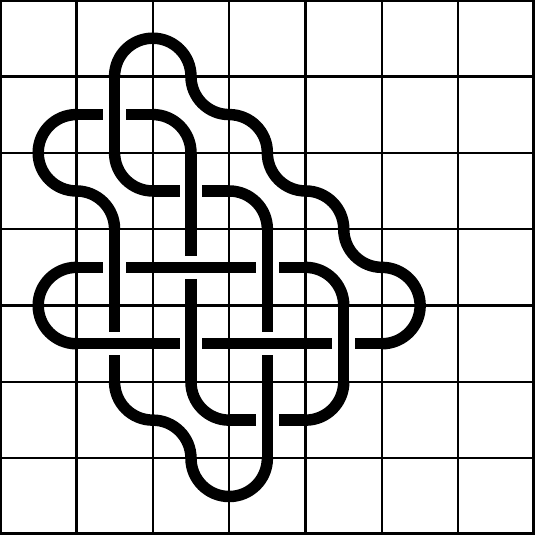}
        \caption*{$10_{117}$ }
    \end{minipage} \hfill
    \begin{minipage}{0.155\linewidth}
        \captionsetup{skip=3pt}
        \centering
        \includegraphics[width=\linewidth]{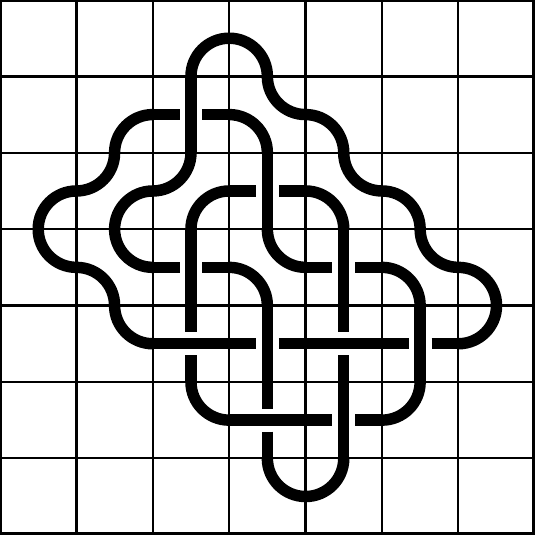}
        \caption*{$10_{118}$}
    \end{minipage} \newline
\end{figure}
\unskip

\begin{figure}[H]
    \centering
    \begin{minipage}{0.155\linewidth}
        \captionsetup{skip=3pt}
        \centering
        \includegraphics[width=\linewidth]{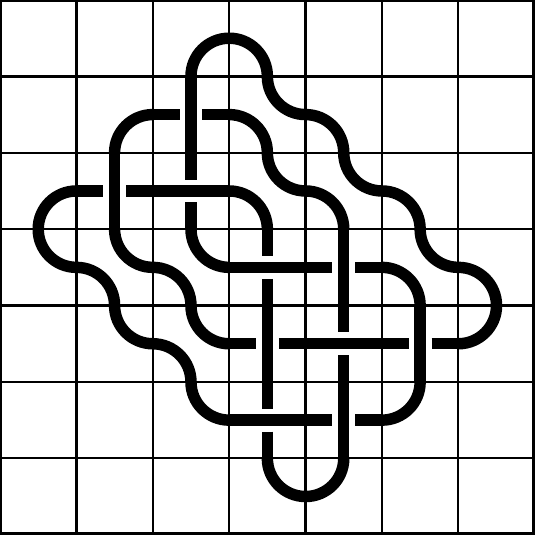}
        \caption*{ $10_{119}$ }
    \end{minipage} \hfill
    \begin{minipage}{0.155\linewidth}
        \captionsetup{skip=3pt}
        \centering
        \includegraphics[width=\linewidth]{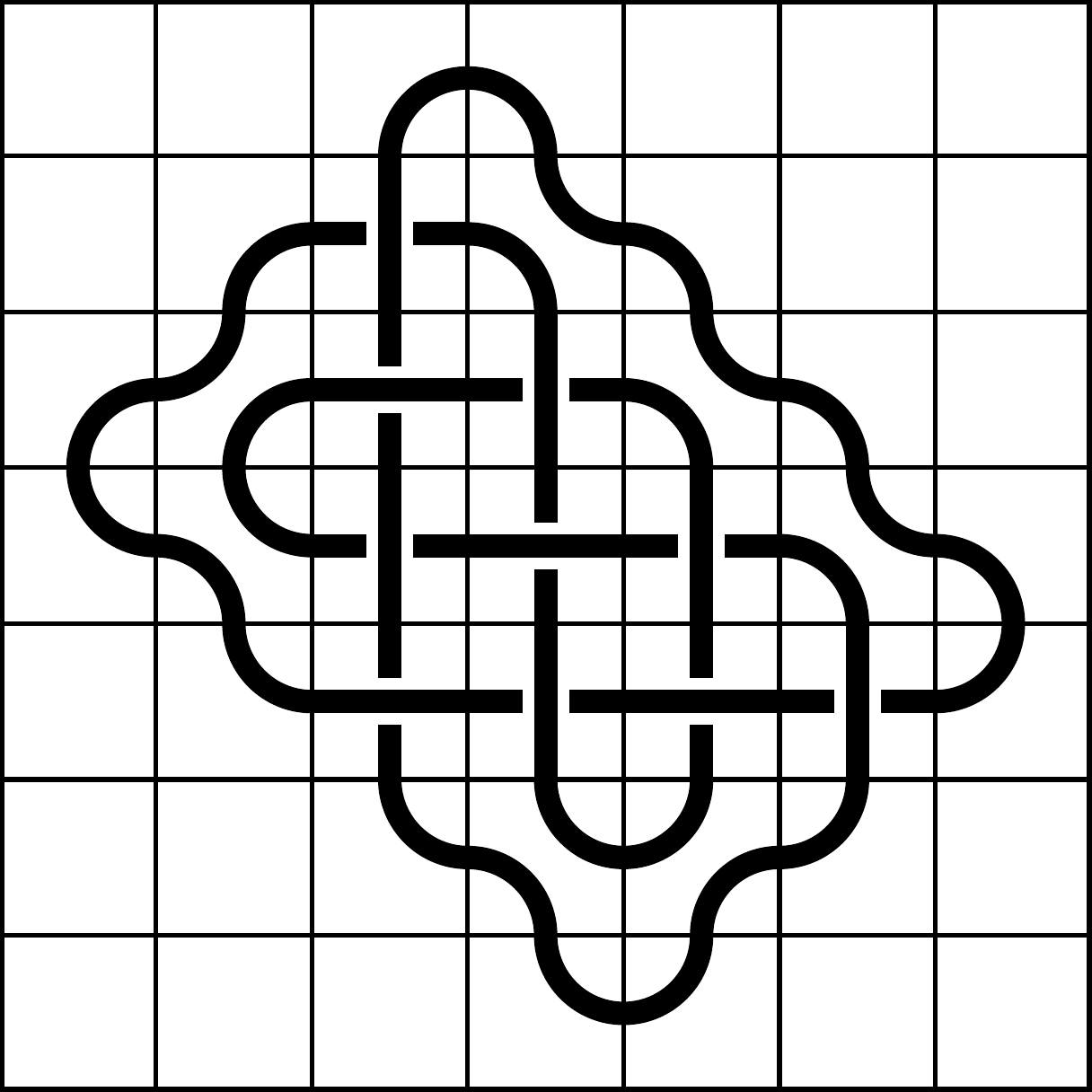}
        \caption*{$10_{120}$}
    \end{minipage} \hfill
    \begin{minipage}{0.155\linewidth}
        \captionsetup{skip=3pt}
        \centering
        \includegraphics[width=\linewidth]{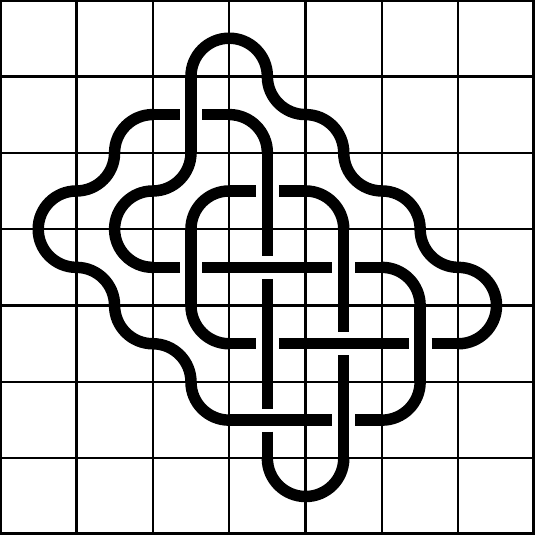}
        \caption*{ $10_{121}$ }
    \end{minipage} \hfill
    \begin{minipage}{0.155\linewidth}
        \captionsetup{skip=3pt}
        \centering
        \includegraphics[width=\linewidth]{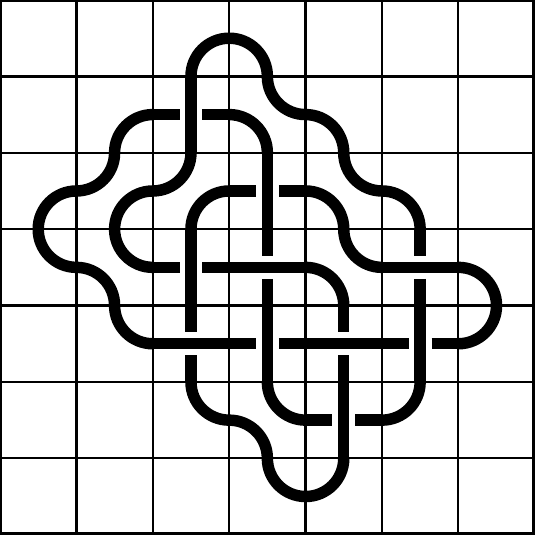}
        \caption*{ $10_{122}$ }
    \end{minipage} \hfill
    \begin{minipage}{0.155\linewidth}
        \captionsetup{skip=3pt}
        \centering
        \includegraphics[width=\linewidth]{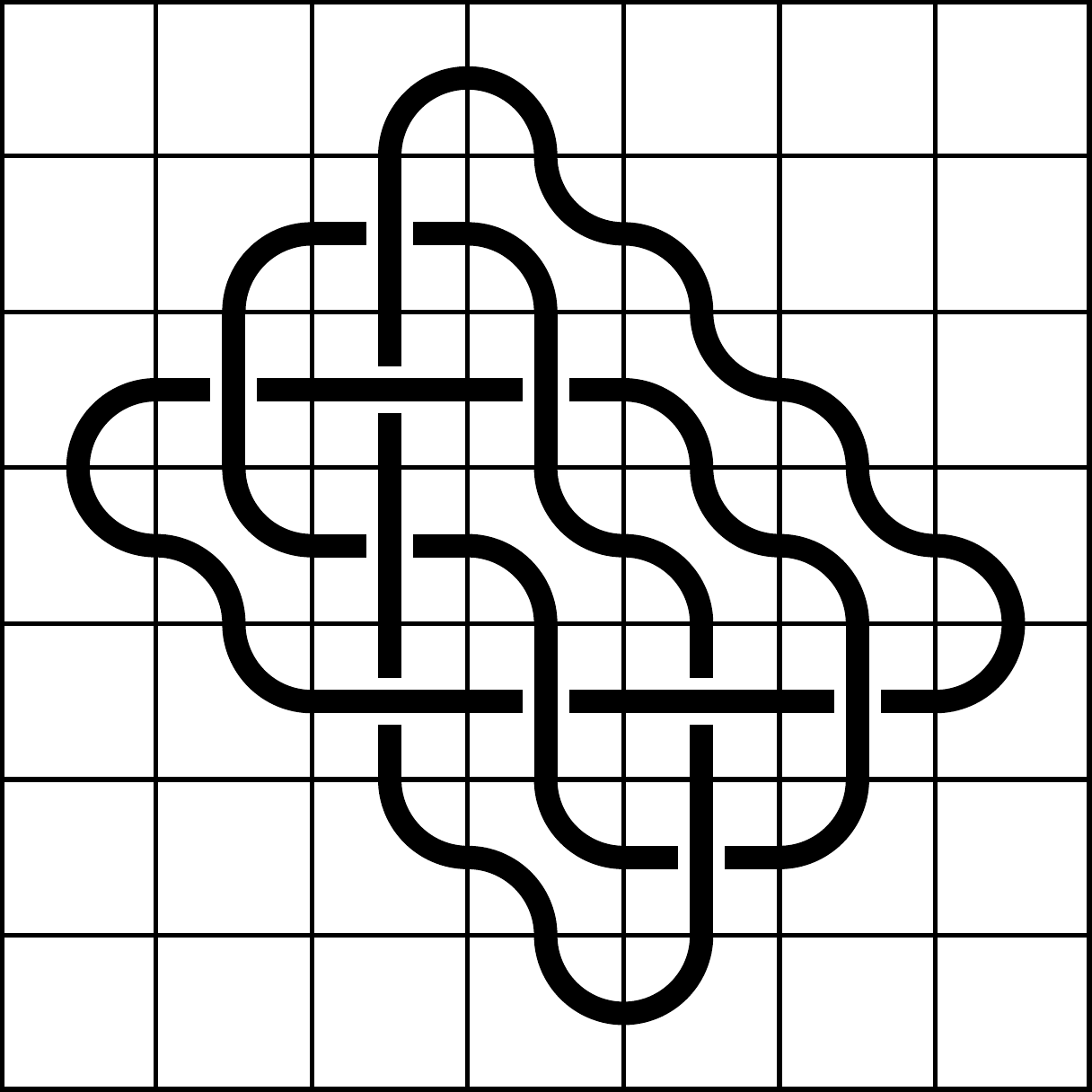}
        \caption*{$10_{123}$}
    \end{minipage} \hfill
    \begin{minipage}{0.155\linewidth}
        \captionsetup{skip=3pt}
        \centering
        \includegraphics[width=\linewidth]{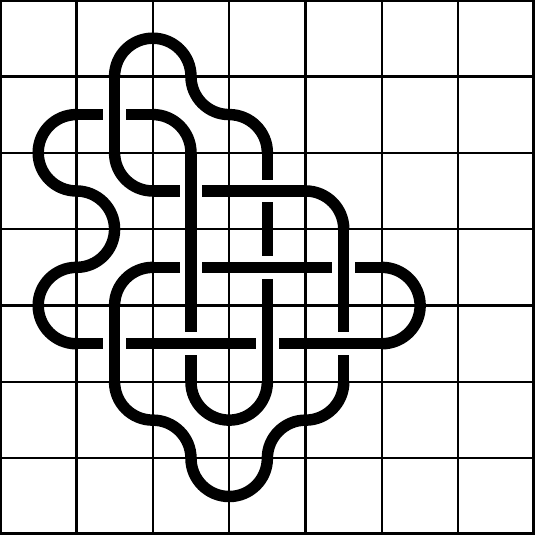}
        \caption*{$10_{128}$ }
    \end{minipage} \newline
\end{figure}
\unskip

\begin{figure}[H]
    \centering
    \begin{minipage}{0.155\linewidth}
        \captionsetup{skip=3pt}
        \centering
        \includegraphics[width=\linewidth]{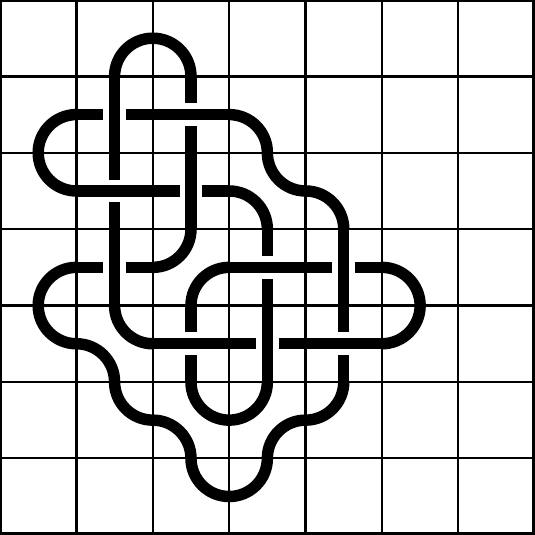}
        \caption*{$10_{129}$ }
    \end{minipage} \hfill
    \begin{minipage}{0.155\linewidth}
        \captionsetup{skip=3pt}
        \centering
        \includegraphics[width=\linewidth]{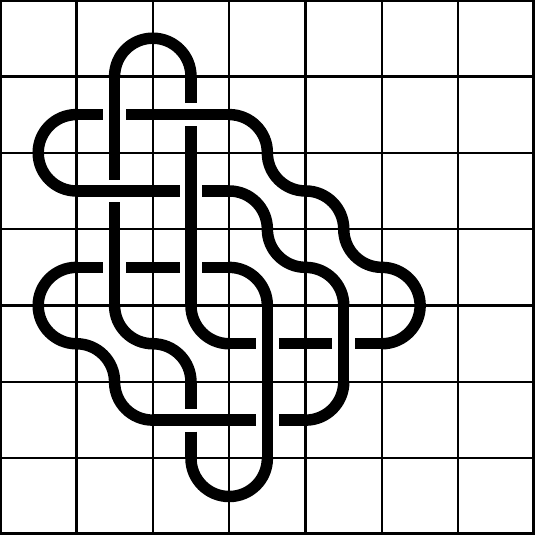}
        \caption*{$10_{130}$ }
    \end{minipage} \hfill
    \begin{minipage}{0.155\linewidth}
        \captionsetup{skip=3pt}
        \centering
        \includegraphics[width=\linewidth]{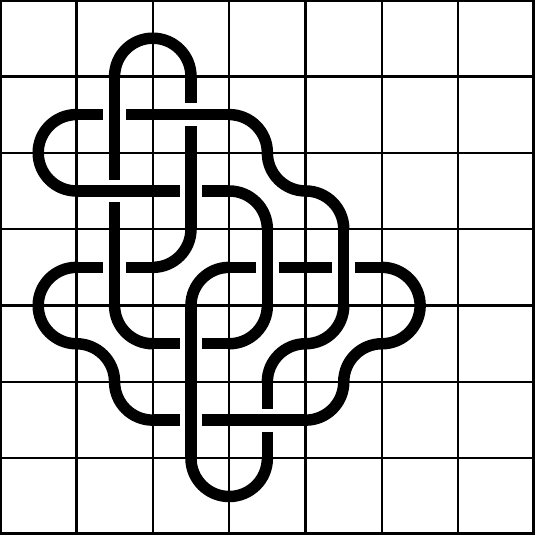}
        \caption*{$10_{131}$ }
    \end{minipage} \hfill
    \begin{minipage}{0.155\linewidth}
        \captionsetup{skip=3pt}
        \centering
        \includegraphics[width=\linewidth]{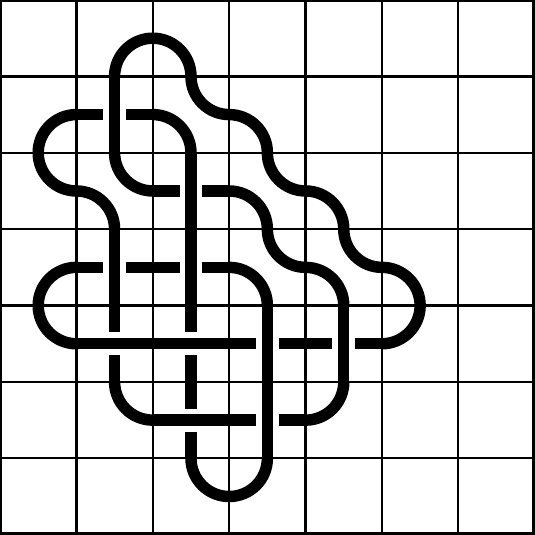}
        \caption*{$10_{132}$ }
    \end{minipage} \hfill
    \begin{minipage}{0.155\linewidth}
        \captionsetup{skip=3pt}
        \centering
        \includegraphics[width=\linewidth]{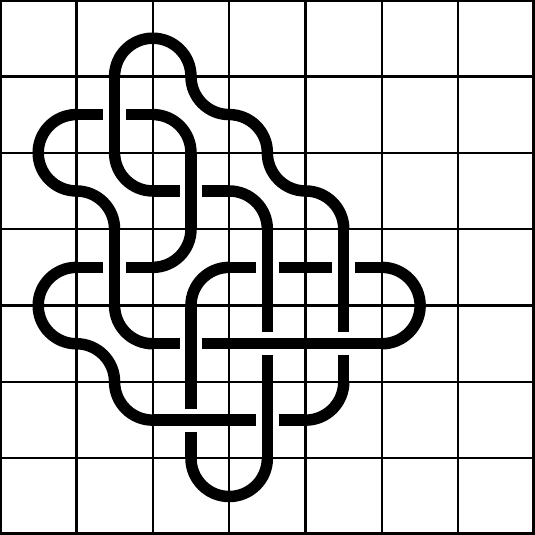}
        \caption*{$10_{133}$ }
    \end{minipage} \hfill
    \begin{minipage}{0.155\linewidth}
        \captionsetup{skip=3pt}
        \centering
        \includegraphics[width=\linewidth]{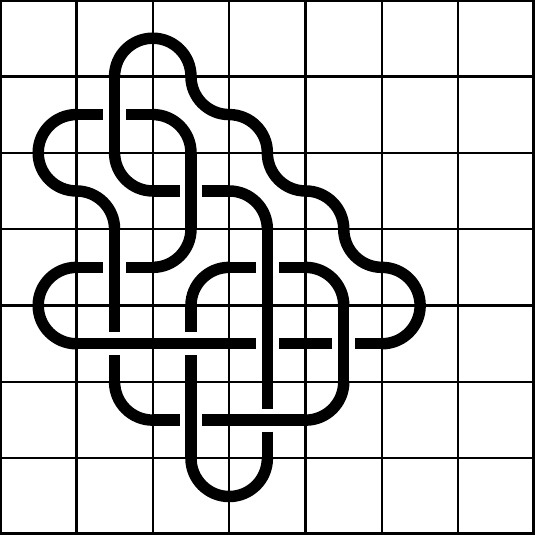}
        \caption*{$10_{134}$ }
    \end{minipage}
\end{figure}
\unskip

\begin{figure}[H]
    \centering
    \begin{minipage}{0.155\linewidth}
        \captionsetup{skip=3pt}
        \centering
        \includegraphics[width=\linewidth]{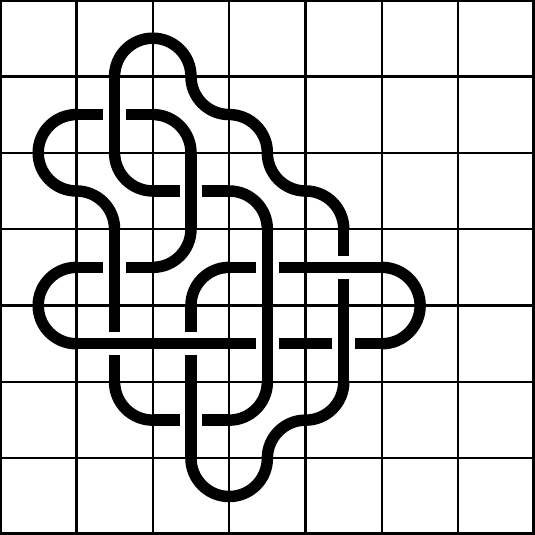}
        \caption*{$10_{135}$ }
    \end{minipage} \hfill
    \begin{minipage}{0.155\linewidth}
        \captionsetup{skip=3pt}
        \centering
        \includegraphics[width=\linewidth]{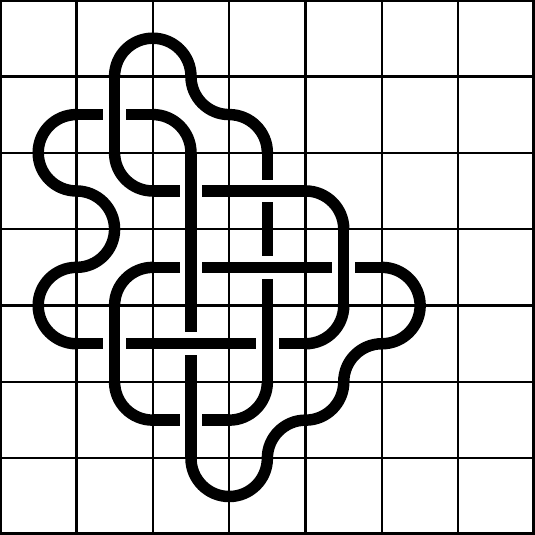}
        \caption*{$10_{136}$ }
    \end{minipage} \hfill
    \begin{minipage}{0.155\linewidth}
        \captionsetup{skip=3pt}
        \centering
        \includegraphics[width=\linewidth]{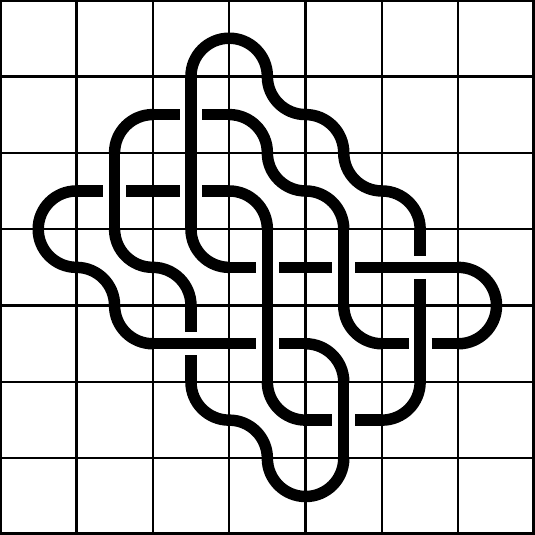}
        \caption*{ $10_{137}$ }
    \end{minipage} \hfill
    \begin{minipage}{0.155\linewidth}
        \captionsetup{skip=3pt}
        \centering
        \includegraphics[width=\linewidth]{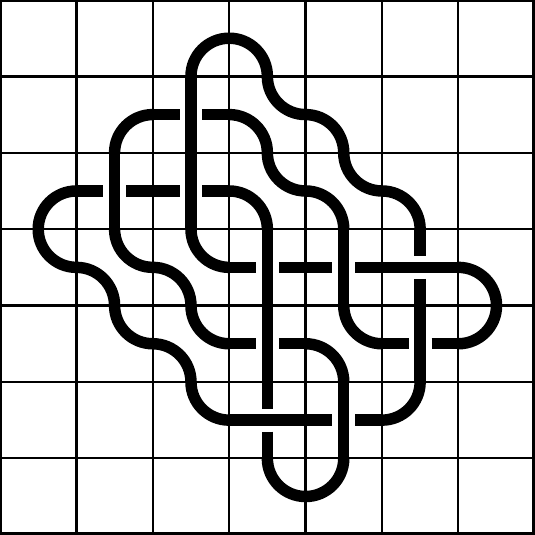}
        \caption*{$10_{138}$}
    \end{minipage} \hfill
    \begin{minipage}{0.155\linewidth}
        \captionsetup{skip=3pt}
        \centering
        \includegraphics[width=\linewidth]{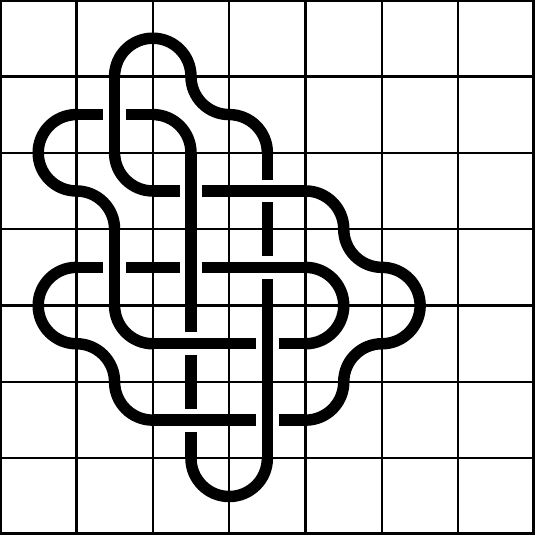}
        \caption*{$10_{145}$ }
    \end{minipage} \hfill
    \begin{minipage}{0.155\linewidth}
        \captionsetup{skip=3pt}
        \centering
        \includegraphics[width=\linewidth]{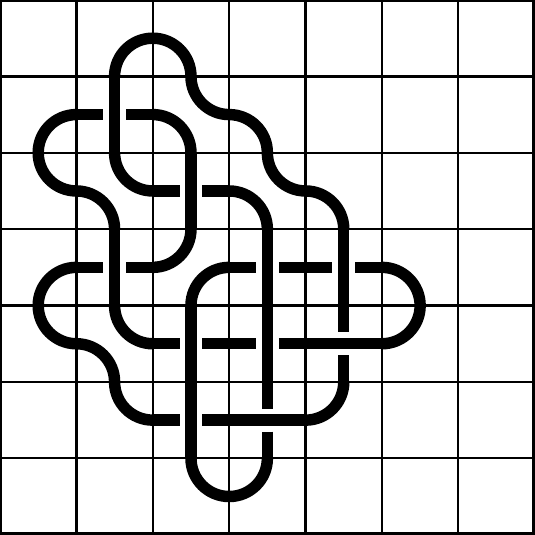}
        \caption*{$10_{146}$ }
    \end{minipage} \newline
\end{figure}
\unskip

\begin{figure}[H]
    \centering
    \begin{minipage}{0.155\linewidth}
        \captionsetup{skip=3pt}
        \centering
        \includegraphics[width=\linewidth]{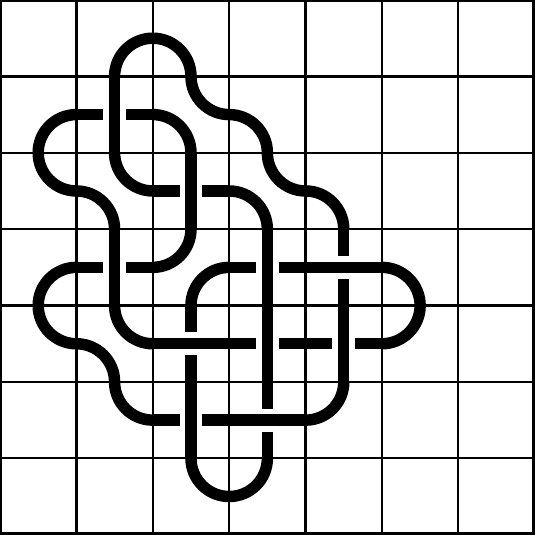}
        \caption*{$10_{147}$ }
    \end{minipage} \hfill
    \begin{minipage}{0.155\linewidth}
        \captionsetup{skip=3pt}
        \centering
        \includegraphics[width=\linewidth]{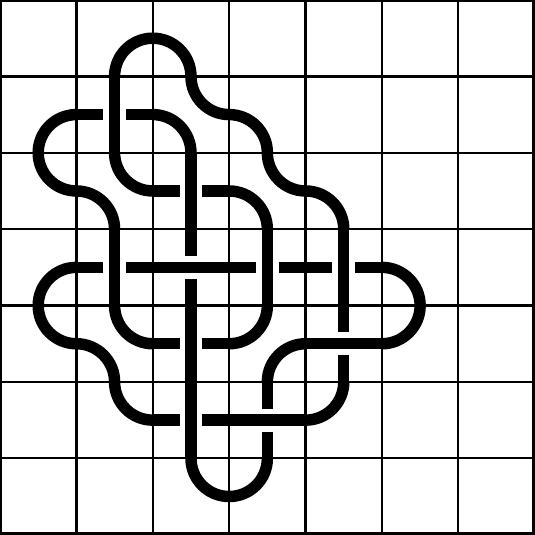}
        \caption*{$10_{149}$ }
    \end{minipage} \hfill
    \begin{minipage}{0.155\linewidth}
        \captionsetup{skip=3pt}
        \centering
        \includegraphics[width=\linewidth]{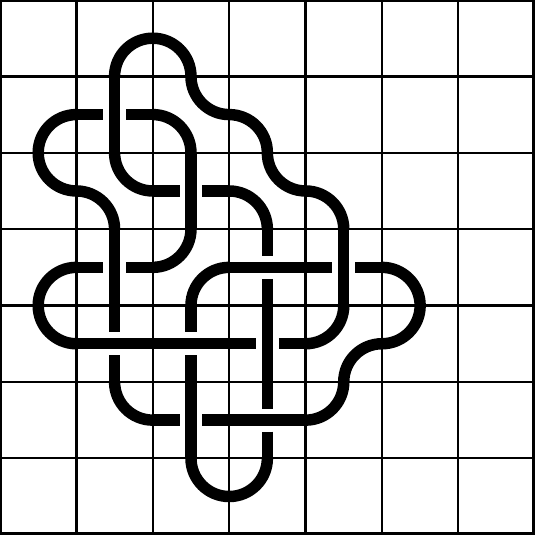}
        \caption*{$10_{150}$ }
    \end{minipage} \hfill
    \begin{minipage}{0.155\linewidth}
        \captionsetup{skip=3pt}
        \centering
        \includegraphics[width=\linewidth]{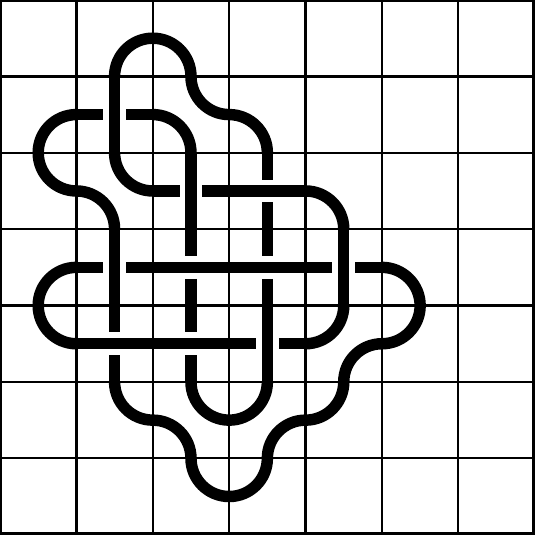}
        \caption*{$10_{151}$ }
    \end{minipage} \hfill
    \begin{minipage}{0.155\linewidth}
        \captionsetup{skip=3pt}
        \centering
        \includegraphics[width=\linewidth]{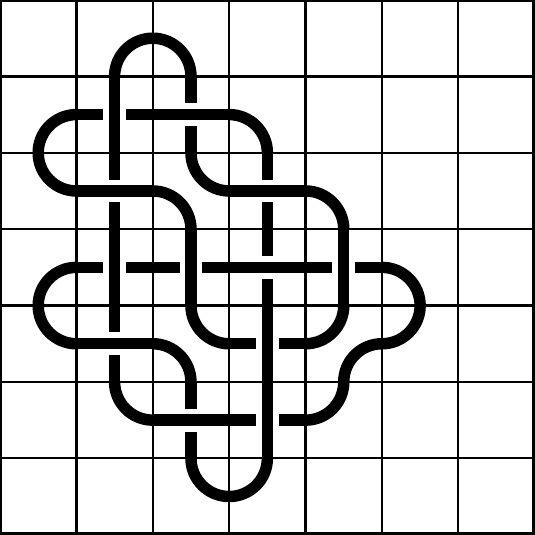}
        \caption*{\phantom{{\Large $\ast$}} $10_{152}$ {\Large $\ast$} }
    \end{minipage} \hfill
    \begin{minipage}{0.155\linewidth}
        \captionsetup{skip=3pt}
        \centering
        \includegraphics[width=\linewidth]{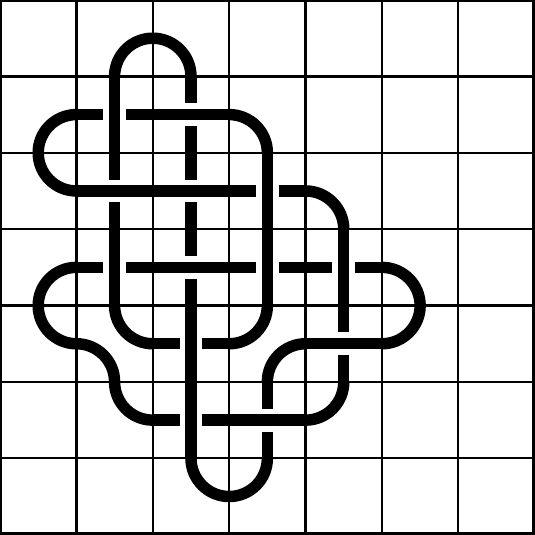}
        \caption*{\phantom{{\Large $\ast$}} $10_{153}$ {\Large $\ast$} }
    \end{minipage}  \newline
\end{figure}
\unskip

\begin{figure}[H]
    \centering
     \begin{minipage}{0.155\linewidth}
        \captionsetup{skip=3pt}
        \centering
        \includegraphics[width=\linewidth]{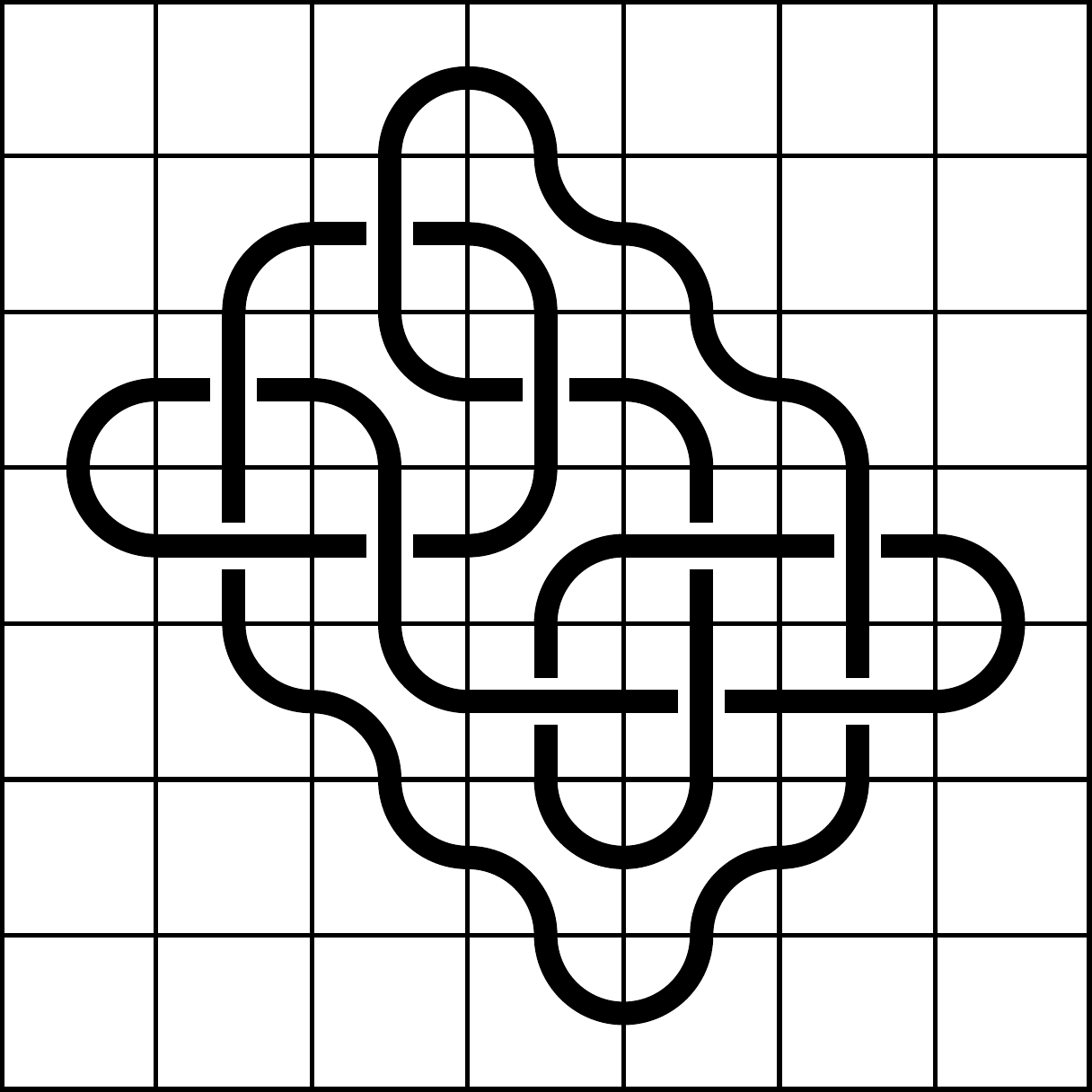}
        \caption*{ $10_{154}$ }
    \end{minipage} \hfill
    \begin{minipage}{0.155\linewidth}
        \captionsetup{skip=3pt}
        \centering
        \includegraphics[width=\linewidth]{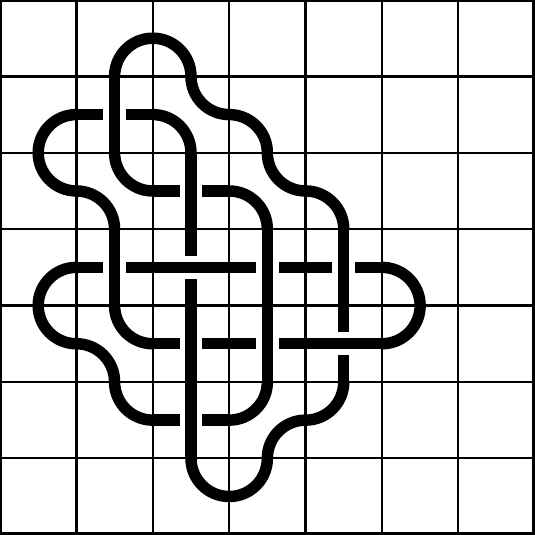}
        \caption*{$10_{156}$ }
    \end{minipage} \hfill
    \begin{minipage}{0.155\linewidth}
        \captionsetup{skip=3pt}
        \centering
        \includegraphics[width=\linewidth]{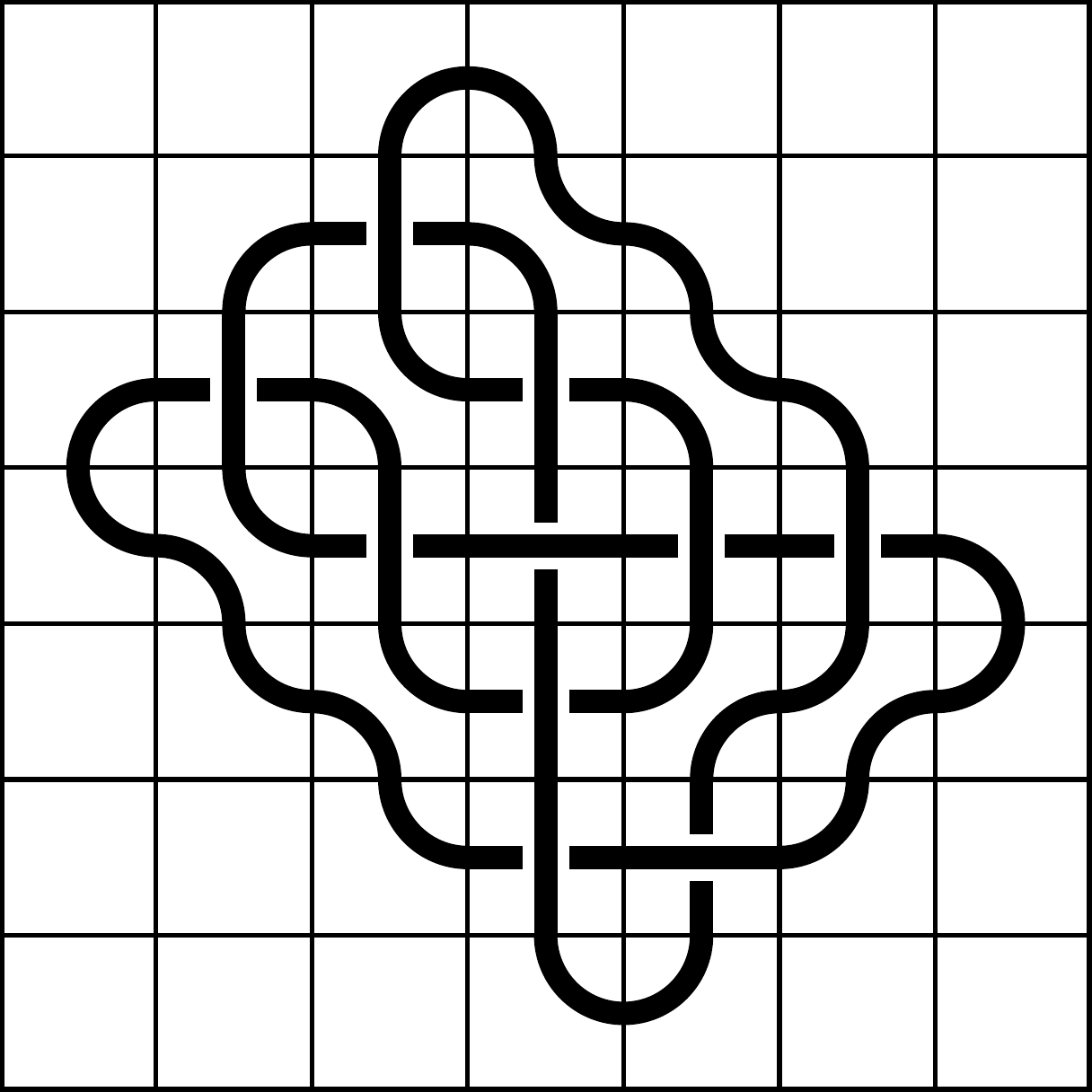}
        \caption*{ $10_{157}$ }
    \end{minipage}  \hfill
    \begin{minipage}{0.155\linewidth}
        \captionsetup{skip=3pt}
        \centering
        \includegraphics[width=\linewidth]{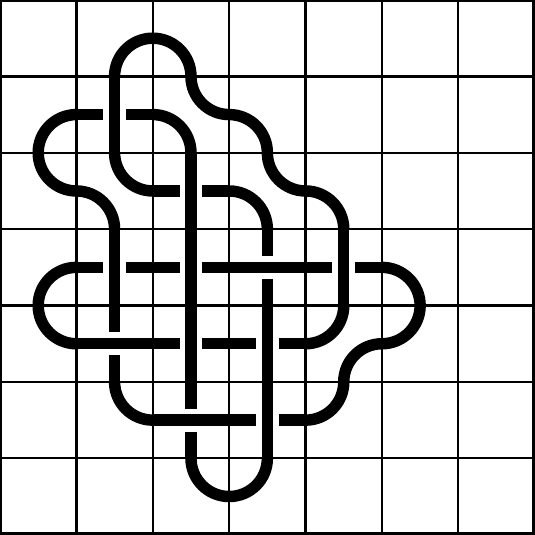}
        \caption*{\phantom{{\Large $\ast$}} $10_{158}$ {\Large $\ast$} }
    \end{minipage} \hfill
    \begin{minipage}{0.155\linewidth}
        \captionsetup{skip=3pt}
        \centering
        \includegraphics[width=\linewidth]{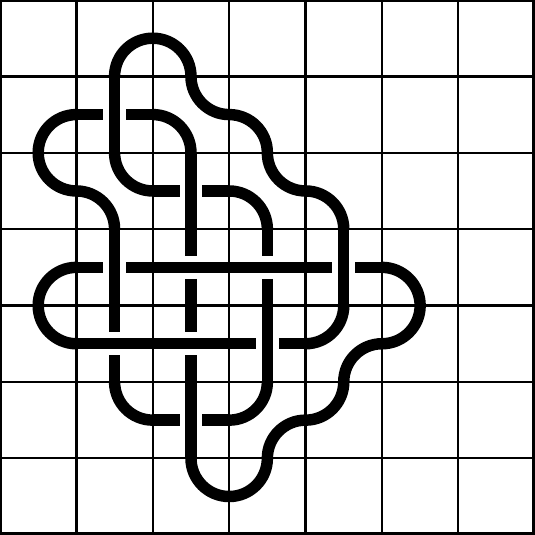}
        \caption*{$10_{160}$ }
    \end{minipage} \hfill
    \begin{minipage}{0.155\linewidth}
        \captionsetup{skip=3pt}
        \centering
        \includegraphics[width=\linewidth]{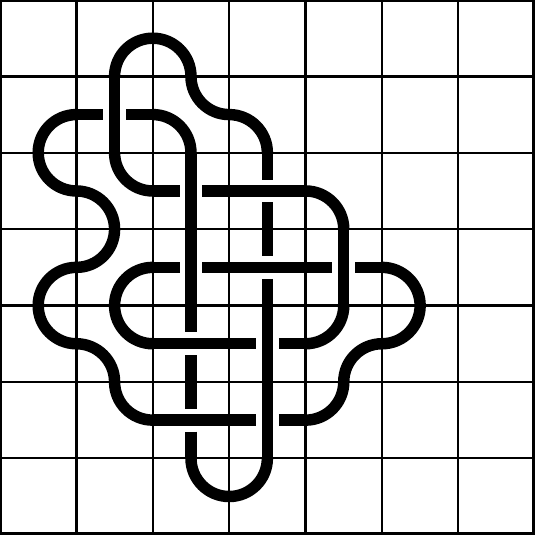}
        \caption*{$10_{161}$ }
    \end{minipage}  \newline
\end{figure}
\unskip

\begin{figure}[H]
    \centering
    \begin{minipage}{0.155\linewidth}
        \captionsetup{skip=3pt}
        \centering
        \includegraphics[width=\linewidth]{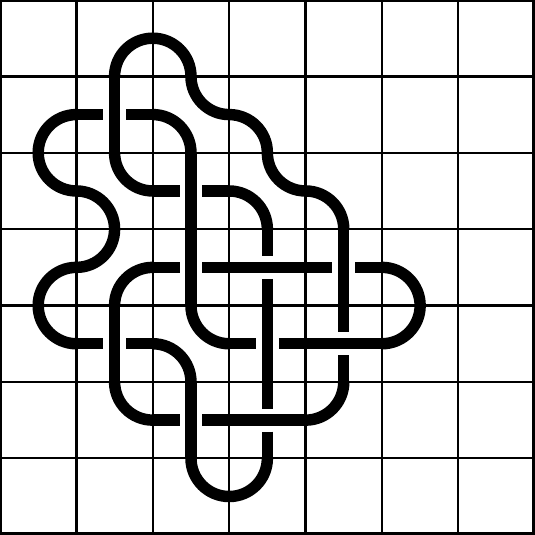}
        \caption*{$10_{162}$ }
    \end{minipage}
    \begin{minipage}{0.155\linewidth}
        \captionsetup{skip=3pt}
        \centering
        \includegraphics[width=\linewidth]{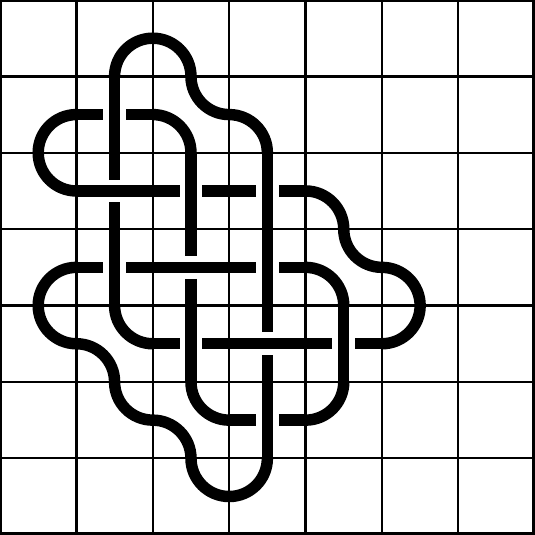}
        \caption*{\phantom{{\Large $\ast$}} $10_{163}$ {\Large $\ast$} }
    \end{minipage}
    \begin{minipage}{0.155\linewidth}
        \captionsetup{skip=3pt}
        \centering
        \includegraphics[width=\linewidth]{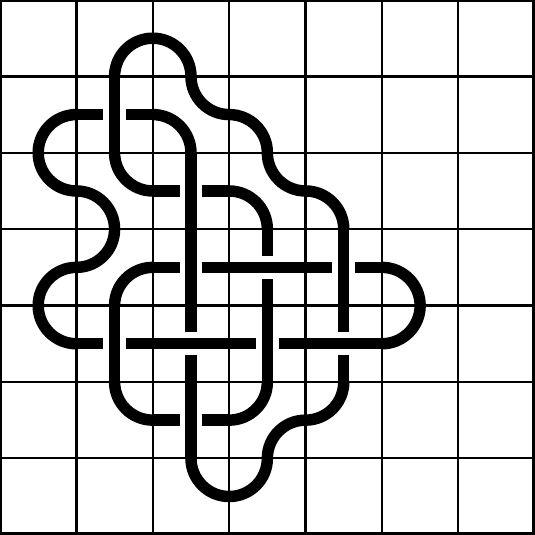}
        \caption*{$10_{164}$ }
    \end{minipage}
    \begin{minipage}{0.155\linewidth}
        \captionsetup{skip=3pt}
        \centering
        \includegraphics[width=\linewidth]{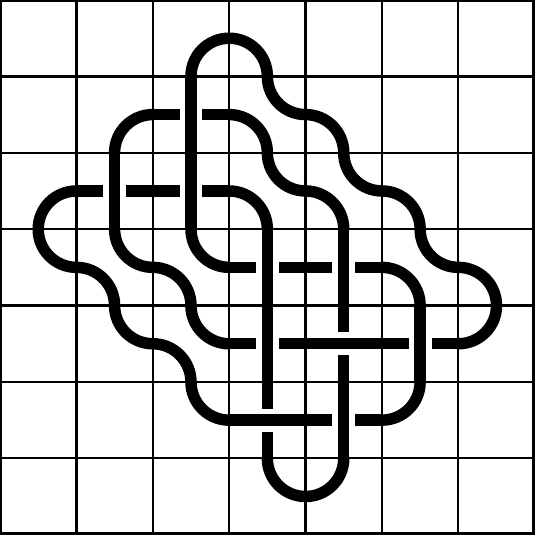}
        \caption*{$10_{165}$}
    \end{minipage}
\end{figure}

\clearpage

\subsection{Mosaics from Theorem \ref{thm:6mosaic7tile}}

These are the prime knots with crossing number 11 or larger, mosaic number 6 needing 32 non-blank tiles, and tile number 27, 29, or 31 realized on a 7-mosaic.

\begin{figure}[H]
    \centering
    \begin{minipage}{0.155\linewidth}
        \captionsetup{skip=3pt}
        \centering
        \includegraphics[width=\linewidth]{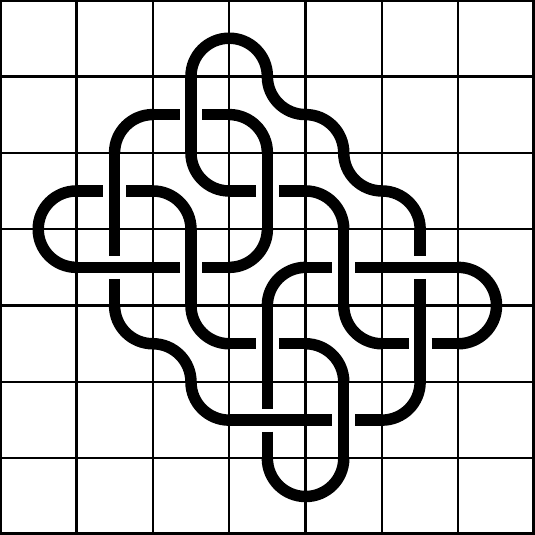}
        \caption*{\ka{11}{43}}
    \end{minipage} \hfill
     \begin{minipage}{0.155\linewidth}
        \captionsetup{skip=3pt}
        \centering
        \includegraphics[width=\linewidth]{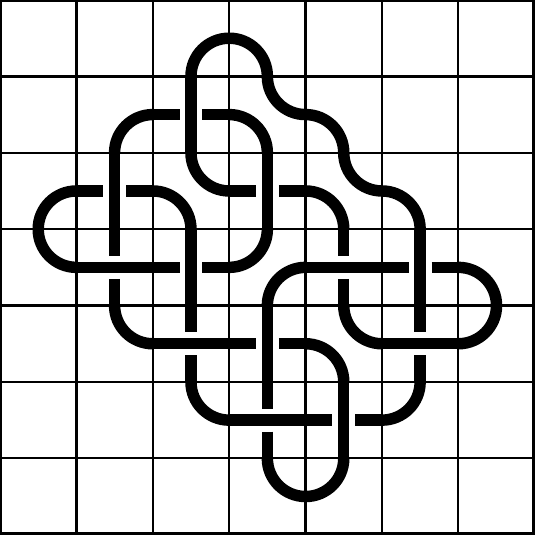}
        \caption*{\phantom{{\Large $\ast$}} \ka{11}{44} {\Large $\ast$}}
    \end{minipage} \hfill
    \begin{minipage}{0.155\linewidth}
        \captionsetup{skip=3pt}
        \centering
        \includegraphics[width=\linewidth]{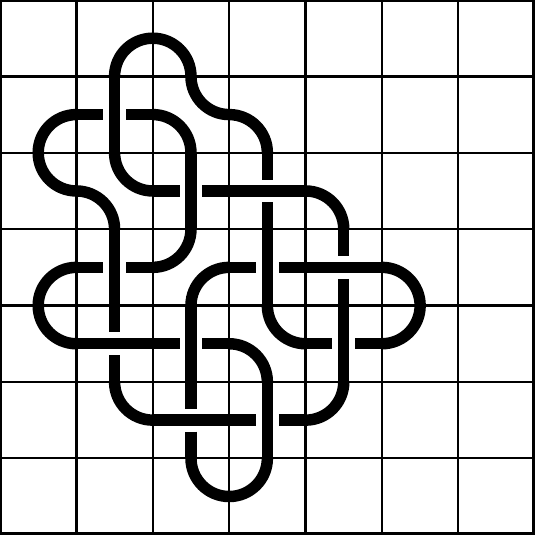}
        \caption*{\ka{11}{46}}
    \end{minipage} \hfill
     \begin{minipage}{0.155\linewidth}
        \captionsetup{skip=3pt}
        \centering
        \includegraphics[width=\linewidth]{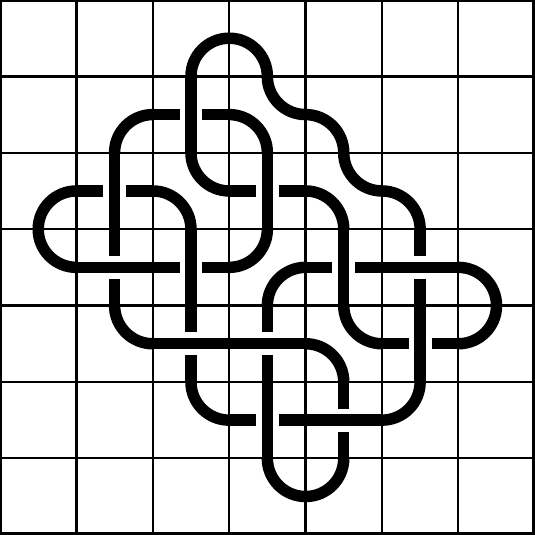}
        \caption*{\phantom{{\Large $\ast$}} \ka{11}{47} {\Large $\ast$}}
    \end{minipage} \hfill
    \begin{minipage}{0.155\linewidth}
        \captionsetup{skip=3pt}
        \centering
        \includegraphics[width=\linewidth]{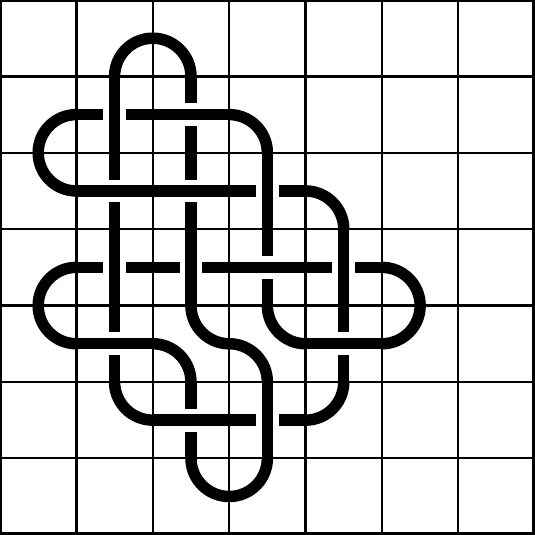}
        \caption*{\phantom{{\Large $\ast$}} \ka{11}{58} {\Large $\ast$}}
    \end{minipage}  \hfill
    \begin{minipage}{0.155\linewidth}
        \captionsetup{skip=3pt}
        \centering
        \includegraphics[width=\linewidth]{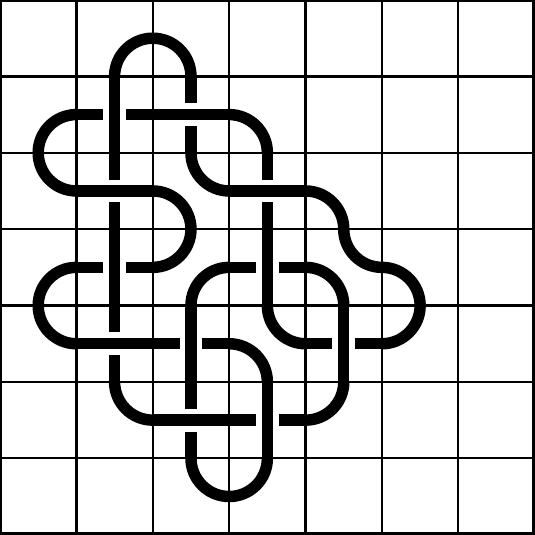}
        \caption*{\ka{11}{59}}
    \end{minipage} \newline
\end{figure}
\unskip

\begin{figure}[H]
    \centering
    \begin{minipage}{0.155\linewidth}
        \captionsetup{skip=3pt}
        \centering
        \includegraphics[width=\linewidth]{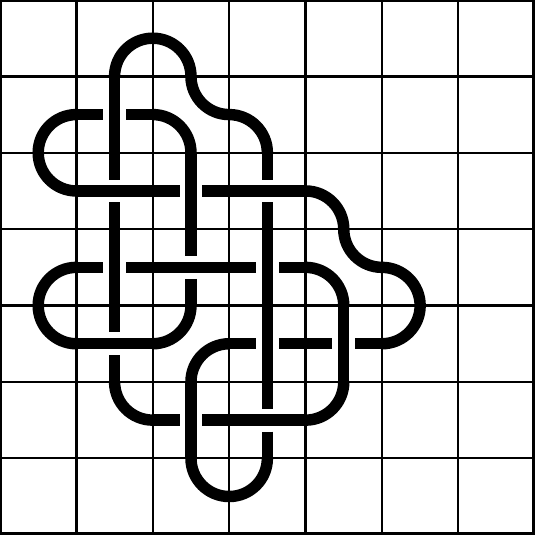}
        \caption*{\phantom{{\Large $\ast$}} \ka{11}{106} {\Large $\ast$}}
    \end{minipage} \hfill
     \begin{minipage}{0.155\linewidth}
        \captionsetup{skip=3pt}
        \centering
        \includegraphics[width=\linewidth]{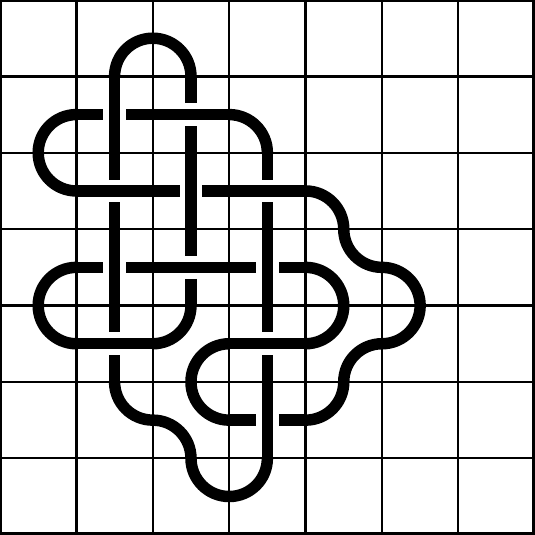}
        \caption*{\ka{11}{139}}
    \end{minipage} \hfill
    \begin{minipage}{0.155\linewidth}
        \captionsetup{skip=3pt}
        \centering
        \includegraphics[width=\linewidth]{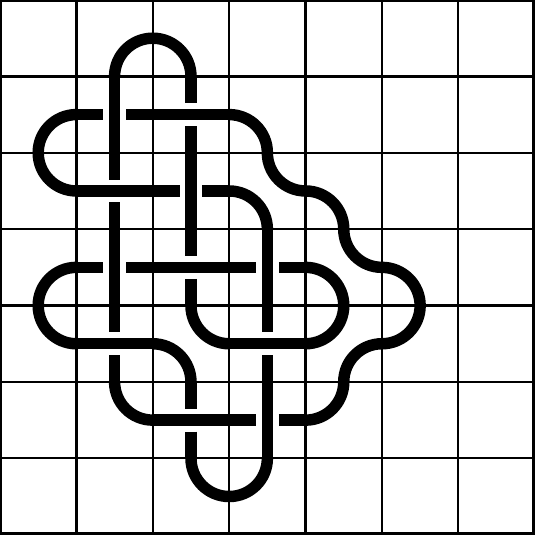}
        \caption*{\ka{11}{165}}
    \end{minipage} \hfill
     \begin{minipage}{0.155\linewidth}
        \captionsetup{skip=3pt}
        \centering
        \includegraphics[width=\linewidth]{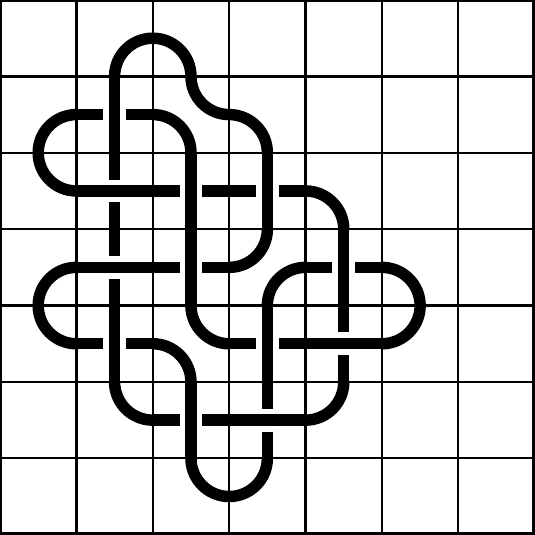}
        \caption*{\phantom{{\Large $\ast$}} \ka{11}{166} {\Large $\ast$}}
    \end{minipage} \hfill
    \begin{minipage}{0.155\linewidth}
        \captionsetup{skip=3pt}
        \centering
        \includegraphics[width=\linewidth]{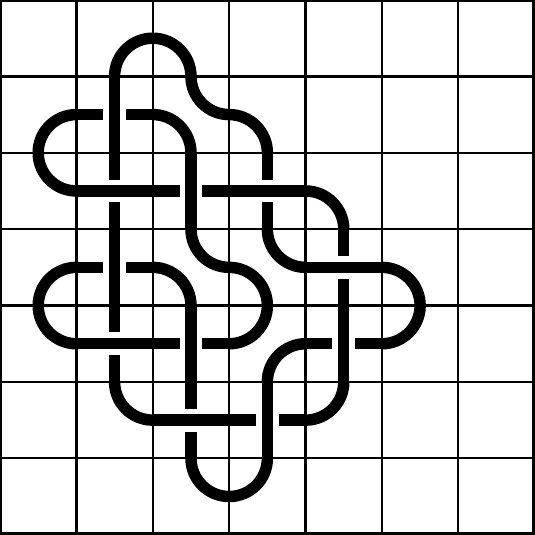}
        \caption*{\ka{11}{179}}
    \end{minipage}  \hfill
    \begin{minipage}{0.155\linewidth}
        \captionsetup{skip=3pt}
        \centering
        \includegraphics[width=\linewidth]{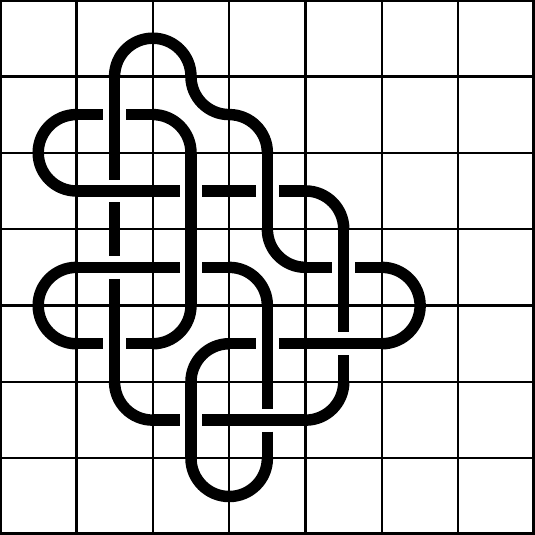}
        \caption*{\phantom{{\Large $\ast$}} \ka{11}{181} {\Large $\ast$}}
    \end{minipage} \newline
\end{figure}
\unskip

\begin{figure}[H]
    \centering
    \begin{minipage}{0.155\linewidth}
        \captionsetup{skip=3pt}
        \centering
        \includegraphics[width=\linewidth]{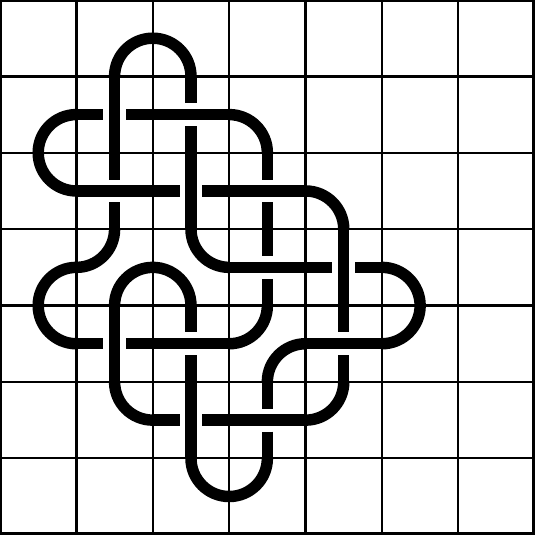}
        \caption*{\phantom{{\Large $\ast$}} \ka{11}{246} {\Large $\ast$}}
    \end{minipage} \hfill
     \begin{minipage}{0.155\linewidth}
        \captionsetup{skip=3pt}
        \centering
        \includegraphics[width=\linewidth]{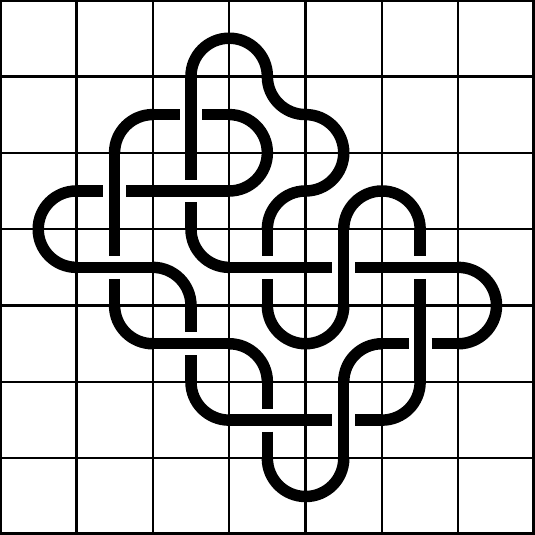}
        \caption*{\ka{11}{247}}
    \end{minipage} \hfill
    \begin{minipage}{0.155\linewidth}
        \captionsetup{skip=3pt}
        \centering
        \includegraphics[width=\linewidth]{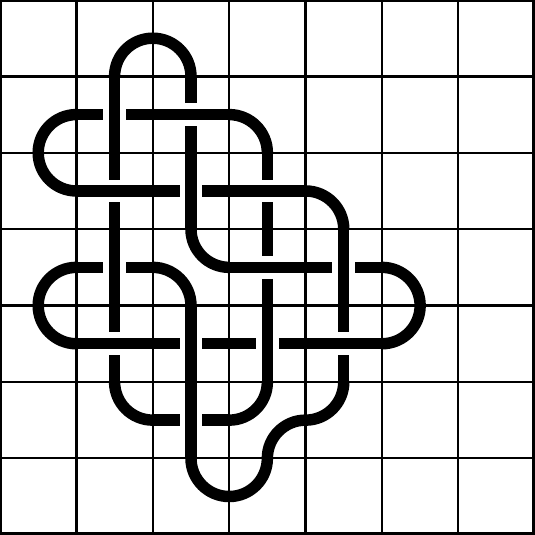}
        \caption*{\phantom{{\Large $\ast$}} \ka{11}{339} {\Large $\ast$}}
    \end{minipage} \hfill
     \begin{minipage}{0.155\linewidth}
        \captionsetup{skip=3pt}
        \centering
        \includegraphics[width=\linewidth]{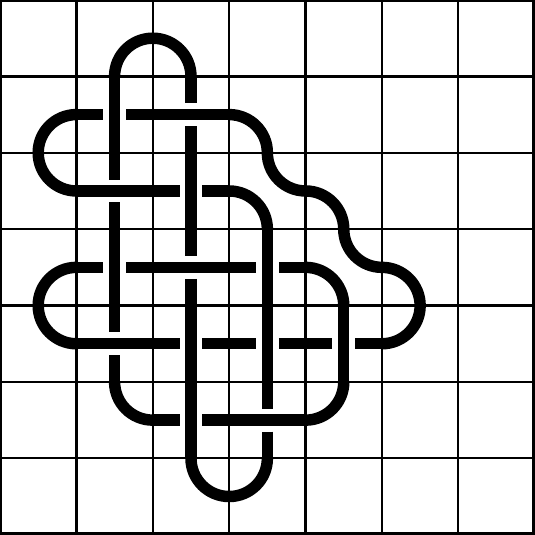}
        \caption*{\phantom{{\Large $\ast$}} \ka{11}{340} {\Large $\ast$}}
    \end{minipage} \hfill
    \begin{minipage}{0.155\linewidth}
        \captionsetup{skip=3pt}
        \centering
        \includegraphics[width=\linewidth]{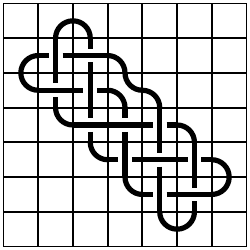}
        \caption*{\phantom{{\Large $\ast$}} \ka{11}{341} {\Large $\ast$}}
    \end{minipage}  \hfill
    \begin{minipage}{0.155\linewidth}
        \captionsetup{skip=3pt}
        \centering
        \includegraphics[width=\linewidth]{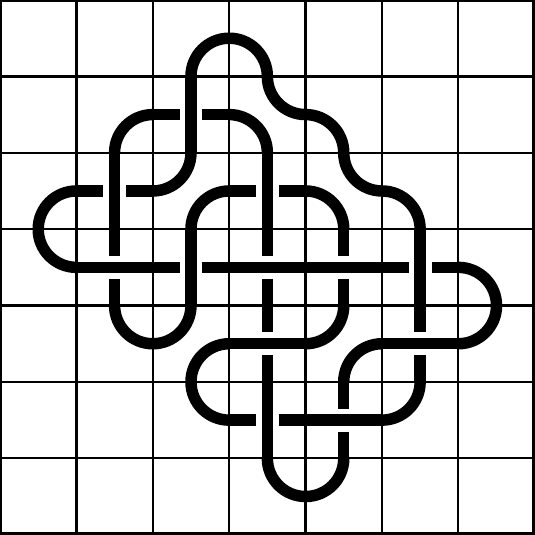}
        \caption*{\phantom{{\Large $\ast$}} \ka{11}{342} {\Large $\ast$}}
    \end{minipage} \newline
\end{figure}
\unskip

\begin{figure}[H]
    \centering
    \begin{minipage}{0.155\linewidth}
        \captionsetup{skip=3pt}
        \centering
        \includegraphics[width=\linewidth]{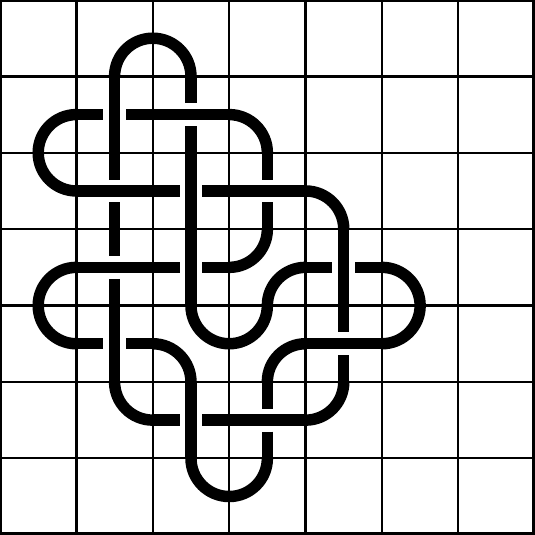}
        \caption*{\phantom{{\Large $\ast$}} \ka{11}{364} {\Large $\ast$}}
    \end{minipage} \hfill
     \begin{minipage}{0.155\linewidth}
        \captionsetup{skip=3pt}
        \centering
        \includegraphics[width=\linewidth]{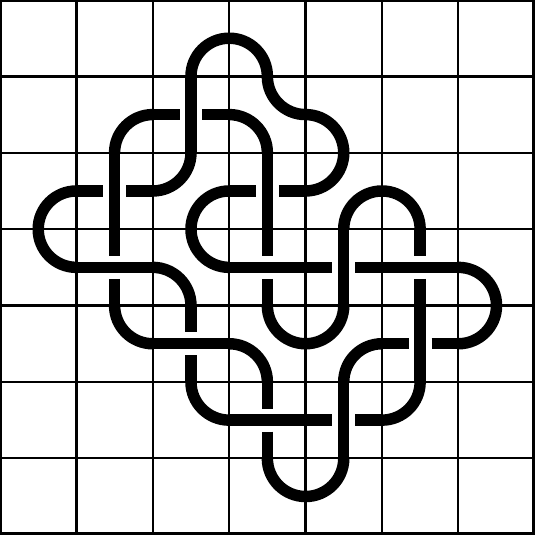}
        \caption*{\ka{11}{367}}
    \end{minipage} \hfill
    \begin{minipage}{0.155\linewidth}
        \captionsetup{skip=3pt}
        \centering
        \includegraphics[width=\linewidth]{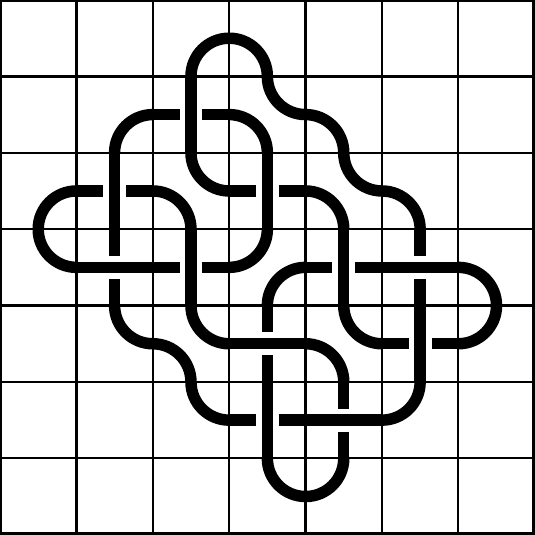}
        \caption*{\kn{11}{71}}
    \end{minipage} \hfill
     \begin{minipage}{0.155\linewidth}
        \captionsetup{skip=3pt}
        \centering
        \includegraphics[width=\linewidth]{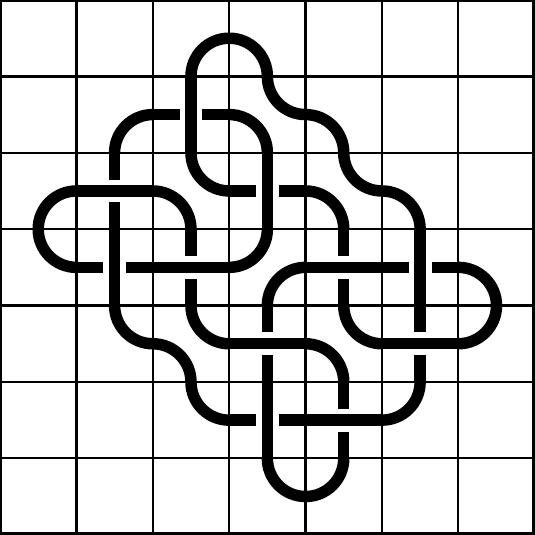}
        \caption*{\kn{11}{72}}
    \end{minipage} \hfill
    \begin{minipage}{0.155\linewidth}
        \captionsetup{skip=3pt}
        \centering
        \includegraphics[width=\linewidth]{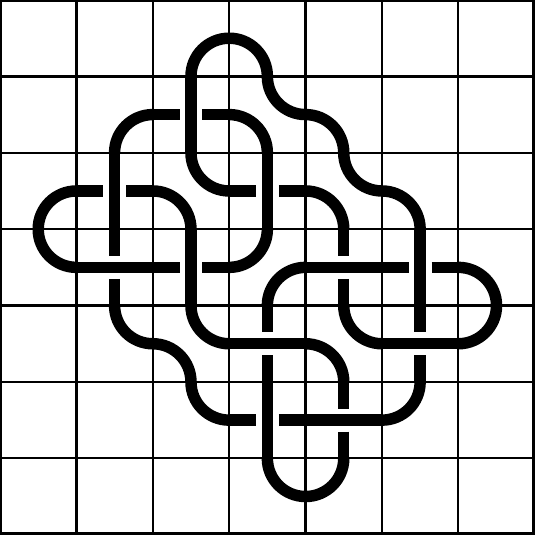}
        \caption*{\kn{11}{73}}
    \end{minipage}  \hfill
    \begin{minipage}{0.155\linewidth}
        \captionsetup{skip=3pt}
        \centering
        \includegraphics[width=\linewidth]{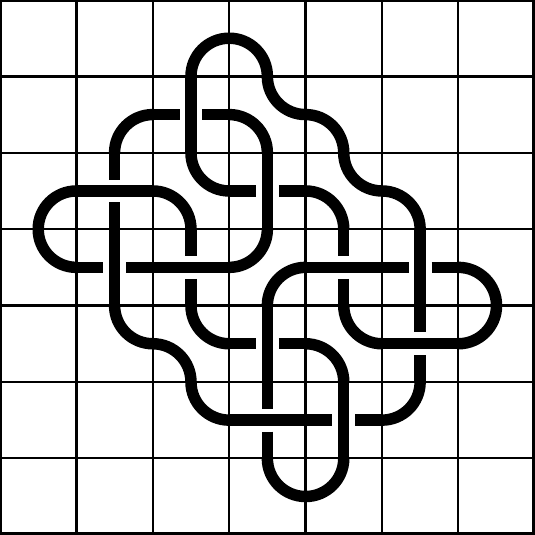}
        \caption*{\kn{11}{74}}
    \end{minipage} \newline
\end{figure}
\unskip

\begin{figure}[H]
    \centering
    \begin{minipage}{0.155\linewidth}
        \captionsetup{skip=3pt}
        \centering
        \includegraphics[width=\linewidth]{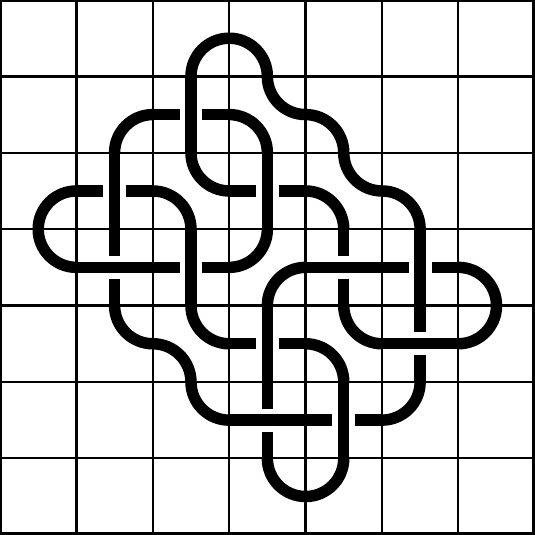}
        \caption*{\kn{11}{75}}
    \end{minipage} \hfill
     \begin{minipage}{0.155\linewidth}
        \captionsetup{skip=3pt}
        \centering
        \includegraphics[width=\linewidth]{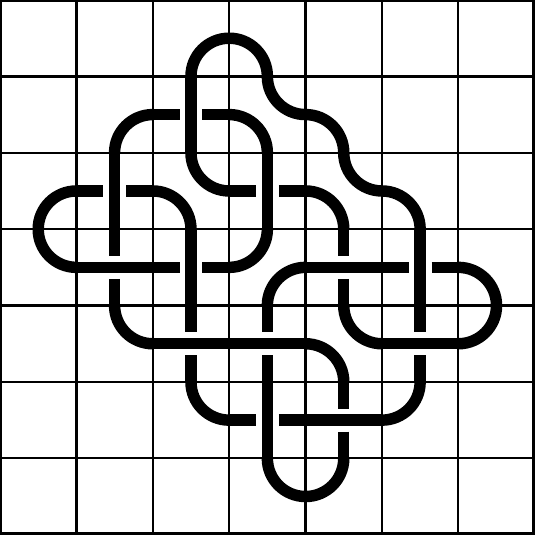}
        \caption*{\phantom{{\Large $\ast$}} \kn{11}{76} {\Large $\ast$}}
    \end{minipage} \hfill
    \begin{minipage}{0.155\linewidth}
        \captionsetup{skip=3pt}
        \centering
        \includegraphics[width=\linewidth]{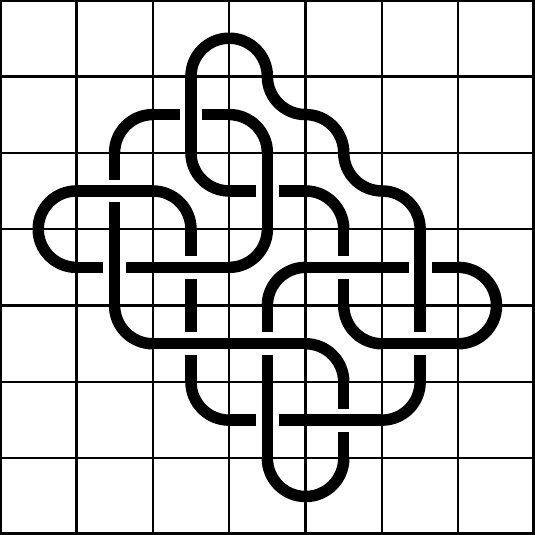}
        \caption*{\phantom{{\Large $\ast$}} \kn{11}{77} {\Large $\ast$}}
    \end{minipage} \hfill
     \begin{minipage}{0.155\linewidth}
        \captionsetup{skip=3pt}
        \centering
        \includegraphics[width=\linewidth]{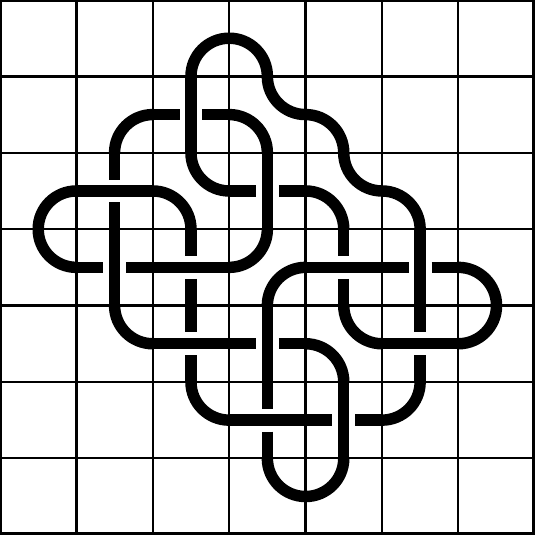}
        \caption*{\phantom{{\Large $\ast$}} \kn{11}{78} {\Large $\ast$}}
    \end{minipage} \hfill
    \begin{minipage}{0.155\linewidth}
        \captionsetup{skip=3pt}
        \centering
        \includegraphics[width=\linewidth]{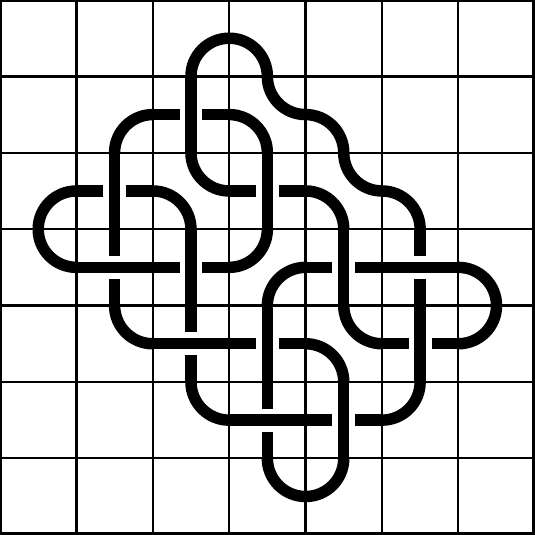}
        \caption*{\ka{12}{119}}
    \end{minipage}  \hfill
    \begin{minipage}{0.155\linewidth}
        \captionsetup{skip=3pt}
        \centering
        \includegraphics[width=\linewidth]{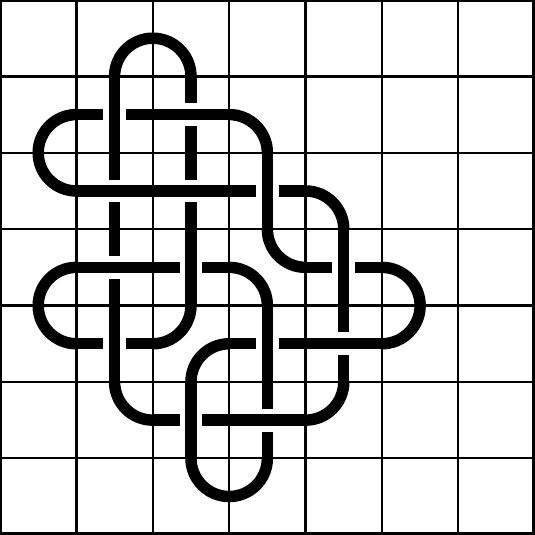}
        \caption*{\phantom{{\Large $\ast$}} \ka{12}{165} {\Large $\ast$}}
    \end{minipage} \newline
\end{figure}
\unskip

\begin{figure}[H]
    \centering
    \begin{minipage}{0.155\linewidth}
        \captionsetup{skip=3pt}
        \centering
        \includegraphics[width=\linewidth]{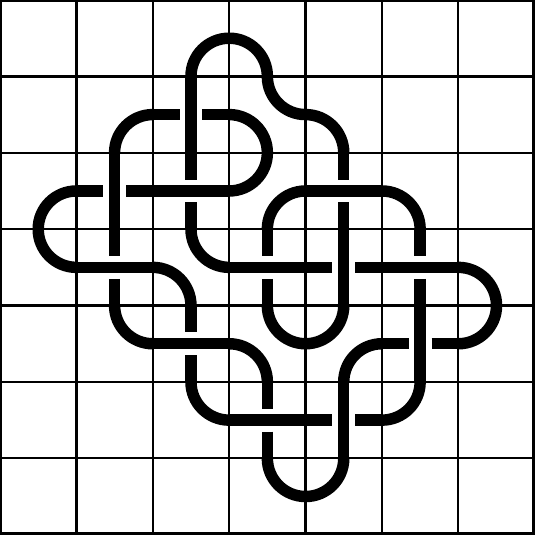}
        \caption*{\ka{12}{169}}
    \end{minipage} \hfill
     \begin{minipage}{0.155\linewidth}
        \captionsetup{skip=3pt}
        \centering
        \includegraphics[width=\linewidth]{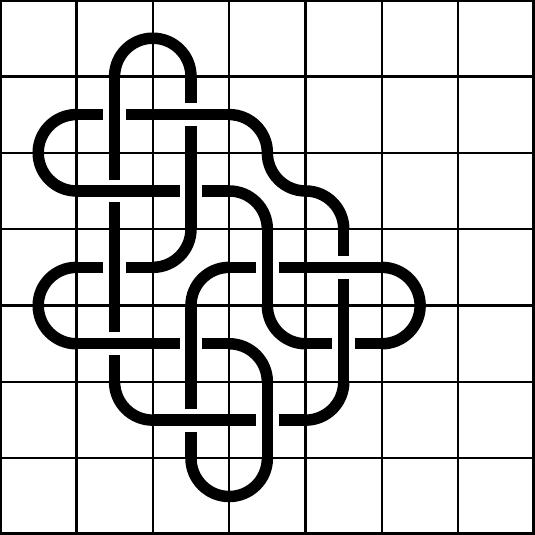}
        \caption*{\ka{12}{373}}
    \end{minipage} \hfill
    \begin{minipage}{0.155\linewidth}
        \captionsetup{skip=3pt}
        \centering
        \includegraphics[width=\linewidth]{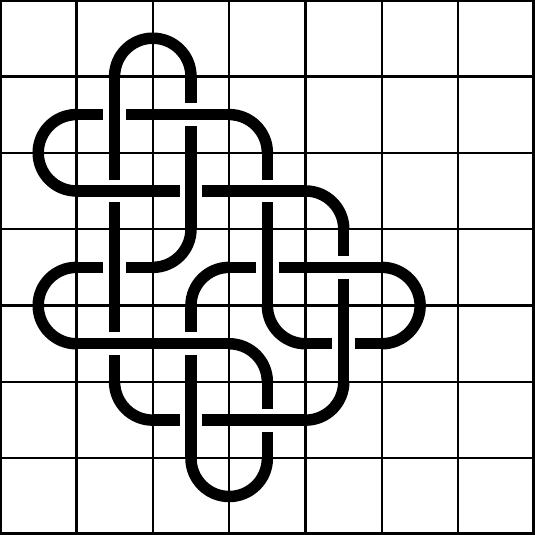}
        \caption*{\phantom{{\Large $\ast$}} \ka{12}{376} {\Large $\ast$}}
    \end{minipage} \hfill
     \begin{minipage}{0.155\linewidth}
        \captionsetup{skip=3pt}
        \centering
        \includegraphics[width=\linewidth]{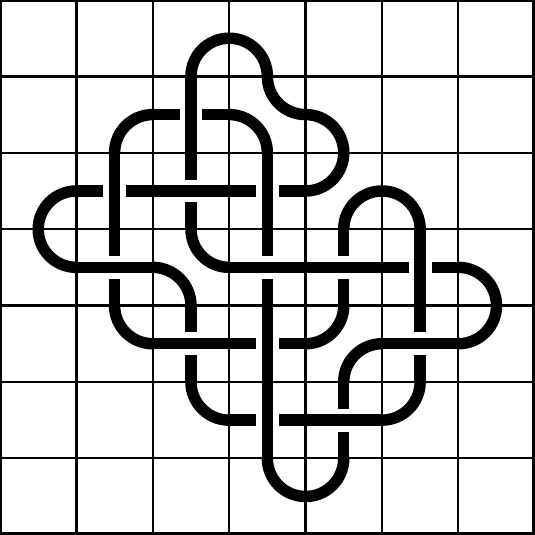}
        \caption*{\phantom{{\Large $\ast$}} \ka{12}{379} {\Large $\ast$}}
    \end{minipage} \hfill
    \begin{minipage}{0.155\linewidth}
        \captionsetup{skip=3pt}
        \centering
        \includegraphics[width=\linewidth]{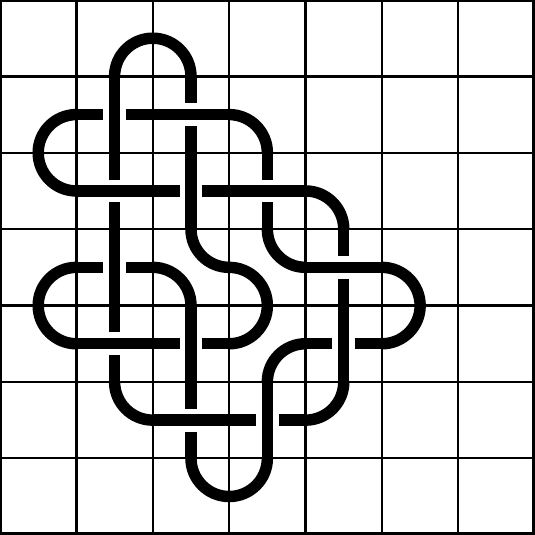}
        \caption*{\ka{12}{380}}
    \end{minipage}  \hfill
    \begin{minipage}{0.155\linewidth}
        \captionsetup{skip=3pt}
        \centering
        \includegraphics[width=\linewidth]{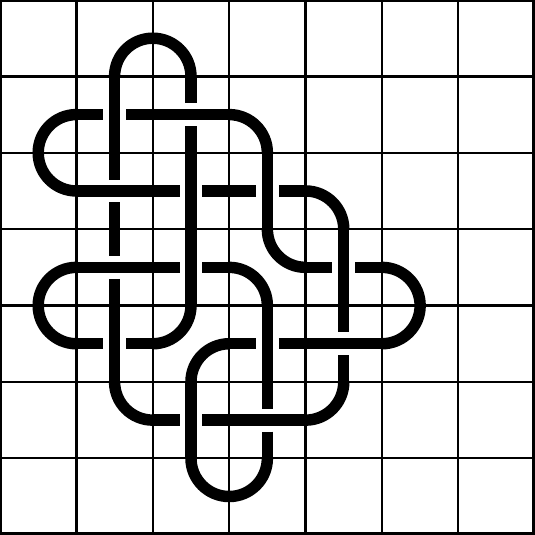}
        \caption*{\phantom{{\Large $\ast$}} \ka{12}{444} {\Large $\ast$}}
    \end{minipage} \newline
\end{figure}
\unskip

\begin{figure}[H]
    \centering
    \begin{minipage}{0.155\linewidth}
        \captionsetup{skip=3pt}
        \centering
        \includegraphics[width=\linewidth]{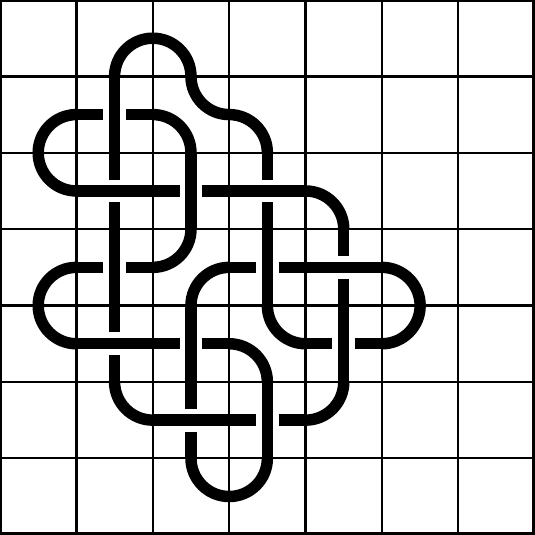}
        \caption*{\ka{12}{503}}
    \end{minipage} \hfill
     \begin{minipage}{0.155\linewidth}
        \captionsetup{skip=3pt}
        \centering
        \includegraphics[width=\linewidth]{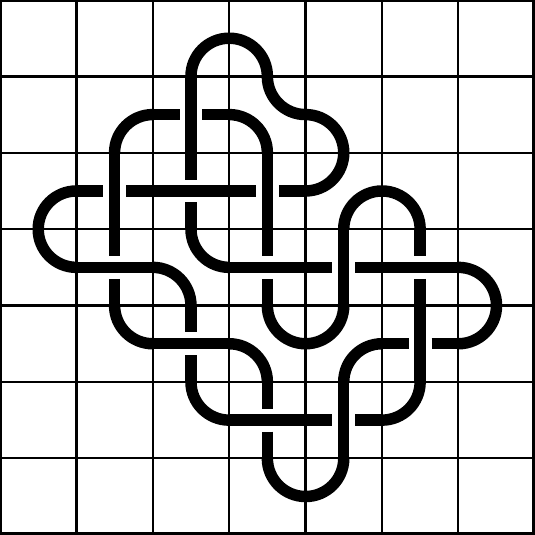}
        \caption*{\ka{12}{722}}
    \end{minipage} \hfill
    \begin{minipage}{0.155\linewidth}
        \captionsetup{skip=3pt}
        \centering
        \includegraphics[width=\linewidth]{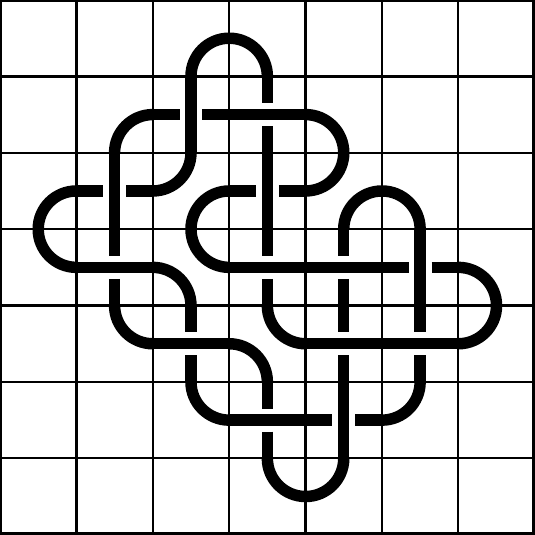}
        \caption*{\phantom{{\Large $\ast$}} \ka{12}{803} {\Large $\ast$}}
    \end{minipage} \hfill
     \begin{minipage}{0.155\linewidth}
        \captionsetup{skip=3pt}
        \centering
        \includegraphics[width=\linewidth]{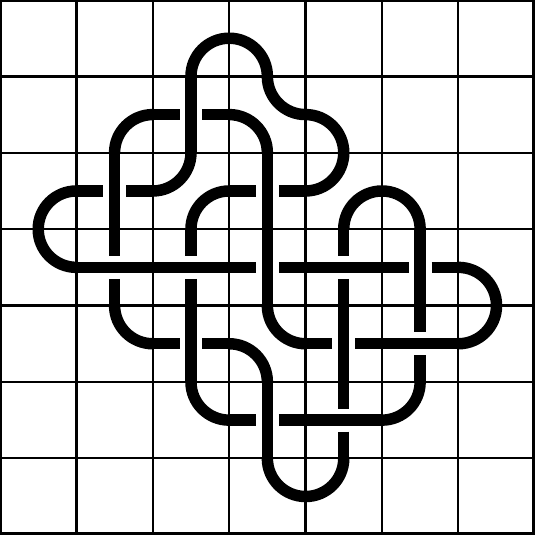}
        \caption*{\phantom{{\Large $\ast$}} \ka{12}{1148} {\Large $\ast$}}
    \end{minipage} \hfill
    \begin{minipage}{0.155\linewidth}
        \captionsetup{skip=3pt}
        \centering
        \includegraphics[width=\linewidth]{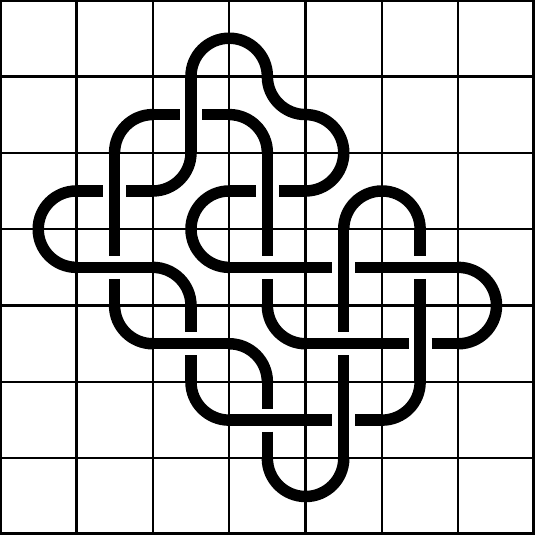}
        \caption*{\ka{12}{1149}}
    \end{minipage}  \hfill
    \begin{minipage}{0.155\linewidth}
        \captionsetup{skip=3pt}
        \centering
        \includegraphics[width=\linewidth]{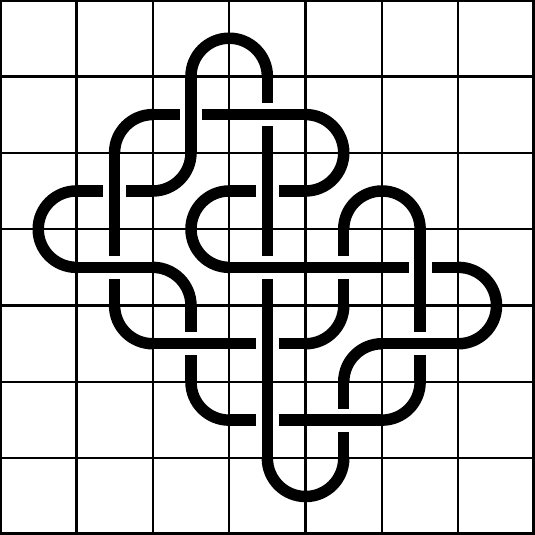}
        \caption*{\phantom{{\Large $\ast$}} \ka{12}{1166} {\Large $\ast$}}
    \end{minipage} \newline
\end{figure}
\unskip

\begin{figure}[H]
    \centering
    \begin{minipage}{0.155\linewidth}
        \captionsetup{skip=3pt}
        \centering
        \includegraphics[width=\linewidth]{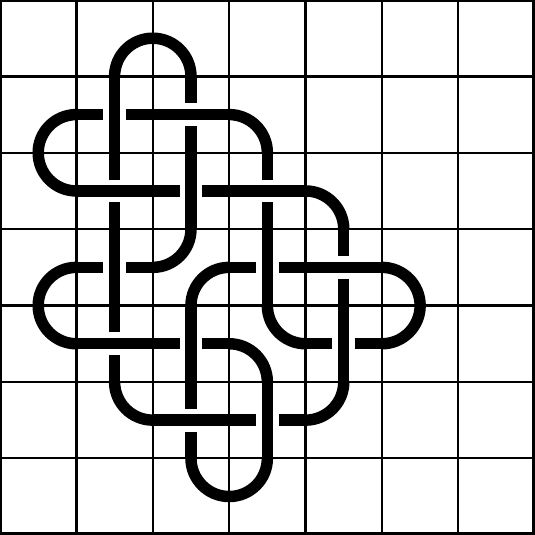}
        \caption*{\ka{13}{1236}}
    \end{minipage}
     \begin{minipage}{0.155\linewidth}
        \captionsetup{skip=3pt}
        \centering
        \includegraphics[width=\linewidth]{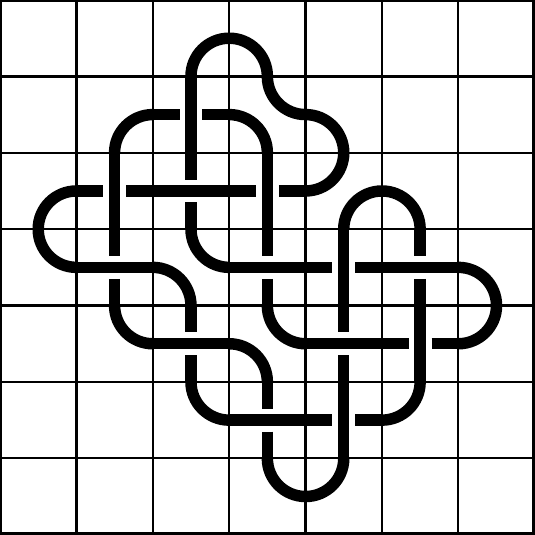}
        \caption*{\ka{13}{1461}}
    \end{minipage}
    \begin{minipage}{0.155\linewidth}
        \captionsetup{skip=3pt}
        \centering
        \includegraphics[width=\linewidth]{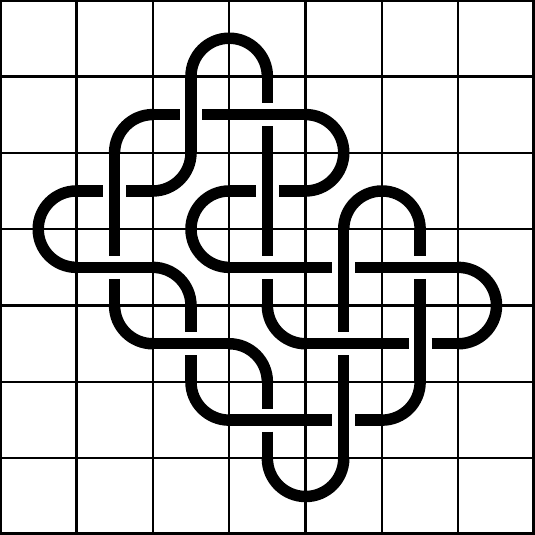}
        \caption*{\ka{13}{4573}}
    \end{minipage}
\end{figure}
\unskip

\subsection{Mosaics from Theorem \ref{thm:CrossingRealized}}

Knot mosaics in which the crossing number is realized have been found for every knot with crossing number 10 or less. Those that were not previously provided on a minimal mosaic or a mosaic in which the tile number was realized are included here.

\begin{figure}[H]
    \centering
     \begin{minipage}{0.155\linewidth}
        \captionsetup{skip=3pt}
        \centering
        \includegraphics[width=\linewidth]{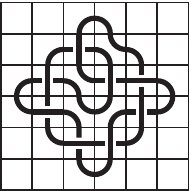}
        \caption*{$8_{1}$}
    \end{minipage} \hfill
    \begin{minipage}{0.155\linewidth}
        \captionsetup{skip=3pt}
        \centering
        \includegraphics[width=\linewidth]{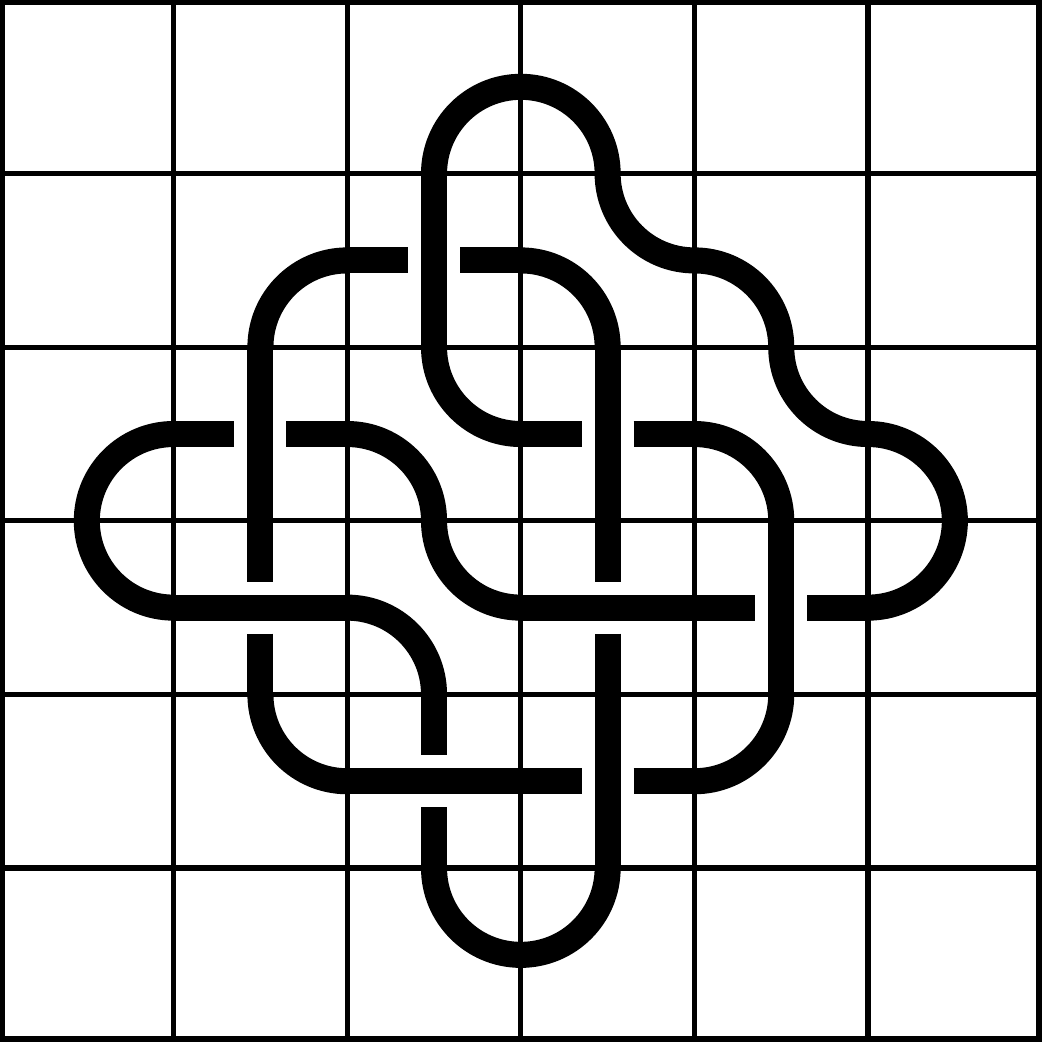}
        \caption*{$8_{7}$}
    \end{minipage} \hfill
    \begin{minipage}{0.155\linewidth}
        \captionsetup{skip=3pt}
        \centering
        \includegraphics[width=\linewidth]{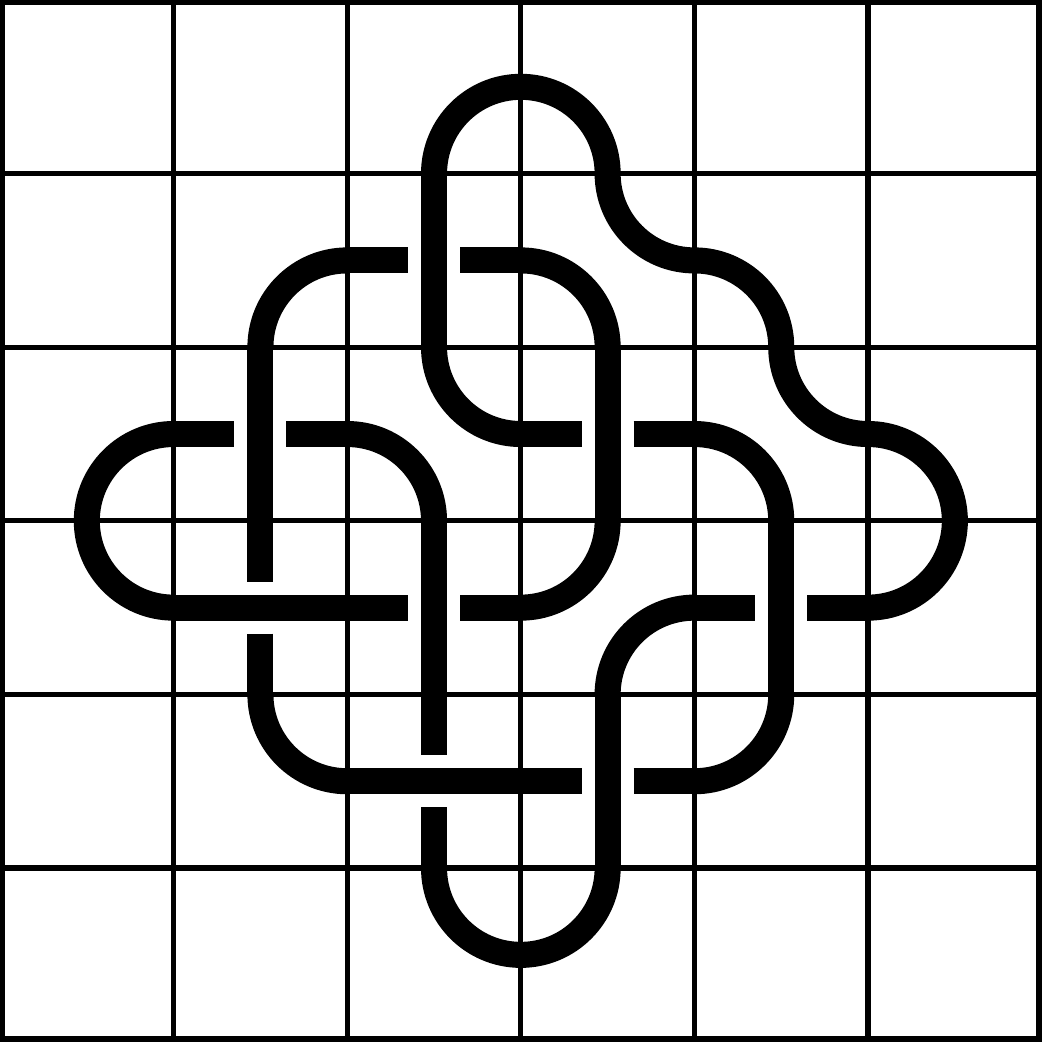}
        \caption*{$8_{8}$}
    \end{minipage}  \hfill
    \begin{minipage}{0.155\linewidth}
        \captionsetup{skip=3pt}
        \centering
        \includegraphics[width=\linewidth]{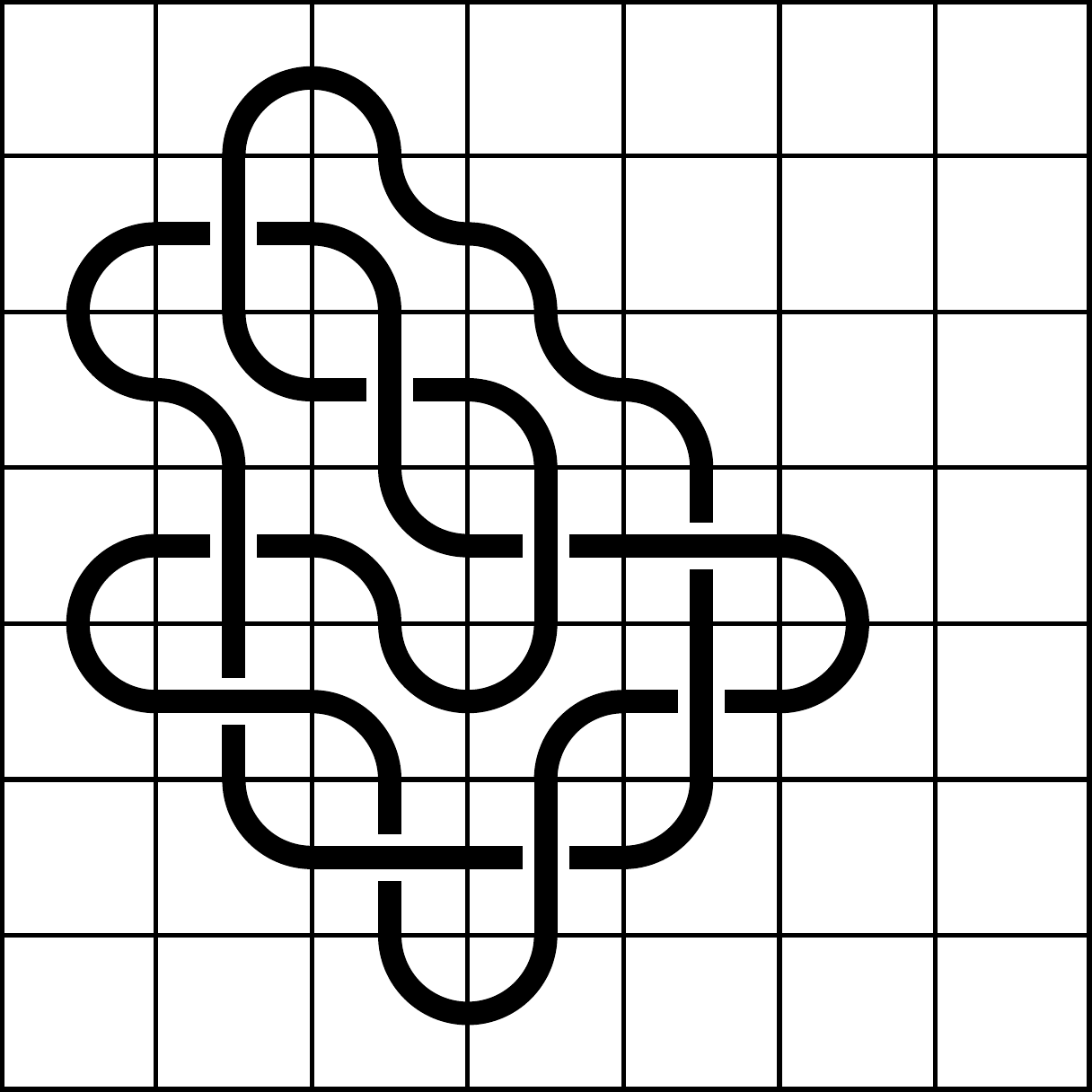}
        \caption*{$9_{3}$}
    \end{minipage} \hfill
    \begin{minipage}{0.155\linewidth}
        \captionsetup{skip=3pt}
        \centering
        \includegraphics[width=\linewidth]{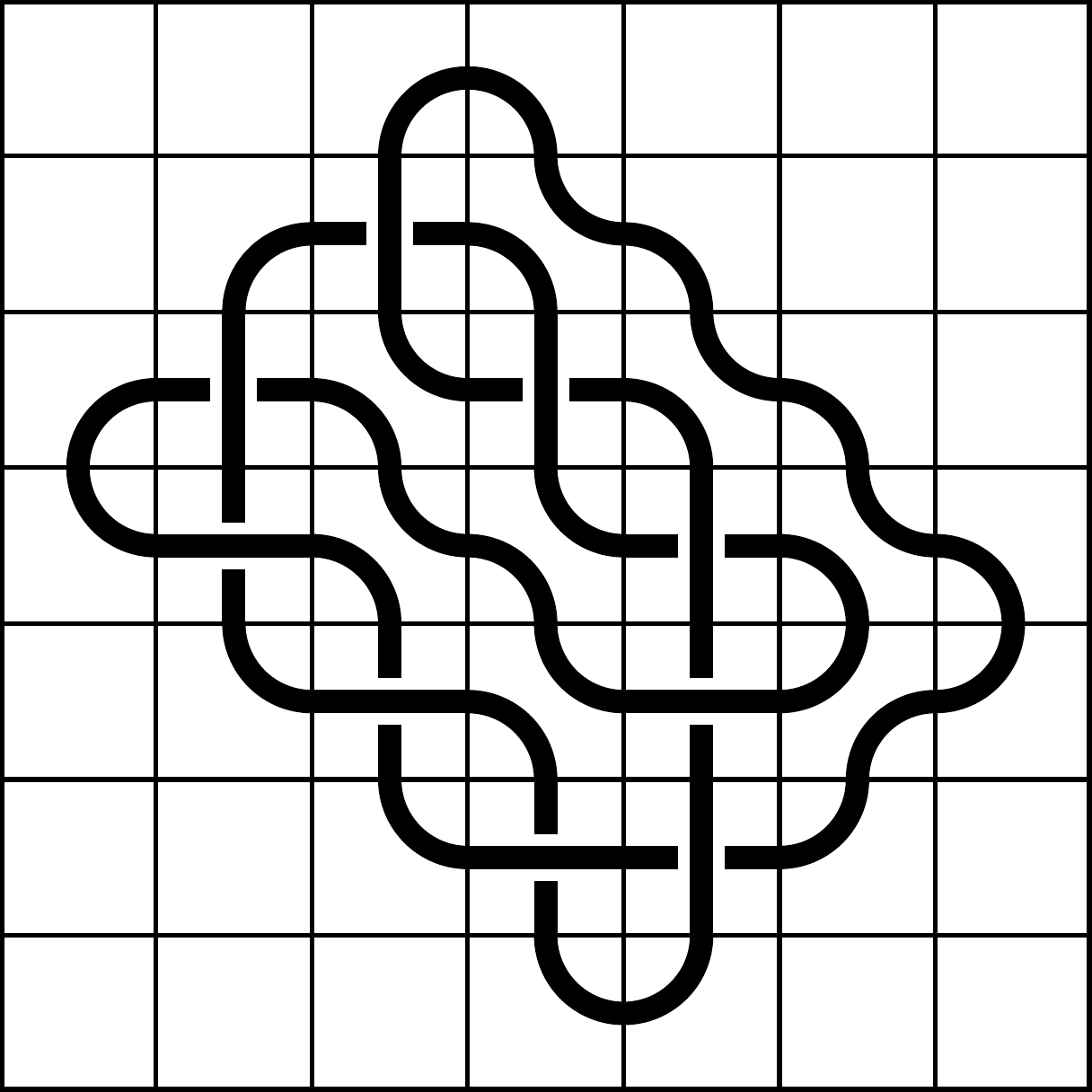}
        \caption*{$9_{4}$}
    \end{minipage} \hfill
    \begin{minipage}{0.155\linewidth}
        \captionsetup{skip=3pt}
        \centering
        \includegraphics[width=\linewidth]{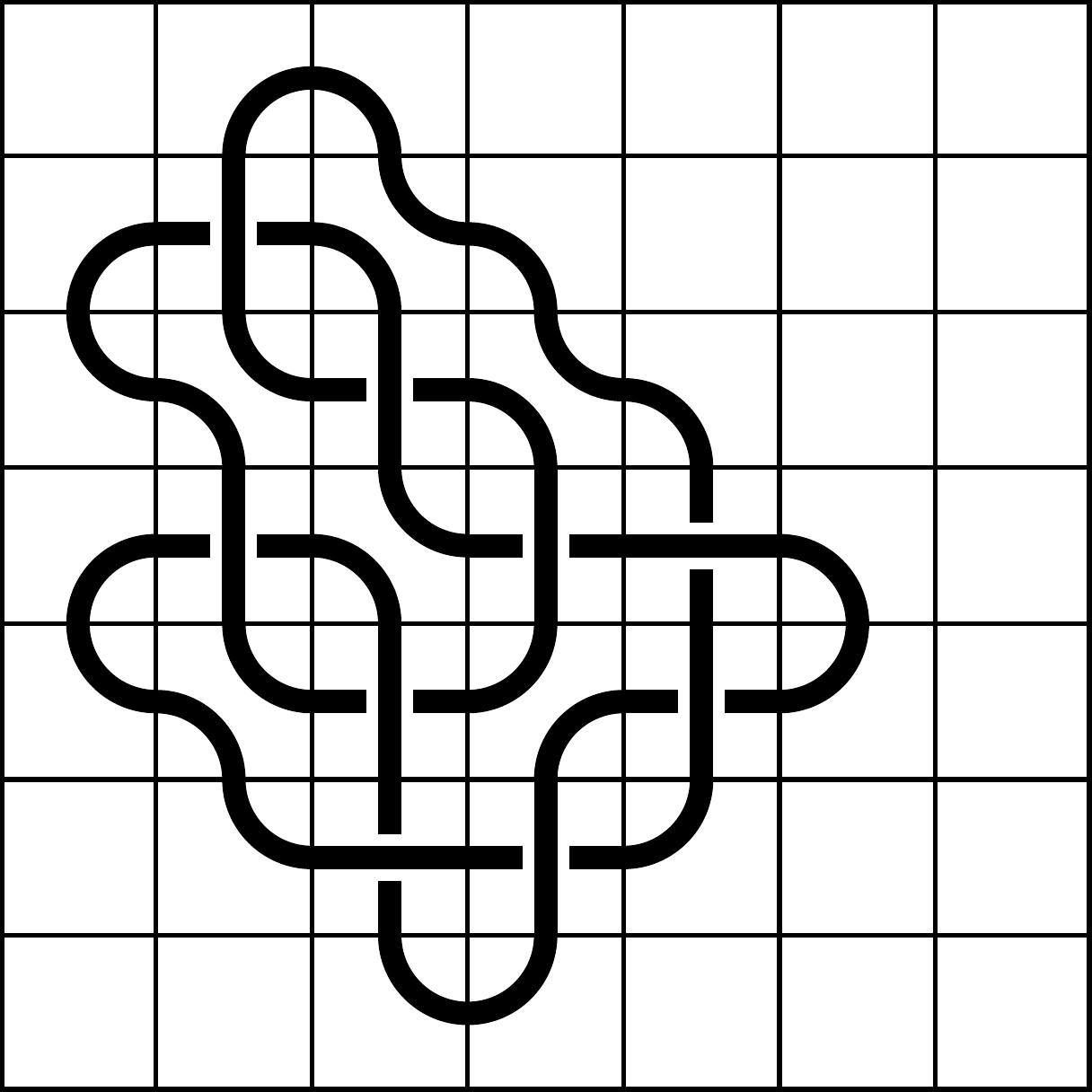}
        \caption*{$9_{7}$ }
    \end{minipage}  \newline
\end{figure}
\unskip

\begin{figure}[H]
    \centering
     \begin{minipage}{0.155\linewidth}
        \captionsetup{skip=3pt}
        \centering
        \includegraphics[width=\linewidth]{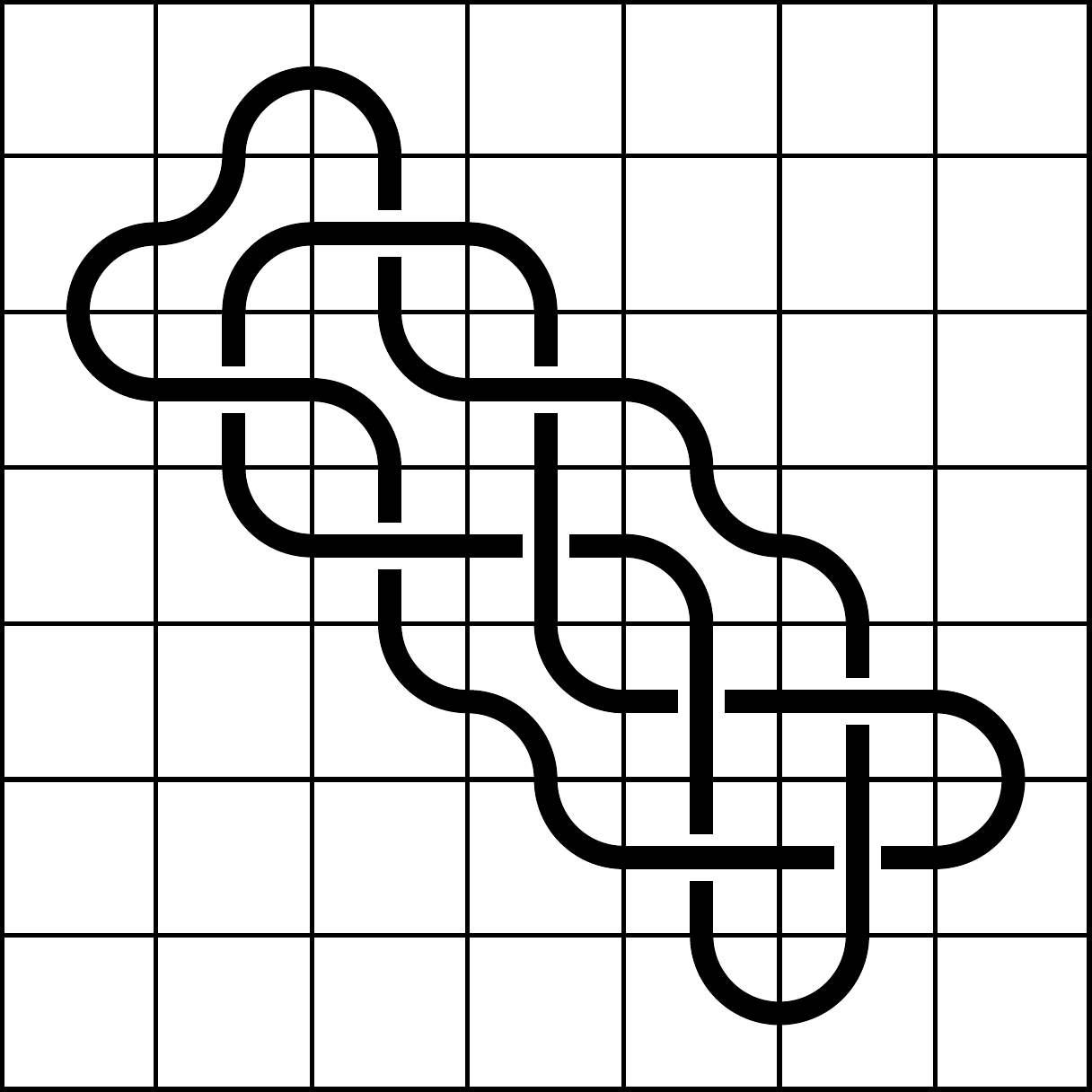}
        \caption*{$9_{9}$}
    \end{minipage} \hfill
    \begin{minipage}{0.155\linewidth}
        \captionsetup{skip=3pt}
        \centering
        \includegraphics[width=\linewidth]{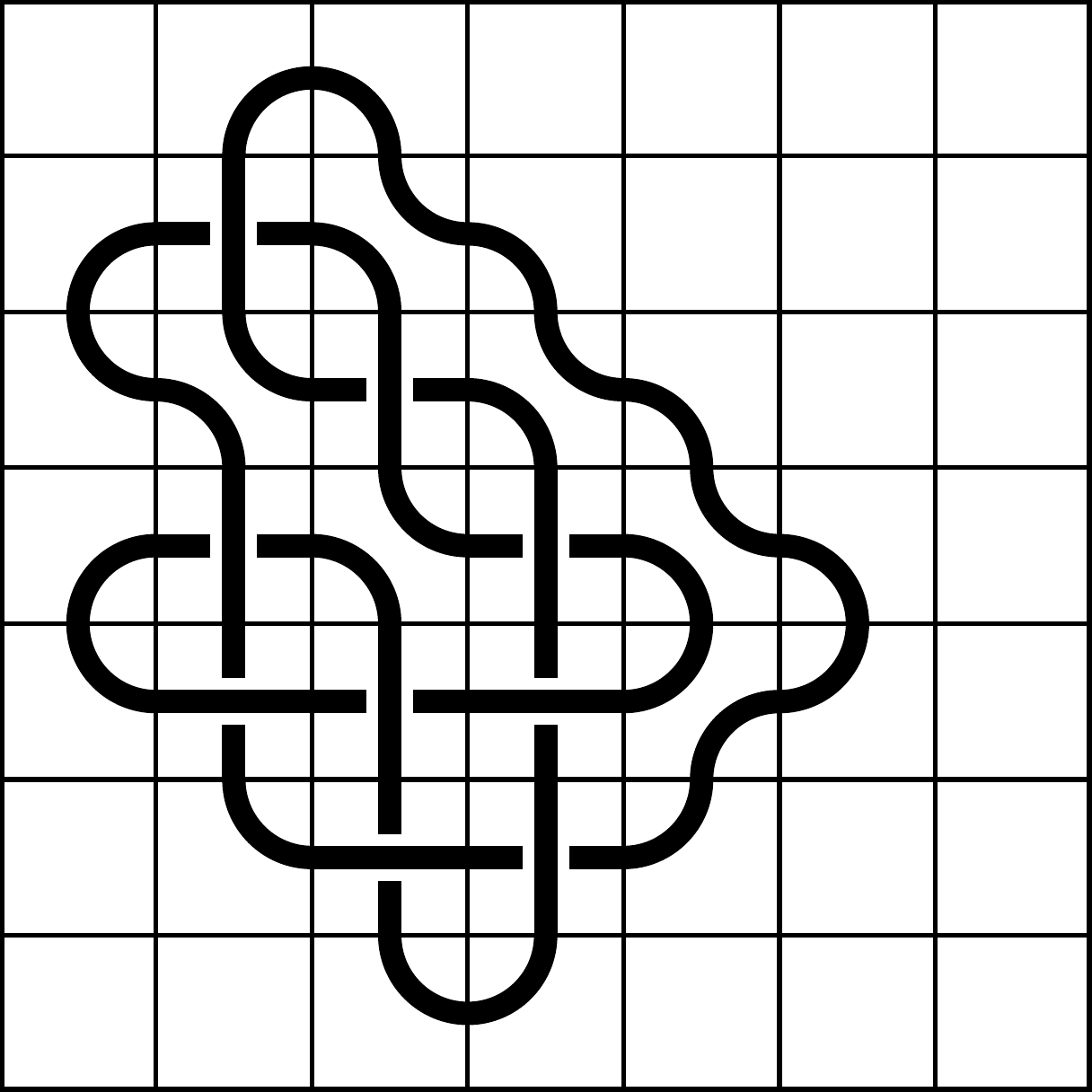}
        \caption*{$9_{12}$}
    \end{minipage} \hfill
    \begin{minipage}{0.155\linewidth}
        \captionsetup{skip=3pt}
        \centering
        \includegraphics[width=\linewidth]{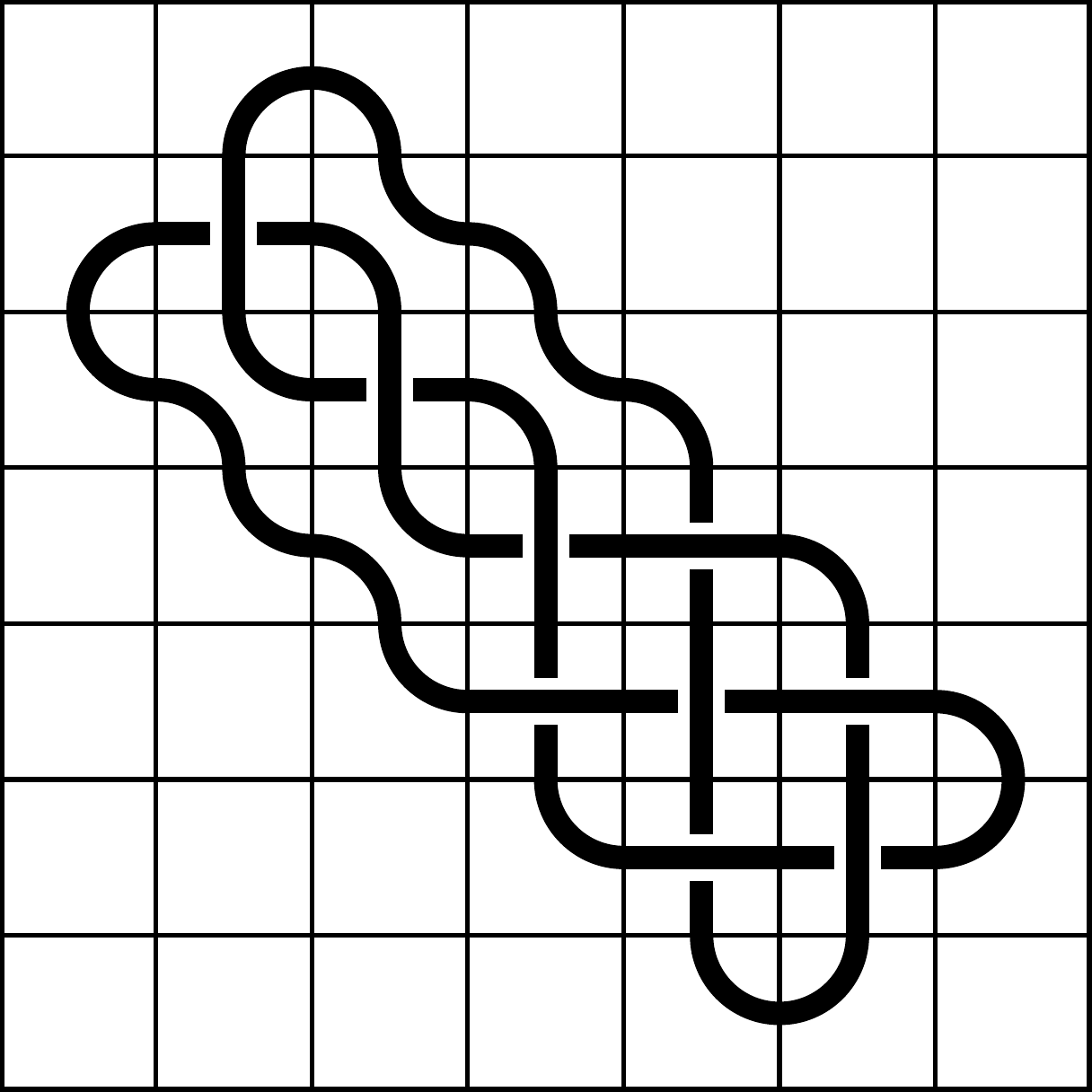}
        \caption*{$9_{13}$}
    \end{minipage}  \hfill
    \begin{minipage}{0.155\linewidth}
        \captionsetup{skip=3pt}
        \centering
        \includegraphics[width=\linewidth]{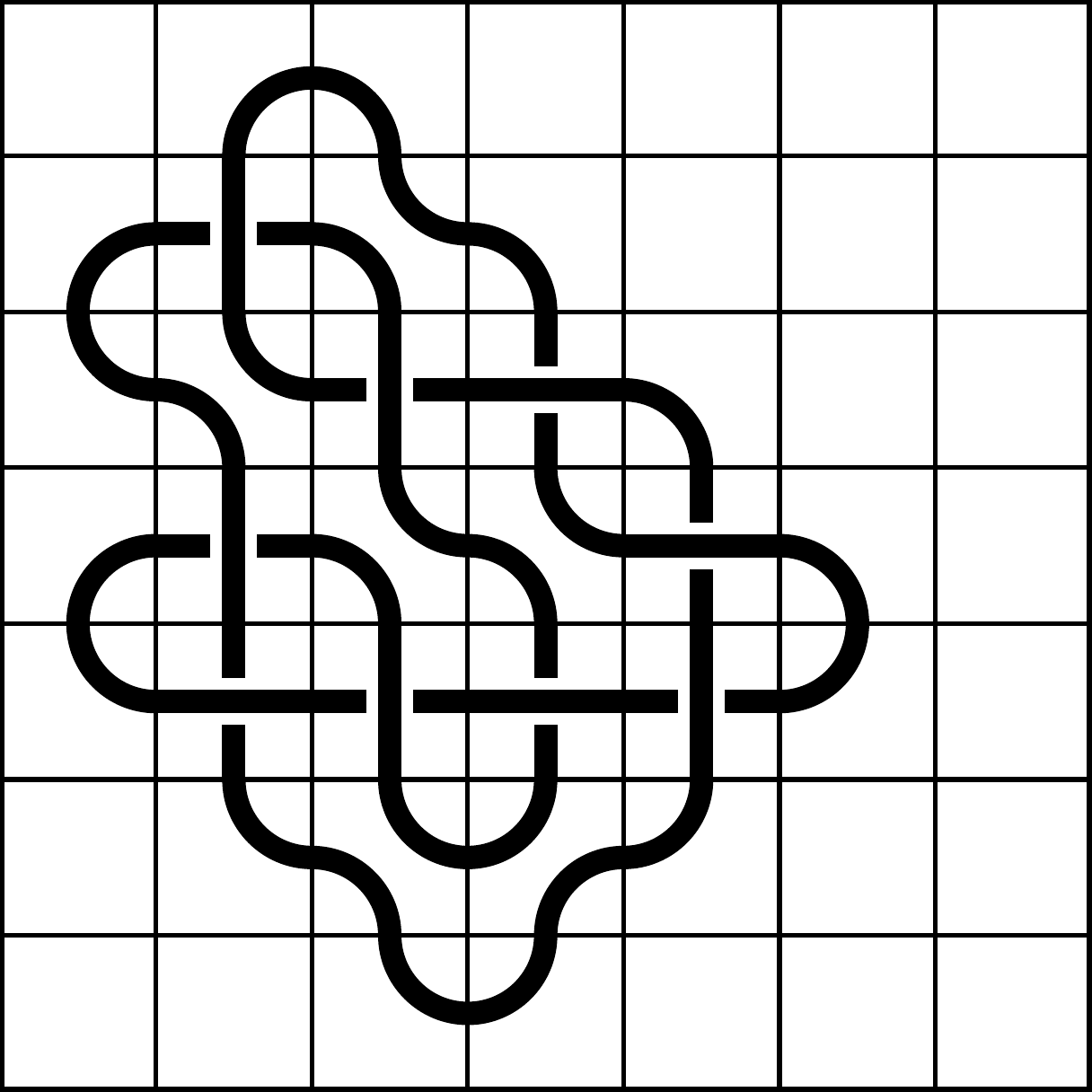}
        \caption*{$9_{15}$}
    \end{minipage} \hfill
    \begin{minipage}{0.155\linewidth}
        \captionsetup{skip=3pt}
        \centering
        \includegraphics[width=\linewidth]{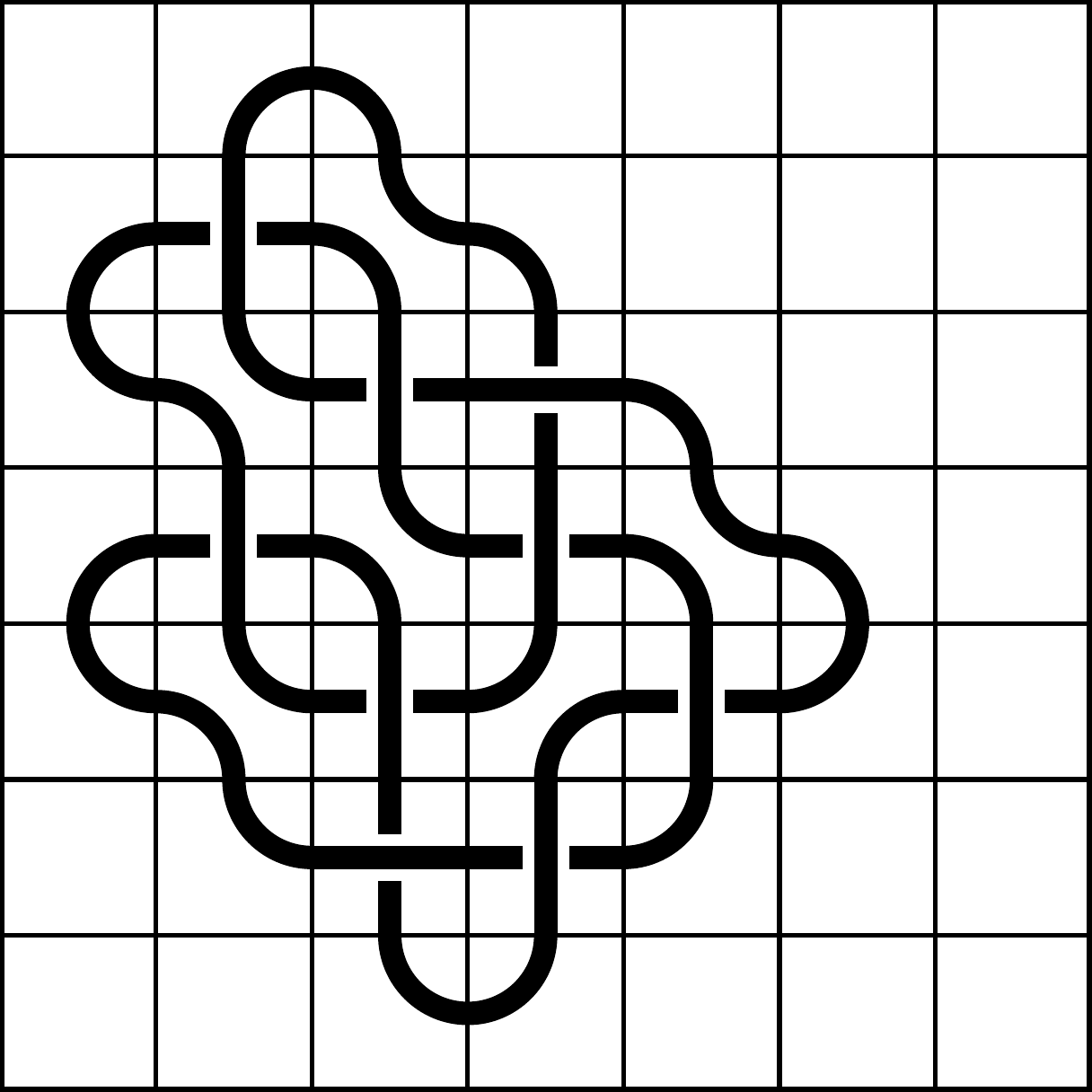}
        \caption*{$9_{19}$}
    \end{minipage} \hfill
    \begin{minipage}{0.155\linewidth}
        \captionsetup{skip=3pt}
        \centering
        \includegraphics[width=\linewidth]{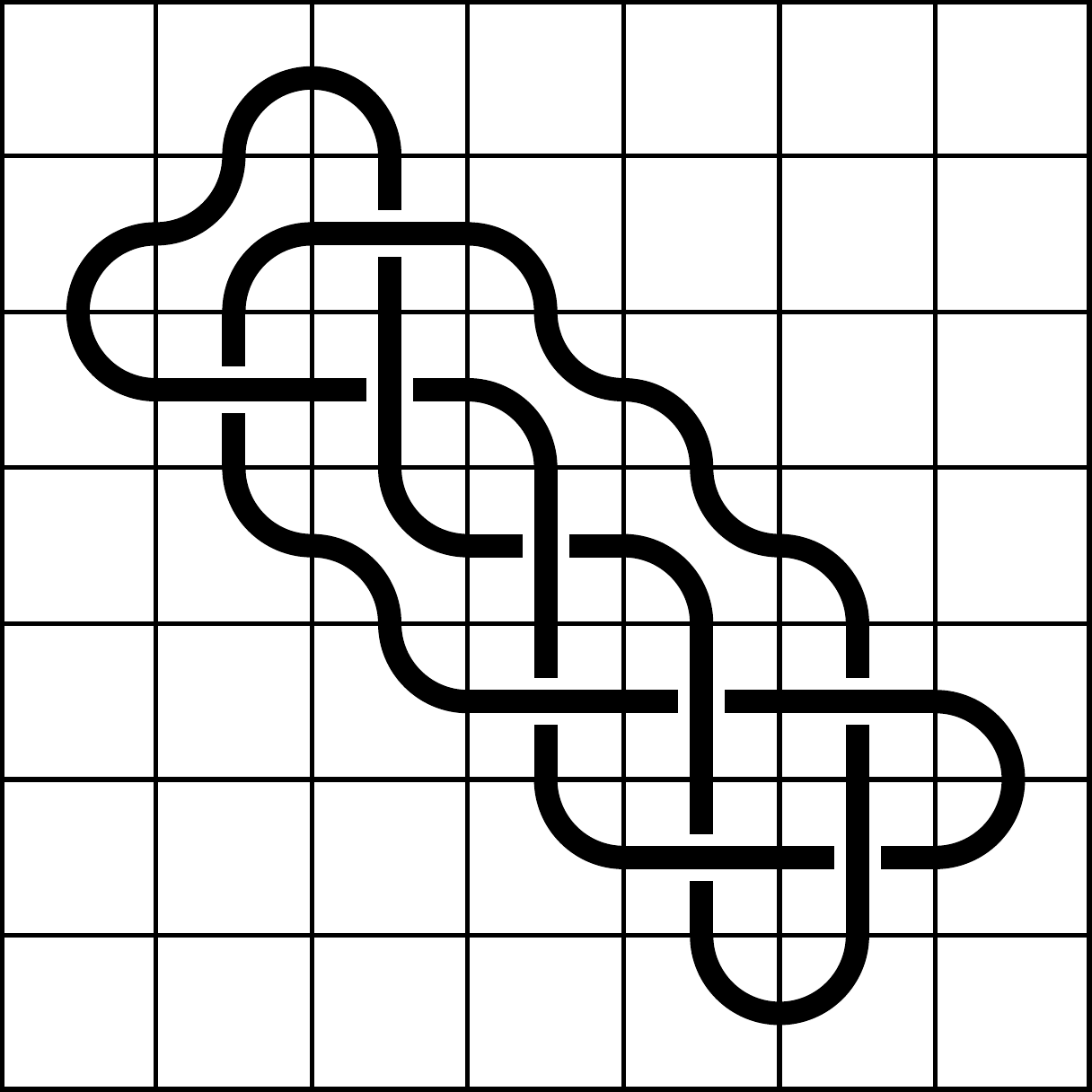}
        \caption*{$9_{21}$ }
    \end{minipage}  \newline
\end{figure}
\unskip

\begin{figure}[H]
    \centering
     \begin{minipage}{0.155\linewidth}
        \captionsetup{skip=3pt}
        \centering
        \includegraphics[width=\linewidth]{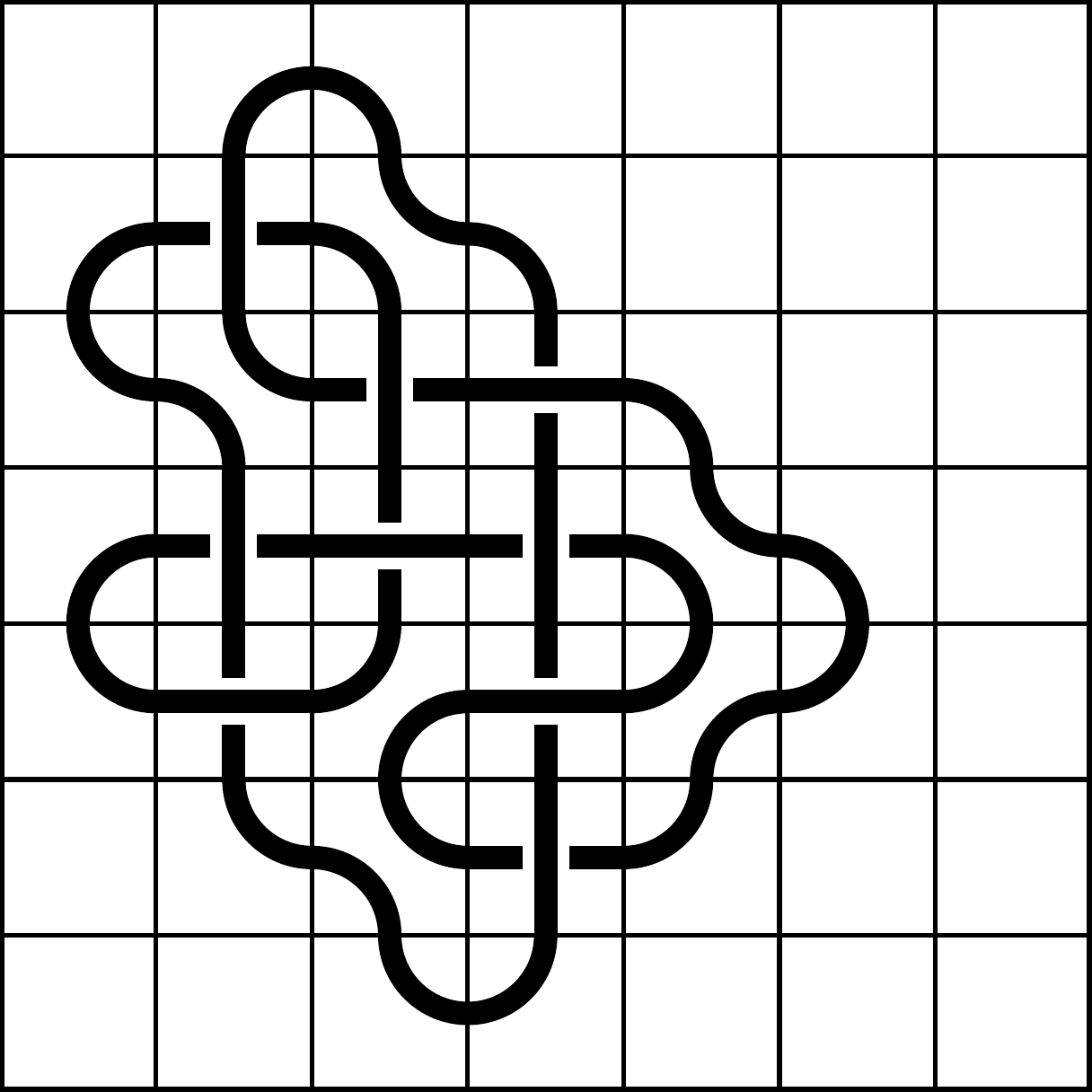}
        \caption*{$9_{24}$}
    \end{minipage} \hfill
    \begin{minipage}{0.155\linewidth}
        \captionsetup{skip=3pt}
        \centering
        \includegraphics[width=\linewidth]{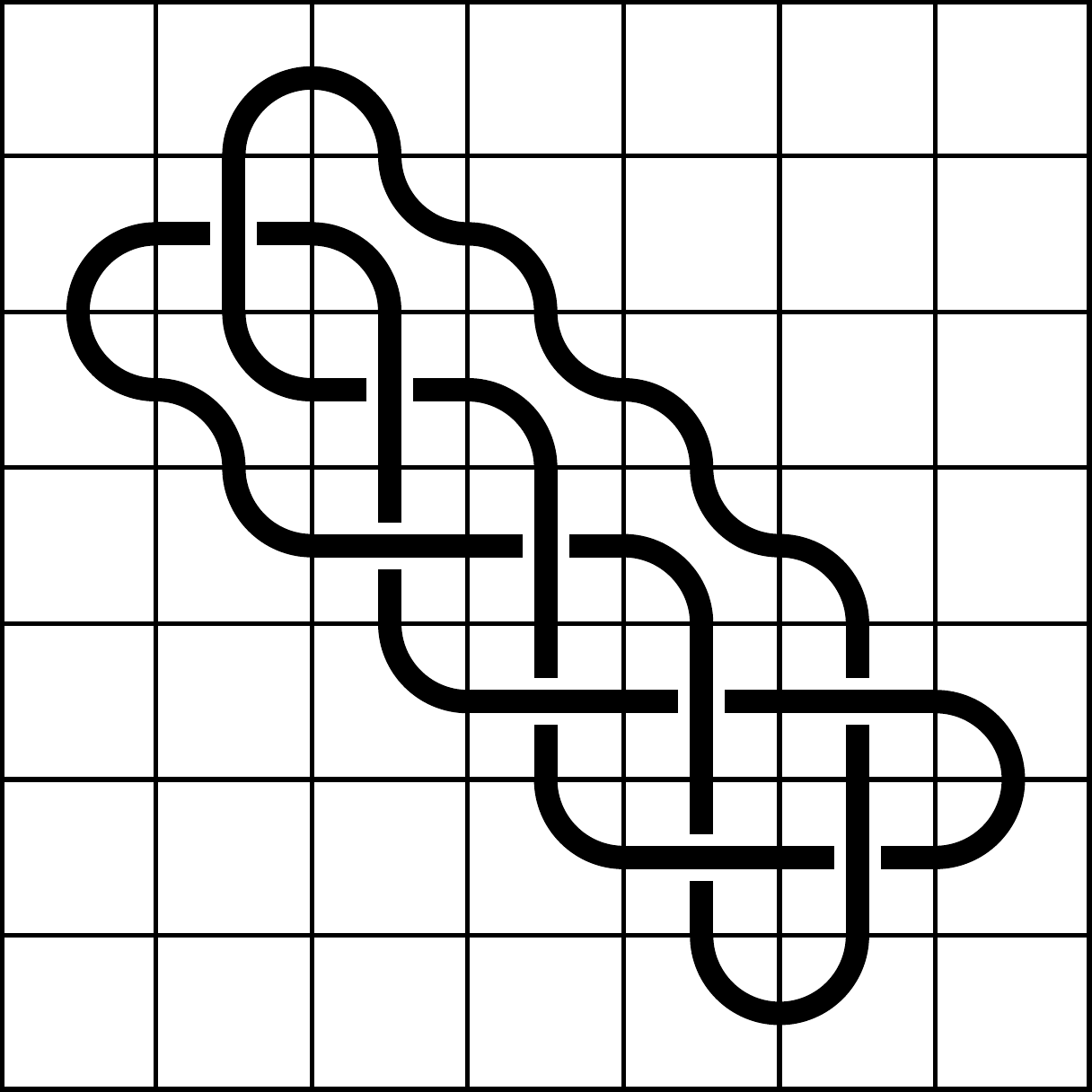}
        \caption*{$9_{26}$}
    \end{minipage} \hfill
    \begin{minipage}{0.155\linewidth}
        \captionsetup{skip=3pt}
        \centering
        \includegraphics[width=\linewidth]{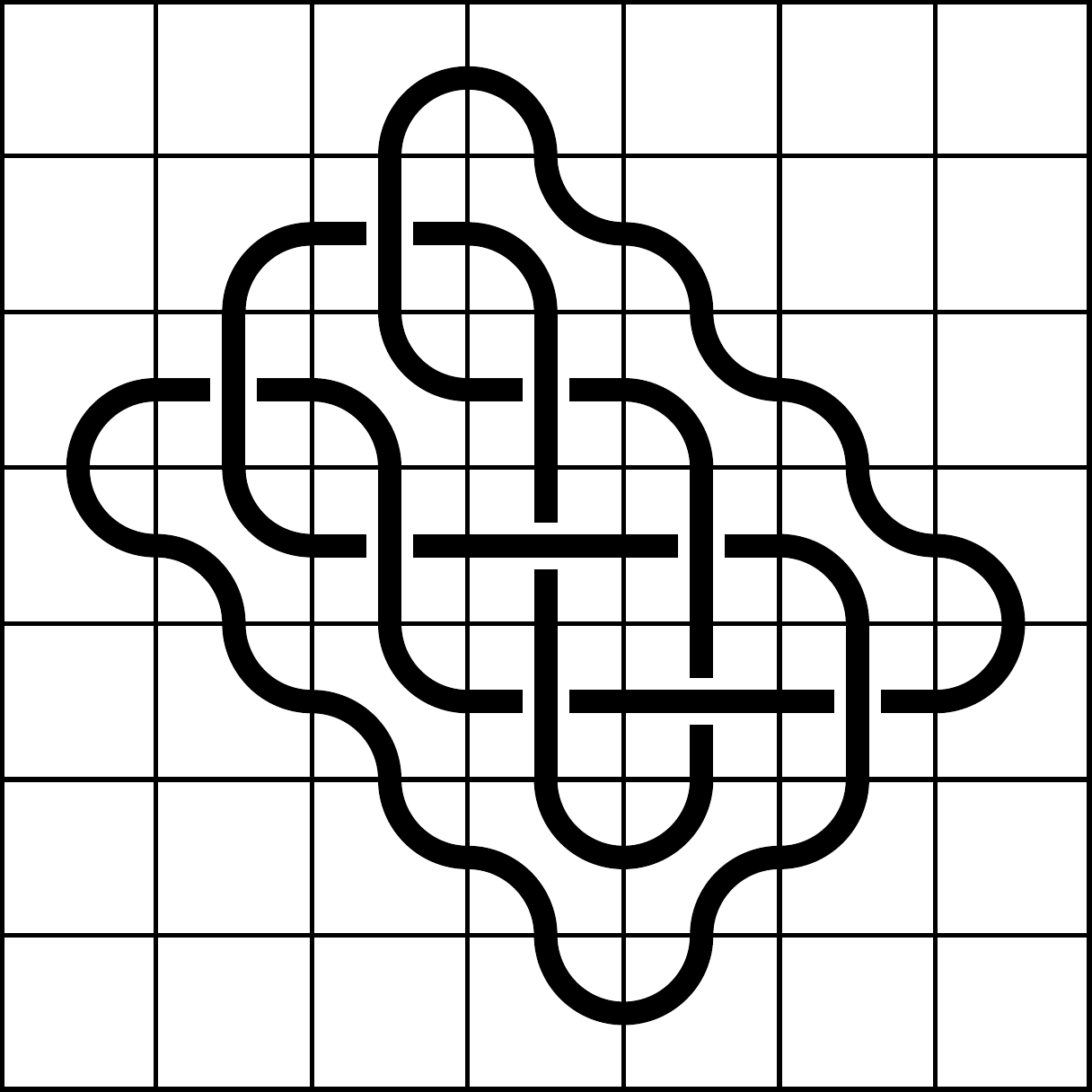}
        \caption*{$9_{29}$}
    \end{minipage}  \hfill
    \begin{minipage}{0.155\linewidth}
        \captionsetup{skip=3pt}
        \centering
        \includegraphics[width=\linewidth]{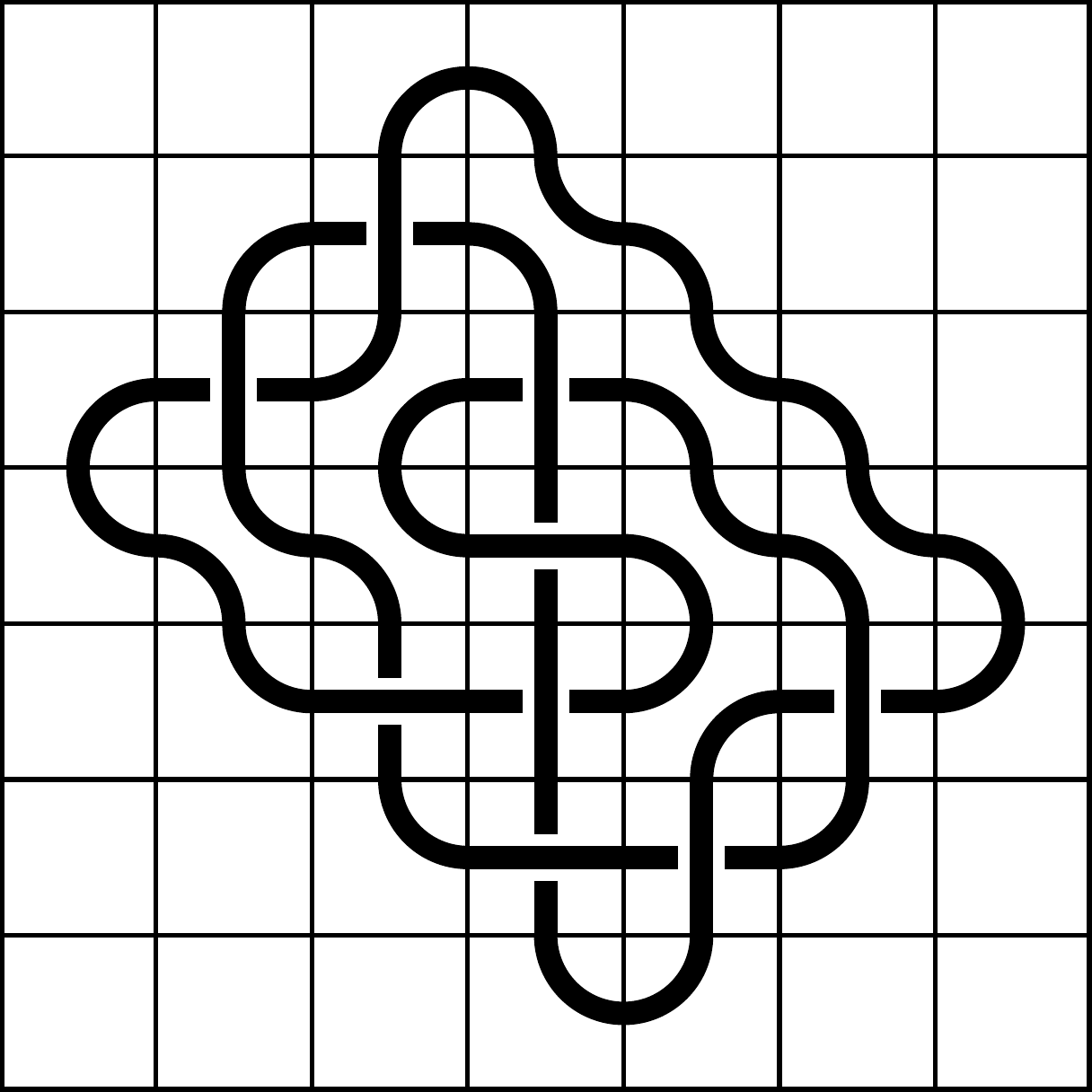}
        \caption*{$9_{35}$}
    \end{minipage} \hfill
    \begin{minipage}{0.155\linewidth}
        \captionsetup{skip=3pt}
        \centering
        \includegraphics[width=\linewidth]{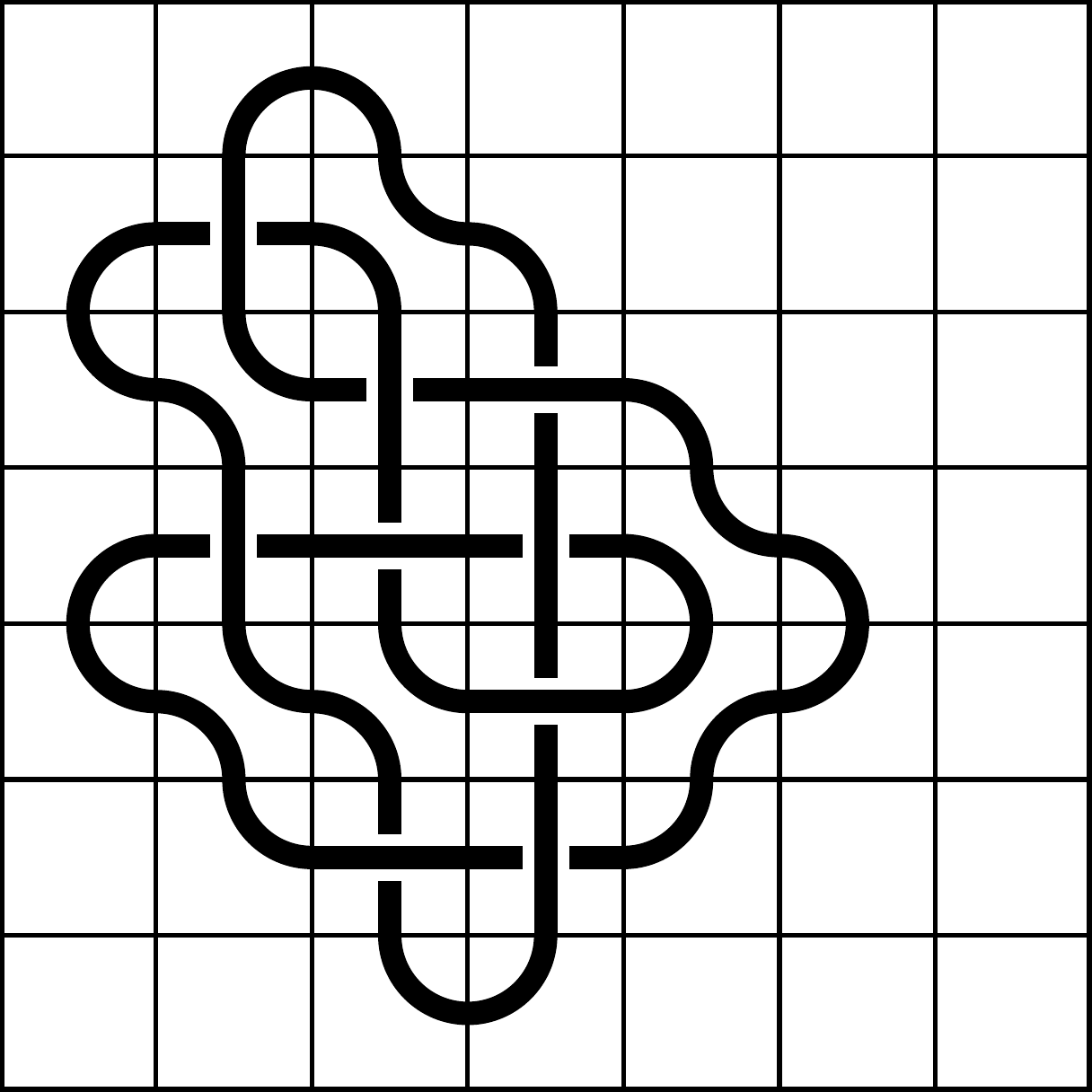}
        \caption*{$9_{37}$}
    \end{minipage} \hfill
    \begin{minipage}{0.155\linewidth}
        \captionsetup{skip=3pt}
        \centering
        \includegraphics[width=\linewidth]{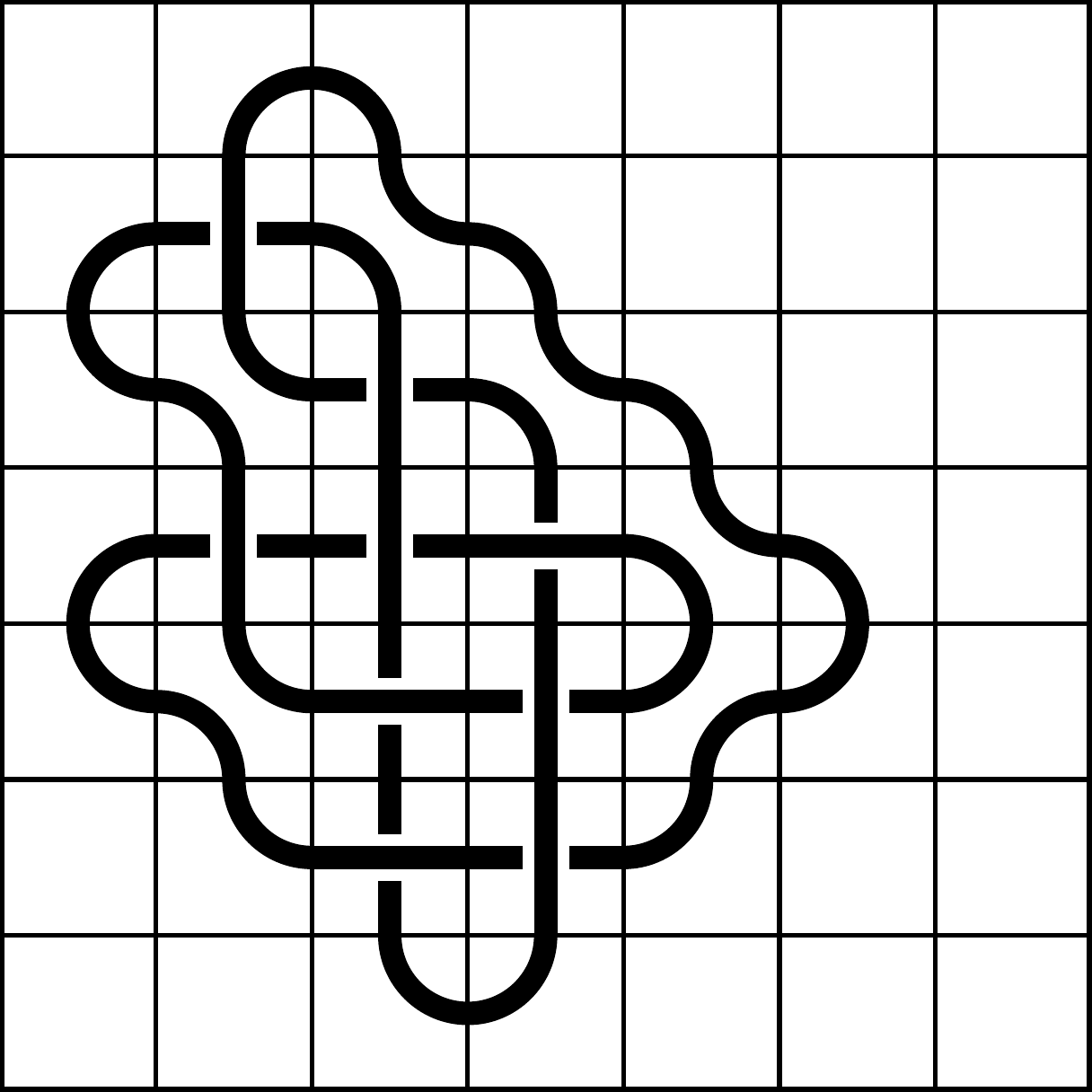}
        \caption*{$9_{46}$ }
    \end{minipage}  \newline
\end{figure}
\unskip

\begin{figure}[H]
    \centering
     \begin{minipage}{0.155\linewidth}
        \captionsetup{skip=3pt}
        \centering
        \includegraphics[width=\linewidth]{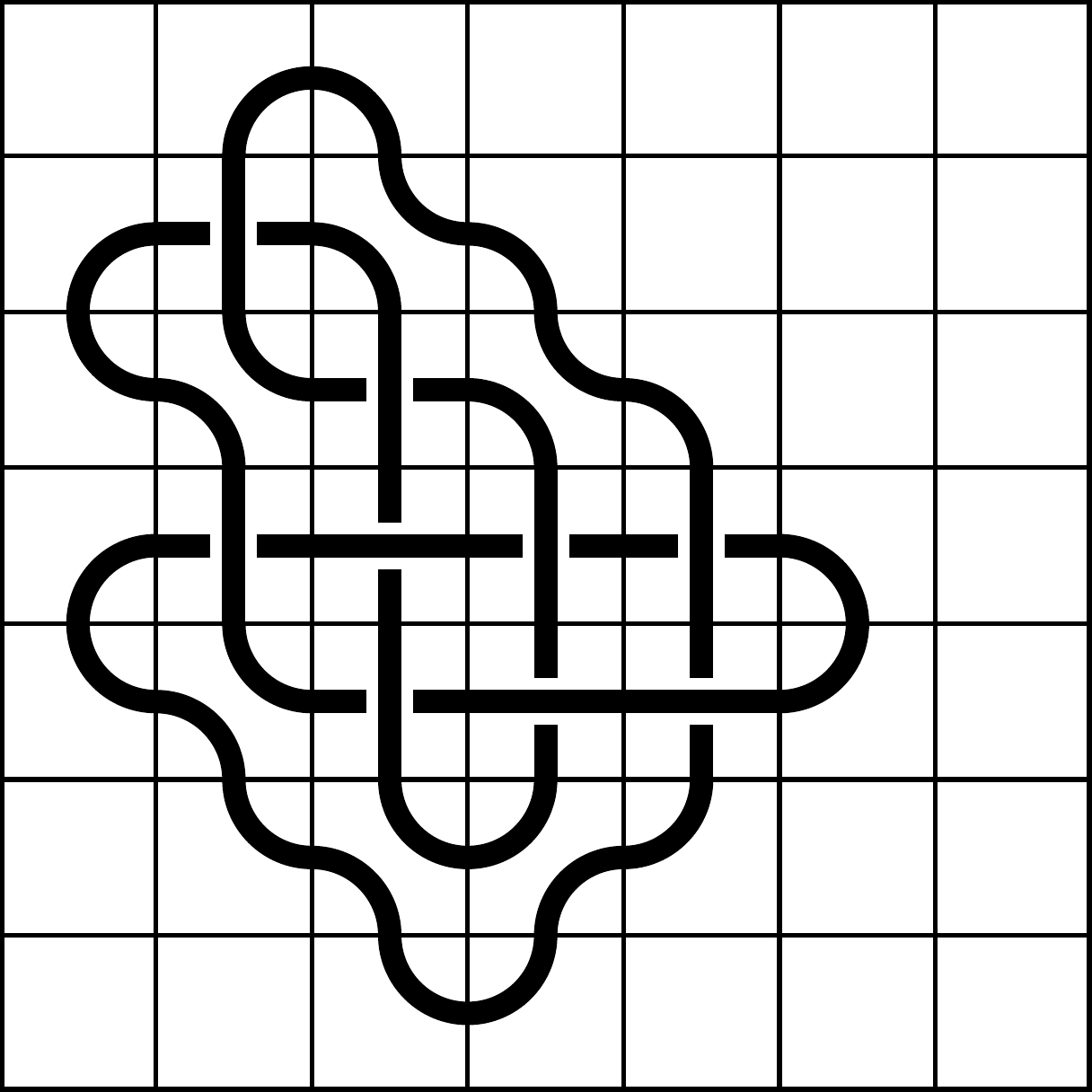}
        \caption*{$9_{48}$}
    \end{minipage} \hfill
    \begin{minipage}{0.155\linewidth}
        \captionsetup{skip=3pt}
        \centering
        \includegraphics[width=\linewidth]{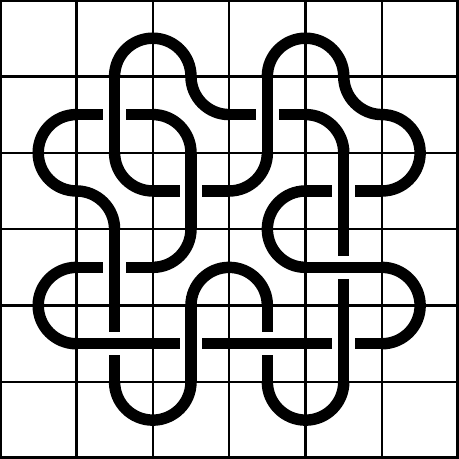}
        \caption*{$10_{1}$ {\footnotesize (m=6)}}
    \end{minipage} \hfill
    \begin{minipage}{0.155\linewidth}
        \captionsetup{skip=3pt}
        \centering
        \includegraphics[width=\linewidth]{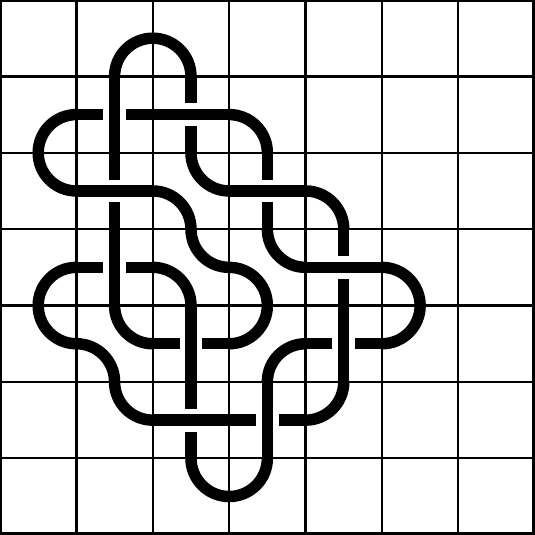}
        \caption*{$10_{1}$ {\footnotesize (less tiles)}}
    \end{minipage}  \hfill
    \begin{minipage}{0.155\linewidth}
        \captionsetup{skip=3pt}
        \centering
        \includegraphics[width=\linewidth]{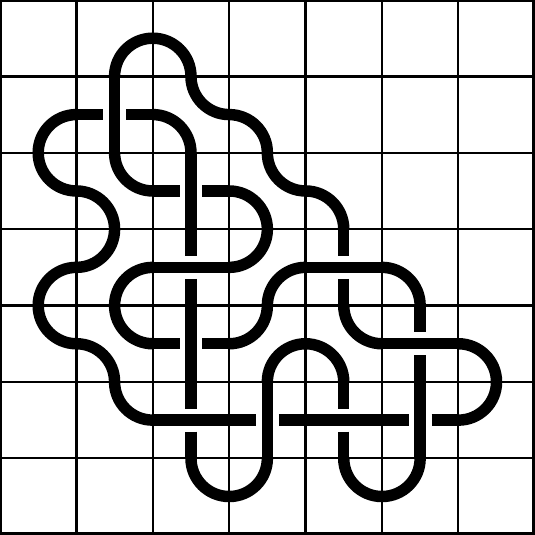}
        \caption*{$10_{3}$}
    \end{minipage} \hfill
    \begin{minipage}{0.155\linewidth}
        \captionsetup{skip=3pt}
        \centering
        \includegraphics[width=\linewidth]{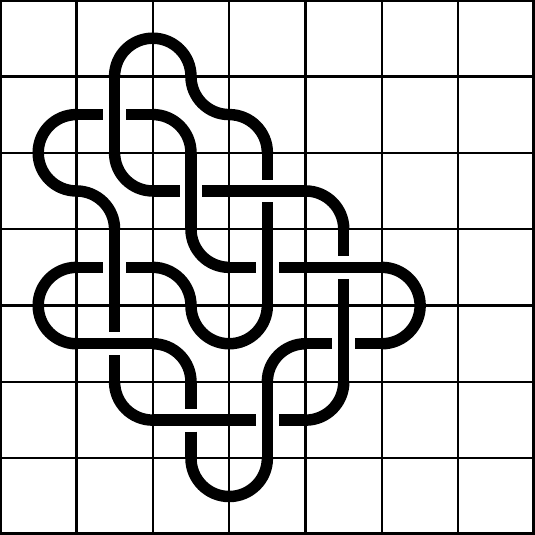}
        \caption*{$10_{5}$}
    \end{minipage} \hfill
    \begin{minipage}{0.155\linewidth}
        \captionsetup{skip=3pt}
        \centering
        \includegraphics[width=\linewidth]{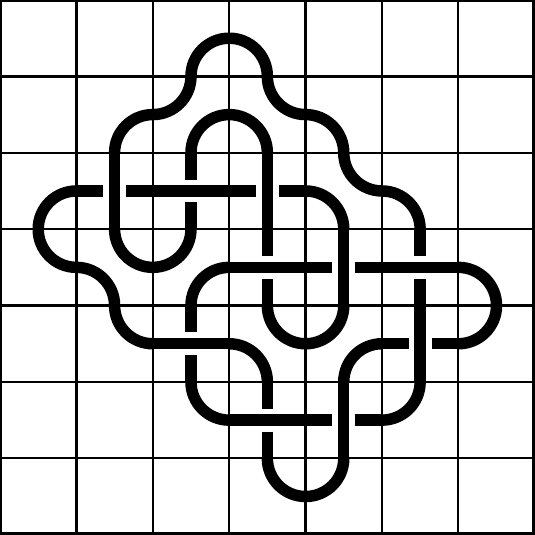}
        \caption*{$10_{6}$ }
    \end{minipage}  \newline
\end{figure}
\unskip

\begin{figure}[H]
    \centering
     \begin{minipage}{0.155\linewidth}
        \captionsetup{skip=3pt}
        \centering
        \includegraphics[width=\linewidth]{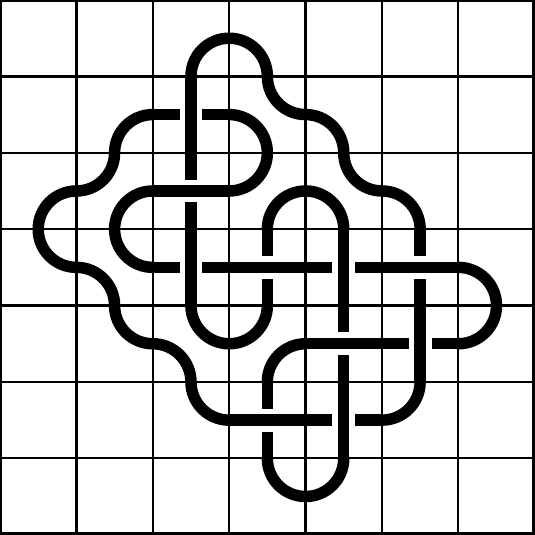}
        \caption*{$10_{7}$}
    \end{minipage} \hfill
    \begin{minipage}{0.155\linewidth}
        \captionsetup{skip=3pt}
        \centering
        \includegraphics[width=\linewidth]{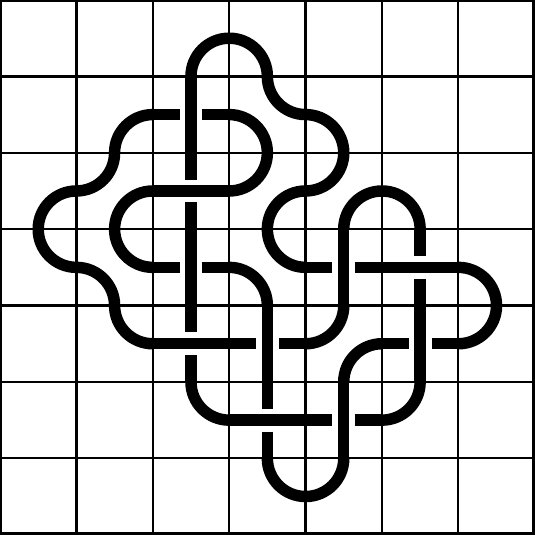}
        \caption*{$10_{9}$}
    \end{minipage} \hfill
    \begin{minipage}{0.155\linewidth}
        \captionsetup{skip=3pt}
        \centering
        \includegraphics[width=\linewidth]{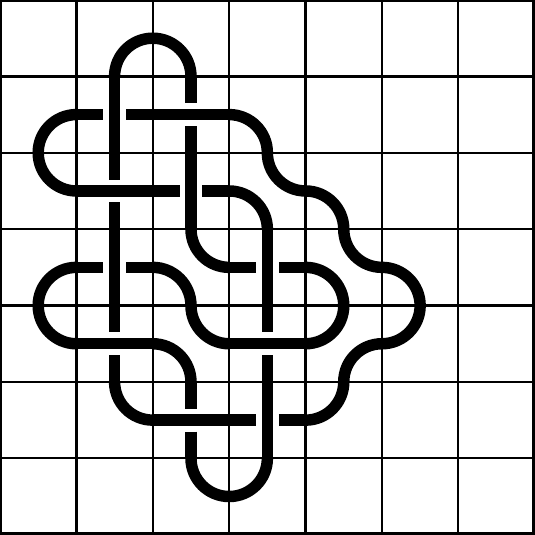}
        \caption*{$10_{11}$}
    \end{minipage}  \hfill
    \begin{minipage}{0.155\linewidth}
        \captionsetup{skip=3pt}
        \centering
        \includegraphics[width=\linewidth]{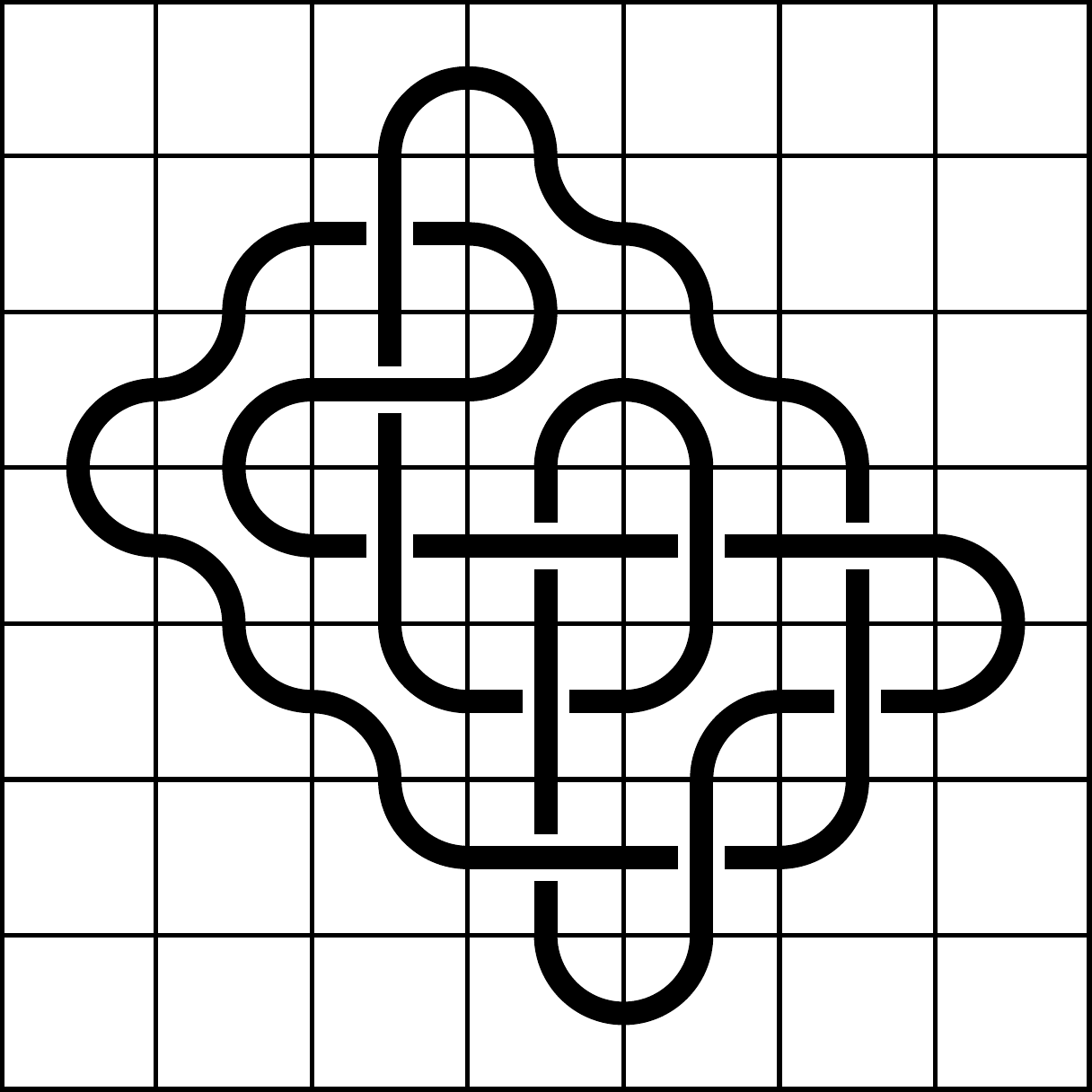}
        \caption*{$10_{12}$}
    \end{minipage} \hfill
    \begin{minipage}{0.155\linewidth}
        \captionsetup{skip=3pt}
        \centering
        \includegraphics[width=\linewidth]{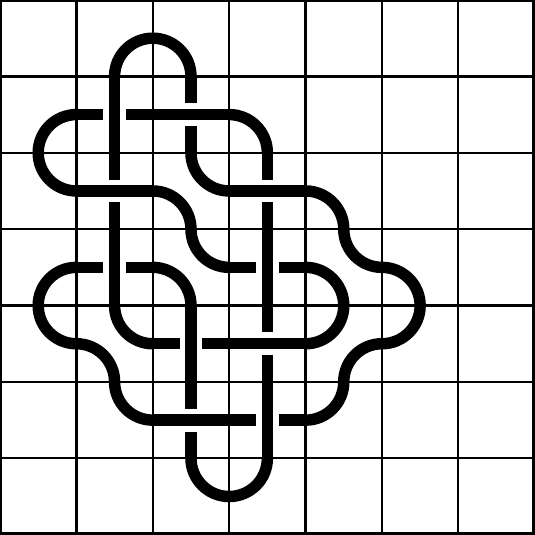}
        \caption*{$10_{13}$}
    \end{minipage} \hfill
    \begin{minipage}{0.155\linewidth}
        \captionsetup{skip=3pt}
        \centering
        \includegraphics[width=\linewidth]{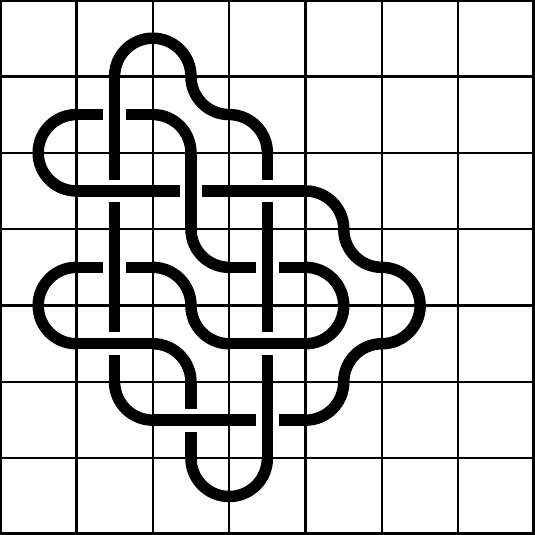}
        \caption*{$10_{14}$}
    \end{minipage}  \newline
\end{figure}
\unskip

\begin{figure}[H]
    \centering
     \begin{minipage}{0.155\linewidth}
        \captionsetup{skip=3pt}
        \centering
        \includegraphics[width=\linewidth]{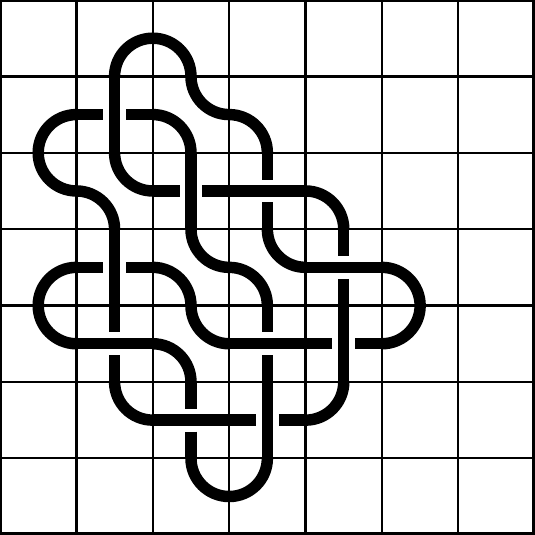}
        \caption*{$10_{15}$}
    \end{minipage} \hfill
    \begin{minipage}{0.155\linewidth}
        \captionsetup{skip=3pt}
        \centering
        \includegraphics[width=\linewidth]{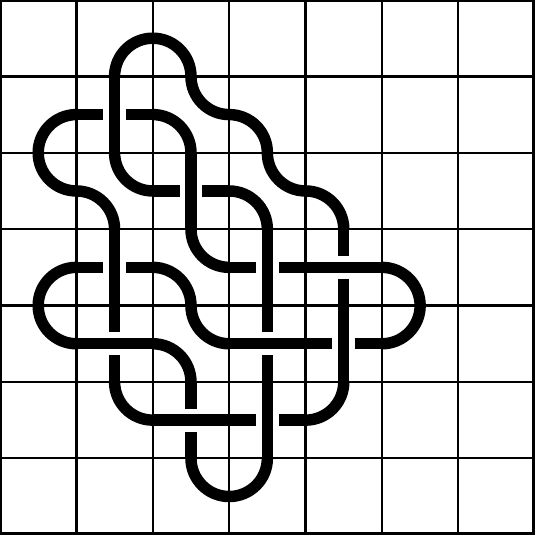}
        \caption*{$10_{16}$}
    \end{minipage} \hfill
    \begin{minipage}{0.155\linewidth}
        \captionsetup{skip=3pt}
        \centering
        \includegraphics[width=\linewidth]{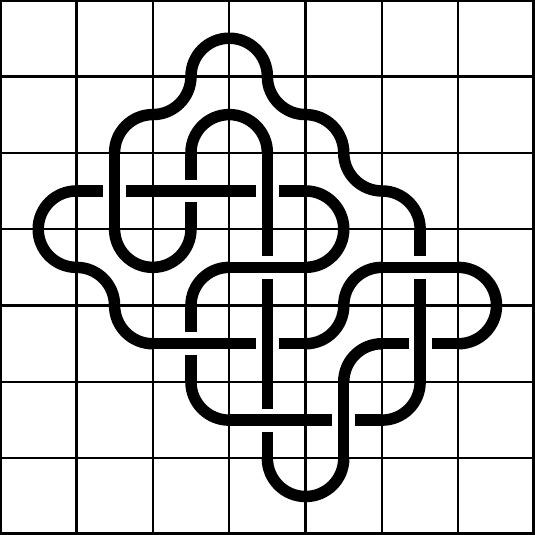}
        \caption*{$10_{17}$}
    \end{minipage}  \hfill
    \begin{minipage}{0.155\linewidth}
        \captionsetup{skip=3pt}
        \centering
        \includegraphics[width=\linewidth]{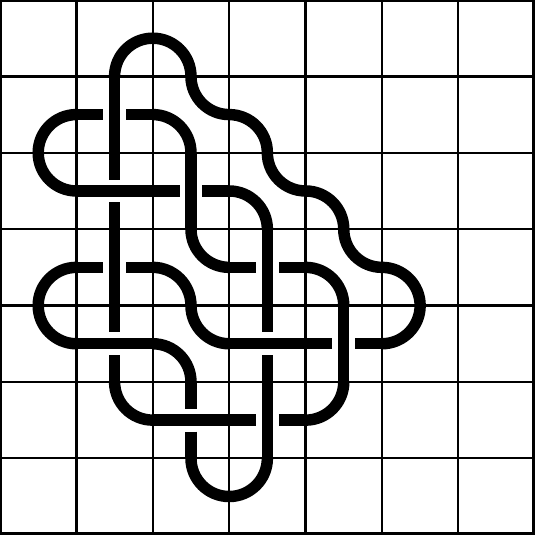}
        \caption*{$10_{18}$}
    \end{minipage} \hfill
    \begin{minipage}{0.155\linewidth}
        \captionsetup{skip=3pt}
        \centering
        \includegraphics[width=\linewidth]{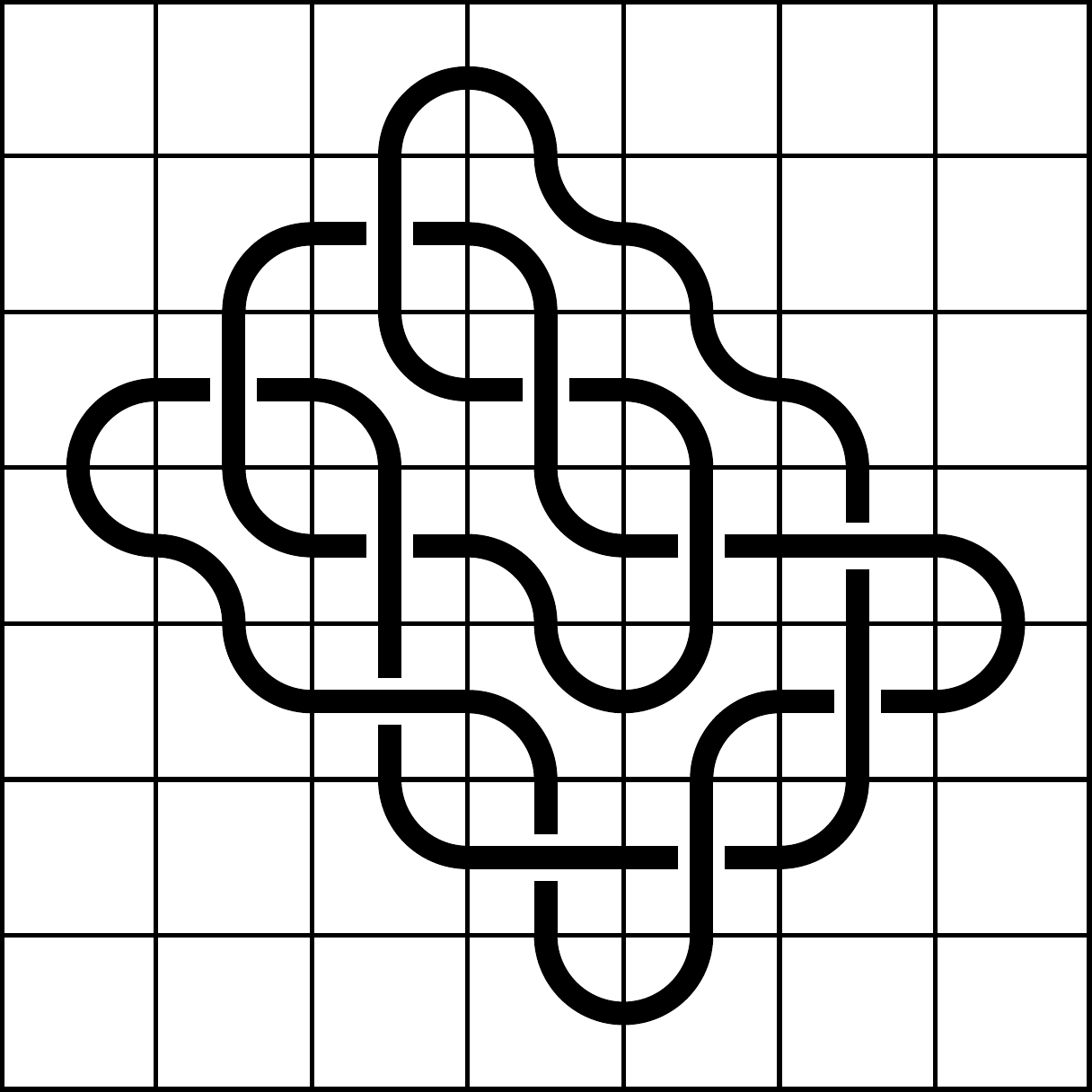}
        \caption*{$10_{20}$}
    \end{minipage} \hfill
    \begin{minipage}{0.155\linewidth}
        \captionsetup{skip=3pt}
        \centering
        \includegraphics[width=\linewidth]{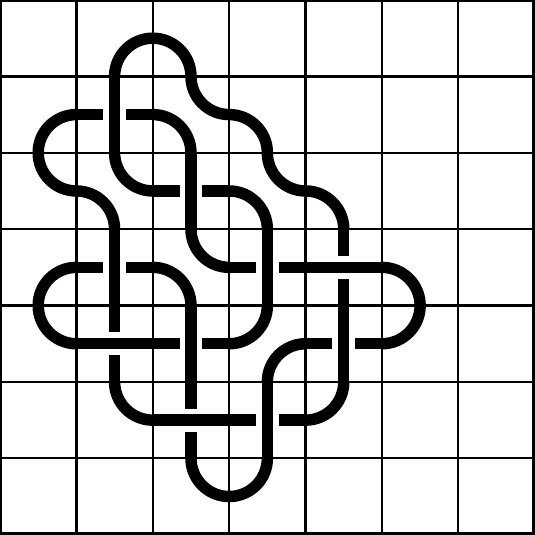}
        \caption*{$10_{21}$}
    \end{minipage}  \newline
\end{figure}
\unskip

\begin{figure}[H]
    \centering
     \begin{minipage}{0.155\linewidth}
        \captionsetup{skip=3pt}
        \centering
        \includegraphics[width=\linewidth]{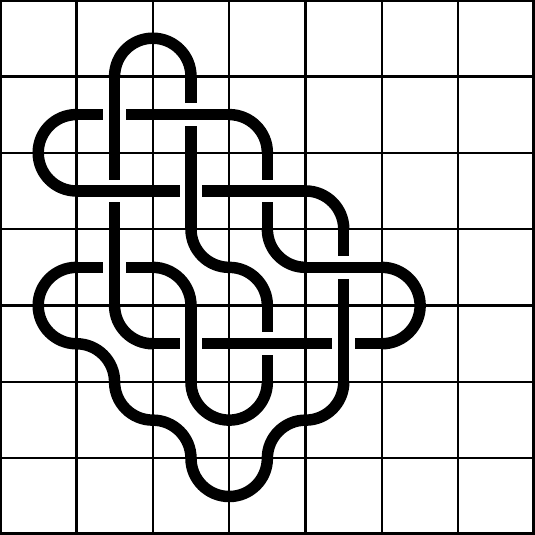}
        \caption*{$10_{22}$}
    \end{minipage} \hfill
    \begin{minipage}{0.155\linewidth}
        \captionsetup{skip=3pt}
        \centering
        \includegraphics[width=\linewidth]{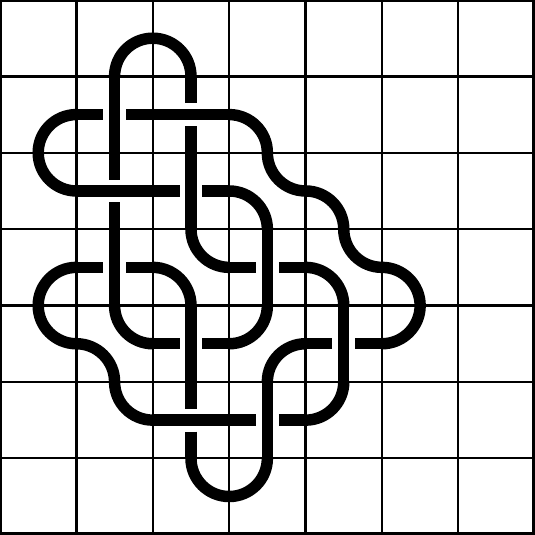}
        \caption*{$10_{24}$}
    \end{minipage} \hfill
    \begin{minipage}{0.155\linewidth}
        \captionsetup{skip=3pt}
        \centering
        \includegraphics[width=\linewidth]{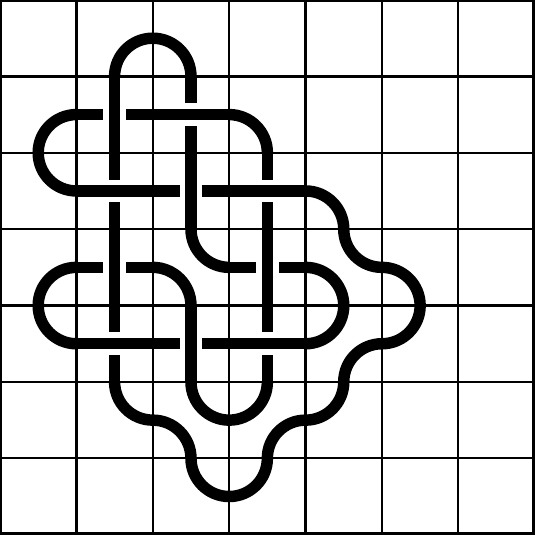}
        \caption*{$10_{31}$}
    \end{minipage}  \hfill
    \begin{minipage}{0.155\linewidth}
        \captionsetup{skip=3pt}
        \centering
        \includegraphics[width=\linewidth]{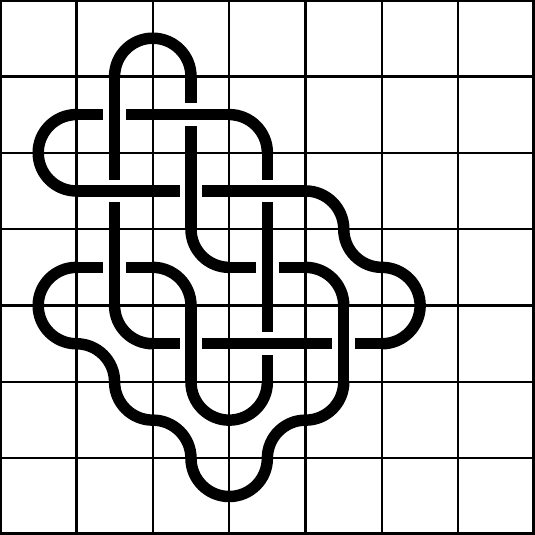}
        \caption*{$10_{33}$}
    \end{minipage} \hfill
    \begin{minipage}{0.155\linewidth}
        \captionsetup{skip=3pt}
        \centering
        \includegraphics[width=\linewidth]{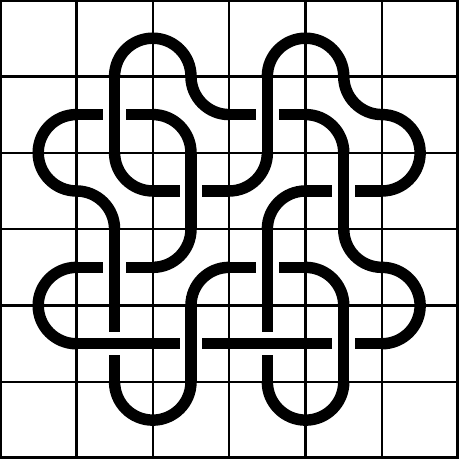}
        \caption*{\small $10_{34}$ {\footnotesize (m=6)}}
    \end{minipage} \hfill
    \begin{minipage}{0.155\linewidth}
        \captionsetup{skip=3pt}
        \centering
        \includegraphics[width=\linewidth]{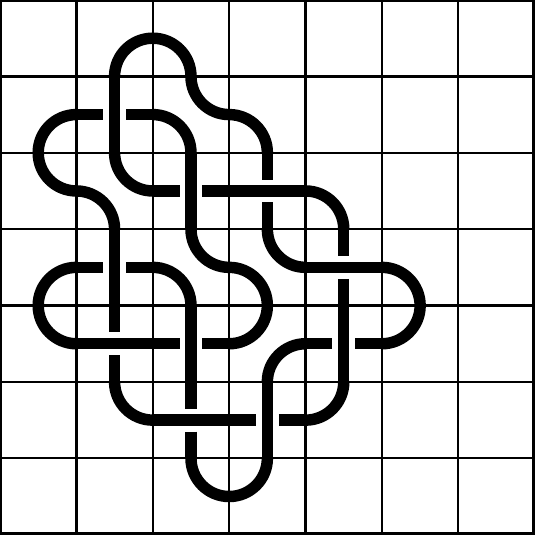}
        \caption*{\small $10_{34}$ {\footnotesize (less tiles)}}
    \end{minipage}  \newline
\end{figure}
\unskip

\begin{figure}[H]
    \centering
     \begin{minipage}{0.155\linewidth}
        \captionsetup{skip=3pt}
        \centering
        \includegraphics[width=\linewidth]{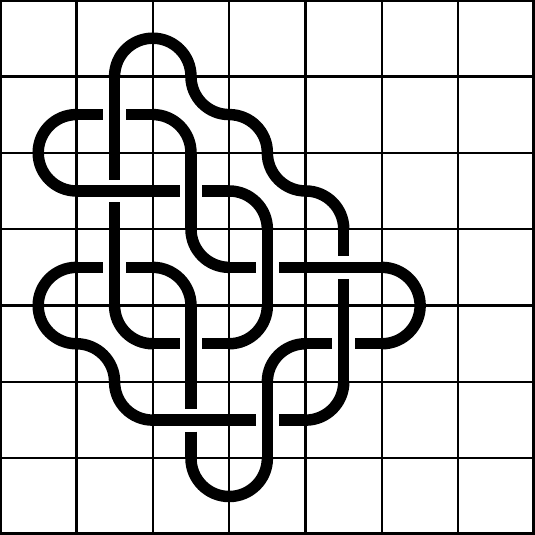}
        \caption*{$10_{35}$}
    \end{minipage} \hfill
    \begin{minipage}{0.155\linewidth}
        \captionsetup{skip=3pt}
        \centering
        \includegraphics[width=\linewidth]{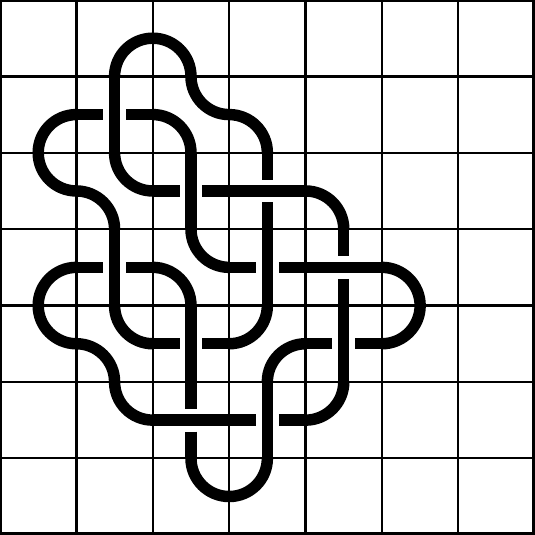}
        \caption*{$10_{36}$}
    \end{minipage} \hfill
    \begin{minipage}{0.155\linewidth}
        \captionsetup{skip=3pt}
        \centering
        \includegraphics[width=\linewidth]{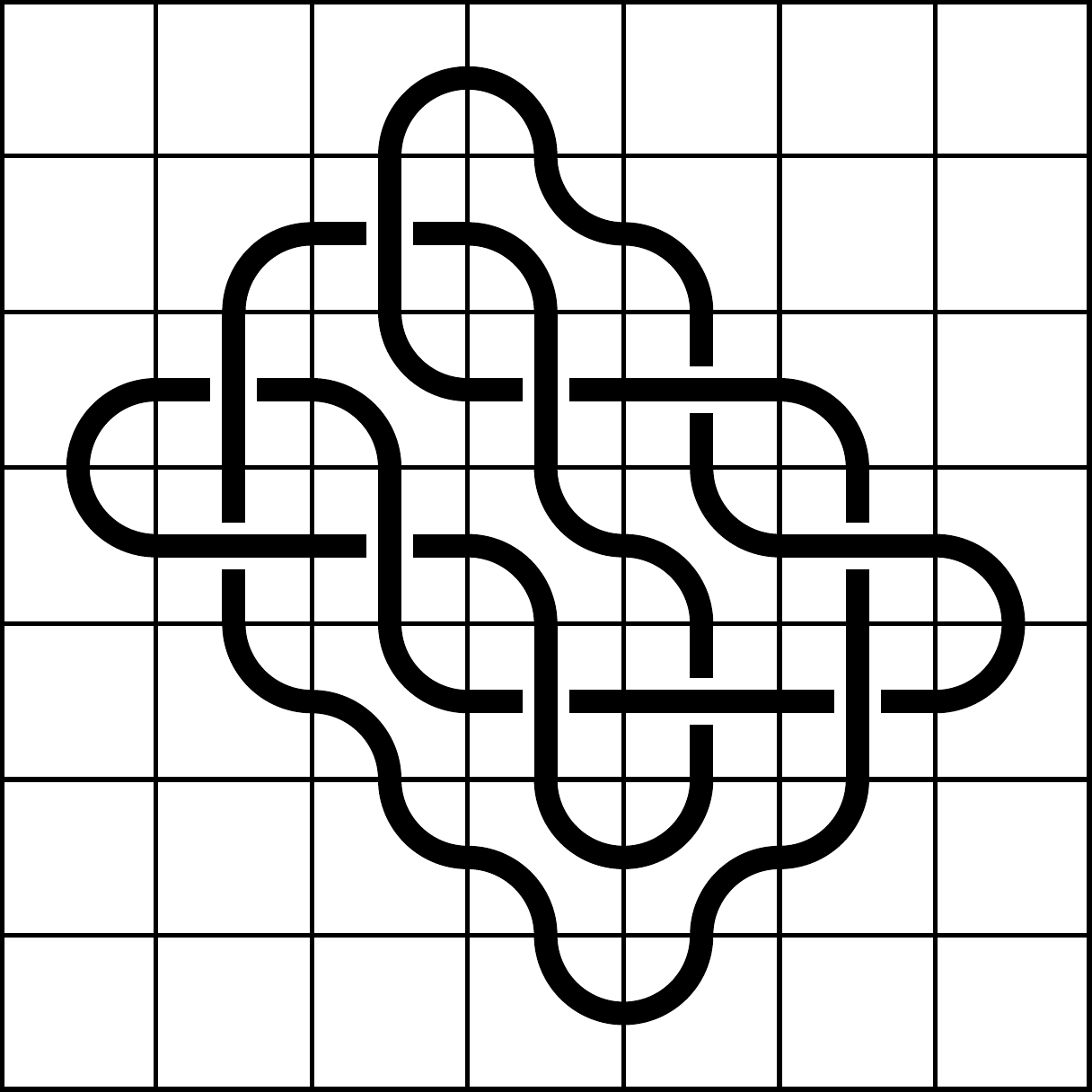}
        \caption*{$10_{37}$}
    \end{minipage}  \hfill
    \begin{minipage}{0.155\linewidth}
        \captionsetup{skip=3pt}
        \centering
        \includegraphics[width=\linewidth]{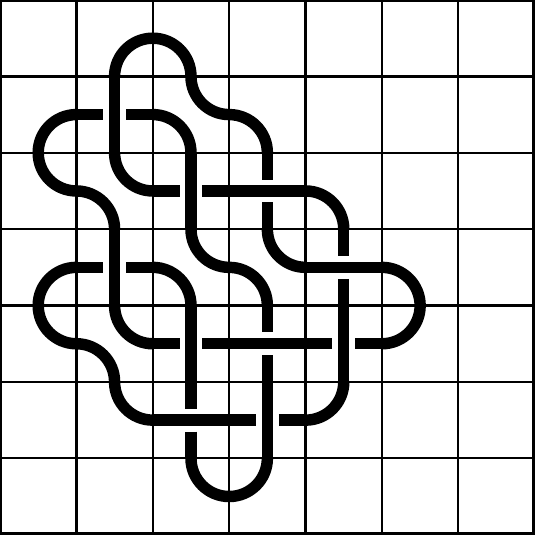}
        \caption*{$10_{38}$}
    \end{minipage} \hfill
    \begin{minipage}{0.155\linewidth}
        \captionsetup{skip=3pt}
        \centering
        \includegraphics[width=\linewidth]{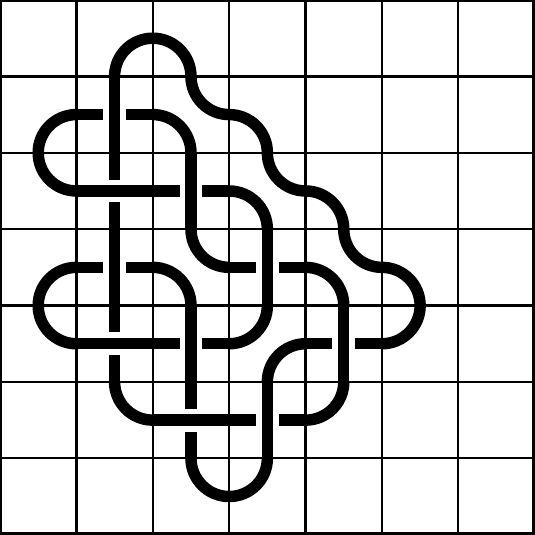}
        \caption*{$10_{39}$}
    \end{minipage} \hfill
    \begin{minipage}{0.155\linewidth}
        \captionsetup{skip=3pt}
        \centering
        \includegraphics[width=\linewidth]{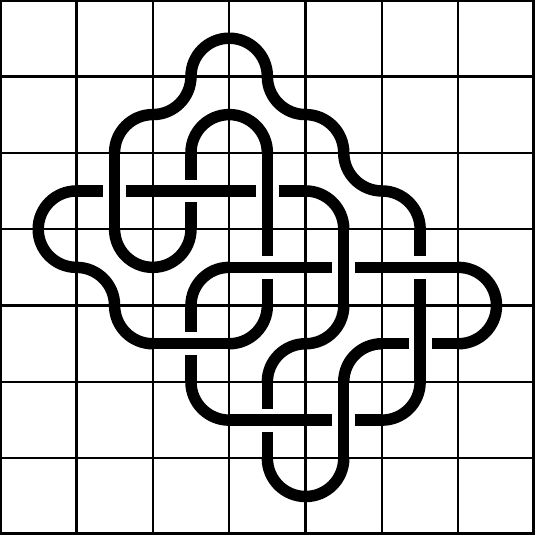}
        \caption*{$10_{48}$}
    \end{minipage}  \newline
\end{figure}
\unskip

\begin{figure}[H]
    \centering
     \begin{minipage}{0.155\linewidth}
        \captionsetup{skip=3pt}
        \centering
        \includegraphics[width=\linewidth]{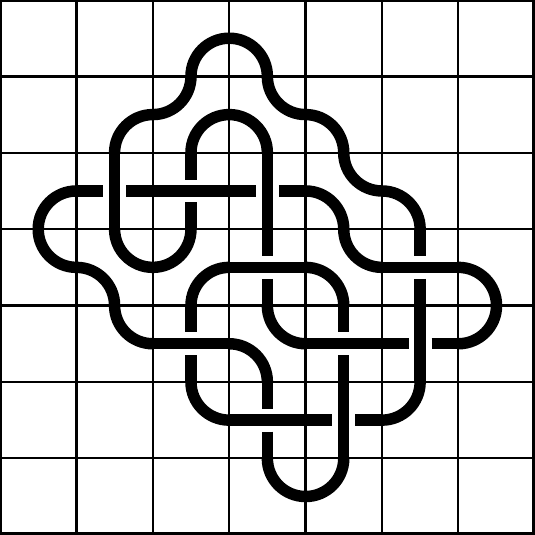}
        \caption*{$10_{50}$}
    \end{minipage} \hfill
    \begin{minipage}{0.155\linewidth}
        \captionsetup{skip=3pt}
        \centering
        \includegraphics[width=\linewidth]{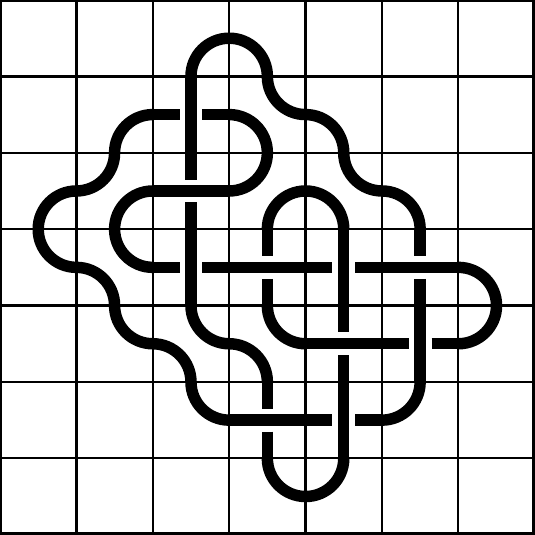}
        \caption*{$10_{51}$}
    \end{minipage} \hfill
    \begin{minipage}{0.155\linewidth}
        \captionsetup{skip=3pt}
        \centering
        \includegraphics[width=\linewidth]{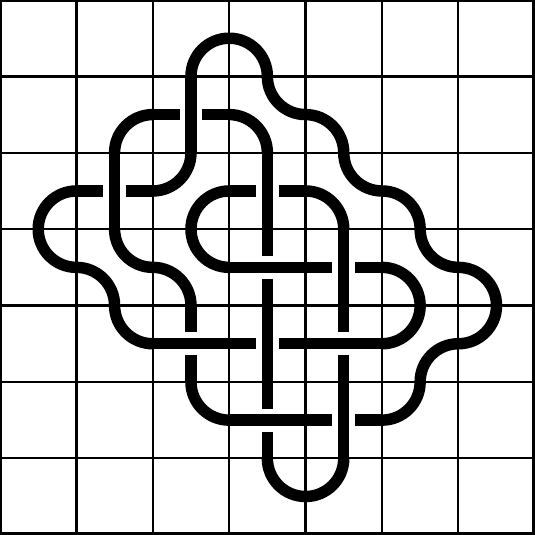}
        \caption*{$10_{56}$}
    \end{minipage}  \hfill
    \begin{minipage}{0.155\linewidth}
        \captionsetup{skip=3pt}
        \centering
        \includegraphics[width=\linewidth]{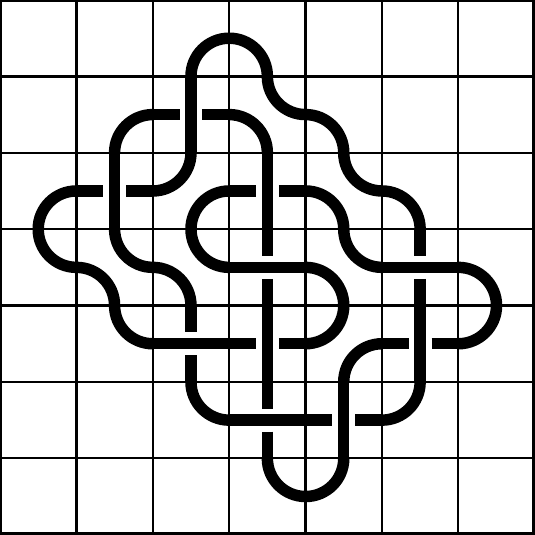}
        \caption*{$10_{61}$}
    \end{minipage} \hfill
    \begin{minipage}{0.155\linewidth}
        \captionsetup{skip=3pt}
        \centering
        \includegraphics[width=\linewidth]{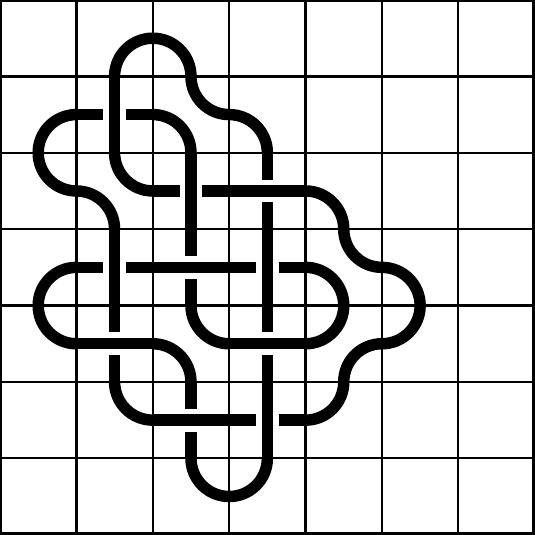}
        \caption*{$10_{63}$}
    \end{minipage} \hfill
    \begin{minipage}{0.155\linewidth}
        \captionsetup{skip=3pt}
        \centering
        \includegraphics[width=\linewidth]{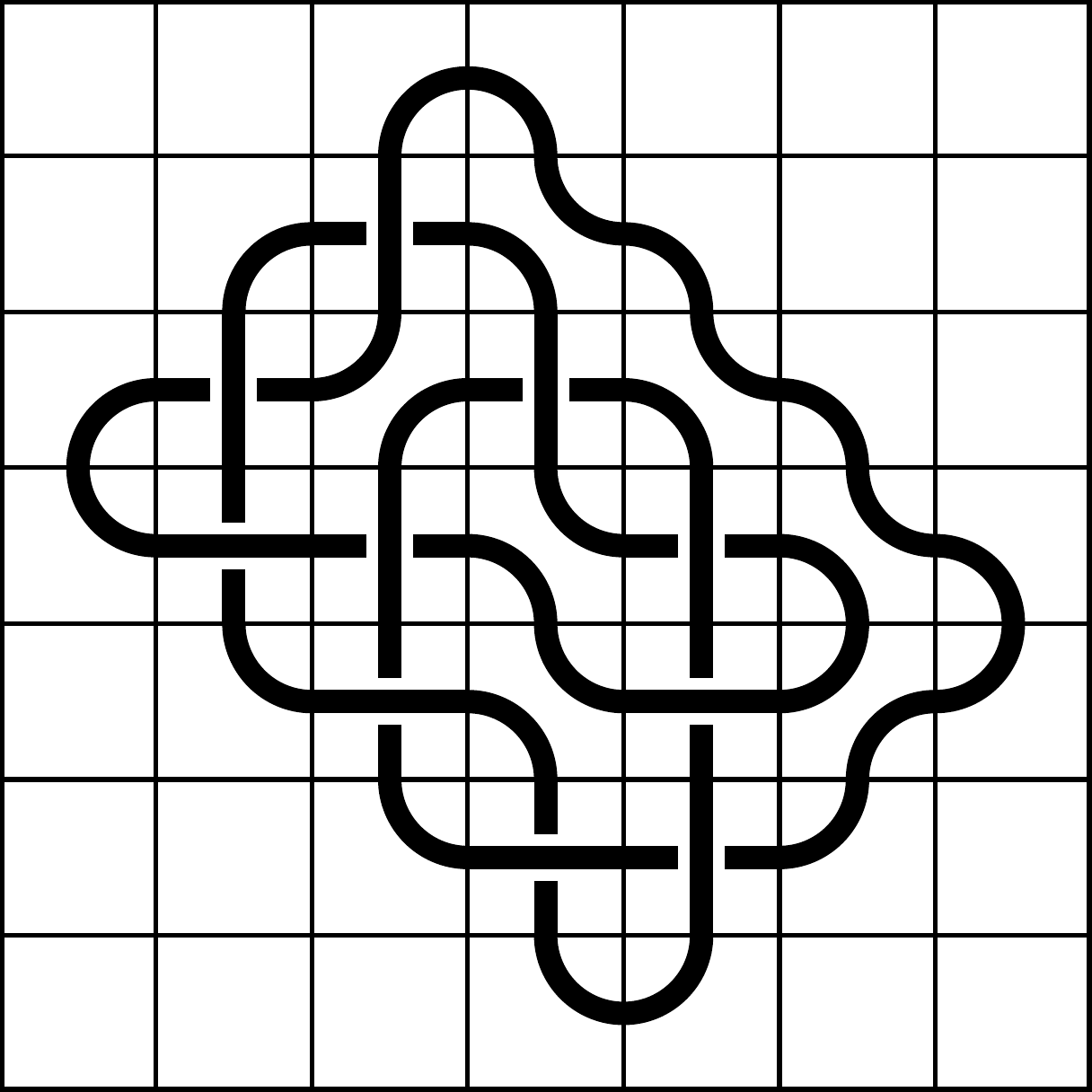}
        \caption*{$10_{64}$}
    \end{minipage}  \newline
\end{figure}
\unskip

\begin{figure}[H]
    \centering
     \begin{minipage}{0.155\linewidth}
        \captionsetup{skip=3pt}
        \centering
        \includegraphics[width=\linewidth]{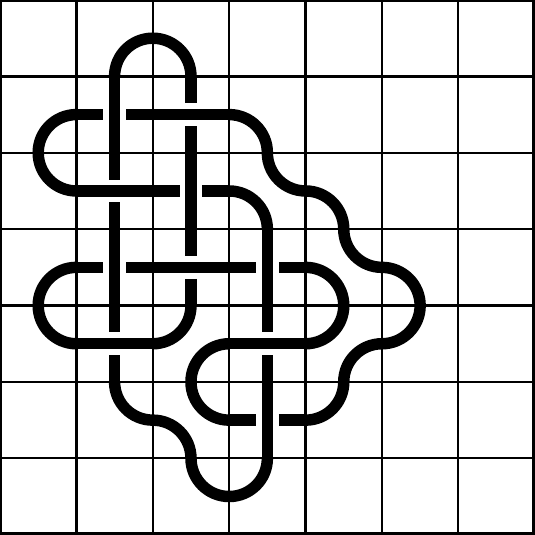}
        \caption*{$10_{65}$}
    \end{minipage} \hfill
    \begin{minipage}{0.155\linewidth}
        \captionsetup{skip=3pt}
        \centering
        \includegraphics[width=\linewidth]{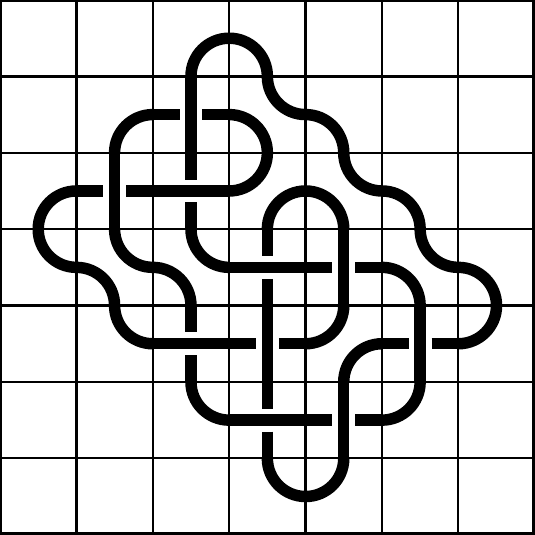}
        \caption*{$10_{67}$}
    \end{minipage} \hfill
    \begin{minipage}{0.155\linewidth}
        \captionsetup{skip=3pt}
        \centering
        \includegraphics[width=\linewidth]{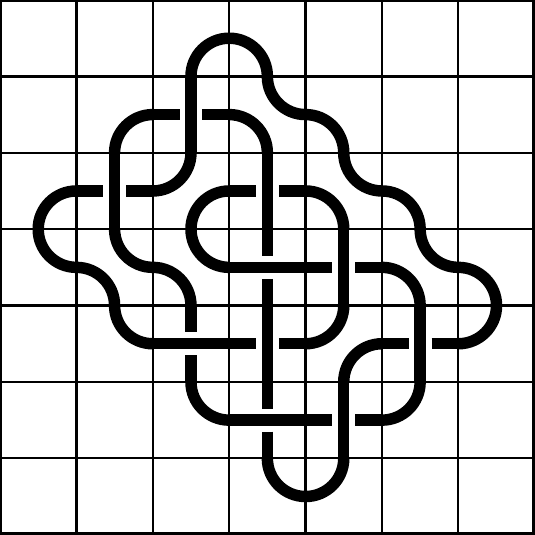}
        \caption*{$10_{68}$}
    \end{minipage}  \hfill
    \begin{minipage}{0.155\linewidth}
        \captionsetup{skip=3pt}
        \centering
        \includegraphics[width=\linewidth]{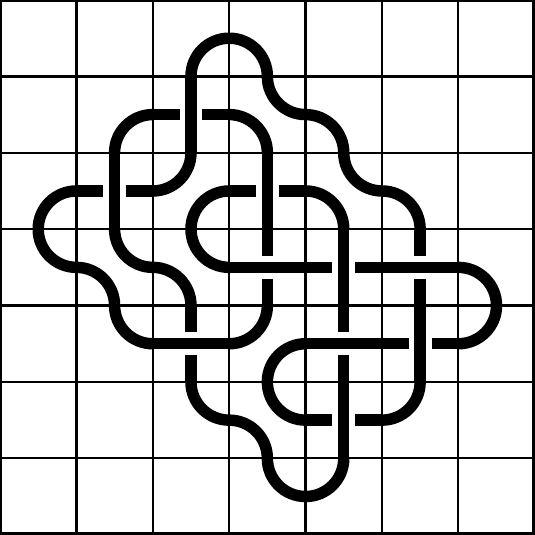}
        \caption*{$10_{70}$}
    \end{minipage} \hfill
    \begin{minipage}{0.155\linewidth}
        \captionsetup{skip=3pt}
        \centering
        \includegraphics[width=\linewidth]{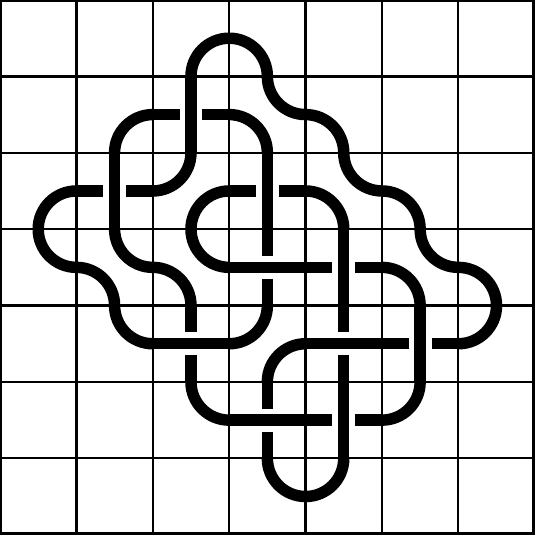}
        \caption*{$10_{72}$}
    \end{minipage} \hfill
    \begin{minipage}{0.155\linewidth}
        \captionsetup{skip=3pt}
        \centering
        \includegraphics[width=\linewidth]{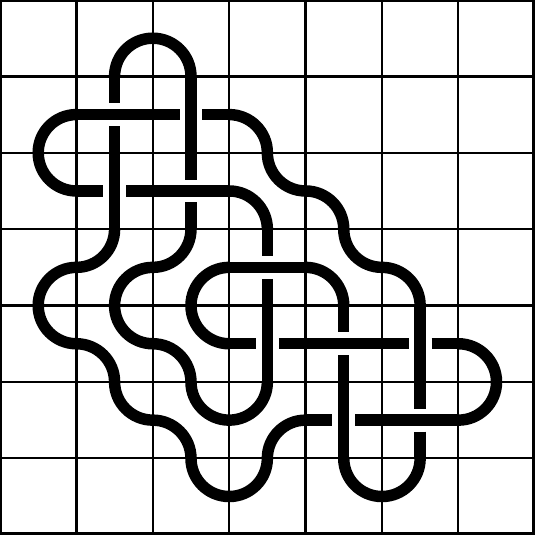}
        \caption*{$10_{76}$}
    \end{minipage}  \newline
\end{figure}
\unskip

\begin{figure}[H]
    \centering
     \begin{minipage}{0.155\linewidth}
        \captionsetup{skip=3pt}
        \centering
        \includegraphics[width=\linewidth]{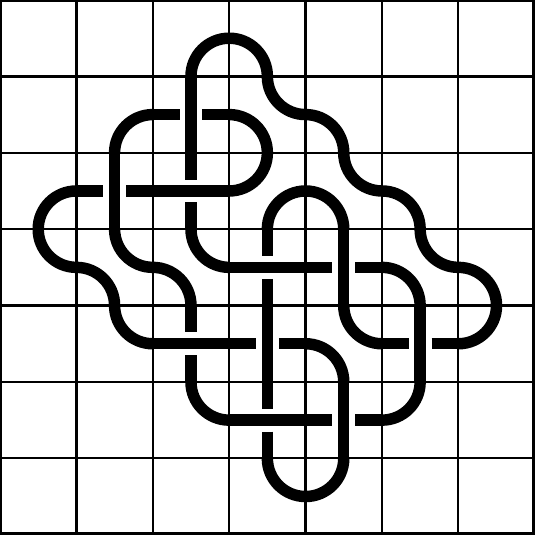}
        \caption*{$10_{77}$}
    \end{minipage} \hfill
    \begin{minipage}{0.155\linewidth}
        \captionsetup{skip=3pt}
        \centering
        \includegraphics[width=\linewidth]{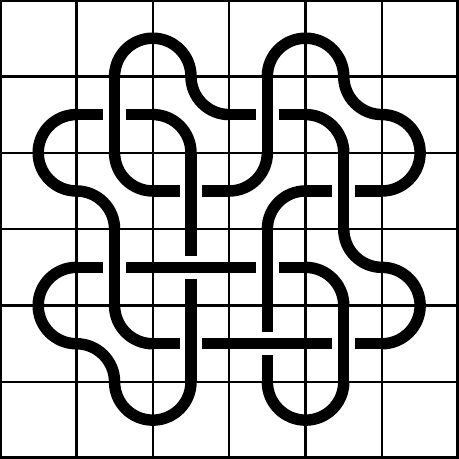}
        \caption*{\small $10_{78}$ {\footnotesize (m=6)}}
    \end{minipage} \hfill
    \begin{minipage}{0.155\linewidth}
        \captionsetup{skip=3pt}
        \centering
        \includegraphics[width=\linewidth]{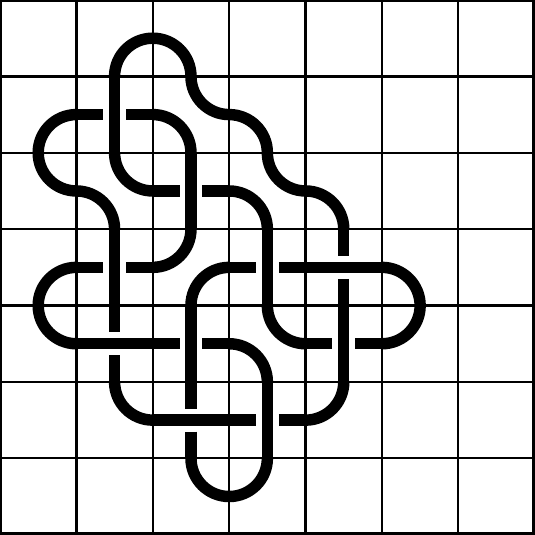}
        \caption*{\small $10_{78}$ {\footnotesize (less tiles)}}
    \end{minipage}  \hfill
    \begin{minipage}{0.155\linewidth}
        \captionsetup{skip=3pt}
        \centering
        \includegraphics[width=\linewidth]{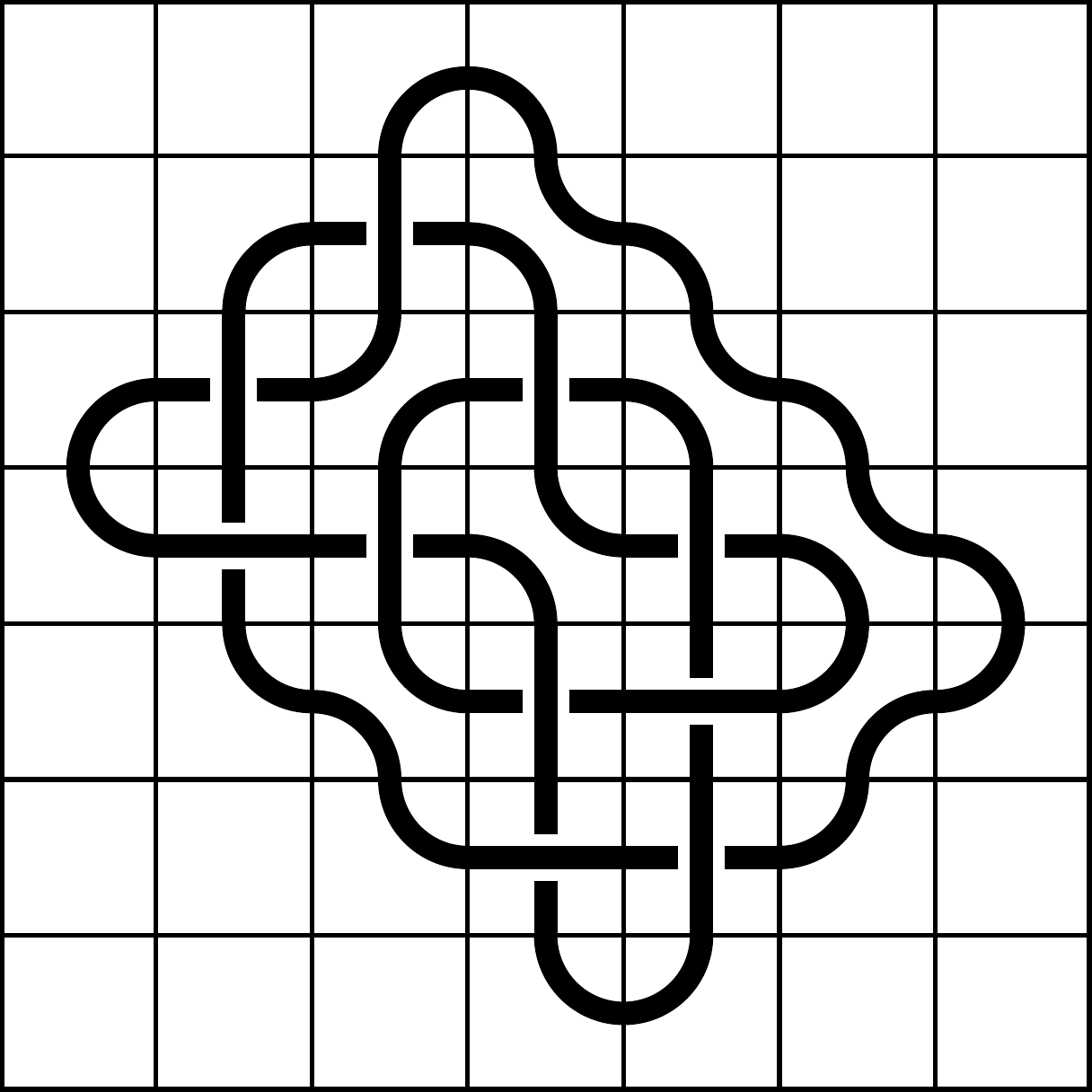}
        \caption*{$10_{79}$}
    \end{minipage} \hfill
    \begin{minipage}{0.155\linewidth}
        \captionsetup{skip=3pt}
        \centering
        \includegraphics[width=\linewidth]{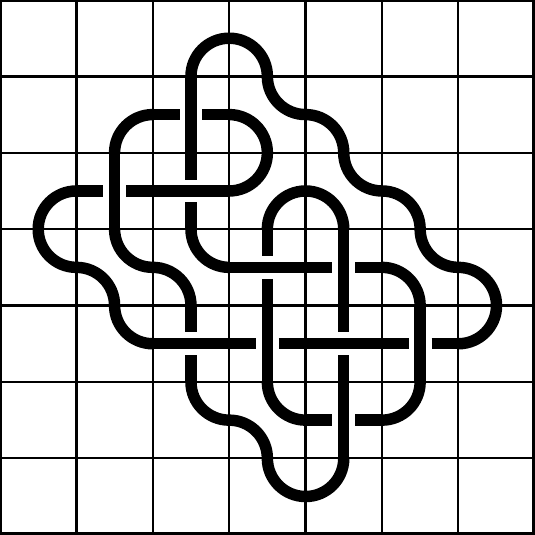}
        \caption*{$10_{84}$}
    \end{minipage} \hfill
    \begin{minipage}{0.155\linewidth}
        \captionsetup{skip=3pt}
        \centering
        \includegraphics[width=\linewidth]{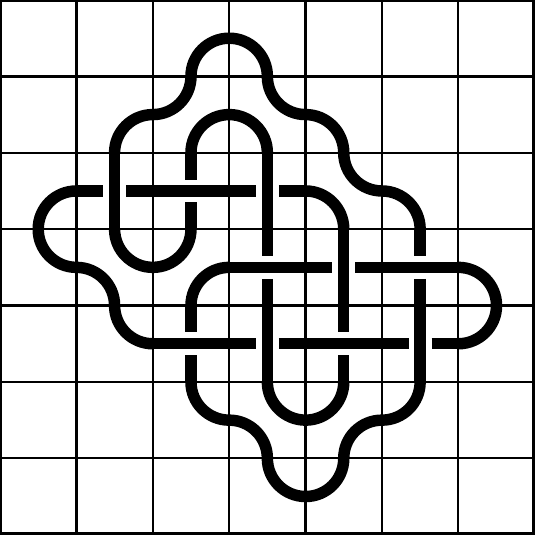}
        \caption*{$10_{90}$}
    \end{minipage}  \newline
\end{figure}
\unskip

\begin{figure}[H]
    \centering
     \begin{minipage}{0.155\linewidth}
        \captionsetup{skip=3pt}
        \centering
        \includegraphics[width=\linewidth]{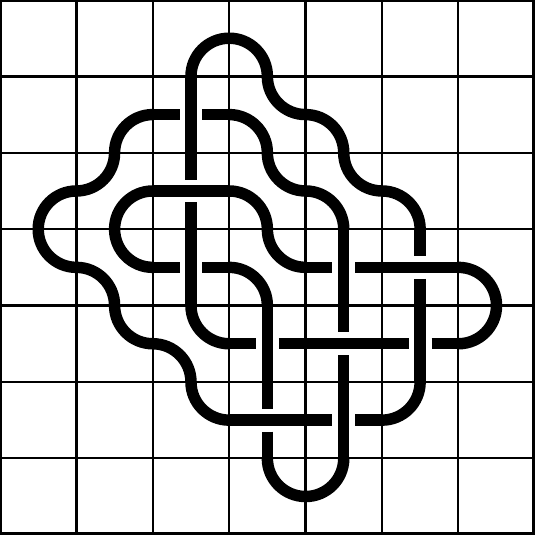}
        \caption*{$10_{91}$}
    \end{minipage} \hfill
    \begin{minipage}{0.155\linewidth}
        \captionsetup{skip=3pt}
        \centering
        \includegraphics[width=\linewidth]{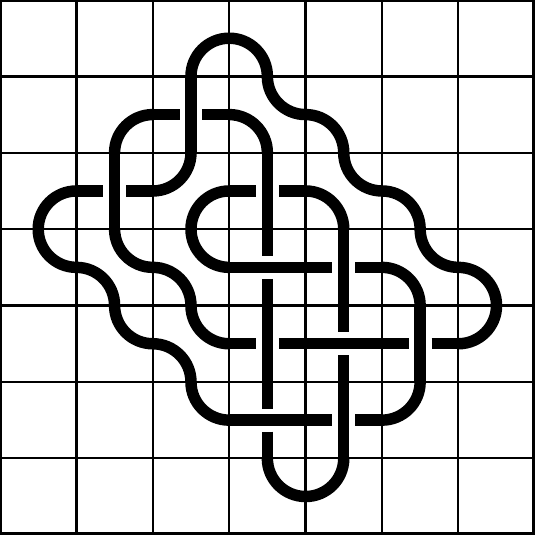}
        \caption*{$10_{92}$}
    \end{minipage} \hfill
    \begin{minipage}{0.155\linewidth}
        \captionsetup{skip=3pt}
        \centering
        \includegraphics[width=\linewidth]{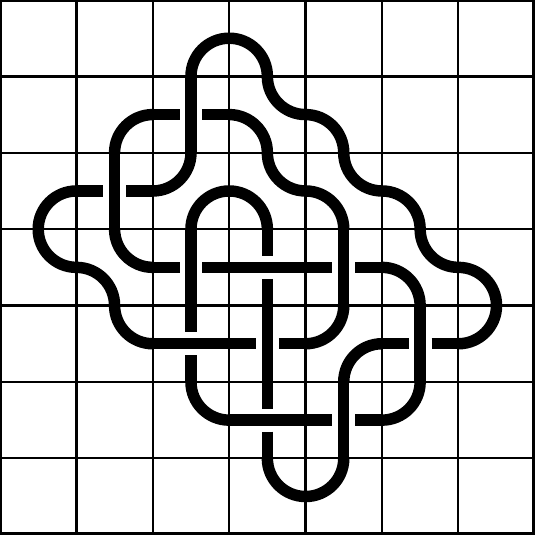}
        \caption*{$10_{93}$}
    \end{minipage}  \hfill
    \begin{minipage}{0.155\linewidth}
        \captionsetup{skip=3pt}
        \centering
        \includegraphics[width=\linewidth]{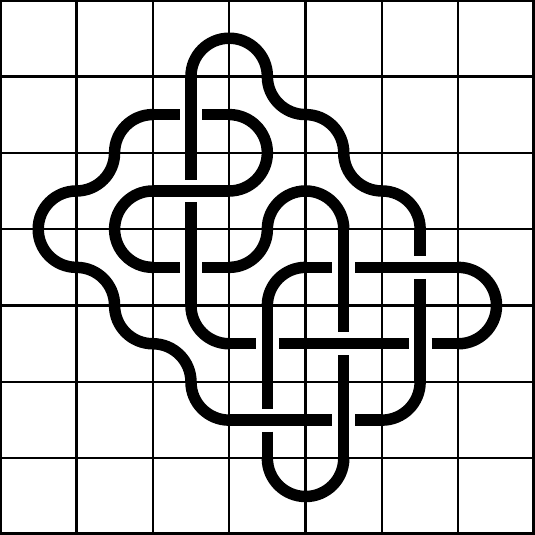}
        \caption*{$10_{103}$}
    \end{minipage} \hfill
    \begin{minipage}{0.155\linewidth}
        \captionsetup{skip=3pt}
        \centering
        \includegraphics[width=\linewidth]{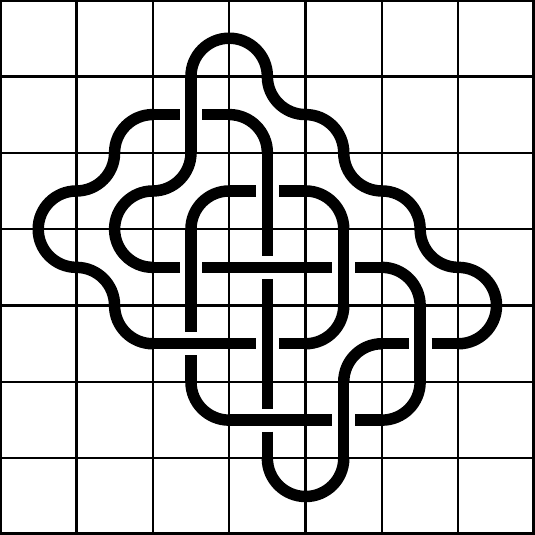}
        \caption*{$10_{114}$}
    \end{minipage} \hfill
    \begin{minipage}{0.155\linewidth}
        \captionsetup{skip=3pt}
        \centering
        \includegraphics[width=\linewidth]{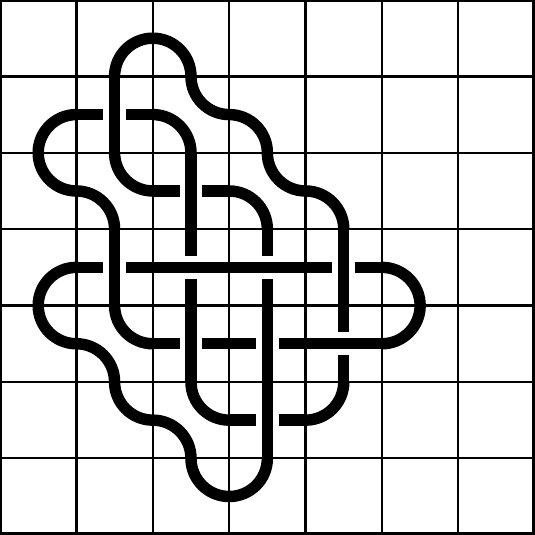}
        \caption*{$10_{140}$}
    \end{minipage}  \newline
\end{figure}
\unskip

\begin{figure}[H]
    \centering
     \begin{minipage}{0.155\linewidth}
        \captionsetup{skip=3pt}
        \centering
        \includegraphics[width=\linewidth]{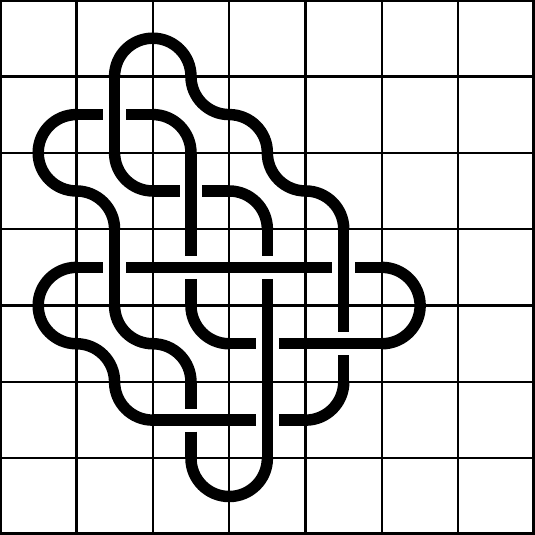}
        \caption*{$10_{142}$}
    \end{minipage} \hfill
    \begin{minipage}{0.155\linewidth}
        \captionsetup{skip=3pt}
        \centering
        \includegraphics[width=\linewidth]{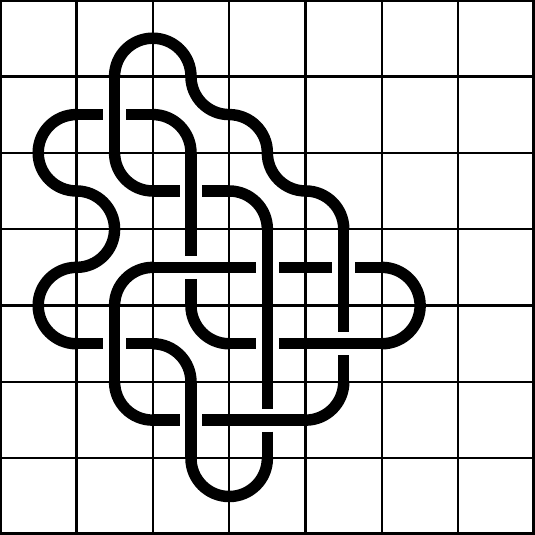}
        \caption*{$10_{144}$}
    \end{minipage} \hfill
    \begin{minipage}{0.155\linewidth}
        \captionsetup{skip=3pt}
        \centering
        \includegraphics[width=\linewidth]{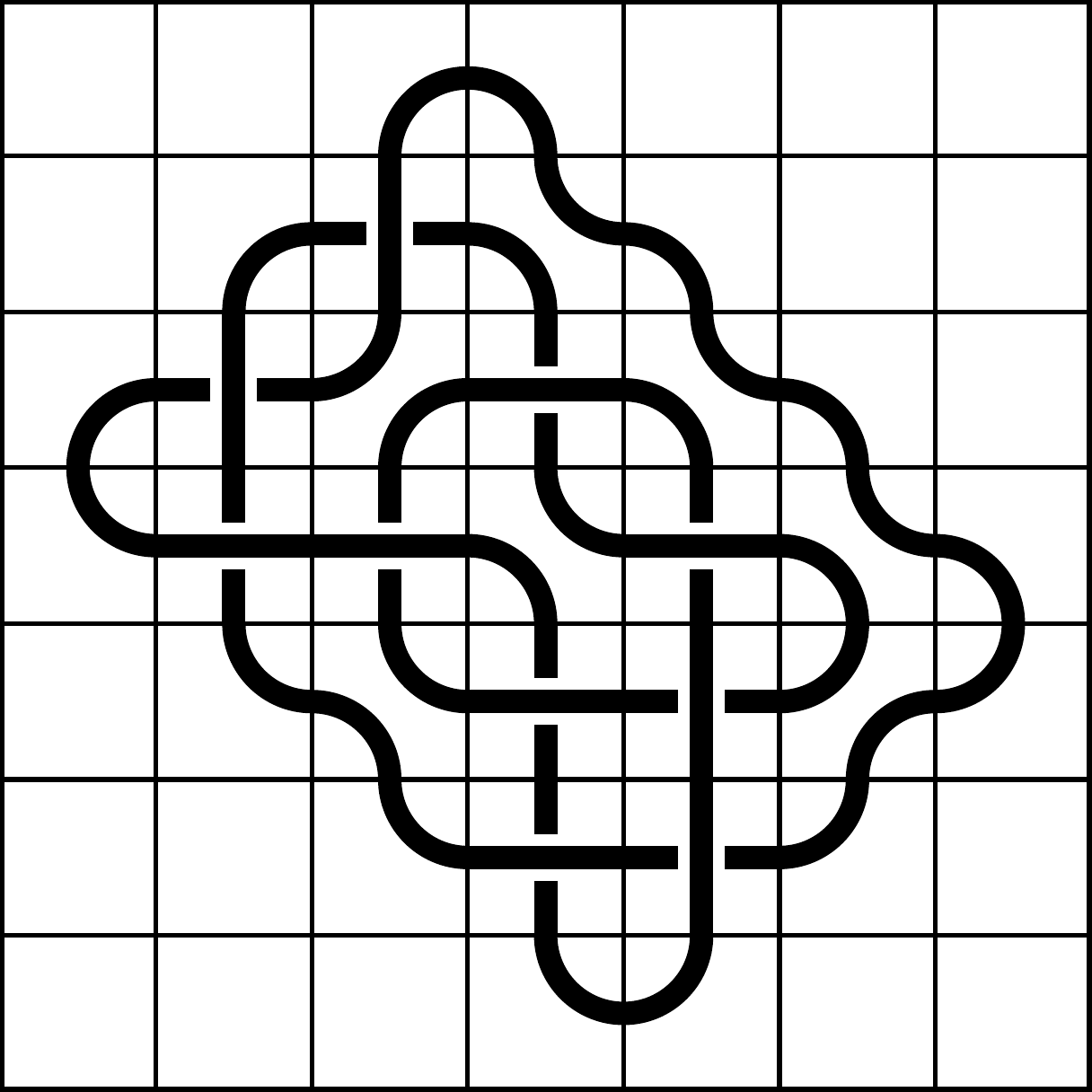}
        \caption*{$10_{152}$}
    \end{minipage}  \hfill
    \begin{minipage}{0.155\linewidth}
        \captionsetup{skip=3pt}
        \centering
        \includegraphics[width=\linewidth]{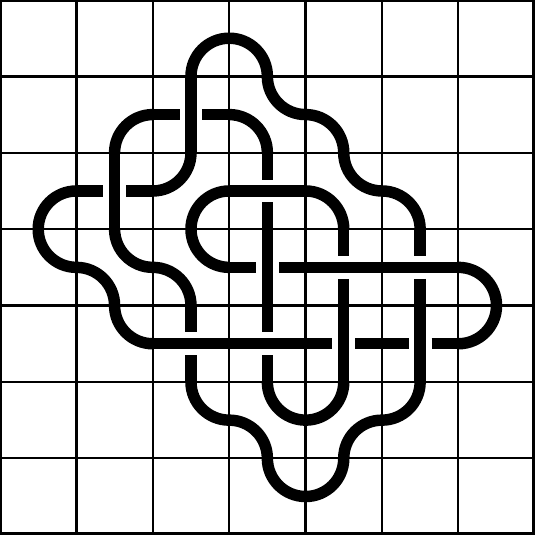}
        \caption*{$10_{153}$}
    \end{minipage} \hfill
    \begin{minipage}{0.155\linewidth}
        \captionsetup{skip=3pt}
        \centering
        \includegraphics[width=\linewidth]{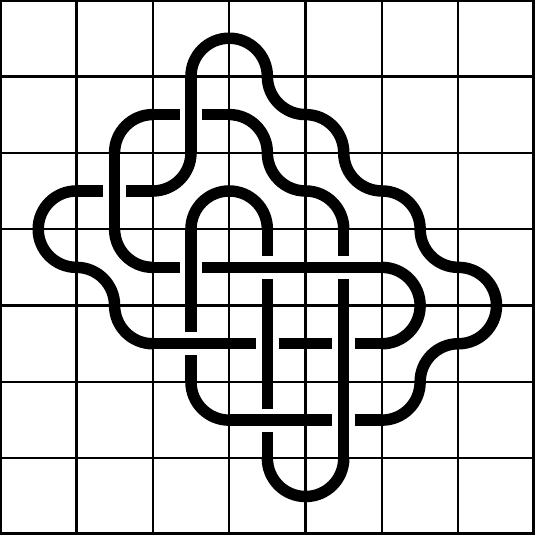}
        \caption*{$10_{158}$}
    \end{minipage} \hfill
    \begin{minipage}{0.155\linewidth}
        \captionsetup{skip=3pt}
        \centering
        \includegraphics[width=\linewidth]{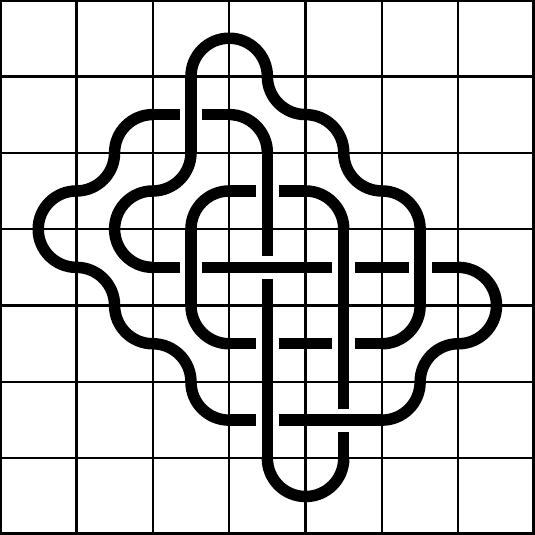}
        \caption*{$10_{163}$}
    \end{minipage}
\end{figure}


\clearpage
\bibliographystyle{amsplain}
\bibliography{bibliography}
\addcontentsline{toc}{section}{\refname}

\end{document}